# Parallel and Distributed Successive Convex Approximation Methods for Big-Data Optimization

Gesualdo Scutari and Ying Sun

January 15, 2018

Gesualdo Scutari
Purdue University, West Lafayette, IN, USA, e-mail: `gscutari@purdue.edu`

Ying Sun
Purdue University, West Lafayette, IN, USA, e-mail: `sun578@purdue.edu`

This work has been partially supported by the USA National Science Foundation under Grants CIF 1564044, CIF 1719205, and CAREER Award 1555850; and in part by the Office of Naval Research under Grant N00014-16-1-2244.

## Introduction

Recent years have witnessed a surge of interest in parallel and distributed optimization methods for large-scale systems. In particular, *nonconvex large-scale* optimization problems have found a wide range of applications in several engineering fields as diverse as (networked) information processing (e.g., parameter estimation, detection, localization, graph signal processing), communication networks (e.g., resource allocation in peer-to-peer/multi-cellular systems), sensor networks, data-based networks (including Facebook, Google, Twitter, and YouTube), swarm robotic, and machine learning (e.g., nonlinear least squares, dictionary learning, matrix completion, tensor factorization), just to name a few–see Figure 1.

The design and the analysis of such complex, large-scale, systems pose several challenges and call for the development of new optimization models and algorithms.

- **Big-Data:** Many of the aforementioned applications lead to *huge-scale* optimization problems (i.e., problems with a very large number of variables). These problems are often referred to as *big-data*. This calls for the development of solution methods that operate *in parallel*, exploiting hierarchical computational architectures (e.g., multicore systems, cluster computers, cloud-based networks), if available, to cope with the curse of dimensionality and accommodate the need of fast (real-time) processing and optimization. The challenge is that such optimization problems are in general not separable in the optimization variables, which makes the design of parallel schemes not a trivial task.

- **In-network optimization:** The networked systems under consideration are typically spatially distributed over a large area (or virtually distributed). Due to the size of such networks (hundreds to millions of agents), and often to the proprietary regulations, these systems do not possess a single central coordinator or access point with the complete system information, which is thus able to solve alone the entire optimization problem. Network/data information is instead distributed among the entities comprising the network. Furthermore, there are some networks such as surveillance networks or some cyber-physical systems, where a centralized architecture is not desirable, as it makes the system prone to central entity fails and external attacks. Additional challenges are encountered from the network topology and connectivity that can be time-varying, due, e.g., to link failures, power outage, and agents' mobility. In this setting, the goal is to develop *distributed* solution methods that operate seamless *in-network*, by leveraging the network connectivity and local information (e.g., neighbor information) to cope with the lack of global knowledge on the optimization problem and offer robustness to possible failures/attacks of central units and/or to time-varying connectivity.

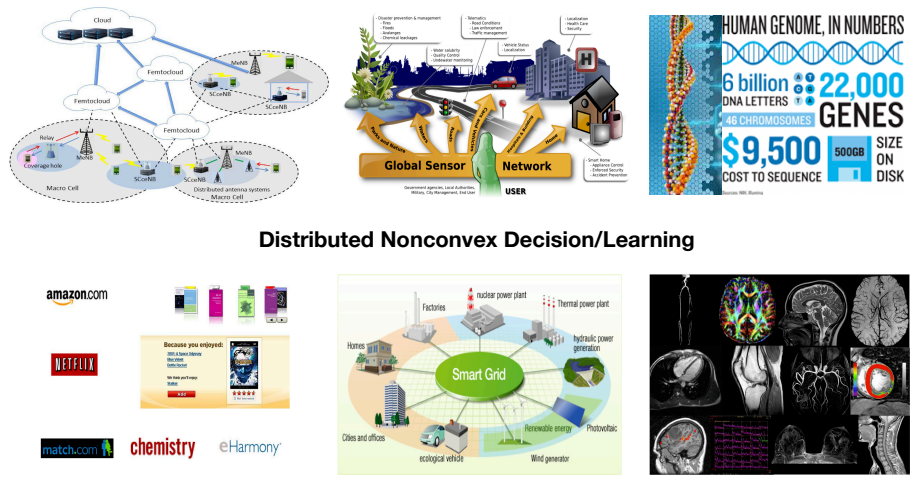

Fig. 1: *A bird's-eye view of some relevant applications generating nonconvex large-scale (networked) optimization problems.*

- **Nonconvexity:** Many formulations of interest are *nonconvex*, with nonconvex objective functions and/or constraints. Except for very special classes of nonconvex problems, whose solution can be obtained in closed form, computing the global optimal solution might be computationally prohibitive in several practical applications. This is the case, for instance, of distributed systems composed of workers with limited computational capabilities and power (e.g., motes or smart dust sensors). The desiderata is designing (parallel/distributed) solution methods that are easy to implement (in the sense that the computations performed by the workers are not expensive), with provable convergence to stationary solutions of the nonconvex problem under consideration (e.g., local optimal solutions). To this regard, a powerful and general tool is offered by the so-called *Successive Convex Approximation* (SCA) techniques: as proxy of the nonconvex problem, a sequence of "more tractable" (possibly convex) subproblems is solved, wherein the original nonconvex functions are replaced by properly chosen "simpler" surrogates. By tailoring the choice of the surrogate functions to the specific structure of the optimization problem under consideration, SCA techniques offer a lot of freedom and flexibility in the algorithmic design.

As a concrete example, consider the emerging field of *in-network big-data analytics*: the goal is to preform some, generally nonconvex, analytic tasks from a sheer volume of data, distributed over a network–see Fig. 2–examples include machine learning problems such as nonlinear least squares, dictionary learning, matrix completion, and tensor factorization, just to name a few. In these data-intensive applications, the huge volume and spatial/temporal disparity of data render centralized processing and storage a formidable task. This happens, for instance, whenever the volume of data overwhelms the storage capacity of a single computing device. Moreover, collecting sensor-network data, which are observed across a large number of spatially scattered centers/servers/agents, and routing all this local information to centralized processors, under energy, privacy constraints and/or link/hardware failures, is often infeasible or inefficient.

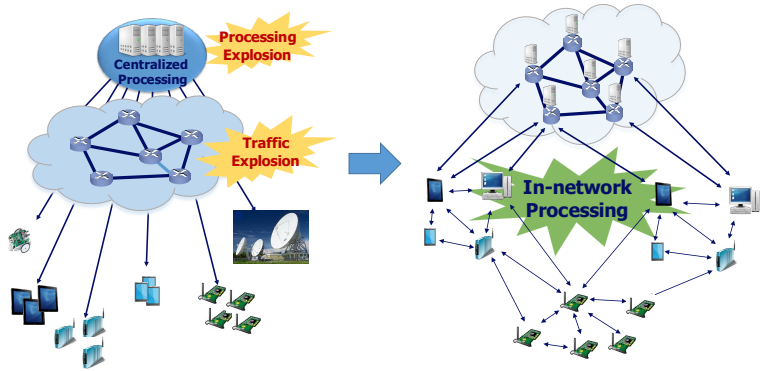

*Fig. 2: In-network big-data analytics: Traditional centralized processing and optimization are often infeasible or inefficient when dealing with large volumes of data distributed over large-scale networks. There is a necessity to develop fully decentralized algorithms that operate seamless in-network.*

The above challenges make the traditional (centralized) optimization and control techniques inapplicable, thus calling for the development of new computational models and algorithms that support efficient, *parallel* and *distributed nonconvex* optimization over networks. The major contribution of this paper is to put forth a general, unified, algorithmic framework, based on *SCA techniques*, for the parallel and distributed solution of a general class of non-convex constrained (non-separable) problems. The presented framework unifies and generalizes several existing SCA methods, making them appealing for a parallel/distributed implementation while offering a flexible selection of function approximants, step size schedules, and control of the computation/communication efficiency.

This paper is organized according to the lectures that one of the authors delivered at the CIME Summer School on *Centralized and Distributed Multi-agent Optimization Models and Algorithms* held in Cetraro, Italy, June 23–27, 2014. These lectures are:

**Lecture I–Successive Convex Approximation Methods: Basics**.

**Lecture II–Parallel Successive Convex Approximation Methods**.

**Lecture III–Distributed Successive Convex Approximation Methods**.

**Omissions:** Consistent with the main theme of the Summer School, the lectures aim at presenting SCA-based algorithms as a powerful framework for parallel and distributed, nonconvex multi-agent optimization. Of course, other algorithms have been proposed in the literature for parallel and distributed optimization. This paper does not cover schemes that are not directly related to SCA-methods or provably applicable to nonconvex problems. Examples of omissions are: primal-dual methods; augmented Lagrangian methods, including the alternating direction methods of multipliers (ADMM); and Newton methods and their inexact versions. When relevant, we provide citations of the omitted algorithms at the end of each lecture, in the section of "Source and Notes".

## Lecture I – Successive Convex Approximation Methods: Basics

This lecture overviews the majorization-minimization (MM) algorithmic framework, a particular instance of Successive Convex Approximation (SCA) Methods. The MM basic principle is introduced along with its convergence properties, which will set the ground for the design and analysis of SCA-based algorithms in the subsequent lectures. Several examples and applications are also discussed.

Consider the following general class of nonconvex optimization problems

$$\underset{\mathbf{x} \in X}{\text{minimize}}\ V(\mathbf{x}), \tag{1}$$

where $X \subseteq \mathbb{R}^m$ is a nonempty closed convex set and $V : O \to \mathbb{R}$ is continuous (possibly nonconvex and nonsmooth) on $O$, an open set containing $X$. Further assumptions on $V$ are introduced as needed.

The MM method applied to Problem (1) is based on the solution of a sequence of "more tractable" subproblems whereby the objective function $V$ is replaced by a "simpler" suitably chosen surrogate function. At each iteration $k$, a subproblem is solved of the type

$$\mathbf{x}^{k+1} \in \underset{\mathbf{x} \in X}{\operatorname{argmin}}\ \widetilde{V}(\mathbf{x}\,|\,\mathbf{x}^k), \tag{2}$$

where $\widetilde{V}(\bullet\,|\,\mathbf{x}^k)$ is a surrogate function (generally dependent on the current iterate $\mathbf{x}^k$) that upperbounds $V$ globally (further assumptions on $\widetilde{V}$ are introduced as needed). The sequence of majorization-minimization steps are pictorially shown in Fig. I.1. The underlying idea of the approach is that the surrogate function $\widetilde{V}$ is chosen so that the resulting subpoblem (2) can be efficiently solved. Roughly speaking, surrogate functions enjoying the following features are desirable:

- (Strongly) Convexity: this would lead to (strongly) convex subproblems (2);
- (Additively) Block-separability in the optimization variables: this is a key enabler for parallel/distributed solution methods, which are desirable to solve large-scale problems;
- Minimizer over $X$ in closed-form: this reduces the cost per iteration of the MM algorithm.

Finding the "right" surrogate function for the problem under consideration (possibly enjoying the properties above) might not be an easy task. A major goal of this section is to put forth general construction techniques for $\widetilde{V}$ and show their application to some representative problems in signal processing, data analysis, and communications. Some instances of $\widetilde{V}$ are drawn from the literature, e.g., [228], while some others are new and introduced for the first time in this chapter. The rest of this lecture is organized as follows. After introducing in Sec. I.1 some basic results which will lay the foundations for the analysis of SCA methods in the subsequent sections, in Sec. I.2 we describe in details the MM framework along with its convergence properties; several examples of valid surrogate functions are also discussed (cf. Sec. I.2.1). When the surrogate function $\widetilde{V}$ is block separable and so are the constraints in (2), subproblems (2) can be solved leveraging parallel algorithms. For unstructured functions $V$, in general separable surrogates are difficult to be found. When dealing with large scale optimization problems, solving (2) with respect to *all* variables might not be efficient or even possible; in all these cases, parallel block schemes are mandatory. This motivates the study of so-called "block MM" algorithms only some blocks of the variables are selected and optimized at a time. Sec. I.3 is devoted to the study of such algorithms. In Sec. I.4 we will present several applications of MM methods to problems in signal processing, machine learning, and communications. Finally, in Sec. I.5 we overview the main literature and highlight some extensions and generalizations of the methods described in this lecture.

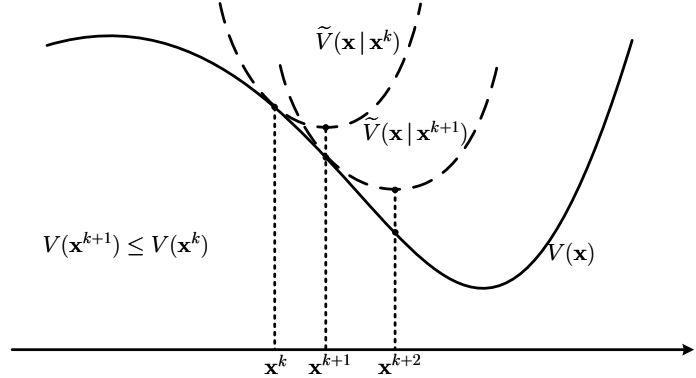

Fig. I.1: Pictorial description of the MM procedure.

## I.1. Preliminaries

We introduce here some preliminary basic results which will be extensively used through the whole paper.

We begin with the definition of directional derivative of a function and some basic properties of directional derivatives.

**Definition I.1 (directional derivative).** *A function $f : \mathbb{R}^m \to (-\infty, \infty]$ is directionally differentiable at $\mathbf{x} \in \text{dom} f \triangleq \{\mathbf{x} \in \mathbb{R}^m : f(\mathbf{x}) < \infty\}$ along a direction $\mathbf{d} \in \mathbb{R}^m$ if the following limit*

$$f'(\mathbf{x};\mathbf{d}) \triangleq \lim_{\lambda \downarrow 0} \frac{f(\mathbf{x}+\lambda \mathbf{d}) - f(\mathbf{x})}{\lambda} \tag{3}$$

*exists; this limit $f'(\mathbf{x};\mathbf{d})$ is called the directional derivative of $f$ at $\mathbf{x}$ along $\mathbf{d}$. If $f$ is directionally differentiable at $\mathbf{x}$ along all directions, then we say that $f$ is directionally differentiable at $\mathbf{x}$.* □

If $f$ is differentiable at $\mathbf{x}$, then $f'(\mathbf{x};\mathbf{d})$ reads: $f'(\mathbf{x};\mathbf{d}) = \nabla f(\mathbf{x})^T \mathbf{d}$, where $\nabla f(\mathbf{x})$ is the gradient of $f$ at $\mathbf{x}$. Some examples of directional derivatives of some structured functions (including convex functions) are discussed next.

- **Case study 1: Convex functions.** Throughout this example, we assume that $f : \mathbb{R}^m \to (-\infty, \infty]$ is a convex, closed, proper function; and $\text{int}(\text{dom} f) \neq \emptyset$ (otherwise, one can work with the relative interior of $\text{dom} f$), with $\text{int}(\text{dom} f)$ denoting the interior of $\text{dom} f$.

We show next that if $\mathbf{x} \in \text{dom} f$, $f'(\mathbf{x};\mathbf{d})$ is well defined, taking values in $[-\infty, +\infty]$. In particular, if $\mathbf{x} \in \text{dom} f$ can be approached by the direction $\mathbf{d} \in \mathbb{R}^m$, then $f'(\mathbf{x};\mathbf{d})$ is finite. For $\mathbf{x} \in \text{dom} f$, $\mathbf{d} \in \mathbb{R}^m$ and nonzero $\lambda \in \mathbb{R}$, define

$$\lambda \mapsto g_\lambda(\mathbf{x};\mathbf{d}) \triangleq \frac{f(\mathbf{x}+\lambda\mathbf{d}) - f(\mathbf{x})}{\lambda}.$$

A simple argument by convexity (increasing slopes) shows that $g(\mathbf{d};\lambda)$ is increasing in $\lambda$. Therefore, the limit in (3) exists in $[-\infty, \infty]$ and can be replaced by

$$f'(\mathbf{x};\mathbf{d}) = \inf_{\lambda > 0} \frac{1}{\lambda} \left[ f(\mathbf{x}+\lambda\mathbf{d}) - f(\mathbf{x}) \right].$$

Moreover, for $0 < \lambda \leq \beta \in \mathbb{R}$, it holds

$$g_{-\beta}(\mathbf{x};\mathbf{d}) \leq g_{-\lambda}(\mathbf{x};\mathbf{d}) \leq g_\lambda(\mathbf{x};\mathbf{d}) \leq g_\beta(\mathbf{x};\mathbf{d}).$$

If $\mathbf{x} \in \text{int}(\text{dom} f)$, both $g_{-\beta}(\mathbf{x};\mathbf{d})$ and $g_\beta(\mathbf{x};\mathbf{d})$ are finite, for sufficiently small $\beta > 0$; therefore, we have

$$-\infty < g_{-\beta}(\mathbf{x};\mathbf{d}) \leq f'(\mathbf{x};\mathbf{d}) = \inf_{\lambda > 0} g_\lambda(\mathbf{d};\mathbf{x}) \leq g_\beta(\mathbf{x};\mathbf{d}) < +\infty.$$

Finally, since $f$ is convex, it is locally Lipschitz continuous: for sufficiently small $\beta > 0$, there exists some finite $L > 0$ such that $g_\beta(\mathbf{x};\mathbf{d}) \leq L\|\mathbf{d}\|$ and $g_{-\beta}(\mathbf{x};\mathbf{d}) \geq -L\|\mathbf{d}\|$. We have proved the following result.

**Proposition I.2.** *For convex functions $f : \mathbb{R}^m \to (-\infty, \infty]$, at any $\mathbf{x} \in \text{dom} f$ and for any $\mathbf{d} \in \mathbb{R}^m$, the directional derivative $f'(\mathbf{x};\mathbf{d})$ exists in $[-\infty, +\infty]$ and it is given by*

$$f'(\mathbf{x};\mathbf{d}) = \inf_{\lambda > 0} \frac{1}{\lambda} \left[ f(\mathbf{x}+\lambda\mathbf{d}) - f(\mathbf{x}) \right].$$

*If $\mathbf{x} \in \text{int}(\text{dom} f)$, there exists a finite constant $L > 0$ such that $|f'(\mathbf{x};\mathbf{d})| \leq L\|\mathbf{d}\|$, for all $\mathbf{d} \in \mathbb{R}^m$.* □

**Directional derivative and subgradients.** The directional derivative of a convex function can be also written in terms of its subgradients, as outlined next. We first introduce the definition of subgradient along with some of its properties.

**Definition I.3 (subgradient).** *A vector $\boldsymbol{\xi} \in \mathbb{R}^m$ is a subgradient of $f$ at a point $\mathbf{x} \in \text{dom} f$ if*

$$f(\mathbf{x}+\mathbf{d}) \geq f(\mathbf{x}) + \boldsymbol{\xi}^T \mathbf{d}, \ \forall \mathbf{d} \in \mathbb{R}^m. \tag{4}$$

The subgradient set (a.k.a. subdifferential) of $f$ at $\mathbf{x} \in \mathrm{dom} f$ is defined as

$$\partial f(\mathbf{x}) \triangleq \left\{ \boldsymbol{\xi} \in \mathbb{R}^m : f(\mathbf{x}+\mathbf{d}) \geq f(\mathbf{x}) + \boldsymbol{\xi}^T \mathbf{d}, \ \forall \mathbf{d} \in \mathbb{R}^m \right\}. \tag{5}$$

Partitioning $\mathbf{x}$ in blocks, $\mathbf{x} = (\mathbf{x}_i)_{i=1}^n$, with $\mathbf{x}_i \in \mathbb{R}^{m_i}$ and $\sum_{i=1}^n m_i = m$, similarly to (5), we can define the *block-subdifferential* with respect to each $\mathbf{x}_i$, as given below, where $(\mathbf{x})_i \triangleq (\mathbf{0}^T, \ldots, \mathbf{x}_i^T, \ldots, \mathbf{0}^T)^T \in \mathbb{R}^m$.

**Definition I.4 (block-subgradient).** *The subgradient set $\partial_i f(\mathbf{x})$ of $f$ at $\mathbf{x} = (\mathbf{x}_i)_{i=1}^n \in \mathrm{dom} f$ with respect to $\mathbf{x}_i$ is defined as*

$$\partial_i f(\mathbf{x}) \triangleq \left\{ \boldsymbol{\xi}_i \in \mathbb{R}^{m_i} : f(\mathbf{x}+(\mathbf{d})_i) \geq f(\mathbf{x}) + \boldsymbol{\xi}_i^T \mathbf{d}_i, \ \forall \mathbf{d}_i \in \mathbb{R}^{m_i} \right\}. \tag{6}$$

Intuitively, when a function $f$ is convex, the subgradient generalizes the derivative of $f$; in fact, $f(\mathbf{x}) + \boldsymbol{\xi}^T \mathbf{d}$ is a global linear underestimator of $f$ at $\mathbf{x}$. Since a convex function has global linear underestimators of itself, the subgradient set $\partial f(\mathbf{x})$ should be non-empty and consist of supporting hyperplanes to the epigraph of $f$. This is formally stated in the next result (see, e.g., [24, 102] for the proof).

**Theorem I.5.** *Let $\mathbf{x} \in \mathrm{int}(\mathrm{dom} f)$. Then, $\partial f(\mathbf{x})$ is nonempty, compact, and convex.*

Note that, in the above theorem, we cannot relax the assumption $\mathbf{x} \in \mathrm{int}(\mathrm{dom} f)$ with $\mathbf{x} \in \mathrm{dom} f$. For instance, consider the function $f(x) = -\sqrt{x}$, with $\mathrm{dom} f = [0, \infty)$. We have $\partial f(0) = \emptyset$.

The subgradient definition describes a global properties of the function whereas the (directional) derivative is a local property. The connection between a directional derivative and the subdifferential of a convex function is contained in the next two results, whose proof can be found in [24, Ch.3].

**Lemma I.6.** *The subgradient set (5) at $\mathbf{x} \in \mathrm{dom} f$ can be equivalently written as*

$$\partial f(\mathbf{x}) \triangleq \left\{ \boldsymbol{\xi} \in \mathbb{R}^m : f'(\mathbf{x};\mathbf{d}) \geq \boldsymbol{\xi}^T \mathbf{d}, \quad \forall \mathbf{d} \in \mathbb{R}^m \right\}. \tag{7}$$

Note that, since $f'(\mathbf{x};\mathbf{d})$ is finite for all $\mathbf{d} \in \mathbb{R}^m$ (cf. Proposition I.2), the above representation readily shows that $\partial f(\mathbf{x})$, $\mathbf{x} \in \mathrm{int}(\mathrm{dom} f)$, is a compact set (as proved already in Theorem I.5). Furthermore, $\boldsymbol{\xi} \in \partial f(\mathbf{x})$ satisfies

$$\|\boldsymbol{\xi}\|_2 = \sup_{\mathbf{d}: \|\mathbf{d}\|_2 \leq 1} \boldsymbol{\xi}^T \mathbf{d} \leq \sup_{\mathbf{d}: \|\mathbf{d}\|_2 \leq 1} f'(\mathbf{x};\mathbf{d}) < \infty.$$

Lemma I.6 above showed how to identify subgradients from directional derivative. Lemma I.7 below shows how to move in the reverse direction.

**Lemma I.7 (max formula).** *At any $\mathbf{x} \in \mathrm{int}(\mathrm{dom} f)$ and all $\mathbf{d} \in \mathbb{R}^m$, it holds*

$$f'(\mathbf{x};\mathbf{d}) = \sup_{\boldsymbol{\xi} \in \partial f(\mathbf{x})} \boldsymbol{\xi}^T \mathbf{d}. \tag{8}$$

Lastly, we recall a straightforward result, stating that the subgradient is simply the gradient of differentiable convex functions. This is a direct consequence of Lemma I.6. Indeed, if $f$ is differentiable at $\mathbf{x}$, we can write [cf. (7)]

$$\boldsymbol{\xi}^T \mathbf{d} \leq f'(\mathbf{x};\mathbf{d}) = \nabla f(\mathbf{x})^T \mathbf{d}, \qquad \forall \boldsymbol{\xi} \in \partial f(\mathbf{x}).$$

Since the above inequality holds for all $\mathbf{d} \in \mathbb{R}^m$, we also have $\boldsymbol{\xi}^T(-\mathbf{d}) \leq f'(\mathbf{x};-\mathbf{d}) = \nabla f(\mathbf{x})^T(-\mathbf{d})$, and thus $\boldsymbol{\xi}^T \mathbf{d} = \nabla f(\mathbf{x})^T \mathbf{d}$, for all $\mathbf{d} \in \mathbb{R}^m$. This proves $\partial f(\mathbf{x}) = \{\nabla f(\mathbf{x})\}$.

The subgradient is also intimately related to optimality conditions for convex minimization. We discuss this relationship in the next subsection. We conclude this brief review with some basic examples of calculus of subgradient.

*Examples of subgradients.* As the first example, consider

$$f(x) = |x|.$$

It is not difficult to check that

$$\partial |x| = \begin{cases} \text{sign}(x), & \text{if } x \neq 0; \\ [-1, 1] & \text{if } x = 0; \end{cases} \qquad (9)$$

where $\text{sign}(x) = 1$, if $x > 0$; $\text{sign}(x) = 0$, if $x = 0$; and $\text{sign}(x) = -1$, if $x < 0$.

Similarly, consider the $\ell_1$ norm function, $f(\mathbf{x}) = \|\mathbf{x}\|_1$. We have

$$\begin{aligned} \partial \|\mathbf{x}\|_1 &= \sum_{i=1}^m \partial |x_i| = \sum_{i=1}^m \begin{cases} \mathbf{e}_i \cdot \text{sign}(x_i), & \text{if } x_i \neq 0; \\ \mathbf{e}_i \cdot [-1, 1], & \text{if } x_i = 0; \end{cases} \\ &= \sum_{x_i > 0} \mathbf{e}_i - \sum_{x_i < 0} \mathbf{e}_i + \sum_{x_i = 0} [-\mathbf{e}_i, \mathbf{e}_i], \end{aligned} \qquad (10)$$

where $\mathbf{e}_i$ denotes the $i$-th standard basis vector of $\mathbb{R}^m$; and the sum for $x_i = 0$ is the Minkowski sum. Therefore,

$$\partial \|\mathbf{0}\|_1 = \sum_{i=1}^m [-\mathbf{e}_i, \mathbf{e}_i] = \{\mathbf{x} \in \mathbb{R}^m : \|\mathbf{x}\|_\infty \leq 1\}.$$

A more complex example is given by considering any norm function $\|\bullet\|$. Introducing the dual norm

$$\|\mathbf{x}\|_* \triangleq \sup_{\mathbf{y}:\|\mathbf{y}\|\leq 1} \mathbf{x}^T \mathbf{y},$$

one can show that

$$\partial \|\mathbf{x}\| = \left\{ \boldsymbol{\xi} \in \mathbb{R}^m : \|\boldsymbol{\xi}\|_* \leq 1, \boldsymbol{\xi}^T \mathbf{x} = \|\mathbf{x}\| \right\}. \qquad (11)$$

As a concrete example, consider the $\ell_2$ norm, $f(\mathbf{x}) = \|\mathbf{x}\|_2$. Observing that

$$\|\mathbf{x}\|_2 = \sup_{\|\mathbf{y}\|_2 \leq 1} \mathbf{x}^T \mathbf{y},$$

a direct application of (11) yields

$$\partial \|\mathbf{x}\|_2 = \{\boldsymbol{\xi} \in \mathbb{R}^m : \|\boldsymbol{\xi}\|_2 \leq 1,\ \boldsymbol{\xi}^T \mathbf{x} = \|\mathbf{x}\|_2\}$$
$$= \begin{cases} \dfrac{\mathbf{x}}{\|\mathbf{x}\|_2}, & \text{if } \mathbf{x} \neq \mathbf{0}; \\ \{\boldsymbol{\xi} \in \mathbb{R}^m : \|\boldsymbol{\xi}\|_2 \leq 1\}, & \text{if } \mathbf{x} = \mathbf{0}. \end{cases}$$

• **Case study 2: Pointwise Max of functions.** Consider the pointwise maximum of (possibly) nonconvex functions

$$g(\mathbf{x}) \triangleq \max_{i=1,\ldots,I} f_i(\mathbf{x}), \tag{12}$$

where each $f_i : \mathbb{R}^m \to (-\infty, +\infty]$ is assumed to be directionally differentiable at a given $\bar{\mathbf{x}}$ along the direction $\mathbf{d}$ (with finite directional derivative). For notational simplicity, we assume that all $f_i$ have the same effective domain. The following lemma shows that $g(\mathbf{x})$ is directional differentiable at $\bar{\mathbf{x}}$ along $\mathbf{d}$ and provides an explicit expression for $g'(\bar{\mathbf{x}}; \mathbf{d})$.

**Lemma I.8.** *In the above setting, the function g defined in* (12) *is directionally differentiable at $\bar{\mathbf{x}}$ along $\mathbf{d}$, with*

$$g'(\bar{\mathbf{x}}; \mathbf{d}) = \max_{i \in A(\bar{\mathbf{x}})} f_i'(\bar{\mathbf{x}}; \mathbf{d}), \tag{13}$$

*where* $A(\mathbf{x}) \triangleq \{i = 1, \ldots, I : f_i(\mathbf{x}) = g(\mathbf{x})\}$.

*Proof.* The proof follows similar steps of that of Danskin's theorem [57]. Let $\{t^k\}$ be a sequence of positive numbers $t^k$ such that $t^k \to 0$ as $k \to \infty$. Define $\mathbf{x}^k = \bar{\mathbf{x}} + t^k \mathbf{d}$; and let $i^k$ and $\bar{i}$ be two indices in $A(\mathbf{x}^k)$ and $A(\bar{\mathbf{x}})$, respectively. We prove (13) by showing that

$$\limsup_{k \to \infty} \frac{g(\mathbf{x}^k) - g(\bar{\mathbf{x}})}{t^k} \leq \max_{i \in A(\bar{\mathbf{x}})} f_i'(\bar{\mathbf{x}}; \mathbf{d}) \leq \liminf_{k \to \infty} \frac{g(\mathbf{x}^k) - g(\bar{\mathbf{x}})}{t^k}. \tag{14}$$

We prove the right inequality first. We have

$$\begin{aligned} \liminf_{k \to \infty} \frac{g(\mathbf{x}^k) - g(\bar{\mathbf{x}})}{t^k} &= \liminf_{k \to \infty} \frac{f_{i^k}(\mathbf{x}^k) - f_{\bar{i}}(\bar{\mathbf{x}})}{t^k} \\ &= \liminf_{k \to \infty} \frac{f_{i^k}(\mathbf{x}^k) - f_{\bar{i}}(\mathbf{x}^k) + f_{\bar{i}}(\mathbf{x}^k) - f_{\bar{i}}(\bar{\mathbf{x}})}{t^k} \\ &\stackrel{(a)}{\geq} \liminf_{k \to \infty} \frac{f_{\bar{i}}(\mathbf{x}^k) - f_{\bar{i}}(\bar{\mathbf{x}})}{t^k} \stackrel{(b)}{=} f_{\bar{i}}'(\bar{\mathbf{x}}; \mathbf{d}), \end{aligned} \tag{15}$$

where (a) follows from $f_{i^k}(\mathbf{x}^k) - f_{\bar{i}}(\mathbf{x}^k) \geq 0$; and in (b) we used the fact that each $f_i$ is directionally differentiable at $\bar{x}$ along $\mathbf{d}$. Since $\bar{i}$ is any arbitrary index in $A(\bar{\mathbf{x}})$, we have

$$\liminf_{k \to \infty} \frac{g(\mathbf{x}^k) - g(\bar{\mathbf{x}})}{t^k} \geq \max_{i \in A(\bar{\mathbf{x}})} f_i'(\bar{\mathbf{x}}; \mathbf{d}). \tag{16}$$

We prove next the inequality on the left of (14). Following similar steps, we have

$$\begin{aligned}
\limsup_{k \to \infty} \frac{g(\mathbf{x}^k) - g(\bar{\mathbf{x}})}{t^k} &= \limsup_{k \to \infty} \frac{f_{i^k}(\mathbf{x}^k) - f_{\bar{i}}(\bar{\mathbf{x}})}{t^k} \\
&= \limsup_{k \to \infty} \frac{f_{i^k}(\mathbf{x}^k) - f_{i^k}(\bar{\mathbf{x}}) + f_{i^k}(\bar{\mathbf{x}}) - f_{\bar{i}}(\bar{\mathbf{x}})}{t^k} \\
&\stackrel{(a)}{\leq} \limsup_{k \to \infty} \frac{f_{i^k}(\mathbf{x}^k) - f_{i^k}(\bar{\mathbf{x}})}{t^k} \\
&\stackrel{(b)}{\leq} \max_{i \in A(\bar{\mathbf{x}})} f_i'(\bar{\mathbf{x}}; \mathbf{d}),
\end{aligned} \tag{17}$$

where (a) comes from $f_{i^k}(\bar{\mathbf{x}}) - f_{\bar{i}}(\bar{\mathbf{x}}) \leq 0$; and (b) is due to $\lim_{k \to \infty} \min_{j \in A(\bar{\mathbf{x}})} |i^k - j| = 0$, which is a consequence of $\lim_{k \to \infty} \|\mathbf{x}^k - \bar{\mathbf{x}}\| = 0$, as showed next.

Suppose that the above statement is not true. Then, there exits a subsequence of $i^k$–say $i^{k_\ell}$, with $i^{k_\ell} \in A(\mathbf{x}^{k_\ell})$–such that $\lim_{\ell \to \infty} i^{k_\ell} = i^\infty \notin A(\bar{\mathbf{x}})$ (note that $A(\bar{\mathbf{x}})$ has a finite cardinality). Therefore, for sufficiently large $\ell$, it holds $g(\mathbf{x}^{k_\ell}) = f_{i^{k_\ell}}(\mathbf{x}^{k_\ell}) = f_{i^\infty}(\mathbf{x}^{k_\ell})$. Letting $\ell \to +\infty$ and invoking continuity of $g$ (it is the point-wise maximum of finitely many continuous functions), we get $g(\bar{\mathbf{x}}) = f_{i^\infty}(\bar{\mathbf{x}})$, which is in contradiction with $i^\infty \notin A(\bar{\mathbf{x}})$. □

**Optimality conditions.** As a non-convex optimization problem, globally optimal solutions of Problem (1) are in general not possible to be computed. Thus, one has to settle for computing a "stationary" solution in practice. Even in this case, there are many kinds of stationary solutions for Problem (1). Ideally, one would like to identify a stationary solution of the sharpest kind. Arguably, for the convex constrained nonconvex program (1), a d(irectional)-stationary solution defined in terms of the directional derivatives of the objective function would qualify for this purpose. For the sake of semplicity, we will make the following blanket assumptions on Problem (1): i) $V$ is directionally differentiable on $X$; and ii) $X$ is closed and convex. We introduce next two concepts of stationarity, namely: d-stationarity and *coordinate-wise* d-stationarity.

**Definition I.9 (d-stationarity).** *Given Problem* (1) *in the above setting,* $\mathbf{x}^* \in X$ *is a d-stationary solution of* (1) *if*

$$V'(\mathbf{x}^*; \mathbf{y} - \mathbf{x}^*) \geq 0, \quad \forall \mathbf{y} \in X. \tag{18}$$

Two remarks are in order. When $V$ is convex, it follows from Lemma I.7 that $\mathbf{x}^*$ is a d-stationary (and thus a global optimal) solution of Problem (1) if there exists a $\xi \in \partial f(\mathbf{x}^*)$ such that $\xi^T(\mathbf{y} - \mathbf{x}^*) \geq 0$, $\forall \mathbf{y} \in X$. Furthermore, if $V$ is differentiable, since $V'(\mathbf{x}; \mathbf{d}) = \nabla V(\mathbf{x})^T \mathbf{d}$, (18) reads $(\mathbf{y} - \mathbf{x}^*)^T \nabla V(\mathbf{x}^*) \geq 0$, $\forall \mathbf{y} \in X$.

**Definition I.10 (coordinate-wise d-stationary).** *Given Problem* (1)*, with* $X = X_1 \times \cdots \times X_n$*,* $X_i \subseteq \mathbb{R}^{m_i}$*, and* $\sum_{i=1}^n m_i = m$*,* $\mathbf{x}^*$ *is a coordinate-wise d-stationary solution of Problem* (1) *if* $V'(\mathbf{x}^*; (\mathbf{y} - \mathbf{x}^*)_i) \geq 0$*,* $\forall \mathbf{y} \in X$ *and all* $i = 1, \ldots, n$.

In words, a coordinate-wise stationary solution is a point for which $\mathbf{x}^*$ is stationary w.r.t. every block of variables. Coordinate-wise stationarity is a weaker form of stationarity. It is the standard property of a limit point of a convergent coordinate-wise scheme (see, e.g., [250]). It is clear that a stationary point is always a coordinate-wise stationary point; the converse however is not always true, unless extra conditions on $V$ are satisfied.

**Definition I.11 (regularity).** *Problem* (1) *is regular at a coordinate-wise d-stationary solution* $\mathbf{x}^* \in X = X_1 \times \cdots \times X_n$*, if* $\mathbf{x}^*$ *is also a d-stationary point of the problem.*

The regularity condition is readily satisfied in the following two simple cases:

(a) $V$ is additively separable (possibly nonsmooth), i.e., $V(\mathbf{x}) = \sum_{i=1}^n V_i(\mathbf{x}_i)$;
(b) $V$ is differentiable around $\mathbf{x}^*$.

Note that (a) is due to the separability of the directional derivative, that is, $V'(\mathbf{x}; \mathbf{d}) = \sum_{i=1}^n V_i'(\mathbf{x}_i; \mathbf{d}_i)$; and so does (b).

Of course the two cases above are not at all inclusive of situations for which regularity holds. As an example of a nonseparable function for which regularity holds at a point at which is not continuously differentiable, consider the function arising in logistic regression problems $V(\mathbf{x}) = \sum_{i=1}^I \log(1 + e^{-a_i \mathbf{y}_i^T \mathbf{x}}) + c \cdot \|\mathbf{x}\|_2$, where $X = \mathbb{R}^m$, and $\mathbf{y}_i \in \mathbb{R}^m$ and $a_i \in \{-1, 1\}$ are given constants. Such a function $V$ is continuosly differentiable, and thus regular, at any stationary point but $\mathbf{x}^* \neq \mathbf{0}$. It is easy to check that $V$ is regular also at $\mathbf{x}^* = \mathbf{0}$, if $c < \log 2$.

Finally, an example of a nonsmooth, nonseparable function that is not regular is $V(\mathbf{x}) = \|\mathbf{A}\mathbf{x}\|_1$, with $\mathbf{A} = [3\ 4; 2\ 1]$ and $X = \mathbb{R}^2$. Point $\mathbf{x}^* = [-4\ 3]^T$ is a coordinate-wise d-stationary point, but not d-stationary [cf. Fig. I.2].

## I.2. The Majorization-Minimization (MM) Algorithm

We study Problem (1) under the following blanket assumptions.

**Assumption I.12.** *Given Problem* (1)*, we assume that:*

1. *$X \neq \emptyset$ is a closed and convex set in $\mathbb{R}^m$;*
2. *$V : O \to \mathbb{R}$ is continuous on the open set $O \supseteq X$;*
3. *$V'(\mathbf{x}; \mathbf{d})$ exists at any $\mathbf{x} \in X$ and for all feasible directions $\mathbf{d} \in \mathbb{R}^m$ at $\mathbf{x}$;*
4. *$V$ is bounded from below.*

Note that the above assumptions are quite standard and are satisfied by most of the problems of practical interest; see Sec. I.4 for some illustrative examples.

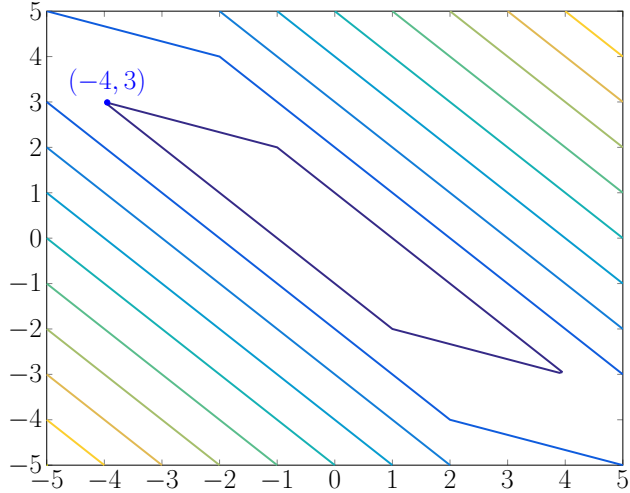

Fig. I.2: Contour of $V(\mathbf{x}) = \|\mathbf{A}\mathbf{x}\|_1$, with $\mathbf{A} = [3\ 4; 2\ 1]$. The function is not regular at $\mathbf{x}^* = [-4,3]^T$ [193].

As already anticipated at the beginning of this section, the idea of the MM procedure is to approximate, at each iteration, the objective function $V$ in (1) by a "simpler" properly chosen surrogate function $\widetilde{V}(\bullet|\mathbf{x}^k)$ and solve the resulting optimization problem (2). Convergence of this iterative method is guaranteed if the following conditions are satisfied in the choice of $\widetilde{V}$. In what follows, we denote by $\widetilde{V}'(\mathbf{y};\mathbf{d}|\mathbf{x})$ the directional derivative of $\widetilde{V}(\bullet|\mathbf{x})$ at $\mathbf{y}$ along the direction $\mathbf{d}$.

**Assumption I.13.** *The surrogate function $\widetilde{V}: O \times O \to \mathbb{R}$ satisfies the following conditions:*

1. *$\widetilde{V}(\bullet|\bullet)$ is continuous on $X \times X$;*
2. *$\mathbf{y} \in \mathrm{argmin}_{\mathbf{x} \in X}\ \widetilde{V}(\mathbf{x}|\mathbf{y}) - V(\mathbf{x})$;*
3. *The directional derivative of $\widetilde{V}$ satisfies $\widetilde{V}'(\mathbf{x};\mathbf{d}|\mathbf{x}) = V'(\mathbf{x};\mathbf{d})$, for all $\mathbf{x} \in X$ and feasible directions $\mathbf{d} \in \mathbb{R}^m$ at $\mathbf{x}$.*

Assumption I.13.2 states that, at any feasible $\mathbf{y}$, $\widetilde{V}(\bullet|\mathbf{y})$ upperbounds $V(\bullet)$ on $X$, in the following sense:
$$\widetilde{V}(\mathbf{x}|\mathbf{y}) \geq V(\mathbf{x}) + c_y, \quad \forall \mathbf{x} \in X, \tag{19}$$
with $c_y \triangleq \widetilde{V}(\mathbf{y}|\mathbf{y}) - V(\mathbf{y})$, where the equality is achieved when $\mathbf{x} = \mathbf{y}$. Assumption I.13.3 is a derivative consistency condition: roughly speaking, it ensures that $\widetilde{V}(\bullet|\mathbf{x})$ has the same first order properties of $V(\bullet)$ at $\mathbf{x} \in X$.

The MM algorithm is summarized in Algorithm 1 and its convergence is stated in Theorem I.14.

**Algorithm 1: The Majorization-Minimization (MM) Algorithm**

**Data** : $\mathbf{x}^0 \in X$. Set $k = 0$.
(S.1) : If $\mathbf{x}^k$ satisfies a termination criterion: STOP;
(S.2) : Update $\mathbf{x}$ as
$$\mathbf{x}^{k+1} \in \underset{\mathbf{x} \in X}{\operatorname{argmin}} \, \widetilde{V}(\mathbf{x} \mid \mathbf{x}^k); \tag{20}$$
(S.3) : $k \leftarrow k+1$, and go to (S.1).

**Theorem I.14.** *Let $\{\mathbf{x}^k\}_{k \in \mathbb{N}_+}$ be the sequence generated by Algorithm 1 under Assumptions I.12 and I.13. Then, every limit point of $\{\mathbf{x}^k\}_{k \in \mathbb{N}_+}$ (if exists) is a d-stationary solution of Problem* (1).

*Proof.* The main properties of Algorithm 1 is to generate a nonincreasing sequence $\{V(\mathbf{x}^k)\}$. Indeed, we have

$$V(\mathbf{x}^{k+1}) \overset{(19)}{\leq} \widetilde{V}(\mathbf{x}^{k+1} \mid \mathbf{x}^k) - c^k \overset{(20)}{\leq} \widetilde{V}(\mathbf{x}^k \mid \mathbf{x}^k) - c^k \overset{(a)}{=} V(\mathbf{x}^k), \tag{21}$$

where $c^k \triangleq \widetilde{V}(\mathbf{x}^k \mid \mathbf{x}^k) - V(\mathbf{x}^k)$, and (a) follows from the fact that the inequality (19) is achieved with equality at $\mathbf{x} = \mathbf{x}^k$.

Let $\mathbf{x}^*$ be a limit point of $\{\mathbf{x}^k\}_{k \in \mathbb{N}_+}$, that is, $\lim_{t \to \infty} \mathbf{x}^{k_t} = \mathbf{x}^* \in X$. We have

$$\widetilde{V}(\mathbf{x}^{k_{t+1}} \mid \mathbf{x}^{k_{t+1}}) - c^{k_{t+1}} = V(\mathbf{x}^{k_{t+1}}) \overset{(21)}{\leq} V(\mathbf{x}^{k_t+1}) \leq \widetilde{V}(\mathbf{x}^{k_t+1} \mid \mathbf{x}^{k_t}) - c^{k_t} \leq \widetilde{V}(\mathbf{x} \mid \mathbf{x}^{k_t}) - c^{k_t},$$

for all $\mathbf{x} \in X$. Let $t \to +\infty$; invoking the continuity of $\widetilde{V}(\bullet \mid \bullet)$ and $V(\bullet)$, we have that the sequence $\{c^{k_t}\}_{t \in \mathbb{N}_+}$ converges (to a finite value). Therefore,

$$\widetilde{V}(\mathbf{x}^* \mid \mathbf{x}^*) \leq \widetilde{V}(\mathbf{x} \mid \mathbf{x}^*), \quad \forall \mathbf{x} \in X, \tag{22}$$

which implies

$$0 \leq \widetilde{V}'(\mathbf{x}^*; \mathbf{d} \mid \mathbf{x}^*) \overset{(a)}{=} V'(\mathbf{x}^*; \mathbf{d}), \quad \forall \mathbf{d} \in \mathbb{R}^m \quad \text{such that} \quad \mathbf{x}^* + \mathbf{d} \in X, \tag{23}$$

where (a) follows from Assumption I.13.3. This shows that $\mathbf{x}^*$ is a d-stationary solution of Problem (1). □

Note that, since the sequence $\{V(\mathbf{x}^k)\}$ is nonincreasing, a sufficient condition for $\{\mathbf{x}^k\}$ to admit a limit point is that the set $\{\mathbf{x} \in X : V(\mathbf{x}) \leq V(\mathbf{x}^0)\}$ is compact. A sufficient condition for that is the coercivity of $V$ on $X$.

*On the termination criterion.* We briefly discuss how to choose the termination criterion in Step 2 of Algorithm 1; we refer the interested reader to [241] for more details. Let $M : X \to X$ be a map such that $M(\mathbf{x}^k) \in \operatorname{argmin}_{\mathbf{x} \in X} \widetilde{V}(\mathbf{x} \mid \mathbf{x}^k)$ and $\mathbf{x}^{k+1} = M(\mathbf{x}^k)$ [cf. (2)]. In words, among all the global minimizers of $\widetilde{V}(\bullet \mid \mathbf{x}^k)$ on $X$, $M$ uniquely selects the one, $\mathbf{x}^{k+1}$, used in Step 2 of the MM algorithm. We shown next that $V(M(\mathbf{x})) = V(\mathbf{x})$ is a sufficient condition of $\mathbf{x}$ being a d-stationary solution of Problem (1). It follows from (21) that $V(M(\mathbf{x})) = V(\mathbf{x})$ forces $\widetilde{V}(M(\mathbf{x})|\mathbf{x}) = \widetilde{V}(\mathbf{x}|\mathbf{x})$, which implies that $\mathbf{x}$ is a minimizer of $\widetilde{V}(\bullet|\mathbf{x})$ on $X$, and thus (by Assumption I.13.3) a d-stationary point of (1). Based on the above observation and assuming that $M(\bullet)$ is *continuous*, the following is a valid merit function to measure distance from stationarity of the iterate $\mathbf{x}^k$:

$$J^{k+1} \triangleq \frac{V(\mathbf{x}^k) - V(\mathbf{x}^{k+1})}{\max(1, |V(\mathbf{x}^k)|)}. \tag{24}$$

The continuity assumption of $M$ maybe hard to check directly. A stronger condition implying continuity of $M$ is that the minimizer of $\widetilde{V}(\bullet|\mathbf{x}^k)$ over $X$ is unique, for all $\mathbf{x}^k \in X$ [241, Lemma 1]. Note that, in such a case, it is not difficult to check that, if $\mathbf{x}$ is a fixed-point of $M$, then it must be a d-stationary point of Problem (1). Therefore, in the aforementioned setting, an alternative merit function is

$$J^{k+1} \triangleq \|\mathbf{x}^{k+1} - \mathbf{x}^k\|. \tag{25}$$

Other termination criteria are discussed in Lecture II for SCA-based algorithms.

### I.2.1 Discussion on Algorithm 1

The following comments on Algorithm 1 are in order.

#### On the choice of surrogate function

The successful application of the MM algorithms relies on the possibility of finding a valid surrogate function $\widetilde{V}$. The critical assumption to be satisfied is undoubtedly the upperbound condition, as stated in Assumption I.13.2. We provide next some systematic rules which help to build a surrogate function that meets this condition (and the other required assumptions); several illustrating examples are also discussed. Although for specific (structured) problems it is possible to find a nonconvex surrogate function whose minimizer can be computed efficiently, a convex surrogate is in general preferred, since it leads to a convex subproblem (2). Therefore, next we mainly focus on convex surrogates.

**1) First order Taylor expansion:** Suppose $V$ is a differentiable *concave* function on $X$. A natural choice for $\widetilde{V}$ satisfying Assumption I.13 is then: given $\mathbf{y} \in X$,

$$\widetilde{V}(\mathbf{x}|\mathbf{y}) = V(\mathbf{y}) + \nabla V(\mathbf{y})^T (\mathbf{x} - \mathbf{y}). \tag{26}$$

More generally, $\widetilde{V}$ can be chosen as any convex differentiable function on $X$, say $\widetilde{V}_{\text{cvx}}$, satisfying the gradient consistency condition $\nabla \widetilde{V}_{\text{cvx}}(\mathbf{x}|\mathbf{x}) = \nabla V(\mathbf{x})$, for all $\mathbf{x} \in X$; this is enough for Assumption I.13.2 (and thus Assumption I.13) to be satisfied, as shown by the following chain of inequalities

$$V(\mathbf{x}) - V(\mathbf{y}) \overset{(a)}{\leq} \nabla V(\mathbf{y})^T (\mathbf{x} - \mathbf{y}) \overset{(b)}{=} \nabla \widetilde{V}_{\text{cvx}}(\mathbf{y}|\mathbf{y})^T (\mathbf{x} - \mathbf{y}) \overset{(c)}{\leq} \widetilde{V}_{\text{cvx}}(\mathbf{x}|\mathbf{y}) - \widetilde{V}_{\text{cvx}}(\mathbf{y}|\mathbf{y}), \tag{27}$$

where (a) follows from the concavity of $V$; in (b) we used $\nabla \widetilde{V}_{\text{cvx}}(\mathbf{y}|\mathbf{y}) = \nabla V(\mathbf{y})$; and (c) follows from the convexity of $\widetilde{V}_{\text{cvx}}$. Eq. (27) shows that (19) (thus Assumption I.13.2) is satisfied by such a $\widetilde{V}_{\text{cvx}}$.

*Example 1.* $V(x) = \log(x)$ is concave on $(0, +\infty]$. Hence, it can be majorized by $\widetilde{V}(x|y) = x/y$, for any given $y > 0$.

*Example 2.* $V(x) = |x|^p$, with $p \in (0,1)$, is concave on $(-\infty, 0)$ and $(0, +\infty)$. It thus can be majorized by the quadratic function $\widetilde{V}(x|y) = \frac{p}{2} |y|^{p-2} x^2$, for any given $y \neq 0$.

**2) Second order Taylor expansion:** Suppose $V$ is $C^1$, with $L$-Lipschitz gradient on $X$. Then $V$ can be majorized by the surrogate function: given $\mathbf{y} \in X$,

$$\widetilde{V}(\mathbf{x}|\mathbf{y}) = V(\mathbf{y}) + \nabla V(\mathbf{y})^T (\mathbf{x} - \mathbf{y}) + \frac{L}{2} \|\mathbf{x} - \mathbf{y}\|^2. \tag{28}$$

Moreover, if $V$ is twice differentiable and there exists a matrix $\mathbf{M} \in \mathbb{R}^{m \times m}$ such that $\mathbf{M} - \nabla^2 V(\mathbf{x}) \succeq \mathbf{0}$, for all $\mathbf{x} \in X$, then $V$ can be majorized by the following valid surrogate function

$$\widetilde{V}(\mathbf{x}|\mathbf{y}) = V(\mathbf{y}) + \nabla V(\mathbf{y})^T (\mathbf{x} - \mathbf{y}) + \frac{1}{2}(\mathbf{x} - \mathbf{y})^T \mathbf{M} (\mathbf{x} - \mathbf{y}). \tag{29}$$

*Example 3 (The proximal gradient algorithm).* Suppose that $V$ admits the structure $V = F + G$, where $F : \mathbb{R}^m \to \mathbb{R}$ is $C^1$, with $L$-Lipschitz gradient on $X$, and $G : \mathbb{R}^m \to \mathbb{R}$ is convex (possibly nonsmooth) on $X$. Using (28) to majorize $F$, a valid surrogate for $V$ is: given $\mathbf{y} \in X$,

$$\widetilde{V}(\mathbf{x}|\mathbf{y}) = F(\mathbf{x}^k) + \nabla F(\mathbf{y})^T (\mathbf{x} - \mathbf{y}) + \frac{L}{2} \|\mathbf{x} - \mathbf{y}\|^2 + G(\mathbf{x}). \tag{30}$$

Quite interestingly, the above choice leads to a strongly convex subproblem (2), whose minimizer has the following closed form:

$$\mathbf{x}^{k+1} = \text{prox}_{1/L, G} \left( \mathbf{x}^k - \frac{1}{L} \nabla F(\mathbf{x}^k) \right), \tag{31}$$

where $\text{prox}_{\gamma, G}(\bullet)$ is the proximal response, defined as

$$\text{prox}_{\gamma, G}(\mathbf{x}) \triangleq \underset{\mathbf{z}}{\text{argmin}} \left\{ G(\mathbf{z}) + \frac{1}{2\gamma} \|\mathbf{z} - \mathbf{x}\|^2 \right\}.$$

The resulting MM algorithm (Algorithm 1) turns out to be the renowned *proximal gradient algorithm*, with step-size $\gamma = 1/L$.

**3) Pointwise maximum:** Suppose $V : \mathbb{R}^m \to \mathbb{R}$ can be written as the pointwise maximum of functions $\{f_i\}_{i=1}^I$, i.e.,

$$V(\mathbf{x}) \triangleq \max_{i=1,\ldots,I} f_i(\mathbf{x}),$$

where each $f_i : \mathbb{R}^m \to \mathbb{R}$ satisfies Assumption I.12.2 & I.12.3. Then $V$ can be majorized by

$$\widetilde{V}(\mathbf{x}\,|\,\mathbf{y}) = \max_{i=1,\ldots,I} \widetilde{f}_i(\mathbf{x}\,|\,\mathbf{y}), \tag{32}$$

for any given $\mathbf{y} \in X$, where $\widetilde{f}_i : X \times X \to \mathbb{R}$ is a surrogate function of $f_i$ satisfying Assumption I.13 and $\widetilde{f}_i(\mathbf{y}\,|\,\mathbf{y}) = f_i(\mathbf{y})$. It is not difficult to verify that $\widetilde{V}$ above satisfies Assumption I.13. Indeed, the continuity of $\widetilde{V}$ follows from that of $\widetilde{f}_i$. Moreover, we have $\widetilde{V}(\mathbf{y}\,|\,\mathbf{y}) = V(\mathbf{y})$. Finally, condition I.13.3 is a direct consequence of Lemma I.8:

$$V'(\mathbf{x};\mathbf{d}) \stackrel{(13)}{=} \max_{i \in A(\mathbf{x})} f_i'(\mathbf{x};\mathbf{d}) = \max_{i \in A(\mathbf{x})} \widetilde{f}_i'(\mathbf{x};\mathbf{d}\,|\,\mathbf{x}) = \widetilde{V}'(\mathbf{x};\mathbf{d}\,|\,\mathbf{x}), \tag{33}$$

where $A(\mathbf{x}) = \{i : f_i(\mathbf{x}) = V(\mathbf{x})\} = \{i : \widetilde{f}_i(\mathbf{x}\,|\,\mathbf{x}) = \widetilde{V}(\mathbf{x})\}$.

**4) Composition by a convex function:** Suppose $V : \mathbb{R}^m \to \mathbb{R}$ can be expressed as $V(\mathbf{x}) \triangleq f(\sum_{i=1}^n \mathbf{A}_i \mathbf{x})$, where $f : \mathbb{R}^m \to \mathbb{R}$ is a convex function and $\mathbf{A}_i \in \mathbb{R}^{m \times m}$ are given matrices. Then, one can construct a surrogate function of $V$ leveraging the following inequality due to convexity of $f$:

$$f\left(\sum_{i=1}^n w_i \mathbf{x}_i\right) \leq \sum_{i=1}^n w_i f(\mathbf{x}_i) \tag{34}$$

for all $\sum_{i=1}^n w_i = 1$ and each $w_i > 0$. Specifically, rewrite first $V$ as: given $\mathbf{y} \in \mathbb{R}^m$,

$$V(\mathbf{x}) = f\left(\sum_{i=1}^n \mathbf{A}_i \mathbf{x}\right) = f\left(\sum_{i=1}^n w_i \left(\frac{\mathbf{A}_i(\mathbf{x}-\mathbf{y})}{w_i} + \sum_{i=1}^n \mathbf{A}_i \mathbf{y}\right)\right).$$

Then, using (34) we can upperbound $V$ as

$$V(\mathbf{x}) \leq \widetilde{V}(\mathbf{x}|\mathbf{y}) = \sum_{i=1}^n w_i f\left(\frac{\mathbf{A}_i(\mathbf{x}-\mathbf{y})}{w_i} + \sum_{i=1}^n \mathbf{A}_i \mathbf{y}\right). \tag{35}$$

It is not difficult to check that $\widetilde{V}$ satisfies Assumption I.13. Eq. (35) is particularly useful to construct surrogate functions that are additively separable in the (block) variables, which opens the way to parallel solution methods wherein (blocks of) variables are updated in parallel. Two examples are discussed next.

*Example 4.* Let $V : \mathbb{R}^m \to \mathbb{R}$ be convex and let $\mathbf{x} \in \mathbb{R}^m$ be partitioned as $\mathbf{x} = (\mathbf{x}_i)_{i=1}^n$, where $\mathbf{x}_i \in \mathbb{R}^{m_i}$ and $\sum_{i=1}^n m_i = m$. Let us rewrite $\mathbf{x}$ in terms of its block as $\mathbf{x} = \sum_{i=1}^n \mathbf{A}_i \mathbf{x}$, where $\mathbf{A}_i \in \mathbb{R}^{m \times m}$ is the block diagonal matrix such that $\mathbf{A}_i \mathbf{x} = (\mathbf{x})_i$, with $(\mathbf{x})_i \triangleq [\mathbf{0}^T, \ldots, \mathbf{0}^T, \mathbf{x}_i^T, \mathbf{0}^T, \ldots, \mathbf{0}^T]^T$ denoting the operator nulling all the blocks of $\mathbf{x}$ except the $i$-th one. Then using (35) one can choose the following surrogate function: given $\mathbf{y} = (\mathbf{y}_i)_{i=1}^n$, with each $\mathbf{y}_i \in \mathbb{R}^{m_i}$,

$$V(\mathbf{x}) = V\left(\sum_{i=1}^{n} \mathbf{A}_i \mathbf{x}\right) \leq \widetilde{V}(\mathbf{x}|\mathbf{y}) \triangleq \sum_{i=1}^{n} w_i V\left(\frac{1}{w_i}((\mathbf{x})_i - (\mathbf{y})_i) + \mathbf{y}\right).$$

It is easy to check that such a $\widetilde{V}$ is separable in the blocks $\mathbf{x}_i$'s.

*Example 5.* Let $V : \mathbb{R} \to \mathbb{R}$ be convex and let vectors $\mathbf{x}, \mathbf{a} \in \mathbb{R}^m$ be partitioned as $\mathbf{x} = (\mathbf{x}_i)_{i=1}^n$ and $\mathbf{a} = (\mathbf{a}_i)_{i=1}^n$, respectively, with $\mathbf{x}_i$ and $\mathbf{a}_i$ having the same size. Then, invoking (35), a valid surrogate function of the composite function $V(\mathbf{a}^T \mathbf{x})$ is: given $\mathbf{y} = (\mathbf{y}_i)_{i=1}^n$, partitioned according to $\mathbf{x}$,

$$\begin{aligned} V(\mathbf{a}^T \mathbf{x}) = V\left(\sum_{i=1}^{n} (\mathbf{a})_i^T \mathbf{x}\right) &\leq \widetilde{V}(\mathbf{x}|\mathbf{y}) \triangleq \sum_{i=1}^{n} w_i V\left(\frac{(\mathbf{a})_i^T(\mathbf{x}-\mathbf{y})}{w_i} + \sum_{i=1}^{n}(\mathbf{a})_i^T \mathbf{y}\right) \\ &= \sum_{i=1}^{n} w_i V\left(\frac{\mathbf{a}_i^T(\mathbf{x}_i - \mathbf{y}_i)}{w_i} + \mathbf{a}^T \mathbf{y}\right). \end{aligned} \quad (36)$$

This is another example of additively (block) separable surrogate function.

**5) Surrogates based on special inequalities:** Other techniques often used to construct valid surrogate functions leverage specific inequalities, such as the Jensen's inequality, the arithmetic-geometric-mean inequality, and the Cauchy-Schwartz inequality. The way of using these inequalities, however, depends on the specific expression of the objective function $V$ under consideration; generalizing these approaches to arbitrary $V$'s seems not possible. We provide next two illustrative (nontrivial) case studies based on the Jensen's inequality and the arithmetic-geometric-mean inequality while we refer the interested reader to [228] for more examples building on this approach.

*Example 6 (The Expectation-Maximization algorithm).* Given a pair of random (vector) variables $(\mathbf{s}, \mathbf{z})$ whose joint probability distribution $p(\mathbf{s}, \mathbf{z}|\mathbf{x})$ is parametrized by $\mathbf{x}$, we consider the maximum likelihood estimation problem of estimating $\mathbf{x}$ *only* from $\mathbf{s}$ while the random variable $\mathbf{z}$ is unobserved/hidden. The problem is formulated as

$$\hat{\mathbf{x}}_{\mathrm{ML}} = \underset{\mathbf{x}}{\operatorname{argmin}} \left\{ V(\mathbf{x}) \triangleq -\log p(\mathbf{s}|\mathbf{x}) \right\}, \quad (37)$$

where $p(\mathbf{s}|\mathbf{x})$ is the (conditional) marginal distribution of $\mathbf{s}$.

In general, the expression of $p(\mathbf{s}|\mathbf{x})$ is not available in closed form; moreover numerical evaluations of the integration of $p(\mathbf{s}, \mathbf{z}|\mathbf{x})$ with respect to $\mathbf{z}$ can be computationally too costly, especially if the dimension of $\mathbf{z}$ is large. In the following we show how to attach Problem (37) using the MM framework. We build next a valid surrogate function for $V$ leading to a simpler optimization problem to solve.

Specifically, we can rewrite $V$ as

$$\begin{aligned}
V(\mathbf{x}) &= -\log p(\mathbf{s}|\mathbf{x}) \\
&= -\log \int p(\mathbf{s}|\mathbf{z},\mathbf{x})p(\mathbf{z}|\mathbf{x})\mathrm{d}\mathbf{z} \\
&= -\log \int \left(\frac{p(\mathbf{s}|\mathbf{z},\mathbf{x})p(\mathbf{z}|\mathbf{s},\mathbf{x}^k)}{p(\mathbf{z}|\mathbf{s},\mathbf{x}^k)}\right) p(\mathbf{z}|\mathbf{x})\mathrm{d}\mathbf{z} \\
&= -\log \int \left(\frac{p(\mathbf{s}|\mathbf{z},\mathbf{x})p(\mathbf{z}|\mathbf{x})}{p(\mathbf{z}|\mathbf{s},\mathbf{x}^k)}\right) p(\mathbf{z}|\mathbf{s},\mathbf{x}^k)\mathrm{d}\mathbf{z} \\
&\stackrel{(a)}{\leq} -\int \log \left(\frac{p(\mathbf{s}|\mathbf{z},\mathbf{x})p(\mathbf{z}|\mathbf{x})}{p(\mathbf{z}|\mathbf{s},\mathbf{x}^k)}\right) p(\mathbf{z}|\mathbf{s},\mathbf{x}^k)\mathrm{d}\mathbf{z} \\
&= -\int \log(p(\mathbf{s},\mathbf{z}|\mathbf{x}))\, p(\mathbf{z}|\mathbf{s},\mathbf{x}^k)\mathrm{d}\mathbf{z} + \underbrace{\int \log\left(p(\mathbf{z}|\mathbf{s},\mathbf{x}^k)\right) p(\mathbf{z}|\mathbf{s},\mathbf{x}^k)\mathrm{d}\mathbf{z}}_{\text{constant}},
\end{aligned}$$

where (a) follows from the Jensen's inequality. This naturally suggests the following surrogate function of $V$: given $\mathbf{y}$,

$$\widetilde{V}(\mathbf{x}|\mathbf{y}) = -\int \log(p(\mathbf{s},\mathbf{z}|\mathbf{x}))\, p(\mathbf{z}|\mathbf{s},\mathbf{y})\mathrm{d}\mathbf{z}. \tag{38}$$

In fact, it is not difficult to check that such a $\widetilde{V}$ satisfies Assumption I.13.

The update of $\mathbf{x}$ resulting from the MM algorithm then reads

$$\mathbf{x}^{k+1} \in \underset{\mathbf{x}}{\operatorname{argmin}}\left\{-\int \log(p(\mathbf{s},\mathbf{z}|\mathbf{x}))\, p(\mathbf{z}|\mathbf{s},\mathbf{x}^k)\mathrm{d}\mathbf{z}\right\}. \tag{39}$$

Problem (39) can be efficiently solved for specific probabilistic models, including those belonging to the exponential family, the Gaussian/multinomial mixture model, and the linear Gaussian latent model.

Quite interestingly, the resulting MM algorithm [Algorithm 1 based on the update (39)] turns out to be the renowned *Expectation-Maximization* (EM) algorithm [66]. The EM algorithm starts from an initial estimate $\mathbf{x}^0$ and generates a sequence $\mathbf{x}^k$ by repeating the following two steps:

- *E-step*: Calculate $\widetilde{V}(\mathbf{x}|\mathbf{x}^k)$ as in (38);
- *M-step*: Update $\mathbf{x}^{k+1}$ as $\mathbf{x}^{k+1} \in \operatorname{argmax}_{\mathbf{x}} -\widetilde{V}(\mathbf{x}|\mathbf{x}^k)$.

Note that these two steps correspond exactly to the majorization and minimization steps (38) and (39), respectively, showing that the EM algorithm belongs to the family of MM schemes, based on the surrogate function (38).

*Example 7 (Geometric programming).* Consider the problem of minimizing a signomial

$$V(\mathbf{x}) = \sum_{j=1}^{J} c_j \prod_{i=1}^{n} x_i^{\alpha_{ij}}$$

on the nonnegative orthant $\mathbb{R}^n_+$, with $c_j, \alpha_{ij} \in \mathbb{R}$. In the following, we assume that $V$ is coercive on $\mathbb{R}^n_+$. A sufficient condition for this to hold is that for all $i = 1, \ldots, n$ there exists at least one $j$ such that $c_j > 0$ and $\alpha_{ij} > 0$, and at least one $j$ such that $c_j > 0$ and $\alpha_{ij} < 0$ [136].

We construct a separable surrogate function of $V$ at given $\mathbf{y} \in \mathbb{R}^n_{++}$. We first derive an upperbound for the summand in $V$ with $c_j > 0$ and a lowerbound for those with $c_j < 0$, using the arithmetic-geometric mean inequality and the concavity of the log function (cf. Example 1), respectively.

Let $z_i$ and $\alpha_i$ be nonnegative scalars, the arithmetic-geometric mean inequality reads

$$\prod_{i=1}^n z_i^{\alpha_i} \leq \sum_{i=1}^n \frac{\alpha_i}{\|\boldsymbol{\alpha}\|_1} z_i^{\|\boldsymbol{\alpha}\|_1}. \tag{40}$$

Since $y_i > 0$ for all $i = 1, \ldots, n$, let $z_i = x_i/y_i$ for $\alpha_i > 0$ and $z_i = y_i/x_i$ for $\alpha_i < 0$. Then (40) implies that the monomial $\prod_{i=1}^n x_i^{\alpha_i}$ can be upperbounded on $\mathbb{R}^n_{++}$ as

$$\prod_{i=1}^n x_i^{\alpha_i} \leq \left(\prod_{i=1}^n (y_i)^{\alpha_i}\right) \sum_{i=1}^n \frac{|\alpha_i|}{\|\boldsymbol{\alpha}\|_1} \left(\frac{x_i}{y_i}\right)^{\|\boldsymbol{\alpha}\|_1 \operatorname{sign}(\alpha_i)}. \tag{41}$$

To upperbound the terms in $V$ with negative $c_j$ on $\mathbb{R}^n_{++}$, which is equivalent to find a lowerbound of $\prod_{i=1}^n x_i^{\alpha_i}$, we use the bound introduced in Example 1, with $x = \prod_{i=1}^n x_i^{\alpha_i}$, which yields

$$\log\left(\prod_{i=1}^n x_i^{\alpha_i}\right) \leq \log\left(\prod_{i=1}^n (y_i)^{\alpha_i}\right) + \left(\prod_{i=1}^n (y_i)^{\alpha_i}\right)^{-1} \left(\prod_{i=1}^n x_i^{\alpha_i} - \prod_{i=1}^n (y_i)^{\alpha_i}\right).$$

Rearranging the terms we have

$$\prod_{i=1}^n x_i^{\alpha_i} \geq \prod_{i=1}^n (y_i)^{\alpha_i} \left(1 + \sum_{i=1}^n \alpha_i \log x_i - \sum_{i=1}^n \alpha_i \log y_i\right). \tag{42}$$

Combining (41) and (42) leads to the separable surrogate function $\widetilde{V}(\mathbf{x}|\mathbf{y}) \triangleq \sum_{i=1}^n \widetilde{V}_i(x_i|\mathbf{y})$ of $V$, with

$$\widetilde{V}_i(x_i|\mathbf{y}) \triangleq \sum_{j: c_j > 0} c_j \left(\prod_{\ell=1}^n (y_\ell)^{\alpha_{\ell j}}\right) \frac{|\alpha_{ij}|}{\|\boldsymbol{\alpha}_j\|_1} \left(\frac{x_i}{y_i}\right)^{\|\boldsymbol{\alpha}_j\|_1 \operatorname{sign}(\alpha_{ij})}$$
$$+ \sum_{j: c_j < 0} c_j \left(\prod_{\ell=1}^n (y_\ell)^{\alpha_{\ell j}}\right) \alpha_{ij} \log x_i. \tag{43}$$

Function $\widetilde{V}(\bullet|\mathbf{y})$ is coercive since it is continuous and upperbounds $V$ on $\mathbb{R}^n_{++}$, therefore $\mathbf{x}^{k+1} \in \mathbb{R}^n_{++}$ (recall that we assumed $\mathbf{x}^k \in \mathbb{R}^n_{++}$). Moreover, the following chain of inequalities holds:

$$V(\mathbf{x}^{k+1}) \leq \widetilde{V}(\mathbf{x}^{k+1}|\mathbf{x}^k) \leq \widetilde{V}(\mathbf{x}^k|\mathbf{x}^k) = V(\mathbf{x}^k) \leq \cdots \leq V(\mathbf{x}^0) < +\infty, \tag{44}$$

which implies that, given $\mathbf{x}^0 \in \mathbb{R}_{++}$, the sequence $\{\mathbf{x}^k\}$ generated by the MM algorithm always stays in the compact set $X^0 \triangleq \{\mathbf{x} \,|\, V(\mathbf{x}) \leq V(\mathbf{x}^0)\} \subseteq \mathbb{R}_{++}^n$.

As $\widetilde{V}(\mathbf{x}\,|\,\mathbf{x}^k)$ is separable, $\mathbf{x}^{k+1}$ thus can be computed in parallel. In addition, focusing on a single component $x_i$ and performing the change of variable $x_i = e^{z_i}$, it is not difficult to check that $\widetilde{V}_i$ is convex in $z_i$; its global minimizer can be then computed efficiently by solving a one dimensional convex optimization problem.

**Parallel implementation**

When dealing with large-scale optimization instances of Problem (1), solving subproblem (2) with respect to the entire $\mathbf{x}$ might be either computationally too costly or practically infeasible. When the feasible set $X$ admits a Cartesian product structure, i.e., $X = X_1 \times \cdots \times X_n$, with $X_i \subseteq \mathbb{R}^{m_i}$, a natural approach to cope with the curse of dimensionally is partitioning the vector of variables $\mathbf{x}$ into $n$ blocks according to $X$, i.e., $\mathbf{x} = \left[\mathbf{x}_1^T, \ldots, \mathbf{x}_n^T\right]^T$ and each $\mathbf{x}_i \in \mathbb{R}^{m_i}$, and leveraging a multi-core computing environment to optimize the blocks *in parallel*. When using the MM algorithm, a natural way to achieve this goal is to construct a surrogate function $\widetilde{V}(\bullet\,|\,\mathbf{y})$ that is *additively separable* in the blocks, i.e., $\widetilde{V}(\mathbf{x}\,|\,\mathbf{y}) = \sum_{i=1}^n \widetilde{V}_i(\mathbf{x}_i\,|\,\mathbf{y})$, so that the subproblem (2) can be decoupled in independent optimization problems, one per block. Some of the examples discussed above [cf. (26), (28), (35)] show cases where this can be readily achieved by exploring the special structure of the objective function $V$. However, when the objective function $V$ does not enjoy any special structure, it seems that there is no general procedure to leverage to obtain a separable surrogate function satisfying Assumption I.13. In such cases, a way to cope with the curse of dimensionality is to solve the subproblem (2) by minimizing the (nonseparable) surrogate function block-wise. Such a version of the MM algorithm, termed *Block-MM*, is discussed in details in the next section.

## I.3. The Block Majorization-Minimization Algorithm

In this section, we introduce the Block MM (BMM) algorithm solving Problem (1), whose feasible set is now assumed to have a Cartesian product structure.

**Assumption I.15.** *The feasible set $X$ of Problem* (1) *admits the Cartesian product structure $X = X_1 \times \cdots \times X_n$, with each $X_i \in \mathbb{R}^{m_i}$ and $\sum_{i=1}^n m_i = m$.*

In the above setting, the optimization variables of (1) can be partitioned in (nonoverlapping) blocks $\mathbf{x}_i$. When judiciously exploited, this block-structure can lead to low-complexity algorithms that are implementable in a parallel/distributed manner. The BMM algorithm is a variant of the MM scheme (cf. Algorithm 1) whereby *only one* block of variables is updated at a time. More specifically, at iteration $k$, a block is selected according to some predetermined rule, say block $i^k$, and updated by solving the subproblem

$$\mathbf{x}_{i^k}^{k+1} \in \operatorname*{argmin}_{\mathbf{x}_{i^k} \in X_{i^k}} \widetilde{V}_{i^k}\left(\mathbf{x}_{i^k}\,|\,\mathbf{x}^k\right), \tag{45}$$

whereas $\mathbf{x}_j^{k+1} = \mathbf{x}_j^k$, for all $j \neq i^k$. In (45), $\widetilde{V}_{i^k}(\bullet | \mathbf{x}^k)$ is a suitably chosen surrogate function of $V$, now function only of the block $\mathbf{x}_{i^k}$. Conditions on $\widetilde{V}_{i^k}$ are similar to those already introduced for the MM algorithm (cf. Assumption I.13) and summarized below. We use the following notation: $O_i \subseteq \mathbb{R}^{m_i}$ is an open set containing $X_i$ and $O = O_1 \times \cdots \times O_n$; $\mathbf{x}_{-i} \triangleq [\mathbf{x}_1^T, \ldots, \mathbf{x}_{i-1}^T, \mathbf{x}_{i+1}^T, \ldots, \mathbf{x}_n^T]^T$ denotes the vector containing all the blocks of $\mathbf{x}$ but the $i$-th one; $(\mathbf{x})_i \triangleq [\mathbf{0}^T, \ldots, \mathbf{0}^T, \mathbf{x}_i^T, \mathbf{0}^T, \ldots, \mathbf{0}^T]^T$ denotes the operator nulling all the blocks of $\mathbf{x}$ except for the $i$-th one; and with a slight abuse of notation, we use $(\mathbf{x}_i, \mathbf{y}_{-i})$ to denote the ordered tuple $(\mathbf{y}_1, \ldots, \mathbf{y}_{i-1}, \mathbf{x}_i, \mathbf{y}_{i+1}, \ldots, \mathbf{y}_n)$.

**Assumption I.16.** *Each surrogate function $\widetilde{V}_i : O_i \times O \to \mathbb{R}$ satisfies the following conditions:*

1. $\widetilde{V}_i(\bullet | \bullet)$ *is continuous on $X_i \times X$;*
2. $\mathbf{y}_i \in \mathrm{argmin}_{\mathbf{x}_i \in X_i} \widetilde{V}_i(\mathbf{x}_i | \mathbf{y}) - V(\mathbf{x}_i, \mathbf{y}_{-i})$, *for all $\mathbf{y} \in X$;*
3. *The directional derivative of $\widetilde{V}_i'$ satisfies $\widetilde{V}_i'(\mathbf{x}_i; \mathbf{d}_i | \mathbf{x}) = V'(\mathbf{x}; (\mathbf{d})_i)$, for all $\mathbf{x}_i \in X_i$ and feasible directions $\mathbf{d}_i \in \mathbb{R}^{m_i}$ at $\mathbf{x}_i$.*

The next question is how to choose the blocks to be updated at each iteration, in order to have convergence. Any of the following rules can be adopted.

**Assumption I.17.** *Blocks in (45) are chosen according to any of the following rules:*

1. **Essentially cyclic rule**: *There exits a finite positive $T$ such that $\bigcup_{t=0}^{T-1} \{i^{k+t}\} = \{1, \ldots, n\}$, for all $k = 0, 1, \ldots$;*
2. **Maximum block improvement rule**: *Choose*

$$i^k \in \underset{i=1,\ldots,n}{\mathrm{argmin}} \min_{\mathbf{x}_i \in X_i} \widetilde{V}_i\left(\mathbf{x}_i | \mathbf{x}^k\right);$$

3. **Random-based rule**: *There exists a $p_{\min} > 0$ such that $\mathbb{P}\left(i^k = j | \mathbf{x}^{k-1}, \ldots, \mathbf{x}^0\right) = p_j^k \geq p_{\min}$, for all $j = 1, \ldots n$ and $k = 0, 1, \ldots$.*

The first two rules above are deterministic rules. Roughly speaking, the essentially cyclic rule (Assumption I.17.1) states that all the blocks must be updated at least once within any $T$ consecutive iterations, with some blocks possibly updated more frequently than others. Of course, a special case of this rule is the simple cyclic rule, i.e., $i^k = (k \bmod n) + 1$, whereby all the blocks are updated once every $T$ iterations. The maximum block improvement rule (Assumption I.17.2) is a greedy-based update: only the block that generates the largest decrease of the surrogates $\widetilde{V}_i$ at $\mathbf{x}^k$ is selected and updated. Finally the random-based selection rule (Assumption I.17.3) selects block randomly and independently, with the condition that all the blocks have a bounded away from zero probability to be chosen. The BMM algorithm is summarized in Algorithm 2.

**Algorithm 2: Block MM Algorithm**

**Data** : $\mathbf{x}^0 \in X$. Set $k = 0$.
(S.1) : If $\mathbf{x}^k$ satisfies a termination criterion: STOP;
(S.2) : Choose an index $i^k \in \{1, \ldots, n\}$;
(S.3) : Update $\mathbf{x}^k$ as
      – Set $\mathbf{x}^{k+1}_{i^k} \in \mathrm{argmin}_{\mathbf{x}_{i^k} \in X_{i^k}} \widetilde{V}_{i^k}\left(\mathbf{x}_{i^k} \,|\, \mathbf{x}^k\right)$;
      – Set $\mathbf{x}^{k+1}_j = \mathbf{x}^k_j$, for all $j \neq i_k$;
(S.4) : $k \leftarrow k+1$, and go to (S.1).

---

Convergence results of Algorithm 2 consist of two major statements. Under the essential cyclic rule (Assumption I.17.1), quasi convexity of the objective function is required [along with the uniqueness of the minimizer in (45)], which also guarantees the existence of the limit points. This is in the same spirit as the classical proof of convergence of Block Coordinate Descent methods; see, e.g., [14, 214, 250]. If the maximum block improvement rule (Assumption I.17.2) or the random-based selection rule (Assumption I.17.3) are used, then a stronger convergence result can be proved by relaxing the quasi-convexity assumption and imposing the compactness of the level sets of $V$. In order to state the convergence result, we introduce the following additional assumptions.

**Assumption I.18.** $\widetilde{V}_i(\bullet \,|\, \mathbf{y})$ *is quasi-convex and subproblem* (45) *has a unique solution, for all* $\mathbf{y} \in X$.

**Assumption I.19.** *The level set* $X^0 \triangleq \{\mathbf{x} \in X : V(\mathbf{x}) \leq V(\mathbf{x}^0)\}$ *is compact and subproblem* (45) *has a unique solution for all* $\mathbf{x}^k \in X$ *and at least* $n-1$ *blocks.*

We are now ready to state the main convergence result of the BMM algorithm (Algorithm 2), as given in Theorem I.20 below (the statement holds almost surely for the random-based selection rule). The proof of this theorem can be found in [193]; see also Sec. I.5 for a detailed discussion on (other) existing convergence results.

**Theorem I.20.** *Given Problem* (1) *under Assumptions I.12 and I.15, let* $\{\mathbf{x}^k\}_{k \in \mathbb{N}_+}$ *be the sequence generated by Algorithm 2, with* $\widetilde{V}$ *chosen according to Assumption I.16. Suppose that, in addition, either one of the following two conditions is satisfied:*

*(a) $i^k$ is chosen according to the essentially cyclic rule (Assumption I.17.1) and $\widetilde{V}$ further satisfies either Assumption I.18 or I.19;*
*(b) $i^k$ is chosen according to the maximum block improvement rule (Assumption I.17.2) or the random-based selection rule (Assumption I.17.3).*

*Then, every limit point of $\{\mathbf{x}^k\}_{k \in \mathbb{N}_+}$ is a coordinate-wise d-stationary solution of* (1). *Furthermore, if $V$ is regular, then every limit point is a d-stationary solution of* (1). □

## I.4. Applications

In this section, we show how to apply the (B)MM algorithm to solve some representative nonconvex problems arising from applications in signal processing, data analysis, and communications. More specifically, we consider the following problems: i) Sparse least squares; ii) Nonnegative least squares; iii) Matrix factorization, including low-rank factorization and dictionary learning; and iv) the multicast beamforming problem. Our list is by no means exhaustive; it just gives a flavor of the kind of structured nonconvexity and applications which (B)MM can be successfully applied to.

### I.4.1 Nonconvex Sparse Least Squares

Retrieving a sparse signal from its linear measurements is a fundamental problem in machine learning, signal processing, bioinformatics, physics, etc.; see [269] and [4, 97] for a recent overview and some books, respectively, on the subject. Consider a linear model, $\mathbf{z} = \mathbf{A}\mathbf{x} + \mathbf{n}$, where $\mathbf{x} \in \mathbb{R}^m$ is the sparse signal to estimate, $\mathbf{z} \in \mathbb{R}^q$ is the vector of available measurements, $\mathbf{A} \in \mathbb{R}^{q \times m}$ is the given measurement matrix, and $\mathbf{n} \in \mathbb{R}^q$ is the observation noise. To estimate the sparse signal $\mathbf{x}$, a mainstream approach in the literature is to solve the following optimization problem

$$\underset{\mathbf{x}}{\text{minimize}} \quad V(\mathbf{x}) \triangleq \|\mathbf{z} - \mathbf{A}\mathbf{x}\|^2 + \lambda G(\mathbf{x}), \tag{46}$$

where the first term in the objective function measures the model fitness whereas the regularization $G$ is used to promote sparsity in the solution, and the regularization parameter $\lambda \geq 0$ is chosen to balance the trade-off between the model fitness and sparsity of the solution.

The ideal choice for $G$ would be the cardinality of $\mathbf{x}$, also referred to as $\ell_0$ "norm" of $\mathbf{x}$. However, its combinatorial nature makes the resulting optimization problem numerically intractable as the variable dimension $m$ becomes large. Due to its favorable theoretical guarantees (under some regularity conditions on $\mathbf{A}$ [33, 34]) and the existence of efficient solution methods for convex instances of (46), the $\ell_1$ norm has been widely adopted in the literature as convex surrogate $G$ of the $\ell_0$ function (in fact, the $\ell_1$ norm is the convex envelop of the $\ell_0$ function on $[-1, 1]^m$) [30, 235]. Yet there is increasing evidences supporting the use of nonconvex formulations to enhance the sparsity of the solution as well as the realism of the models [2, 36, 94, 152, 223]. For instance, it is well documented that *nonconvex* surrogates of the $\ell_0$ function, such as the SCAD [80], the "transformed" $\ell_1$, the logarithmic, the exponential, and the $\ell_p$ penalty [266], outperform the $\ell_1$ norm in enhancing the sparsity of the solution. Table I.1 summarizes these nonconvex surrogates whereas Fig. I.3 shows their graph.

Quite interestingly, it has been recently shown that the aforementioned nonconvex surrogates of the $\ell_0$ function enjoy a separable DC (Difference of Convex) structure (see, e.g., [2, 233] and references therein); specifically, we have the following

Table I.1: Examples of nonconvex surrogates of the $\ell_0$ function having a DC structure [cf. (47)]

| Penalty function | Expression |
|---|---|
| Exp [29] | $g_{\exp}(x) = 1 - e^{-\theta\|x\|}$ |
| $\ell_p (0 < p < 1)$ [88] | $g_{\ell_p^+}(x) = (\|x\| + \varepsilon)^{1/\theta}$, |
| $\ell_p (p < 0)$ [192] | $g_{\ell_p^-}(x) = 1 - (\theta\|x\| + 1)^p$ |
| SCAD [80] | $g_{\text{scad}}(x) = \begin{cases} \frac{2\theta}{a+1}\|x\|, & 0 \leq \|x\| \leq \frac{1}{\theta} \\ \frac{-\theta^2\|x\|^2 + 2a\theta\|x\| - 1}{a^2 - 1}, & \frac{1}{\theta} < \|x\| \leq \frac{a}{\theta} \\ 1, & \|x\| > \frac{a}{\theta} \end{cases}$ |
| Log [248] | $g_{\log}(x) = \frac{\log(1+\theta\|x\|)}{\log(1+\theta)}$ |

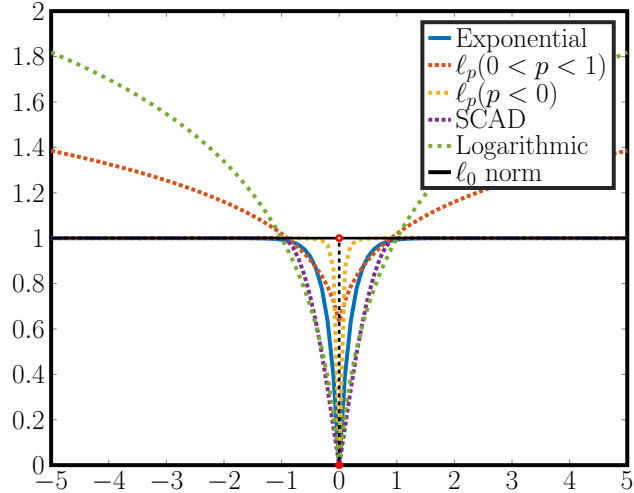

Fig. I.3: Nonconvex surrogate functions of the $\ell_0$ function given in Table I.1.

$$G(\mathbf{x}) = \sum_{i=1}^{m} g(x_i), \quad \text{with} \quad g(x_i) = \underbrace{\eta(\theta)|x_i|}_{\triangleq g^+(x_i)} - \underbrace{(\eta(\theta)|x_i| - g(x_i))}_{\triangleq g^-(x_i)}, \quad (47)$$

where the specific expression of $g : \mathbb{R} \to \mathbb{R}$ is given in Table I.1; and $\eta(\theta)$ is a fixed given function, whose expression depends on the surrogate $g$ under consideration, see Table I.2. Note that the parameter $\theta$ controls the tightness of the approxima-tion of the $\ell_0$ function: in fact, it holds that $\lim_{\theta \to +\infty} g(x_i) = 1$ if $x_i \neq 0$, otherwise $\lim_{\theta \to +\infty} g(x_i) = 0$. Moreover, it can be shown that for all the functions in Table I.1, $g^-$ is convex and has Lipschitz continuous first derivative $dg^-/dx$ [233], whose closed form is given in Table I.2.

*Table I.2: Explicit expression of $\eta(\theta)$ and $dg^-/dx$ [cf. (47)]*

| $g$ | $\eta(\theta)$ | $dg_\theta^-/dx$ |
|---|---|---|
| $g_{\exp}$ | $\theta$ | $\text{sign}(x) \cdot \theta \cdot (1 - e^{-\theta|x|})$ |
| $g_{\ell_p^+}$ | $\frac{1}{\theta} \varepsilon^{1/\theta - 1}$ | $\frac{1}{\theta} \text{sign}(x) \cdot [\varepsilon^{\frac{1}{\theta}-1} - (|x|+\varepsilon)^{\frac{1}{\theta}-1}]$ |
| $g_{\ell_p^-}$ | $-p \cdot \theta$ | $-\text{sign}(x) \cdot p \cdot \theta \cdot [1 - (1+\theta|x|)^{p-1}]$ |
| $g_{\text{scad}}$ | $\frac{2\theta}{a+1}$ | $\begin{cases} 0, & \|x\| \le \frac{1}{\theta} \\ \text{sign}(x) \cdot \frac{2\theta(\theta\|x\|-1)}{a^2-1}, & \frac{1}{\theta} < \|x\| \le \frac{a}{\theta} \\ \text{sign}(x) \cdot \frac{2\theta}{a+1}, & \text{otherwise} \end{cases}$ |
| $g_{\log}$ | $\frac{\theta}{\log(1+\theta)}$ | $\text{sign}(x) \cdot \frac{\theta^2 \|x\|}{\log(1+\theta)(1+\theta\|x\|)}$ |

Motivated by the effectiveness of the aforementioned nonconvex surrogates of the $\ell_0$ function and recent works [2, 149, 208, 261, 266], in this section, we show how to use the MM framework to design efficient algorithms for the solution of Problem (46), where $G$ is assumed to have the DC structure (47), capturing thus in a unified way all the nonconvex $\ell_0$ surrogates reported in Table I.1. The key question is how to construct a valid surrogate function $\widetilde{V}$ of $V$ in (46). We address this issue following two steps: 1) We first find a surrogate $\widetilde{G}$ for the nonconvex $G$ in (47) satisfying Assumption I.16; and 2) then we construct the overall surrogate $\widetilde{V}$ of $V$, building on $\widetilde{G}$.

**Step** 1**: Surrogate $\widetilde{G}$ of** $G$**.** There are two ways to construct $\widetilde{G}$, namely: i) tailoring $\widetilde{G}$ to the specific structure of the function $g$ under consideration (cf. Table I.1); or ii) leveraging the DC structure of $g$ in (47) and obtain a unified expression for $\widetilde{G}$, valid for all the DC functions in Table I.1. Few examples based on the approach i) are shown first, followed by the general design as in ii).

*Example 8 (The log $\ell_0$ surrogate).* Let $G$ be the "log" $\ell_0$ surrogate, i.e., $G(\mathbf{x}) = \sum_{i=1}^m g_{\log}(x_i)$, where $g_{\log}$ is defined in Table I.1. A valid surrogate is obtained majorizing the log function $g_{\log}$ (cf. Example 1), which leads to $\widetilde{G}(\mathbf{x}\,|\,\mathbf{y}) = \sum_{i=1}^m \widetilde{g}_{\log}(x_i\,|\,y_i)$, with

$$\widetilde{g}_{\log}(x_i\,|\,y_i) = \frac{\theta}{\log(1+\theta)} \cdot \frac{1}{1+\theta|y_i|} \, |x_i|. \tag{48}$$

*Example 9 (The $\ell_p(0 < p < 1)$ surrogate).* Let $G$ be the $\ell_p(0 < p < 1)$ function, i.e., $G(\mathbf{x}) = \sum_{i=1}^m g_{\ell_p^+}(x_i)$, where $g_{\ell_p^+}$ is defined in Table I.1. Similar to the Log surrogate, we can derive a majorizer of such a $G(\mathbf{x})$ by exploiting the concavity of $g_{\ell_p^+}$. We have

$$\widetilde{g}_{\ell_p^+}(x_i | y_i) = \frac{1}{\theta}(|y_i| + \varepsilon)^{1/\theta - 1}|x_i|. \tag{49}$$

The desired valid surrogate is then given by $\widetilde{G}(\mathbf{x} | \mathbf{y}) = \sum_{i=1}^m \widetilde{g}_{\ell_p^+}(x_i | y_i)$.

*Example 10 (DC surrogates).* We consider now nonconvex regularizers having the DC structure (47). A natural approach is to keep the convex component $g^+$ in (47) while linearizing the differentiable concave part $-g^-$, which leads to the following convex majorizer:

$$\widetilde{g}(x_i | y_i) = \eta(\theta) |x_i| - \left. \frac{dg^-(x)}{dx} \right|_{x = y_i} \cdot (x_i - y_i), \tag{50}$$

where the expression of $dg(x)^-/dx$ is given in Table I.2. The desired majorization then reads $\widetilde{G}(\mathbf{x} | \mathbf{y}) = \sum_{i=1}^m \widetilde{g}(x_i | y_i)$.

Note that although the log and $\ell_p$ surrogates provided in Example 8 and 9 are special cases of DC surrogates, the majorizer constructed using (50) is different from the ad-hoc surrogates $\widetilde{g}_{\log}$ and $\widetilde{g}_{\ell_p^+}$ in (48) and (49), respectively.

**Step** 2**: Surrogate $\widetilde{V}$ of $V$.** We derive now the surrogate of $V$, when $G$ is given by (47). Since the loss function $\|\mathbf{z} - \mathbf{A}\mathbf{x}\|^2$ is convex, two natural options for $\widetilde{V}$ are: 1) keeping $\|\mathbf{z} - \mathbf{A}\mathbf{x}\|^2$ unaltered while replacing $G$ with the surrogate $\widetilde{G}$ discussed in Step 1; or 2) majorizing also $\|\mathbf{z} - \mathbf{A}\mathbf{x}\|^2$. The former approach "better" preserves the structure of the objective function but at the price of a higher cost of each iteration [cf. (2)]: the minimizer of the resulting $\widetilde{V}$ does not have a closed form expression; the overall algorithm is thus a double-loop scheme. The latter approach is instead motivated by the goal of obtaining low-cost iterations. We discuss next both options and establish an interesting connection between the two resulting MM algorithms.

• **Option 1.** Keeping $\|\mathbf{z} - \mathbf{A}\mathbf{x}\|^2$ unaltered while using (50), leads to the following update in Algorithm 1:

$$\mathbf{x}^{k+1} \in \underset{\mathbf{x}}{\operatorname{argmin}} \left\{ \widetilde{V}(\mathbf{x} | \mathbf{x}^k) \triangleq \|\mathbf{A}\mathbf{x} - \mathbf{z}\|^2 + \lambda \sum_{i=1}^m \widetilde{g}(x_i | x_i^k) \right\}. \tag{51}$$

Problem (51) is convex but does not have a closed form solution. To solve (51), we develop next an ad-hoc iterative soft-thresholding-based algorithm invoking again the MM framework on $\widetilde{V}(\bullet | \mathbf{x}^k)$ in (51).

We denote by $\mathbf{x}^{k,r}$ the $r$-th iterate of the (inner loop) MM algorithm (Algorithm 1) used to solve (51); we initialize the inner algorithm by $\mathbf{x}^{k,0} = \mathbf{x}^k$. Since the quadratic term $\|\mathbf{A}\mathbf{x} - \mathbf{z}\|^2$ in (51) has Lipschitz gradient, a natural surrogate function for $\widetilde{V}(\mathbf{x} | \mathbf{x}^k)$, is [cf. (28)]: given $\mathbf{x}^{k,r}$,

$$\widetilde{V}^k(\mathbf{x} | \mathbf{x}^{k,r}) = 2\mathbf{x}^T \mathbf{A}^T (\mathbf{A}\mathbf{x}^{k,r} - \mathbf{z}) + \frac{L}{2} \|\mathbf{x} - \mathbf{x}^{k,r}\|^2 + \lambda \sum_{i=1}^m \widetilde{g}(x_i | x_i^k), \tag{52}$$

where $L = 2\lambda_{\max}(\mathbf{A}^T\mathbf{A})$. Denoting

$$\mathbf{b}^{k,r} \triangleq \mathbf{x}^{k,r} - \frac{2}{L}\mathbf{A}^T(\mathbf{A}\mathbf{x}^{k,r} - \mathbf{z}) + \frac{\lambda}{L} \cdot \left( \left. \frac{dg^-(x)}{dx} \right|_{x=x_i^k} \right)_{i=1}^m,$$

the main update of the inner MM algorithm minimizing $\widetilde{V}^k(\mathbf{x} \mid \mathbf{x}^{k,r})$ reads

$$\mathbf{x}^{k,r+1} = \underset{\mathbf{x}}{\arg\min} \frac{L}{2}\|\mathbf{x} - \mathbf{b}^{k,r}\|^2 + \lambda\, \eta(\theta)\, \|\mathbf{x}\|_1. \tag{53}$$

Quite interestingly, the solution of Problem (53) can be obtained in closed form. Writing the first order optimality condition

$$\mathbf{0} \in (\mathbf{x} - \mathbf{b}^{k,r}) + \frac{\lambda\, \eta(\theta)}{L} \partial \|\mathbf{x}\|_1,$$

and recall that the subgradient of $\|\mathbf{x}\|_1$ takes the following form [cf. (10)]:

$$\partial \|\mathbf{x}\|_1 = \{\boldsymbol{\zeta} : \boldsymbol{\zeta}^T \mathbf{x} = \|\mathbf{x}\|_1,\, \|\boldsymbol{\zeta}\|_\infty \leq 1\},$$

we have the following expression for $\mathbf{x}^{k,r+1}$:

$$\mathbf{x}^{k,r+1} = \mathrm{sign}(\mathbf{b}^{k,r}) \cdot \max\left\{ |\mathbf{b}^{k,r}| - \frac{\lambda\, \eta(\theta)}{L}, 0 \right\},$$

where the sign and max operators are applied component-wise. Introducing the *soft-thresholding operator*

$$S_\alpha(x) \triangleq \mathrm{sign}(x) \cdot \max\{|x| - \alpha, 0\}, \tag{54}$$

$\mathbf{x}^{k,r+1}$ can be rewritten succinctly as

$$\mathbf{x}^{k,r+1} = S_{\frac{\lambda\, \eta(\theta)}{L}}(\mathbf{b}^{k,r}), \tag{55}$$

where the soft-thresholding operator is applied component-wise.

Overall the double loop algorithm, based on the MM outer updates (51) and MM inner iterates (55) is summarized in Algorithm 3.

**Algorithm 3: MM Algorithm for Nonconvex Sparse Least Squares**

**Data** : $\mathbf{x}^0 \in X$. Set $k = 0$.
(S.1) : If $\mathbf{x}^k$ satisfies a termination criterion: STOP;
(S.2) : Set $r = 0$. Initialize $\mathbf{x}^{k,0}$ as $\mathbf{x}^{k,0} = \mathbf{x}^k$;
(S.3) : If $\mathbf{x}^{k,r}$ satisfies a termination criterion: STOP;
       (a) : Update $\mathbf{x}^{k,r}$ as

$$\mathbf{x}^{k,r+1} = S_{\frac{\lambda \eta(\theta)}{L}}\left(\mathbf{b}^{k,r}\right);$$

       (b) : $r \leftarrow r + 1$, and go to (S.3).
(S.4) : $\mathbf{x}^{k+1} = \mathbf{x}^{k,r}$, $k \leftarrow k+1$, and go to (S.1).

*Termination criteria:* As the termination criterion of Step 1 and Step 3 in Algorithm 3, one can use any valid merit function measuring the distance of the iterates from stationarity of (46) and optimality of (53), respectively. For both the inner and outer loop, it is not difficult to check that the objective functions of the associated optimization problems−(53) and (51), respectively−are strictly convex if $\lambda > 0$; therefore both optimization problem have unique minimizers. It turns out that both functions defined in (24) and (25) can be adopted as valid merit functions. The loop can be then terminated once the value of the chosen function goes below the desired threshold.

• **Option 2.** Algorithm 3 is a double loop MM-based algorithm: in the outer loop, the surrogate function $\widetilde{V}(\bullet \mid \mathbf{x}^k)$ [cf. (51)] is iteratively minimized by means of an inner MM algorithm based on the surrogate function $\widetilde{V}^k(\mathbf{x} \mid \mathbf{x}^{k,r})$ [cf. (52)]. A closed look at (51) and (52) shows that the following relationship holds between $\widetilde{V}$ and $\widetilde{V}^k$:

$$\widetilde{V}^k(\mathbf{x} \mid \mathbf{x}^{k,0}) \geq \widetilde{V}(\mathbf{x} \mid \mathbf{x}^k) \geq V(\mathbf{x}).$$

The above inequality shows that $\widetilde{V}^{k,0}(\mathbf{x} \mid \mathbf{x}^k)$ is in fact a valid surrogate function of $V$ at $\mathbf{x} = \mathbf{x}^k$. This means that the inner loop of Algorithm 3 can be terminated after one iteration without affecting the convergence of the scheme. Specifically, Step 3 of Algorithm 3 can be replaced with the following iterate

$$\mathbf{x}^{k+1} = S_{\frac{\lambda \eta(\theta)}{L}}\left(\mathbf{b}^{k,0}\right). \tag{56}$$

The resulting algorithm is in fact an MM scheme minimizing $\widetilde{V}^{k,0}(\bullet \mid \mathbf{x}^k)$, whose convergence is guaranteed by Theorem I.14.

### I.4.2 Nonnegative Least Squares

Finding a nonnegative solution $\mathbf{x} \in \mathbb{R}^m$ of the linear model $\mathbf{z} = \mathbf{A}\mathbf{x} + \mathbf{n}$ has attracted significant attention in the literature. This problem arises in applications where the measured data is nonnegative; examples include image pixel intensity, economical quantities such as stock price and trading volume, biomedical records such as weight, height, blood pressure, etc [73,139,212]. It is also one of the key ingredients in nonnegative matrix/tensor factorization problem for analysing structured data set. The nonnegative least squares (NNLS) problem consists in finding a nonnegative **x** that minimizes the residual error between the data and the model in least square sense:

$$\begin{aligned}\underset{\mathbf{x}}{\text{minimize}} \quad & V(\mathbf{x}) \triangleq \|\mathbf{z} - \mathbf{A}\mathbf{x}\|^2 \\ \text{subject to} \quad & \mathbf{x} \geq \mathbf{0},\end{aligned} \quad (57)$$

where $\mathbf{z} \in \mathbb{R}^q$ and $\mathbf{A}^{q \times m}$ are given.

Note that Problem (57) is convex. We show next how to construct a surrogate function satisfying Assumption I.13 which is additively separable in the components of **x**, so that the resulting subproblems (2) can be solved in parallel. To this end, we expand the square in the objective function and write

$$V(\mathbf{x}) = \mathbf{x}^T \mathbf{A}^T \mathbf{A} \mathbf{x} - 2 \mathbf{z}^T \mathbf{A} \mathbf{x} + \mathbf{z}^T \mathbf{z}. \quad (58)$$

Let $\mathbf{M} \succeq \mathbf{A}\mathbf{A}^T$. Using (29), we can majorize $V$ by

$$\widetilde{V}(\mathbf{x}|\mathbf{y}) = V(\mathbf{y}) + 2(\mathbf{A}^T \mathbf{A} \mathbf{y} - \mathbf{A}^T \mathbf{z})^T (\mathbf{x} - \mathbf{y}) + (\mathbf{x} - \mathbf{y})^T \mathbf{M} (\mathbf{x} - \mathbf{y}), \quad (59)$$

for any given $\mathbf{y} \in \mathbb{R}^m$. The goal is then finding a matrix $\mathbf{M} \succeq \mathbf{A}\mathbf{A}^T$ that is diagonal, so that $\widetilde{V}(\mathbf{x}|\mathbf{y})$ becomes additively separable. We provide next two alternative expressions for **M**.

• **Option 1:** Since $\nabla V$ is Lipschitz continuous on $\mathbb{R}^m$, we can use the same upperbound as in Example 3. This corresponds to choose $\mathbf{M} = \lambda \mathbf{I}$, with $\lambda$ such that $\lambda \geq \lambda_{\max}(\mathbf{A}^T \mathbf{A})$. This leads overall to the following surrogate of $V$:

$$\widetilde{V}(\mathbf{x}|\mathbf{y}) = V(\mathbf{y}) + 2(\mathbf{A}^T \mathbf{A} \mathbf{y} - \mathbf{A}^T \mathbf{z})^T (\mathbf{x} - \mathbf{y}) + \lambda \|\mathbf{x} - \mathbf{y}\|^2.$$

The above choice leads to a strongly convex subproblem (2), whose minimizer has the following closed form:

$$\mathbf{x}^{k+1} = \left[\mathbf{x}^k - \frac{1}{\lambda}\left(\mathbf{A}^T \mathbf{A} \mathbf{x}^k - \mathbf{A}^T \mathbf{z}\right)\right]_+, \quad (60)$$

where $[\bullet]_+$ denotes the Euclidean projection onto the nonnegative orthant. The resulting MM algorithm (Algorithm 1) based on the update (60) turns out to be the renowned gradient projection algorithm with constant step-size $1/\lambda$.

• **Option 2:** If Problem (57) has some extra structure, the surrogate $\widetilde{V}$ can be tailored to $V$ even further. For instance, suppose that *in addition* to the structure above, there hold $\mathbf{A} \in \mathbb{R}_{++}^{q \times m}$, $\mathbf{z} \in \mathbb{R}_+^q$ and $\mathbf{z} \neq \mathbf{0}$. It has been shown in [58,139] that the following diagonal matrix

$$\mathbf{M} \triangleq \text{Diag}\left(\frac{(\mathbf{A}^T \mathbf{A} \mathbf{y})_1}{y_1}, \cdots, \frac{(\mathbf{A}^T \mathbf{A} \mathbf{y})_m}{y_m}\right) \quad (61)$$

satisfies $\mathbf{M} \succeq \mathbf{A}^T\mathbf{A}$. Substituting (61) in (59), one obtains the following closed form solution of the resulting subproblem (2):

$$\mathbf{x}^{k+1} = (\mathbf{A}^T\mathbf{z}/\mathbf{A}^T\mathbf{A}\mathbf{x}^k) \cdot \mathbf{x}^k,$$

wherein both division and multiplication are intended to be applied element-wise.

### I.4.3 Sparse plus Low-Rank Matrix Decomposition

Another useful paradigm is to decompose a partly or fully observed data matrix into the sum of a low rank and (bilinear) sparse term; the low-rank component captures correlations and periodic trends in the data whereas the bilinear term explains parsimoniously data patterns, (co-)clusters, innovations or outliers.

Let $\mathbf{Y} \in \mathbb{R}^{m \times t}(m \leq t)$ be the data matrix. The goal is to find a low rank matrix $\mathbf{L} \in \mathbb{R}^{m \times t}$ with rank $r_0 \triangleq \text{rank}(\mathbf{L}) \ll m$, and a sparse matrix $\mathbf{S} \in \mathbb{R}^{m \times t}$ such that $\mathbf{Y} = \mathbf{L} + \mathbf{S} + \mathbf{V}$, where $\mathbf{V} \in \mathbb{R}^{m \times t}$ accounts for measurement errors. To cope with the missing data in $\mathbf{Y}$, we introduce i) the set $\Omega \subseteq M \times T$ of index pairs $(i, j)$, where $M \triangleq \{1, \ldots, m\}$ and $T \triangleq \{1, \ldots, t\}$; and ii) the sampling operator $P_\Omega(\cdot)$, which nulls the entries of its matrix argument not in $\Omega$, leaving the rest unchanged. In this way, one can express incomplete and (possibly noise-)corrupted data as

$$P_\Omega(\mathbf{Y}) \approx P_\Omega(\mathbf{L} + \mathbf{S}), \tag{62}$$

where $\approx$ is quantified by a specific loss function (and regularization) [219].

Model (62) subsumes a variety of statistical learning paradigms including: (robust) principal component analysis [32,41], compressive sampling [35], dictionary learning (DL) [75, 178], non-negative matrix factorization [107, 122, 138], matrix completion, and their robust counterparts [112]. Task (62) also emerges in various applications such as (i) network anomaly detection [130, 154, 155]; (i) distributed acoustic signal processing [71, 72]; (iii) distributed localization and sensor identification [205]; (iv) distributed seismic forward modeling in geological applications [171, 271] (e.g., finding the Green's function of some model of a portion of the earth's surface); (v) topic modeling for text corpus from social media [198, 267]; (vi) data and graph clustering [98, 126, 127, 236]; and (vii) power grid state estimation [93, 120].

In the following, we study task (62) adopting the least-squares (LS) error as loss function for $\approx$. We show in detail how to design an MM algorithm for the solution of two classes of problems under (62), namely: i) the low-rank matrix completion problem; and ii) the dictionary learning problem. Similar techniques can be used also to solve other tasks modeled by (62).

**1) Low-rank matrix completion.** The low-rank matrix completion problem arises frequently in learning problems where the task is to fill in the blanks in the partially observed collinear data. For instance, in movie rating problems, $\mathbf{Y}$ is the rating matrix whose entries $y_{ij}$ represent the score of movie $j$ given by individual $i$ if he/she has watched it, and considered missing otherwise. Despite of being highly incomplete, such a data set is also rank deficient as individuals sharing similar interests may give similar ratings (the corresponding rows of **Y** are collinear), which makes the matrix completion task possible. Considering model (62), the question becomes how to impose a low-rank structure on **L** and sparse structure on **S**. We describe next two widely used approaches.

A first approach is enforcing a low-rank structure on **L** by promoting sparsity on the singular values of **L** [denoted by $\sigma_i(\mathbf{L})$, $i = 1, \ldots, m$] as well as on the elements of **S** via regularization. This leads to the following formulation

$$\underset{\mathbf{L},\mathbf{S}}{\text{minimize}} \quad V(\mathbf{L},\mathbf{S}) \triangleq \|P_\Omega(\mathbf{Y} - \mathbf{L} - \mathbf{S})\|_F^2 + \lambda_r \cdot G_r(\mathbf{L}) + \lambda_s \cdot G_s(\mathbf{S}), \tag{63}$$

where $G_r(\mathbf{L}) \triangleq \sum_{i=1}^m g_r(\sigma_i(\mathbf{L}))$ and $G_s(\mathbf{S}) \triangleq \sum_{i=1}^r \sum_{j=1}^t g_s(s_{ij})$ are sparsity promoting regularizers, and $\lambda_r$ and $\lambda_s$ positive coefficients. Since $G_r(\mathbf{L})$ promotes sparsity on the singular values of **L**, this will induce a low-rank structure on **L**. Note that the general formulation (63) contains many popular choices of low-rank inducing penalties. For instance, choosing $g_r(x) = g_s(x) = \text{card}(x)$, one gets $G_r(\mathbf{L}) = \text{rank}(\mathbf{L})$ and $G_s(\mathbf{S}) = \|\mathbf{S}\|_0 \triangleq \|\text{vec}(\mathbf{S})\|_0$, which become the exact (nonconvex) rank penalty on **L** and the cardinality penalty on **S**, respectively. Another popular choice is $g_r(x) = |x|$, which leads to the convex nuclear norm penalty $G_r(\mathbf{L}) = \|\mathbf{L}\|_*$; another example is $g_r(x) = \log(x)$, which yields the nonconvex logdet penalty. To keep the analysis general, in the following we tacitly assume that $g_r$ and $g_s$ are any of the DC surrogate functions of the $\ell_0$ function introduced in Sec. I.4.1 [cf. (47)]. Note that, since $g_r$ and $g_s$ are DC and $\sigma_i(\mathbf{L})$ is a convex function of **L**, they are all directionally differentiable (cf. Sec. I.1). Therefore, $V(\mathbf{L},\mathbf{S})$ in (63) is directionally differentiable.

A second approach to enforce a low-rank structure on **L** is to "hard-wire" a rank (at most) $r$ into the structure of **L** by decomposing **L** as $\mathbf{L} = \mathbf{DX}$, where $\mathbf{D} \in \mathbb{R}^{m \times r}$ and $\mathbf{X} \in \mathbb{R}^{r \times t}$ are two "thin" matrices. The problem then reads

$$\underset{\mathbf{D},\mathbf{X}}{\text{minimize}} \quad \|P_\Omega(\mathbf{Y} - \mathbf{DX} - \mathbf{S})\|_F^2 + \lambda_r \cdot G_r(\mathbf{D},\mathbf{X}) + \lambda_s \cdot G_s(\mathbf{S}), \tag{64}$$

where $G_r$ and $G_s$ promote low-rank and sparsity structures, respectively. While $G_s$ can be chosen as in (63), the choice of $G_r$ acting on the two factors **D** and **X** while imposing the low-rankness of **L** is less obvious; two alternative choices are the following. Since

$$\|\mathbf{L}\|_* = \inf_{\mathbf{DX}=\mathbf{L}} \frac{1}{2} \left( \|\mathbf{D}\|_F^2 + \|\mathbf{X}\|_F^2 \right),$$

an option is choosing

$$G_r(\mathbf{D},\mathbf{X}) = \frac{1}{2} \left( \|\mathbf{D}\|_F^2 + \|\mathbf{X}\|_F^2 \right).$$

Another low-rank inducing regularizer is the max-norm penalty

$$\|\mathbf{L}\|_{\max} \triangleq \inf_{\mathbf{DX}=\mathbf{L}} \|\mathbf{D}\|_{2,\infty} + \|\mathbf{X}\|_{2,\infty},$$

where $\|\bullet\|_{2,\infty}$ denotes the maximum $\ell_2$ row norm of a matrix; this leads to

$$G_r(\mathbf{D}, \mathbf{X}) = \|\mathbf{D}\|_{2,\infty} + \|\mathbf{X}\|_{2,\infty}.$$

Here we focus only on the first formulation, Problem (63); the algorithmic design for (64) will be addressed within the context of the dictionary learning problem, which is the subject of the next section.

To deal with Problem (63), the first step is to rewrite the objective function in a more convenient form, by getting rid of the projection operator $P_\Omega$. Define $\mathbf{Q} \triangleq \mathrm{diag}(\mathbf{q})$, with $q_i = 1$ if $\left(\mathrm{vec}(P_\Omega(\mathbf{Y}))\right)_i \neq 0$; and $q_i = 0$ otherwise. Then, $V(\mathbf{L},\mathbf{S})$ in (63) can be rewritten as

$$V(\mathbf{L},\mathbf{S}) = \underbrace{\|\mathrm{vec}(\mathbf{Y}) - \mathrm{vec}(\mathbf{L}+\mathbf{S})\|_\mathbf{Q}^2}_{Q(\mathbf{L},\mathbf{S})} + \lambda_r \cdot G_r(\mathbf{L}) + \lambda_s \cdot G_s(\mathbf{S}), \qquad (65)$$

where $\|\mathbf{x}\|_\mathbf{Q}^2 \triangleq \mathbf{x}^T \mathbf{Q} \mathbf{x}$. Note that since $q_i$ is equal to either 0 or 1, we have $\lambda_{\max}(\mathbf{Q}) = 1$. In the following, we derive an algorithm that alternately optimizes $\mathbf{L}$ and $\mathbf{S}$ based on the block MM algorithm described in Algorithm 2. In the following, we denote by $m_{ij}$ the $(i,j)$-th entry of a generic matrix $\mathbf{M}$.

We start with the optimization of $\mathbf{S}$, given $\mathbf{L} = \mathbf{L}^k$. One can easily see that $V(\mathbf{L}^k, \mathbf{S})$ is of the same form as the objective function in Problem (46). Therefore, two valid surrogate functions can be readily constructed using the same technique already introduced in Sec. I.4.1. Specifically, two alternative surrogates of $V(\mathbf{L}^k, \mathbf{S})$ are [cf. (51)]: given $\mathbf{S}^k$,

$$\widetilde{V}^{(1)}(\mathbf{L}^k, \mathbf{S} \,|\, \mathbf{S}^k) = Q(\mathbf{L}^k, \mathbf{S}) + \lambda_s \sum_{i=1}^r \sum_{j=1}^t \widetilde{g}_s(s_{ij} \,|\, s_{ij}^k) \qquad (66)$$

and [cf. (52)]

$$\begin{aligned}
\widetilde{V}^{(2)}&(\mathbf{L}^k, \mathbf{S} \,|\, \mathbf{S}^k) \\
&= Q(\mathbf{L}^k, \mathbf{S}^k) + 2\left(\mathrm{vec}(\mathbf{Y}-\mathbf{L}^k) - \mathrm{vec}(\mathbf{S}^k)\right)^T \mathbf{Q} \left(\mathrm{vec}(\mathbf{Y}-\mathbf{L}^k) - \mathrm{vec}(\mathbf{S})\right) \\
&\quad + \left\|\mathrm{vec}(\mathbf{S}) - \mathrm{vec}(\mathbf{S}^k)\right\|^2 + \lambda_s \sum_{i=1}^r \sum_{j=1}^t \widetilde{g}_s(s_{ij} \,|\, s_{ij}^k) \\
&= \left\|\mathbf{S} - \widetilde{\mathbf{Y}}^k - (\lambda_s/2) \cdot \mathbf{W}^k\right\|_F^2 + \lambda_s \eta(\theta) \|\mathbf{S}\|_1 + \mathrm{const.},
\end{aligned} \qquad (67)$$

where $\mathbf{W}^k$ and $\widetilde{\mathbf{Y}}^k$ are matrices of the same size of $\mathbf{S}$, with $(i,j)$-th entries defined as

$$w_{ij}^k \triangleq \left.\frac{dg_s^-(x)}{dx}\right|_{x=s_{ij}^k} \quad \text{and} \quad \widetilde{y}_{ij}^k \triangleq \begin{cases} y_{ij} - \ell_{ij}^k, & \text{if } (i,j) \in \Omega, \\ s_{ij}^k, & \text{otherwise}, \end{cases} \qquad (68)$$

respectively; and in const. we absorbed irrelevant constant terms.

The minimizer of $V^{(1)}(\mathbf{L}^k, \bullet \,|\, \mathbf{S}^k)$ and $V^{(2)}(\mathbf{L}^k, \bullet \,|\, \mathbf{S}^k)$ can be computed following the same steps as described in Option 1 and Option 2 in Sec. I.4.2, respectively.

Next, we only provide the update of $\mathbf{S}$ based on minimizing $V^{(2)}(\mathbf{L}^k, \bullet | \mathbf{S}^k)$, which is given by

$$s_{ij}^{k+1} = \text{sign}\left(\widetilde{y}_{ij} + (\lambda_s/2) \cdot w_{ij}^k\right) \cdot \max\left\{\left|\widetilde{y}_{ij} + (\lambda_s/2) \cdot w_{ij}^k\right| - \eta(\theta)/2, 0\right\}. \quad (69)$$

Next, we fix $\mathbf{S} = \mathbf{S}^{k+1}$ and optimize $\mathbf{L}$. In order to obtain a closed form update of $\mathbf{L}$, we upperbound $Q(\bullet, \mathbf{S}^{k+1})$ and $G_r(\bullet)$ in (65) using (28) and (26), respectively, and following similar steps as to obtain (67). Specifically, a surrogate of $Q(\bullet, \mathbf{S}^{k+1})$ is: given $\mathbf{L}^k$,

$$\widetilde{Q}(\mathbf{L}, \mathbf{S}^{k+1} | \mathbf{L}^k) = \|\mathbf{L} - \mathbf{X}^k\|_F^2 + \text{const.}, \quad (70)$$

where $\mathbf{X}^k$ is a matrix having the same size of $\mathbf{Y}$, with entries defined as

$$x_{ij}^k \triangleq \begin{cases} y_{ij} - s_{ij}^{k+1}, & (i,j) \in \Omega; \\ \ell_{ij}^k, & \text{otherwise.} \end{cases} \quad (71)$$

To upperbound the nonconvex regularizer $G_r(\mathbf{L}) = \sum_{i=1}^m g(\sigma_i(\mathbf{L}))$, we invoke (50) and obtain

$$\widetilde{g}_r\left(\sigma_i(\mathbf{L}) | \sigma_i(\mathbf{L}^k)\right) = \eta(\theta)|\sigma_i(\mathbf{L})| - w_i^k \cdot \left(\sigma_i(\mathbf{L}) - \sigma_i(\mathbf{L}^k)\right), \quad (72)$$

with

$$w_i^k \triangleq \left.\frac{dg_r^-(x)}{dx}\right|_{x=\sigma_i(\mathbf{L}^k)}.$$

Using the directional differentiability of the singular values $\sigma_i(\mathbf{L}^k)$ (see, e.g., [95]) and the chain rule, it is not difficult to check that Assumption I.13 (in particular the directional derivative consistency condition I.13.3) is satisfied for $\widetilde{g}_r$; therefore $\widetilde{g}_r\left(\bullet | \sigma_i(\mathbf{L}^k)\right)$ is a valid surrogate function of $g_r$.

Combining (70) and (72) yields the following surrogate function of $V(\mathbf{L}, \mathbf{S}^{k+1})$:

$$\widetilde{V}\left(\mathbf{L}, \mathbf{S}^{k+1} | \mathbf{L}^k\right) = \|\mathbf{L} - \mathbf{X}^k\|_F^2 + \lambda_r \sum_{i=1}^m \left(\eta|\sigma_i(\mathbf{L})| - w_i^k \sigma_i(\mathbf{L})\right). \quad (73)$$

The final step is computing the minimizer of $\widetilde{V}(\mathbf{L} | \mathbf{L}^k)$. To this end, we first introduce the following lemma [244].

**Lemma I.21 (von Neumann's trace inequality).** *Let $\mathbf{A}$ and $\mathbf{B}$ be two $m \times m$ complex-valued matrices with singular values $\sigma_1(\mathbf{A}) \geq \cdots \geq \sigma_m(\mathbf{A})$ and $\sigma_1(\mathbf{B}) \geq \cdots \geq \sigma_m(\mathbf{B})$, respectively. Then,*

$$|\text{Tr}(\mathbf{AB})| \leq \sum_{i=1}^m \sigma_i(\mathbf{A})\sigma_i(\mathbf{B}). \quad (74)$$

Note that Lemma I.21 can be readily generalized to rectangular matrices. Specifically, given $\mathbf{A}, \mathbf{B}^T \in \mathbb{R}^{m \times t}$, define the augmented square matrices $\widetilde{\mathbf{A}} \triangleq [\mathbf{A}; \mathbf{0}_{(t-m) \times t}]$ and $\widetilde{\mathbf{B}} \triangleq [\mathbf{B}, \mathbf{0}_{t \times (t-m)}]$, respectively. Applying Lemma I.21 to $\widetilde{\mathbf{A}}$ and $\widetilde{\mathbf{B}}$, we get

$$\text{Tr}(\mathbf{A}\mathbf{B}) = \text{Tr}(\widetilde{\mathbf{A}}\widetilde{\mathbf{B}}) \leq \sum_{i=1}^{t} \sigma_i(\widetilde{\mathbf{A}})\sigma_i(\widetilde{\mathbf{B}}) = \sum_{i=1}^{m} \sigma_i(\mathbf{A})\sigma_i(\mathbf{B}), \tag{75}$$

where equality is achieved when $\mathbf{A}$ and $\mathbf{B}^T$ share the same singular vectors. Using (75), we can now derive the closed form of the minimizer of (73).

**Proposition I.22.** *Let $\mathbf{X}^k = \mathbf{U}_X \mathbf{\Sigma}_X \mathbf{V}_X^T$ be the singular value decomposition (SVD) of $\mathbf{X}^k$. The minimizer of $\widetilde{V}\left(\mathbf{L}, \mathbf{S}^{k+1} \,|\, \mathbf{L}^k\right)$ in (73) (and thus the update of $\mathbf{L}$) is given by*

$$\mathbf{L}^{k+1} = \mathbf{U}_X D_{\frac{\eta \lambda_r}{2}}\left(\mathbf{\Sigma}_X + \text{Diag}(\{w_i^k \lambda_r/2\}_{i=1}^m)\right) \mathbf{V}_X^T, \tag{76}$$

*where $D_\alpha(\mathbf{\Sigma})$ denotes a diagonal matrix with the i-th diagonal element equal to $(\mathbf{\Sigma}_{ii} - \alpha)_+$, and $(x)_+ \triangleq \max(0, x)$.*

*Proof.* Expanding the squares we rewrite $\widetilde{V}\left(\mathbf{L}, \mathbf{S}^{k+1} \,|\, \mathbf{L}^k\right)$ as

$$\begin{aligned}
&\widetilde{V}(\mathbf{L}, \mathbf{S}^{k+1} \,|\, \mathbf{L}^k) \\
&= \text{Tr}(\mathbf{L}\mathbf{L}^T) - 2\text{Tr}(\mathbf{L}(\mathbf{X}^k)^T) + \text{Tr}(\mathbf{X}^k(\mathbf{X}^k)^T) + \lambda_r \sum_{i=1}^{m} \left(\eta |\sigma_i(\mathbf{L})| - w_i^k \sigma_i(\mathbf{L})\right) \\
&= \sum_{i=1}^{m} \sigma_i^2(\mathbf{L}) - 2\text{Tr}(\mathbf{L}(\mathbf{X}^k)^T) + \lambda_r \sum_{i=1}^{m} \left(\eta |\sigma_i(\mathbf{L})| - w_i^k \sigma_i(\mathbf{L})\right) + \text{Tr}(\mathbf{X}^k(\mathbf{X}^k)^T).
\end{aligned} \tag{77}$$

To find a minimizer of $\widetilde{V}(\bullet, \mathbf{S}^{k+1} \,|\, \mathbf{L}^k)$, we introduce the SVD of $\mathbf{L} = \mathbf{U}_L \mathbf{\Sigma}_L \mathbf{V}_L^T$, and optimize separately on $\mathbf{U}_L$, $\mathbf{V}_L$, and $\mathbf{\Sigma}_L$, with $\mathbf{\Sigma}_L = \text{diag}(\sigma_1(\mathbf{L}), \ldots, \sigma_m(\mathbf{L}))$.

From (75) we have $\text{Tr}(\mathbf{L}(\mathbf{X}^k)^T) \leq \sum_{i=1}^{m} \sigma_i(\mathbf{L})\sigma_i(\mathbf{X}^k)$, and equality is reached if $\mathbf{U}_L = \mathbf{U}_X$ and $\mathbf{V}_L = \mathbf{V}_X$, respectively; which are thus optimal. To compute the optimal $\mathbf{\Sigma}_L$ let us substitute $\mathbf{U}_L = \mathbf{U}_X$ and $\mathbf{V}_L = \mathbf{V}_X$ in (77) and solve the resulting minimization problem with respect to $\mathbf{\Sigma}_L$:

$$\begin{aligned}
\underset{\{\sigma_i(\mathbf{L})\}_{i=1}^m}{\text{minimize}} \quad & \sum_{i=1}^{m} \sigma_i(\mathbf{L})^2 - 2\sum_{i=1}^{m} \sigma_i(\mathbf{L})\sigma_i(\mathbf{X}^k) + \lambda_r \sum_{i=1}^{m} \left(\eta \sigma_i(\mathbf{L}) - w_i^k \sigma_i(\mathbf{L})\right) \\
\text{subject to} \quad & \sigma_i(\mathbf{L}) \geq 0, \; \forall i = 1, \ldots, m.
\end{aligned} \tag{78}$$

Problem (78) is additively separable. The optimal value of each $\sigma_i(\mathbf{L})$ is

$$\begin{aligned}
\sigma_i^\star(\mathbf{L}) &= \underset{\sigma_i(\mathbf{L}) \geq 0}{\text{argmin}} \; (\sigma_i(\mathbf{L}) - \sigma_i(\mathbf{X}^k) - w_i^k \lambda_r/2 + \eta \lambda_r/2)^2 \\
&= (\sigma_i(\mathbf{X}^k) + w_i^k \lambda_r/2 - \eta \lambda_r/2)_+,
\end{aligned} \tag{79}$$

which completes the proof.

The block MM algorithm solving the matrix completion Problem (63), based on the S-updates (69) and L-update (76), is summarized in Algorithm 4.

**Algorithm 4: Block MM Algorithm for Matrix Completion [cf. (63)]**

**Data** : $\mathbf{L}^0, \mathbf{S}^0 \in \mathbb{R}^{m \times t}$. Set $k = 0$.
(S.1) : If $\mathbf{L}^k$ and $\mathbf{S}^k$ satisfy a termination criterion: STOP;
(S.2) : Alternately optimize $\mathbf{S}$ and $\mathbf{L}$:
    (a) : Update $\mathbf{S}^{k+1}$ as according to (69);
    (b) : Update $\mathbf{L}^{k+1}$ according to (76);
(S.3) : $k \leftarrow k+1$, and go to (S.1).

---

**2) Dictionary Learning.** Given the data matrix $\mathbf{Y} \in \mathbb{R}^{m \times t}$, the dictionary learning (DL) problem consists in finding a basis, the dictionary $\mathbf{D} \in \mathbb{R}^{m \times r}$ (with $r \ll t$), over which $\mathbf{Y}$ can be sparsely represented throughout the coefficients $\mathbf{X} \in \mathbb{R}^{r \times t}$. This problem appears in a wide range of machine learning applications, such as image denoising, video surveillance, face recognition, and unsupervised clustering. We consider the following formulation for the DL problem:

$$\underset{\mathbf{X}, \mathbf{D} \in D}{\text{minimize}} \quad V(\mathbf{D}, \mathbf{X}) \triangleq \underbrace{\|\mathbf{Y} - \mathbf{D}\mathbf{X}\|_F^2}_{F(\mathbf{D}, \mathbf{X})} + \lambda_s G(\mathbf{X}), \tag{80}$$

where $D$ is a convex compact set, bounding the elements of the dictionary so that the optimal solution will not go to infinity due to scaling ambiguity; and $G$ aims at promoting sparsity on $\mathbf{X}$, with $\lambda_s$ being a positive given constant. In the following, we assume that $G(\mathbf{X}) = \sum_{i=1}^r \sum_{j=1}^t g(x_{ij})$, with $g$ being any of the DC functions introduced in (47).

Since $F(\mathbf{D}, \mathbf{X})$ in (80) is biconvex, we can derive an algorithm for Problem (80) based on the block MM algorithm by updating $\mathbf{D}$ and $\mathbf{X}$ alternately.

Given $\mathbf{X} = \mathbf{X}^k$, $F(\mathbf{D}, \mathbf{X}^k)$ is convex in $\mathbf{D}$. A natural choice for a surrogate of $F(\mathbf{D}, \mathbf{X}^k)$ is $F(\mathbf{D}, \mathbf{X}^k)$ itself, that is,

$$\widetilde{F}^{(1)}(\mathbf{D} \,|\, \mathbf{X}^k) \triangleq F(\mathbf{D}, \mathbf{X}^k) = \|\mathbf{Y} - \mathbf{D}\mathbf{X}^k\|_F^2, \tag{81}$$

and update $\mathbf{D}$ solving

$$\mathbf{D}^{k+1} = \underset{\mathbf{D} \in D}{\text{argmin}} \, \|\mathbf{Y} - \mathbf{D}\mathbf{X}^k\|_F^2. \tag{82}$$

Problem (82) is convex, but does not have a closed form solution for a general constraint set $D$. In some special cases, efficient iterative method can be derived to solve (82) by exploiting its structure (see, e.g., [194]). For instance, consider the constraint set $D = \{\mathbf{D} \in \mathbb{R}^{m \times r} : \|\mathbf{D}\|_F^2 \leq \alpha\}$. Writing the KKT conditions of (82) (note that Slater's constraint qualification holds), we get

$$0 \leq \alpha - \|\mathbf{D}\|_F^2 \perp \mu \geq 0, \tag{83a}$$

$$\nabla_{\mathbf{D}} L(\mathbf{D}, \mu) = 0, \tag{83b}$$

where $L(\mathbf{D}, \mu)$ is the Lagragian function, defined as

$$L(\mathbf{D},\mu) = \|\mathbf{Y} - \mathbf{D}\mathbf{X}^k\|_F^2 + \mu(\|\mathbf{D}\|_F^2 - \alpha). \tag{84}$$

For any given $\mu \geq 0$, the solution of (83b) is given by

$$\mathbf{D}(\mu) = \mathbf{Y}(\mathbf{X}^k)^T (\mathbf{X}^k(\mathbf{X}^k)^T + \mu \mathbf{I})^{-1}, \tag{85}$$

where $\mu$ needs to be chosen in order to satisfy the complementarity condition in (83)

$$0 \leq h(\mu) \perp \mu \geq 0, \tag{86}$$

with $h(\mu) \triangleq \alpha - \|\mathbf{D}(\mu)\|_F^2$. Since $h(\bullet)$ is monotone, (86) can be efficiently solved using bisection.

An alternative surrogate function of $F(\mathbf{D},\mathbf{X}^k)$ leading to a closed form solution of the resulting minimization problem can be readily obtained leveraging the Lipschitz continuity of $\nabla_\mathbf{D} F(\mathbf{D},\mathbf{X}^k)$ and using (28), which yields

$$\widetilde{F}^{(2)}(\mathbf{D},\mathbf{X}^k \,|\, \mathbf{D}^k) = 2\operatorname{Tr}\left\{(\mathbf{D}^k \mathbf{X}^k \mathbf{X}^{kT} - \mathbf{Y}\mathbf{X}^{kT})^T (\mathbf{D} - \mathbf{D}^k)\right\} + L\|\mathbf{D} - \mathbf{D}^k\|_F^2 + \text{const.}, \tag{87}$$

where $L \triangleq \lambda_{\max}(\mathbf{S}^k \mathbf{S}^{kT})$ and const. is an irrelevant constant. The update of $\mathbf{D}$ is then given by

$$\begin{aligned}
\mathbf{D}^{k+1} &= P_D\left(\frac{1}{L}(\mathbf{D}^k \mathbf{X}^k \mathbf{X}^{kT} - \mathbf{Y}\mathbf{X}^{kT})\right) \\
&\triangleq \operatorname*{argmin}_{\mathbf{D} \in D} \left\|\mathbf{D} - \left(\mathbf{D}^k - \frac{1}{L}(\mathbf{D}^k \mathbf{X}^k \mathbf{X}^{kT} - \mathbf{Y}\mathbf{X}^{kT})\right)\right\|^2,
\end{aligned} \tag{88}$$

which has a closed form expression for simple constraint sets such as $D = \mathbb{R}_+^{m \times r}$, $D = \{\mathbf{D} \,|\, \mathbf{D} \in \mathbb{R}_+^{m \times r}, \|\mathbf{D}\|_F^2 \leq \alpha\}$, and $D = \{\mathbf{D} \,|\, \|\mathbf{d}_i\|_2^2 \leq \alpha_i, \,\forall i = 1,\ldots,m\}$.

We fix now $\mathbf{D} = \mathbf{D}^{k+1}$, and update $\mathbf{X}$. Problem (80) is separable in the columns of $\mathbf{X}$; the subproblem associated with the $j$-th column of $\mathbf{X}$, denoted by $\mathbf{x}_j$, reads

$$\mathbf{x}_j^{k+1} \in \operatorname*{argmin}_{\mathbf{x}_j} \frac{1}{2}\|\mathbf{y}_j - \mathbf{D}^{k+1}\mathbf{x}_j\|^2 + \lambda_s \sum_{j=1}^t g_s(x_{ij}). \tag{89}$$

Note that Problem (89) is of the same form of (46); therefore, it can be solved using the MM algorithm, based on the surrogate functions derived in Sec. I.4.1 [cf. (51) and (52)].

### I.4.4 Multicast Beamforming

We study the Max-Min Fair (MMF) beamforming problem for single group multicasting [218], where a single base station (BS) equipped with $m$ antennas wishes to transmit common information to a group of $q$ single-antenna users over the same frequency band. The goal of multicast beamforming is to exploit the channels and the spatial diversity offered by the multiple transmit antennas to steer transmitted power towards the group of desired users while limiting interference (leakage) to nearby co-channel users and systems.

Denoting by $\mathbf{w} \in \mathbb{C}^m$ the beamforming vector, the Max-Min beamforming problem reads

$$\begin{aligned}
\underset{\mathbf{w} \in \mathbb{C}^m}{\text{maximize}} \quad & \min_{i=1,\ldots,q} \mathbf{w}^H \mathbf{R}_i \mathbf{w} \\
\text{subject to} \quad & \|\mathbf{w}\|^2 \leq P_T,
\end{aligned} \quad (90)$$

where $P_T$ is the power budget of the BS and $\mathbf{R}_i \in \mathbb{C}^{m \times m}$ is a positive semidefinite matrix modeling the channel between the BS and user $i$. Specifically, $\mathbf{R}_i = \mathbf{h}_i \mathbf{h}_i^H / \sigma_i^2$ if instantaneous Channel State Information (CSI) is assumed, where $\mathbf{h}_i$ is the frequency-flat quasi-static channel vector from the BS to user $i$ and $\sigma_i^2$ is the variance of the zero-mean, wide-sense stationary additive noise at the $i$-th receiver; and $\mathbf{R}_i = \mathbb{E}\{\mathbf{h}_i \mathbf{h}_i^H\} / \sigma_i^2$ represents the spatial correlation matrix if only long-term CSI is available (in the latter case, no special structure for $\mathbf{R}_i$ is assumed).

Problem (90) contains complex variables. One could reformulate the problem into the real domain by using separate variables for the real and imaginary parts of the complex variables, but this approach is not advisable because it does not take advantage of the structure of (real) functions of the complex variables. Following a well-established path in the signal processing and communication communities, here we work directly with complex variables by means of "Wirtinger derivatives". The main advantage of this approach is that we can use the so-called "Wirtinger calculus" to easily compute in practice derivatives of the functions in (90) directly in the complex domain. It can be shown that all the results in this chapter extend to the complex domain when using Wirtinger derivatives instead of classical gradients. Throughout the chapter we will freely use the Wirtinger calculus and refer the reader to [103, 124, 209] for more information on this topic.

We derive now an MM algorithm to solve Problem (90). To do so, we first rewrite (90) in the equivalent minimization form

$$\begin{aligned}
\underset{\mathbf{w} \in \mathbb{C}^m}{\text{minimize}} \quad & \max_{i=1,\ldots,q} -\mathbf{w}^H \mathbf{R}_i \mathbf{w} \\
\text{subject to} \quad & \|\mathbf{w}\|^2 \leq P_T.
\end{aligned}$$

Since $f_i(\mathbf{w}) \triangleq -\mathbf{w}^H \mathbf{R}_i \mathbf{w}$ is a concave function on $\mathbb{C}^m$, it can be majorized by its first order approximation: given $\mathbf{y} \in \mathbb{C}^m$,

$$\widetilde{f}_i(\mathbf{w} \mid \mathbf{y}) = f_i(\mathbf{y}) + 2\operatorname{Re}\left\{(\mathbf{y})^H \mathbf{R}_i (\mathbf{w} - \mathbf{y})\right\}. \quad (91)$$

Using (91) and (32), it is easy to check that the following convex function is a valid surrogate of $V(\mathbf{w}) \triangleq \max_{i=1,\ldots,q} -\mathbf{w}^H \mathbf{R}_i \mathbf{w}$:

$$\widetilde{V}(\mathbf{w} \mid \mathbf{y}) = \max_{i=1,\ldots,q} \widetilde{f}_i(\mathbf{w}; \mathbf{y}). \quad (92)$$

The main iterate of the MM algorithm based on (92) is then given by: given $\mathbf{w}^k$,

$$\mathbf{x}^{k+1} \in \underset{\|\mathbf{w}\|^2 \leq P_T}{\operatorname{argmin}} \left\{ \max_{i=1,\ldots,q} -\mathbf{w}^{kH}\mathbf{R}_i\mathbf{w}^k + 2\operatorname{Re}\{\mathbf{w}^{kH}\mathbf{R}_i(\mathbf{w}-\mathbf{w}^k)\} \right\}. \tag{93}$$

Convergence to d-stationary solutions of (90) is guaranteed by Theorem I.14.

## I.5. Sources and Notes

**A bit of history.** The MM algorithmic framework has a long history that traces back to 1970, when the majorization principle was introduced [179] for descent-based algorithms, using line search. In 1977, the MM principle was applied to multidimensional scaling [62, 65]. Concurrently, its close relative, the (generalized) EM algorithm [66], was proposed by Dempster *et. al.* in 1977 in the context of maximum likelihood estimation with missing observations. The authors proved the monotonicity of the EM scheme showing that the E-step generates an upperbound of the objective function. The idea of successively minimizing an upperbound of the original function appeared in subsequent works, including [63, 90, 91, 134, 188], to name a few; and it was introduced as a general algorithmic framework in [64, 99, 135]. The connection between EM and MM was clarified in [10], where the authors argued that the majorization step (also referred to as "optimization transfer") rather than the missing data is the key ingredient of EM. The EM/MM algorithm has gained significant attention and applied in various fields ever since [108, 159, 254]. A recent tutorial on the MM algorithm along with its applications to problems in signal processing, communications, and machine learning can be found in [228].

Building on the plain MM/EM, as described in Algorithm 1, several generalizations of the algorithm have been developed to improve its convergence speed, practical implementability, as well as scalability. Some representative examples are the following. Both the majorization and minimization steps can be performed inexactly: the global upperbound condition of the surrogate function can be relaxed to be just a local upperbound; and the exact minimizer of the surrogate function can be replaced by one that only decreases the value of the surrogate with respect to the current iterate [66, 131]. MM can also be coupled with line-search schemes to accelerate its convergence [6, 100, 141, 200]. Furthermore, instead of majorizing the objective function on the whole space, the "subspace MM" algorithm constructs majorizers on an iteration-dependent subspace [47–49, 129]. For structured problems whose variables are naturally partitioned in blocks, majorization can be done block-wise to reduce the scale of the subproblems and achieving tighter upperbounds [81]. Sweeping rules of the blocks such as the (essential-)cyclic rule, random-based rule, Gauss-Southwell rule, maximum improvement rule, have been studied in [81, 111, 193]. An incremental MM was proposed in [151] to minimize sum-utilities composed of a large number of cost functions.

**On the convergence of MM/EM.** Due to the intimacy between MM and EM, convergence results of MM cannot be summarized independently from those of EM. Therefore, in the following we will not distinguish between EM and MM. Earlier studies including the proof of monotonicity of the sequence of the objective values along with the characterization of the limit points of the sequence generated my the EM/MM algorithm were presented in [66]. Results were refined in [27, 252], under the assumption that the objective function is differentiable and the iterates lie in the interior of the constraint set: it was shown that, if the surrogate function satisfies some mild conditions, all the limit points of the sequence generated by the EM/MM algorithm are stationary points of the problem. Conditions for the convergence of the whole sequence were given in [241]. A more comprehensive study of MM convergence with extensions to constrained problems and block-wise updates, can be found in [82, 111], where the surrogate only needs to upperbound the objective function locally. Convergence of (block-)MM applied to problems with non-smooth objective functions was studied in [193] (cf. Theorem I.14 and I.20). All the above results considered convex constraints. Convergence of MM under non-convex constraint were only partially investigated; examples include [220] and [59, 157, 263], the latter focusing on specific problems. In many applications, the EM/MM has been observed to be slow [108, 254]. Accelerated versions of EM/MM include: i) [6, 100, 141, 200], based on modifying the step-size; ii) [26, 115, 131, 132, 161], based on adjusting the search direction or inexact computation of the M-step; iii) and [150, 242, 273] based on finding a fixed-point of the algorithmic mapping. We refer the readers to [116, 159, 228] for a comprehensive overview.

**On the choice of surrogate function and related algorithms.** The performance of MM depends crucially on the choice of the surrogate function. On one hand, a desirable surrogate function should be sufficiently "simple" so that the resulting minimization step can be carried out efficiently; on the other hand, it should also preserve the structure of the objective function, which is expected to enhance the practical convergence. Achieving a trade-off between these two goals is often a non-trivial task. Some guidelines on how to construct valid surrogate functions along with several examples were provided in [108, 151, 254]. Quite interestingly, under specific choices of the surrogate function, the resulting (block-)MM algorithm becomes an instance of well-knowns schemes, widely studied in the literature. Examples include the EM algorithm [66], the convex-concave procedure (CCCP) [145, 191, 264], the proximal algorithms [8, 15, 50, 53, 182], the cyclic minimization algorithm [226], and block coordinate descent-based schemes [250]. Finally, the idea of approximating a function using an upperbound has also been adopted to convexify nonconvex constraint functions; examples include the inner approximation algorithm [158], the CCCP procedure, and SCA-based algorithms [77, 207, 208].

**Applications of MM/EM.** In the last few years there has been a growing interest in using the MM framework to solve a gamut of problems in several fields, including signal/image processing, machine learning, communications, and bioinformatics, just to name a few. A non-exhaustive list of specific applications include sparse linear/logistic regression [7, 21, 23, 36, 59, 85, 86, 125, 156], sparse (generalized) principal component analysis [118, 221, 263], matrix factorization/completion [83, 84, 117, 194], phase retrieval [173, 189], edge preserving regularizations in image processing [3, 90, 91], covariance estimation [12, 225, 227, 249, 274], sequence design [220, 253, 270], nonnegative quadratic programming [133, 139, 212], signomial programming [136], and sensor network localization [8, 9, 54, 177].

In the era of big data, the desiderata of MM has steered to low computational complexity and parallel/online computing. This raises new questions and challenges, including i) how to design MM schemes with better convergence rate; ii) how to extend the MM framework to stochastic/online learning problems; and iii) how to design surrogate functions exploiting problem structure and leading to closed form solutions and/or parallel/distributed updates.

## Lecture II – Parallel Successive Convex Approximation Methods

This lecture goes beyond MM-based methods, addressing some limitations and challenges of the MM design. The MM approach calls for the surrogate function $\widetilde{V}$ to be a *global upperbound* of the objective function $V$; this requirement might limit the applicability and the effectiveness of the MM method, for several reasons.

First of all, when $V$ does not have a favorable structure to exploit, it is not easy to build a surrogate $\widetilde{V}$ (possibly convex) that is a (tight) global upper bound of $V$. For instance, suppose that $V$ is twice continuously differentiable and $\widetilde{V}$ is given by (28). A valid choice for $L$ in (28) to meet the upperbound requirement is $L \geq \sup_{\mathbf{x} \in X} \|\nabla^2 V(\mathbf{x})\|_2$. However, computing $\sup_{\mathbf{x} \in X} \|\nabla^2 V(\mathbf{x})\|_2$ for unstructured $V$ is not in general easy. In all such cases, a natural option is leveraging some upper bound of $\sup_{\mathbf{x} \in X} \|\nabla^2 V(\mathbf{x})\|_2$ (e.g., by uniformly bounding the largest eigenvalue of $\nabla^2 V$ on $X$). In practice, however, these bounds can be quite loose (much larger than $\sup_{\mathbf{x} \in X} \|\nabla^2 V(\mathbf{x})\|_2$), resulting in very slow instances of the MM algorithm.

Second, an upper approximation of $V$ might be "too conservative" and not capturing well the "global behaviour" of $V$; this may affect the guarantees of the resulting MM algorithm. Fig. II.1 depicts such a situation: in Fig. II.1a an upper convex approximation is chosen for $\widetilde{V}$ whereas in Fig. II.1b the surrogate function is not an upper bound of $V$ but it shares with $V$ the same gradient at the base point while preserving the "low frequency component" of $V$. As shown in the figure, the two surrogates have different minimizers.

Third, building upper approximations of $V$ that are also (additively) *block separable* is in general a challenging task, making the MM method not suitable for a parallel implementation, which instead is a desirable feature when dealing with large-scale problems. Block updates in MM schemes are possible, but in a sequential (e.g., cyclic) form, as discussed in Sec. I.3 (Lecture I) for the block alternating MM algorithm (Algorithm 2). This contrasts with the intrinsic parallel nature of other algorithms for nonconvex problems, like the (proximal) gradient algorithm.

In this lecture we present a flexible SCA-based algorithmic framework that addresses the above issues. The method hinges on surrogate functions that i) need not be a global upper bound of the objective function but only preserve locally its first order information; ii) are much easier to be constructed than upper approximations; and iii) lead to subproblems that, if constraints permit, are block separable and thus can be solved in parallel. However, the aforementioned surrogates need to be strongly convex, a property that is not required by MM algorithms. Furthermore, to guarantee convergence when the surrogate function is not a global upper bound of the objective function, a step-size is employed in the update of the variables. We begin by first describing a vanilla SCA algorithm in Sec. II.2, where all the blocks are updated in parallel; several choices of the surrogate functions are discussed and convergence of the scheme under different step-size rules is provided. In Sec. II.3, we extend the vanilla algorithm to the case where i) a parallel selective update of the block variables is performed at each iteration–several deterministic and random-based selection rules will be considered–and ii) inexact solutions of the block-subproblems are used. This is motivated by applications where the parallel

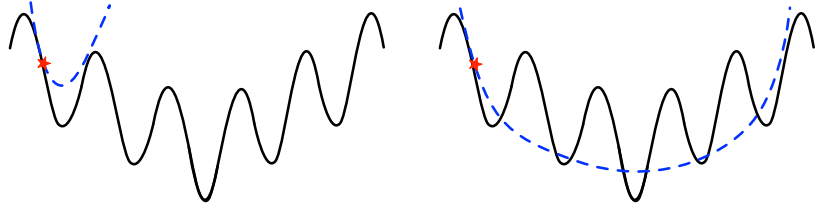

*(a) MM approach: Upper approximation (dotted blue line) of the original function (solid black line) at the base point (red star).*

*(b) SCA approach: Local approximation (dotted blue line) of the original function (solid black line) at the base point (red star).*

Fig. II.1: Upper versus local approximation of the objective function.

update of all blocks and/or the computation of the exact solutions of the subproblems at each iteration is not beneficial or affordable. In Sec. II.4, "hybrid" parallel SCA methods are introduced, which combine deterministic and random-based block selection rules. These schemes have been shown to be very effective when dealing with huge-scale optimization problems. In Sec. II.5, we apply the proposed (parallel) SCA methods to a variety of problems arising from applications in signal processing, machine learning, and communications, and compare their performance with those of existing MM methods. Finally, in Sec. II.7 we overview the main literature and briefly discuss some extensions of the methods described in this lecture.

## II.1. Problem Formulation

We study Problem (1), assuming the following structure for $V$:

$$\underset{\mathbf{x}\in X}{\text{minimize}} \quad V(\mathbf{x}) \triangleq F(\mathbf{x}) + G(\mathbf{x}). \tag{94}$$

**Assumption II.1.** *Given Problem* (94), *we assume that*

1. $X = X_1 \times \cdots X_n$, with each $\emptyset \neq X_i \subseteq \mathbb{R}^{m_i}$ closed and convex;
2. $F: O \to \mathbb{R}$ is $C^1$ on the open set $O \supseteq X$, and $\nabla F$ is L-Lipschitz on $X$;
3. $G: O \to \mathbb{R}$ is convex, possibly nonsmooth;
4. V is bounded from below on $X$.

While Assumption II.1 is slightly more restrictive than Assumption I.12 (it requires $G$ to be convex), it is general enough to cover a gamut of formulations, arising from applications in several fields; some examples are discussed in Sec. II.1.1. On the other hand, under Assumption II.1, we have more flexibility in the choice of the surrogate function and design of parallel algorithms, as it will be discussed shortly.

### II.1.1 Some motivating applications

Many problems in fields as diverse as sensor networks, imaging, machine learning, data analysis, genomics, and geophysics, can be formulated as Problem (94) and satisfy Assumption II.1. Some illustrative examples are documented next; see Sec. II.5 for more details and some numerical results.

**Example #1−LASSO:** Consider a linear regression model with $q$ predictor/feature-response pairs $\{(z_i, \mathbf{a}_i)\}_{i=1}^q$, where $\mathbf{a}_i \in \mathbb{R}^m$ is a $m$-dimensional vector of features or predictors, and $z_i$ is the associated response variable. Let $\mathbf{z} = (z_1, \ldots, z_q)^T$ denote the $q$-dimensional vector of response, and $\mathbf{A} \in \mathbb{R}^{q \times m}$ be the matrix with $\mathbf{a}_i$ in its $i$-th row. Then, the LASSO problem in the so-called Lagrangian form [235], aiming at finding a sparse vector of regression weights $\mathbf{x} \in \mathbb{R}^m$, is an instance of Problem (94), with $F(\mathbf{x}) = \|\mathbf{z} - \mathbf{A}\mathbf{x}\|^2$ and $G(\mathbf{x}) = \lambda \|\mathbf{x}\|_1$, where $\lambda$ is a positive given constant.

**Example #2−Group LASSO:** There are many regression problems wherein the covariates within a group become nonzero (or zero) simultaneously. In such settings, it is natural to select or omit all the coefficient within a group together. The *group* LASSO promotes this structure by using sums of (un-squared) $\ell_2$ penalties. Specifically, consider the regression vector $\mathbf{x} \in \mathbb{R}^m$ possessing the group sparse pattern $\mathbf{x} = [\mathbf{x}_1^T, \ldots, \mathbf{x}_n^T]^T$, i.e., all the elements of $\mathbf{x}_i$ either take large values or close to zero [262]. A widely used group LASSO formulation is the instance of Problem (94), with $F(\mathbf{x}) = \|\mathbf{z} - \mathbf{A}\mathbf{x}\|^2$ and $G(\mathbf{x}) = \lambda \sum_{i=1}^n \|\mathbf{x}_i\|_2$, and $\lambda > 0$.

**Example #3−Sparse logistic regression:** Logistic regression has been popular in biomedical research for half a century, and has recently gained popularity to model a wider range of data. Given the training set $\{\mathbf{z}_i, w_i\}_{i=1}^q$, where $\mathbf{z}_i \in \mathbb{R}^m$ is the feature vector and $w_i \in \{-1, 1\}$ is the label of the $i$-th sample, the logistic regression problem based on a linear logistic model consists in minimizing the negative log likelihood [160, 213] $F(\mathbf{x}) = (1/q) \cdot \sum_{i=1}^q \log(1 + e^{-w_i \cdot \mathbf{z}_i^T \mathbf{x}})$; regularizations can be introduced, e.g., in the form $G(\mathbf{x}) = \lambda \|\mathbf{x}\|_1$ (or $G(\mathbf{x}) = \lambda \sum_{i=1}^n \|\mathbf{x}_i\|_2$), with $\lambda > 0$. Clearly, this is an instance of Problem (94).

**Example #4−Dictionary learning:** Dictionary learning is an unsupervised learning problem that consists in finding a basis $\mathbf{D} \in \mathbb{R}^{q \times m}$–called dictionary–whereby data $\mathbf{z}_i \in \mathbb{R}^q$ can be sparsely represented by coefficients $\mathbf{x}_i \in \mathbb{R}^m$. Let $\mathbf{Z} \in \mathbb{R}^{q \times I}$ and $\mathbf{Z} \in \mathbb{R}^{m \times I}$ be the data and representation matrix whose columns are the data vectors $\mathbf{z}_i$ and coefficients $\mathbf{x}_i$, respectively. The DL problem is the instance of Problem (94), with $F(\mathbf{D}, \mathbf{X}) = \|\mathbf{Z} - \mathbf{D}\mathbf{X}\|_F^2$ and $G(\mathbf{X}) = \lambda \|\mathbf{X}\|_1$, $X = \{(\mathbf{D}, \mathbf{X}) \in \mathbb{R}^{q \times m} \times \mathbb{R}^{m \times I} : \|\mathbf{D}\mathbf{e}_i\|^2 \leq \alpha_i, \forall i = 1, \ldots, m\}$, where $\mathbf{e}_i$ is the $i$-th canonical vector, and $\|\mathbf{X}\|_F$ and $\|\mathbf{X}\|_1$ denote the Frobenius norm and the $\ell_1$ matrix norm of $\mathbf{X}$, respectively. Note that this is an example of $F(\mathbf{D}, \mathbf{X})$ that is not jointly convex in $(\mathbf{D}, \mathbf{X})$, but bi-convex (i.e., convex in $\mathbf{D}$ and $\mathbf{X}$ separately).

**Example #5−(Sparse) empirical risk minimization:** Given a training set $\{D_i\}_{i=1}^I$, the parametric empirical risk minimization problem aims at finding the model $h : \mathbb{R}^m \ni X \to \mathbb{R}^q$, parameterized by $\mathbf{x}$, that minimizes the risk function $\sum_{i=1}^I \ell(h(\mathbf{x}; D_i))$, where $\ell : \mathbb{R}^q \to \mathbb{R}$ is a loss function, measuring the mismatch between the model and the data. This optimization problem is a special case of Problem (94), with $F(\mathbf{x}) = \sum_{i=1}^{I} f_i(\mathbf{x})$ and $f_i(\mathbf{x}) \triangleq \ell\big(h(\mathbf{x}; D_i)\big)$. To promote sparsity, one can add in the objective function a regularizer $G$ using, e.g., any of the surrogates of the $\ell_0$ cardinality function listed in Table I.1 (cf. Sec. I.4.1). By absorbing the smooth $G^-$ part in $F$, the resulting regularized empirical risk minimization problem is still written in the form (94). Note that this general problem contains the previous examples as special cases, and generalizes them by incorporating also nonconvex regularizers.

All the above examples contain separable $G$. Some applications involving non-separable $G$ are discussed next.

**Example #6−Robust linear regression:** Linear least-squares estimates can behave badly when the error distribution is not normal, particularly when the errors are heavy-tailed. One remedy is to remove influential observations from the least-squares fit. Another approach, termed *robust regression*, is to use a fitting criterion that is not as vulnerable as least squares to unusual (outliers) data. Consider the system model as in Example #1; a simple example of robustification is replacing the $\ell_2$ norm loss function with the $\ell_1$ norm, which leads to the instance of Problem (94), with $F(\mathbf{x}) = 0$ and $G(\mathbf{x}) = \|\mathbf{A}\mathbf{x} - \mathbf{v}\|_1$ [96].

**Example #7−The Fermat-Weber problem:** This problem consists in finding $\mathbf{x} \in \mathbb{R}^n$ such that the weighted sum of distances between $\mathbf{x}$ and the $I$ anchors $\mathbf{v}_1, \mathbf{v}_2, \ldots, \mathbf{v}_I$ is minimized [74]. It can be formulated as Problem (94), with $F(\mathbf{x}) = 0$ and $G(\mathbf{x}) = \sum_{i=1}^{I} \omega_i \|\mathbf{A}_i \mathbf{x} - \mathbf{v}_i\|_2$, $X = \mathbb{R}^n$, where $\mathbf{A}_i \in \mathbb{R}^{q \times n}$, $\mathbf{v}_i \in \mathbb{R}^q$, and $\omega_i > 0$ are given constants, for all $i$.

**Example #8−The Total Variation (TV) image reconstruction:** TV minimizing models have become a successful methodology for image processing, including denoising, deconvolution, and restoration, to name a few [40]. The noise-free discrete TV reconstruction problem can be formulated as Problem (94) with $F(\mathbf{X}) = \|\mathbf{Z} - \mathbf{A}\mathbf{X}\|^2$ and $G(\mathbf{X}) = \lambda \cdot \mathrm{TV}(\mathbf{X})$, $X \equiv \mathbb{R}^{m \times m}$, where $\mathbf{A} \in \mathbb{R}^{q \times m}$, $\mathbf{X} \in \mathbb{R}^{m \times m}$, $\mathbf{Z} \in \mathbb{R}^{q \times m}$, and $\mathrm{TV}(\mathbf{X}) \triangleq \sum_{i,j=1}^{m} \|\nabla_{ij}(\mathbf{X})\|_p$ is the discrete total variational seminorm of $\mathbf{X}$, with $p = 1$ or $2$ and $\nabla_{ij}(\mathbf{X})$ being the discrete gradient of $\mathbf{X}$ defined as $\nabla_{ij}(\mathbf{X}) \triangleq [(\nabla_{ij}^{(1)}(\mathbf{X})), (\nabla_{ij}^{(2)}(\mathbf{X}))]$, with

$$\nabla_{ij}^{(1)}(\mathbf{X}) \triangleq \begin{cases} X_{i+1,j} - X_{i,j}, & \text{if } i < m, \\ 0, & i = m; \end{cases} \quad \nabla_{ij}^{(2)}(\mathbf{X})^{(2)} \triangleq \begin{cases} X_{i,j+1} - X_{i,j}, & \text{if } j < m, \\ 0, & j = m; \end{cases}$$

for all $i, j = 1, \ldots, m$.

## II.2. Parallel SCA: Vanilla Version

We begin introducing a vanilla version of the parallel SCA framework wherein all the block variables are updated in parallel; generalizations of this scheme will be considered in Sec. II.3.

The most natural parallel (Jacobi-type) solution method one can employ is solving (94) *blockwise* and *in parallel*: given $\mathbf{x}^k$, all the (block) variables $\mathbf{x}_i$ are updated *simultaneously* by solving the following subproblems

$$\mathbf{x}_i^{k+1} \in \underset{\mathbf{x}_i \in X_i}{\operatorname{argmin}} \left\{ F(\mathbf{x}_i, \mathbf{x}_{-i}^k) + G(\mathbf{x}_i, \mathbf{x}_{-i}^k) \right\}, \qquad \forall i \in N \triangleq \{1, \ldots, n\}. \tag{95}$$

Unfortunately, this method converges only under very restrictive conditions [16] that are seldom verified in practice (even in the absence of the nonsmooth function $G$). Furthermore, the exact computation of $\mathbf{x}_i^{k+1}$ is in general difficult, due to the nonconvexity of $F$. To cope with these two issues, the proposed approach consists in solving for each block instead

$$\widehat{\mathbf{x}}_i(\mathbf{x}^k) \triangleq \underset{\mathbf{x}_i \in X_i}{\operatorname{argmin}} \ \widetilde{F}_i(\mathbf{x}_i \,|\, \mathbf{x}^k) + G(\mathbf{x}_i, \mathbf{x}_{-i}^k), \tag{96}$$

and then setting

$$\mathbf{x}_i^{k+1} = \mathbf{x}_i^k + \gamma^k \left( \widehat{\mathbf{x}}_i(\mathbf{x}^k) - \mathbf{x}_i^k \right). \tag{97}$$

In (96), $\widetilde{F}_i\left(\bullet \,|\, \mathbf{x}^k\right)$ represents a strongly convex surrogate replacing $F(\bullet, \mathbf{x}_{-i}^k)$, and in (97) a step-size $\gamma^k \in (0,1]$ is introduced to control the "length" of the update along the direction $\widehat{\mathbf{x}}_i(\mathbf{x}^k) - \mathbf{x}_i^k$. The step-size is needed if one does not require that the surrogate $\widetilde{F}_i\left(\bullet \,|\, \mathbf{x}^k\right)$ is a global upper bound of $F(\bullet, \mathbf{x}_{-i}^k)$ (as in the MM algorithm).

The surrogate function $\widetilde{F}_i$ has the following properties ($\nabla \widetilde{F}_i$ denotes the partial gradient of $\widetilde{F}_i$ with respect to the first argument).

**Assumption II.2.** *Each function* $\widetilde{F}_i : O_i \times O \to \mathbb{R}$ *satisfies the following conditions:*

1. $\widetilde{F}_i(\bullet \,|\, \mathbf{y})$ *is* $\tau_i$-*strongly convex on* $X_i$, *for all* $\mathbf{y} \in X$;
2. $\widetilde{F}_i(\bullet \,|\, \mathbf{y})$ *is differentiable on* $O_i$ *and* $\nabla_{\mathbf{y}_i} F(\mathbf{y}) = \nabla \widetilde{F}_i(\mathbf{y}_i \,|\, \mathbf{y})$, *for all* $\mathbf{y} \in X$.

Stronger convergence results can be obtained under the following additional assumptions.

**Assumption II.3.** $\nabla \widetilde{F}_i(\mathbf{x}_i \,|\, \bullet)$ *is* $\widetilde{L}_i$-*Lipschitz on* $X$, *for all* $\mathbf{x}_i \in X_i$.

**Assumption II.3\*.** $\widetilde{F}_i(\bullet \,|\, \bullet)$ *is continuous on* $O_i \times O$.

Assumption II.2 states that $\widetilde{F}_i$ should be regarded as a (simple) convex approximation of $F$ at the current iterate $\mathbf{x}^k$ that preserves the first order properties of $F$. Note that, as anticipated, $\widetilde{F}_i$ need not be a global upper bound of $F(\bullet, \mathbf{x}_{-i})$. Furthermore, the above assumptions guarantee that the mapping $\widehat{\mathbf{x}}(\mathbf{x}) \triangleq (\widehat{\mathbf{x}}_i(\mathbf{x}))_{i=1}^n$, with $\widehat{\mathbf{x}}_i : X \to X_i$ defined in (96), enjoys the following properties that are instrumental to prove convergence of the algorithm to stationary solutions of Problem (94).

**Lemma II.4 (Continuity of $\widehat{\mathbf{x}}$).** *Consider Problem* (94) *under Assumption II.1. The following hold:*

(a) *If Assumptions II.2 and II.3\* are satisfied,* $\widehat{\mathbf{x}}(\bullet)$ *is continuous on* $X$;

(b) *If Assumptions II.2 and II.3 are satisfied, and* $G$ *is separable,* $\widehat{\mathbf{x}}(\bullet)$ *is Lipschitz continuous on* $X$.

*Proof.* See Appendix–Sec. II.6.1.

Other properties of the best-response (e.g., in the presence of nonconvex constraints) can be found in [77, 207].

The described algorithm is based on solving in parallel the subproblems in (96), converging thus to fixed-points of the mapping $\widehat{\mathbf{x}}(\bullet)$. It is then natural to ask which relation exists between these fixed points and the (d-)stationary solutions of Problem (94). The following lemma addresses this question.

**Lemma II.5 (On the fixed-points of $\widehat{\mathbf{x}}$).** *Given Problem* (94) *under Assumption II.1, let each $\widetilde{F}_i$ in* (96) *be chosen according to Assumption II.2. The following hold.*

(a) *The set of fixed-points of $\widehat{\mathbf{x}}(\bullet)$ coincides with that of the coordinate-wise d-stationary solutions of* (94);

(b) *If, in addition, G is separable–$G(\mathbf{x}) = \sum_{i=1}^{n} g_i(\mathbf{x}_i)$–then the set of fixed-points coincides with that of d-stationary solutions of* (94).

*Proof.* The proof of statement (b) can be found in [79, Proposition 8]. The proof of statement (a) follows similar steps and thus is omitted. □

To complete the description of the algorithm, we need to specify how to choose the step-size $\gamma^k \in (0,1]$ in (97). Any of the following standard rules can be used.

**Assumption II.6.** *The step-size sequence $\{\gamma^k \in (0,1]\}_{k \in \mathbb{N}_+}$ satisfies any of the following rules:*

1. **Bounded step-size:** $0 < \liminf_{k \to \infty} \gamma^k \leq \limsup_{k \to \infty} \gamma^k < 2c_\tau/L$, where $c_\tau \triangleq \min_{i=1,\ldots,n} \tau_i$ (cf. Assumption II.2.1);
2. **Diminishing step-size:** $\sum_{k=0}^{\infty} \gamma^k = +\infty$ and $\sum_{k=0}^{\infty} (\gamma^k)^2 < +\infty$;
3. **Line search:** let $\alpha, \delta \in (0,1)$, choose $\gamma^k = \delta^{t_k}$, where $t_k$ is the smallest nonnegative integer such that

$$V\left(\mathbf{x}^k + \gamma^k \Delta \widehat{\mathbf{x}}^k\right)$$
$$\leq V(\mathbf{x}^k) + \alpha \cdot \gamma^k \left( \nabla F(\mathbf{x}^k)^T \Delta \widehat{\mathbf{x}}^k + \sum_{i=1}^{n} \left( G(\widehat{\mathbf{x}}_i(\mathbf{x}^k), \mathbf{x}^k_{-i}) - G(\mathbf{x}^k) \right) \right) \quad (98)$$

with $\Delta \widehat{\mathbf{x}}^k \triangleq \widehat{\mathbf{x}}(\mathbf{x}^k) - \mathbf{x}^k$.

The parallel SCA procedure is summarized in Algorithm 5.

**Algorithm 5: Parallel Successive Convex Approximation (p-SCA)**

**Data** : $\mathbf{x}^0 \in X$, $\{\gamma^k \in (0,1]\}_{k \in \mathbb{N}_+}$.
    Set $k = 0$.
(S.1) : If $\mathbf{x}^k$ satisfies a termination criterion: STOP;
(S.2) : For all $i \in N$, solve in parallel

$$\widehat{\mathbf{x}}_i(\mathbf{x}^k) \triangleq \operatorname*{argmin}_{\mathbf{x}_i \in X_i} \widetilde{F}_i(\mathbf{x}_i \,|\, \mathbf{x}^k) + G(\mathbf{x}_i, \mathbf{x}^k_{-i}); \tag{99}$$

(S.3) : Set $\mathbf{x}^{k+1} \triangleq \mathbf{x}^k + \gamma^k \left( \widehat{\mathbf{x}}(\mathbf{x}^k) - \mathbf{x}^k \right)$;
(S.4) : $k \leftarrow k+1$, and go to (S.1).

Convergence of Algorithm 5 is stated below; Theorem II.7 deals with nonseparable $G$ while Theorem II.8 specializes the results to the case of (block) separable $G$. The proof of the theorems is omitted, because they are special cases of more general results (Theorem II.12 and Theorem II.13) that will be introduced in Sec. 6.

**Theorem II.7.** *Consider Problem* (94) *under Assumption II.1. Let* $\{\mathbf{x}^k\}_{k \in \mathbb{N}_+}$ *be the sequence generated by Algorithm 5, with each $\widetilde{F}_i$ chosen according to Assumptions II.2 and II.3\*; and let the step-size $\gamma^k \in (0, 1/n]$, for all $k \in \mathbb{N}_+$. Then, there hold:*

*(a) If $\{\gamma^k\}_{k \in \mathbb{N}_+}$ is chosen according to Assumption II.6.2 (diminishing rule), then*

$$\liminf_{k \to \infty} \|\widehat{\mathbf{x}}(\mathbf{x}^k) - \mathbf{x}^k\| = 0; \tag{100}$$

*(b) If $\{\gamma^k\}_{k \in \mathbb{N}_+}$ is chosen according to Assumption II.6.1 (bounded condition) or Assumption II.6.3 (line-search),*

$$\lim_{k \to \infty} \|\widehat{\mathbf{x}}(\mathbf{x}^k) - \mathbf{x}^k\| = 0. \tag{101}$$

**Theorem II.8.** *Consider Problem* (94) *under Assumption II.1 and $G(\mathbf{x}) = \sum_{i=1}^n g_i(\mathbf{x}_i)$, with each $g_i : O_i \to \mathbb{R}$ being convex (possibly nonsmooth). Let* $\{\mathbf{x}^k\}_{k \in \mathbb{N}_+}$ *be the sequence generated by Algorithm 5, with each $\widetilde{F}_i$ chosen according to Assumption II.2; and let the step-size $\gamma^k \in (0, 1]$, for all $k \in \mathbb{N}_+$. Then, there hold:*

*(a) If $\widetilde{F}_i$ satisfies Assumption II.3\* and $\{\gamma^k\}_{k \in \mathbb{N}_+}$ is chosen according to Assumption II.6.2 (diminishing rule), then (100) holds.*

*(b) Suppose that either one of the following conditions is satisfied:*

  *(i) Each $\widetilde{F}_i$ satisfies Assumption II.3\* and $\{\gamma^k\}_{k \in \mathbb{N}_+}$ is chosen according to Assumption II.6.1 (bounded condition) or Assumption II.6.3 (line-search);*

  *(ii) Each $\widetilde{F}_i$ satisfies Assumption II.3 and $\{\gamma^k\}_{k \in \mathbb{N}_+}$ is chosen according to Assumption II.6.2 (diminishing rule).*

*Then, (101) holds.*

The above theorems establish the following connection between the limit points of the sequence $\{\mathbf{x}^k\}_{k\in\mathbb{N}_+}$ and the stationary points of Problem (94). By (101) and the continuity of $\widehat{\mathbf{x}}(\bullet)$ (cf. Lemma II.4), we infer that every limit point $\mathbf{x}^\infty$ of $\{\mathbf{x}^k\}_{k\in\mathbb{N}_+}$ (if exists) is a fixed point of $\widehat{\mathbf{x}}(\bullet)$ and thus, by Lemma II.5, it is a coordinate-wise d-stationary solution of Problem (94). If, in addition, $G$ is separable, $\mathbf{x}^\infty$ is a d-stationary solution of (94). When (100) holds instead, there exists a subsequence $\{\mathbf{x}^{k_t}\}_{t\in\mathbb{N}_+}$ of $\{\mathbf{x}^k\}_{k\in\mathbb{N}_+}$ such that $\lim_{t\to\infty}\|\widehat{\mathbf{x}}(\mathbf{x}^{k_t})-\mathbf{x}^{k_t}\|=0$, and the aforementioned connection with the (coordinate-wise) stationary solutions of (94) holds for every limit point of such a subsequence. The existence of a limit point of $\{\mathbf{x}^k\}_{k\in\mathbb{N}_+}$ is guaranteed under standard extra conditions on the feasible set $X$–e.g., boundedness–or on the objective function $V$–e.g., coercivity on $X$.

### II.2.1 Discussion on Algorithm 5

Algorithm 5 represents a gamut of parallel solution methods, each of them corresponding to a specific choice of the surrogate functions $\widetilde{F}_i$ and step-size rule. Theorems above provide a unified set of conditions for the convergence of all such schemes. Some representative choices for $\widetilde{F}_i$ and $\gamma^k$ are discussed next.

**On the choice of the surrogate $\widetilde{F}_i$.** Some examples of surrogate functions satisfying Assumption II.2 for specific $F$ are the following.

*1) Block-wise convexity.* Suppose $F(\mathbf{x}_1,\ldots,\mathbf{x}_n)$ is convex in each block $\mathbf{x}_i$ separately (but not necessarily jointly). A natural approximation for such an $F$ exploring its "partial" convexity is: given $\mathbf{y}=(\mathbf{y}_i)_{i=1}^n\in X$,

$$\widetilde{F}(\mathbf{x}\,|\,\mathbf{y})=\sum_{i=1}^n\widetilde{F}_i(\mathbf{x}_i\,|\,\mathbf{y}), \qquad (102)$$

with each $\widetilde{F}_i(\mathbf{x}_i\,|\,\mathbf{y})$ defined as

$$\widetilde{F}_i(\mathbf{x}_i\,|\,\mathbf{y})\triangleq F(\mathbf{x}_i,\mathbf{y}_{-i})+\frac{\tau_i}{2}(\mathbf{x}_i-\mathbf{y}_i)^T\mathbf{H}_i(\mathbf{x}_i-\mathbf{y}_i), \qquad (103)$$

where $\tau_i$ is any positive constant, and $\mathbf{H}_i$ is any $m_i\times m_i$ positive definite matrix (of course, one can always choose $\mathbf{H}_i=\mathbf{I}$). The quadratic term in (103) can be set to zero if $F(\bullet,\mathbf{y}_{-i})$ is strongly convex on $X_i$, for all $\mathbf{y}_{-i}\in X_{-i}\triangleq X_1\times\cdots\times X_{i-1}\times X_{i+1}\times\cdots\times X_n$.

*2) (Proximal) gradient-like approximations.* If no convexity is present in $F$, mimicking proximal-gradient methods, a valid choice of $\widetilde{F}$ is the first order approximation of $F$ (plus a quadratic regularization), that is, $\widetilde{F}$ is given by (102), with each

$$\widetilde{F}_i(\mathbf{x}_i\,|\,\mathbf{y})\triangleq \nabla_{\mathbf{x}_i}F(\mathbf{y})^T(\mathbf{x}_i-\mathbf{y}_i)+\frac{\tau_i}{2}\|\mathbf{x}_i-\mathbf{y}_i\|^2. \qquad (104)$$

Note that the above approximation has the same form of the one used in the MM algorithm [cf. (28)], with the difference that in (104) $\tau_i$ can be any positive number (and not necessarily larger than $L_i$). When $\tau_i<L_i$, $\widetilde{F}_i$ in (104) is no longer a global upper bound of $F(\bullet, \mathbf{x}_{-i})$. In such cases, differently from the MM algorithm, stepsize $\gamma^k = 1$ in Algorithm 5 may not be used at each $k$.

*3) Sum-utility function.* Suppose that $F(\mathbf{x}) \triangleq \sum_{i=1}^{I} f_i(\mathbf{x}_1, \ldots, \mathbf{x}_n)$. This structure arises, e.g., in multi-agent systems wherein $f_i$ is the cost function of agent $i$ that controls its own block variables $\mathbf{x}_i \in X_i$. In many application it is common that the cost functions $f_i$ are convex in some agents' variables (cf. Sec. II.5). To exploit this partial convexity, let us introduce the following set

$$\tilde{C}_i \triangleq \left\{ j : f_j(\bullet, \mathbf{x}_{-i}) \text{ is convex}, \forall \mathbf{x}_{-i} \in X_{-i} \right\}, \tag{105}$$

which represents the set of indices of all the functions $f_j$ that are convex in $\mathbf{x}_i$, for any feasible $\mathbf{x}_{-i}$; and let $C_i \subseteq \tilde{C}_i$ be any subset of $\tilde{C}_i$. Then, the following surrogate function satisfies Assumption II.2 while exploiting the partial convexity of $F$ (if any): given $\mathbf{y} = (\mathbf{y}_i)_{i=1}^{n} \in X$,

$$\widetilde{F}(\mathbf{x}|\mathbf{y}) = \sum_{i=1}^{n} \widetilde{F}_{C_i}(\mathbf{x}_i|\mathbf{y}),$$

with each $\widetilde{F}_{C_i}$ defined as

$$\widetilde{F}_{C_i}(\mathbf{x}_i|\mathbf{y}) \triangleq \sum_{j \in C_i} f_j(\mathbf{x}_i, \mathbf{y}_{-i}) + \sum_{\ell \notin C_i} \nabla_{\mathbf{x}_i} f_\ell(\mathbf{y})^T (\mathbf{x}_i - \mathbf{y}_i) \\ + \frac{\tau_i}{2} (\mathbf{x}_i - \mathbf{y}_i)^T \mathbf{H}_i (\mathbf{x}_i - \mathbf{y}_i), \tag{106}$$

where $\mathbf{H}_i$ is any $m_i \times m_i$ positive definite matrix. Roughly speaking, for each agent $i$, the above approximation function preserves the convex part of $F$ w.r.t. $\mathbf{x}_i$ while it linearizes the nonconvex part.

*4) Product of functions.* Consider an $F$ that is written as the product of functions (see [208] for some examples); without loss of generality, here we study only the case of product of two functions. Let $F(\mathbf{x}) = f_1(\mathbf{x}) f_2(\mathbf{x})$, with $f_1$ and $f_2$ convex and non-negative on $X$; if the functions are just block-wise convex, the proposed approach can be readily extended. In view of the expression of the gradient of $F$ and Assumption II.2.2, $\nabla_{\mathbf{x}} F = f_2 \nabla_{\mathbf{x}} f_1 + f_1 \nabla_{\mathbf{x}} f_2$, it seems natural to consider the following approximation: given $\mathbf{y} \in X$,

$$\widetilde{F}(\mathbf{x}|\mathbf{y}) = f_1(\mathbf{x}) f_2(\mathbf{y}) + f_1(\mathbf{y}) f_2(\mathbf{x}) + \frac{\tau_i}{2} (\mathbf{x} - \mathbf{y})^T \mathbf{H} (\mathbf{x} - \mathbf{y}),$$

where, as usual, $\mathbf{H}$ is a positive definite matrix; this term can be omitted if $f_1$ and $f_2$ are positive on the feasible set and $f_1 + f_2$ is strongly convex (for example if one of the two functions is strongly convex). In case $f_1$ and $f_2$ are still positive but not necessarily convex, we can use the expression

$$\widetilde{F}(\mathbf{x}|\mathbf{y}) = \widetilde{f}_1(\mathbf{x}|\mathbf{y}) f_2(\mathbf{y}) + f_1(\mathbf{y}) \widetilde{f}_2(\mathbf{x}|\mathbf{y}),$$

where $\widetilde{f}_1$ and $\widetilde{f}_2$ are any legitimate surrogates of $f_1$ and $f_2$, respectively. Finally, if $f_1$ and $f_2$ can take nonpositive values, introducing $h_1(\mathbf{x}|\mathbf{y}) \triangleq \widetilde{f}_1(\mathbf{x}|\mathbf{y}) f_2(\mathbf{y})$ and $h_2(\mathbf{x}|\mathbf{y}) \triangleq f_1(\mathbf{y}) \widetilde{f}_2(\mathbf{x}|\mathbf{y})$, one can write

$$\widetilde{F}(\mathbf{x}|\mathbf{y}) = \widetilde{h}_1(\mathbf{x}|\mathbf{y}) + \widetilde{h}_2(\mathbf{x}|\mathbf{y}),$$

where $\widetilde{h}_1$ (resp. $\widetilde{f}_1$) and $\widetilde{h}_2$ (resp. $\widetilde{f}_2$) are legitimate surrogates of $h_1$ (resp. $f_1$) and $h_2$ (resp. $f_2$), respectively. Note that in this last case, we no longer need the quadratic term because it is already included in the approximations $\widetilde{f}_1$ and $\widetilde{f}_2$, and $\widetilde{h}_1$ and $\widetilde{h}_2$, respectively. As the final remark, note that the functions $F$ discussed above belong to a class of nonconvex functions for which it does not seem possible to find a global convex upper bound; therefore, the MM techniques introduced in Lecture I are not readily applicable.

*5) Composition of functions*. Let $F(\mathbf{x}) = h(\mathbf{f}(\mathbf{x}))$, where $h : \mathbb{R}^q \to \mathbb{R}$ is a finite convex smooth function such that $h(u_1, \ldots, u_q)$ is nondecreasing in each component, and $\mathbf{f} : \mathbb{R}^m \to \mathbb{R}^q$ is a smooth mapping, with $\mathbf{f}(\mathbf{x}) = (f_1(\mathbf{x}), \ldots, f_q(\mathbf{x}))^T$ and $f_i$ not necessarily convex. Examples of functions $F$ belonging to such a class are those arising from nonlinear least square-based problems, that is, $F(\mathbf{x}) = \|\mathbf{f}(\mathbf{x})\|^2$, where $\mathbf{f}(\mathbf{x})$ is a smooth nonlinear (possibly) nonconvex map. A convex approximation satisfying Assumption II.2 is: given $\mathbf{y} \in \mathbb{R}^m$,

$$\widetilde{F}(\mathbf{x}|\mathbf{y}) \triangleq h\left(\mathbf{f}(\mathbf{y}) + \nabla \mathbf{f}(\mathbf{y})(\mathbf{x}-\mathbf{y})\right) + \frac{\tau}{2}\|\mathbf{x}-\mathbf{y}\|^2, \tag{107}$$

where $\nabla \mathbf{f}(\mathbf{y})$ denotes the Jacobian of $\mathbf{f}$ at $\mathbf{y}$.

**On the choice of the step-size $\gamma^k$.** Some possible choices for the step-size satisfying Assumption II.6 are the following.

*1) Bounded step-size:* Assumption II.6.1 requires the step-size to be eventually in the interval $[\delta, 2c_\tau/L)$, for some $\delta > 0$. A simple choice is $\gamma^k = \gamma > 0$, with $2c_\tau/\gamma > L$ and for all $k$. This simple (but conservative) condition imposes a constraint only on the ratio $c_\tau/\gamma$, leaving free the choice of one of the two parameters. An interesting case is when the proximal gradient-like approximation in (104) is used. Setting therein each $\tau_i > L/2$ allows one to use step-size $\gamma = 1$, obtaining thus the MM algorithm as a special case.

*2) Diminishing step-size:* In scenarios where the knowledge of system parameters, e.g., $L$, is not available, one can use a diminishing step-size satisfying Assumption II.6.2. Two examples of diminishing step-size rules are:

$$\gamma^k = \gamma^{k-1}\left(1 - \varepsilon \gamma^{k-1}\right), \quad k = 1, \ldots, \quad \gamma^0 < 1/\varepsilon; \tag{108}$$

$$\gamma^k = \frac{\gamma^{k-1} + \alpha(k)}{1 + \beta(k)}, \quad k = 1, \ldots, \quad \gamma^0 = 1; \tag{109}$$

where in (108) $\varepsilon \in (0,1)$ is a given constant, whereas in (109) $\alpha(k)$ and $\beta(k)$ are two nonnegative real functions of $k \geq 1$ such that: i) $0 \leq \alpha(k) \leq \beta(k)$; and ii) $\alpha(k)/\beta(k) \to 0$ as $k \to \infty$ while $\sum_k (\alpha(k)/\beta(k)) = \infty$. Examples of such $\alpha(k)$ and $\beta(k)$ are: $\alpha(k) = \alpha$ or $\alpha(k) = \log(k)^\alpha$, and $\beta(k) = \beta \cdot k$ or $\beta(k) = \beta \cdot \sqrt{k}$, where $\alpha, \beta$ are given constants satisfying $\alpha \in (0,1)$, $\beta \in (0,1)$, and $\alpha \leq \beta$.

*3) Line search:* Assumption II.6.3 is an Armijo-like line-search that employs a backtracking procedure to find the largest $\gamma^k$ generating sufficient descent of the objective function at $\mathbf{x}^k$ along the direction $\Delta \widehat{\mathbf{x}}(\mathbf{x}^k)$. Of course, using a step-size generated by line-search will likely be more efficient in terms of iterations than the one based on diminishing step-size rules. However, as a trade-off, performing line-search requires evaluating the objective function multiple times per iteration; resulting thus in more costly iterations. Furthermore, performing a line-search on a multicore architecture requires some shared memory and coordination among the cores/processors.

## II.3. Parallel SCA: Selective updates

The parallel SCA algorithm introduced in Sec. II.2 consists in updating at each iteration *all* the block variables by computing the solutions of subproblems in the form (96). In this section, we generalize the algorithm by i) unlocking parallel updates of a *subset* of all the blocks at a time, and ii) allowing inexact computations of the solution of each subproblem. This is motivated by applications where computing the exact solutions of large-scale subproblems is computationally demanding and/or updating all the block variables at each iteration is not beneficial; see Sec. II.5 for some examples.

**Inexact solutions:** Subproblems in (96) are solved inexactly by computing $\mathbf{z}^k$ satisfying $\|\mathbf{z}_i^k - \widehat{\mathbf{x}}_i(\mathbf{x}^k)\| \leq \varepsilon_i^k$, where $\varepsilon_i^k$ is the desired accuracy (to be properly chosen). Some conditions on the inexact solution (and thus associated error) are needed to guarantee convergence on the resulting algorithm, as stated next.

**Assumption II.9.** *Given $\widehat{\mathbf{x}}_i(\mathbf{x}^k)$ as defined in* (96)*, the inexact solutions $\mathbf{z}_i^k$ satisfies: for all $i = 1, \ldots, n$,*

1. $\|\mathbf{z}_i^k - \widehat{\mathbf{x}}_i(\mathbf{x}^k)\| \leq \varepsilon_i^k$ *and* $\lim_{k \to \infty} \varepsilon_i^k = 0$;
2. $\widetilde{F}_i(\mathbf{z}_i^k | \mathbf{x}^k) + G(\mathbf{z}_i^k, \mathbf{x}_{-i}^k) \leq \widetilde{F}_i(\mathbf{x}_i^k | \mathbf{x}^k) + G(\mathbf{x}^k)$.

The above conditions are quite natural: Assumption II.9.1 states that the error must asymptotically vanish [subproblems (96) need to be solved with increasing accuracy] while Assumption II.9.2 requires that $\mathbf{z}_i^k$ generates a decrease in the objective function of subproblem (96) at iteration $k$ [$\mathbf{z}_i^k$ need not be a minimizer of (96)].

**Updating only some blocks:** At each iteration $k$ a suitable chosen subset of blocks– say $S^k \subseteq N$ [recall that $N \triangleq \{1, \ldots, n\}$]–is selected and updated by computing for each block $i \in S^k$ an inexact solution $\mathbf{z}_i^k$ of the associated subproblem (96): given $\mathbf{x}^k$ and $S^k$, let

$$\mathbf{x}_i^{k+1} = \begin{cases} \mathbf{x}_i^k + \gamma^k (\mathbf{z}_i^k - \mathbf{x}_i^k), & \text{if } i \in S^k, \\ \mathbf{x}_i^k & \text{if } i \notin S^k. \end{cases}$$

Several options are possible for the block selection rule $S^k$. For instance, one can choose the blocks to update according to some deterministic (cyclic) or random-based rule. Greedy-like schemes–updating at each iteration only the blocks that are "far away" from the optimum–have been shown to be quite effective in some applications. Finally, one can also adopt hybrid rules that properly combine the aforementioned selection methods. For instance, one can first select a subset of blocks uniformly at random, and then within such a pool updating only the blocks resulting from a greedy rule. Of course some minimal conditions on the updating rule are necessary to guarantee convergence, as stated below.

**Assumption II.10.** *The block selection satisfies one of the following rules:*

1. **Essentially cyclic rule**: $S^k$ *is selected so that* $\bigcup_{s=0}^{T-1} S^{k+s} = N$, *for all* $k \in \mathbb{N}_+$ *and some finite* $T > 0$;
2. **Greedy rule**: *Each* $S^k$ *contains at least one index i such that*

$$E_i(\mathbf{x}^k) \geq \rho \max_{j \in N} E_j(\mathbf{x}^k),$$

   *where* $\rho \in (0,1]$ *and* $E_i(\mathbf{x}^k)$ *is an error bound function satisfying*

$$\underline{s}_i \cdot \|\widehat{\mathbf{x}}_i(\mathbf{x}^k) - \mathbf{x}_i^k\| \leq E_i(\mathbf{x}^k) \leq \bar{s}_i \cdot \|\widehat{\mathbf{x}}_i(\mathbf{x}^k) - \mathbf{x}_i^k\|, \tag{110}$$

   *for some* $0 < \underline{s}_i \leq \bar{s}_i < +\infty$;

3. **Random-based rule**: *The sets* $S^k$ *are realizations of independent random sets* $\boldsymbol{S}^k$ *taking value in the power set of N, such that* $\mathbb{P}(i \in \boldsymbol{S}^k) \geq p$, *for all* $i = 1, \ldots, n$ *and* $k \in \mathbb{N}_+$, *and some* $p > 0$.

The above selection rules are quite general and have a natural interpretation. The cyclic rule [Assumption II.10.1] requires that all the blocks are updated (at least) once within $T$ consecutive iterations, where $T$ is an arbitrary (finite) integer. Assumption II.10.2 is a greedy-based rule: only the blocks that are "far" from the optimum need to be updated at each iteration; $E_i(\mathbf{x}^k)$ can be viewed as a local measure of the distance of block $i$ from optimality. The greedy rule in Assumption II.10.2 thus calls for the updates of one block that is within a fraction $\rho$ from the largest distance $E_i(\mathbf{x}^k)$. Some examples of valid error bound functions $E_i$ are discussed in Sec. II.3.1. Finally, Assumption II.10.3 is a random selection rule: blocks can be selected according to any probability distribution as long as they have a positive probability to be picked. Specific rules satisfying Assumption II.10 are discussed in Sec. II.3.1.

The described parallel selective SCA method is summarized in Algorithm 6, and termed "inexact FLEXible parallel sca Algorithm" (FLEXA).

**Algorithm 6: Inexact Flexible Parallel SCA Algorithm (FLEXA)**

**Data** : $\mathbf{x}^0 \in X$, $\{\gamma^k \in (0,1]\}_{k \in \mathbb{N}_+}$, $\varepsilon_i^k \geq 0$, for all $i \in N$ and $k \in \mathbb{N}_+$, $\rho \in (0,1]$.
   Set $k = 0$.
(S.1) : If $\mathbf{x}^k$ satisfies a termination criterion: STOP;
(S.2) : Choose a set $S^k$ according to any of the rules in Assumption II.10;
(S.3) : For all $i \in S^k$, solve (96) with accuracy $\varepsilon_i^k$ :
   find $\mathbf{z}_i^k \in X_i$ s.t. $\|\mathbf{z}_i^k - \widehat{\mathbf{x}}_i(\mathbf{x}^k)\| \leq \varepsilon_i^k$;
   Set $\widehat{\mathbf{z}}_i^k = \mathbf{z}_i^k$ for $i \in S^k$, and $\widehat{\mathbf{z}}_i^k = \mathbf{x}_i^k$ for $i \notin S^k$;
(S.4) : Set $\mathbf{x}^{k+1} \triangleq \mathbf{x}^k + \gamma^k (\widehat{\mathbf{z}}^k - \mathbf{x}^k)$;
(S.5) : $k \leftarrow k+1$, and go to (S.1).

To complete the description of the algorithm, we need to specify how to choose the step-size $\gamma^k$ in Step 4. Assumption II.11 below provides some standard rules.

**Assumption II.11.** *The step-size sequence $\{\gamma^k \in (0,1]\}_{k \in \mathbb{N}_+}$ satisfies any of the following rules:*

1. **Bounded step-size**: $0 < \liminf_{k \to \infty} \gamma^k \leq \limsup_{k \to \infty} \gamma^k < c_\tau / L$, *where* $c_\tau \triangleq \min_{i=1,\ldots n} \tau_i$;
2. **Diminishing step-size**: $\sum_{k=0}^\infty \gamma^k = +\infty$ *and* $\sum_{k=0}^\infty (\gamma^k)^2 < +\infty$. *In addition, if $S^k$ is chosen according to the cyclic rule [Assumption II.10.1], $\gamma^k$ further satisfies $0 < \eta_1 \leq \gamma^{k+1}/\gamma^k \leq \eta_2 < +\infty$, for sufficiently large k, and some $\eta_1 \in (0,1)$ and $\eta_2 \geq 1$;*
3. **Line-search**: *Let $\alpha, \delta \in (0,1)$, choose $\gamma^k = \delta^{t_k}$, where $t_k$ is the smallest nonnegative integer such that*

$$V(\mathbf{x}^k + \gamma^k \Delta \widehat{\mathbf{x}}^k) \leq$$
$$V(\mathbf{x}^k) + \alpha \cdot \gamma^k \left( \nabla F(\mathbf{x}^k)^T \Delta \widehat{\mathbf{x}}^k + \sum_{i \in S^k} \left( G(\mathbf{z}_i^k, \mathbf{x}_{-i}^k) - G(\mathbf{x}^k) \right) \right), \quad (111)$$

*where $\Delta \widehat{\mathbf{x}}^k \triangleq (\widehat{\mathbf{z}}^k - \mathbf{x}^k)$.*

Convergence of Algorithm 6 is stated below and summarized in the flow chart in Fig. II.2. Theorem II.12 applies to settings where the step-size is chosen according to the bounded rule or line-search while Theorem II.13 states convergence under the diminishing step-size rule.

**Theorem II.12.** *Consider Problem* (94) *under Assumption II.1. Let $\{\mathbf{x}^k\}_{k \in \mathbb{N}_+}$ be the sequence generated by Algorithm 6, under the following conditions:*
 *(i) Each surrogate function $\widetilde{F}_i$ satisfies Assumptions II.2-II.3 or II.2-II.3*;*
*(ii) $S^k$ is chosen according to any of the rules in Assumption II.10;*

*(iii) Each inexact solution $\mathbf{z}_i^k$ satisfies Assumption II.9;*

*(iv) $\{\gamma^k\}_{k\in\mathbb{N}_+}$ is chosen according to either Assumption II.11.1 (bounded rule) or Assumption II.11.3 (line-search); in addition, if G is nonseparable, $\{\gamma^k\}_{k\in\mathbb{N}_+}$ also satisfies $\gamma^k \in (0, 1/n]$, for all $k \in \mathbb{N}_+$.*

*Then (101) holds [almost surely if $S^k$ is chosen according to Assumption II.10.3 (random-based rule)].*

**Theorem II.13.** *Consider Problem* (94) *under Assumption II.1. Let $\{\mathbf{x}^k\}_{k\in\mathbb{N}_+}$ be the sequence generated by Algorithm 6, under conditions (i), (ii) and (iii) of Theorem II.12. Suppose that $\{\gamma^k\}_{k\in\mathbb{N}_+}$ is chosen according to Assumption II.11.2 (diminishing rule).*

*Then, (100) holds [almost surely if $S^k$ is chosen according to Assumption II.10.3 (random-based rule)]. Furthermore, if G is separable and the surrogate functions $\widetilde{F}_i$ satisfy Assumption II.3, then also (101) holds (almost surely under Assumption II.10.3).*

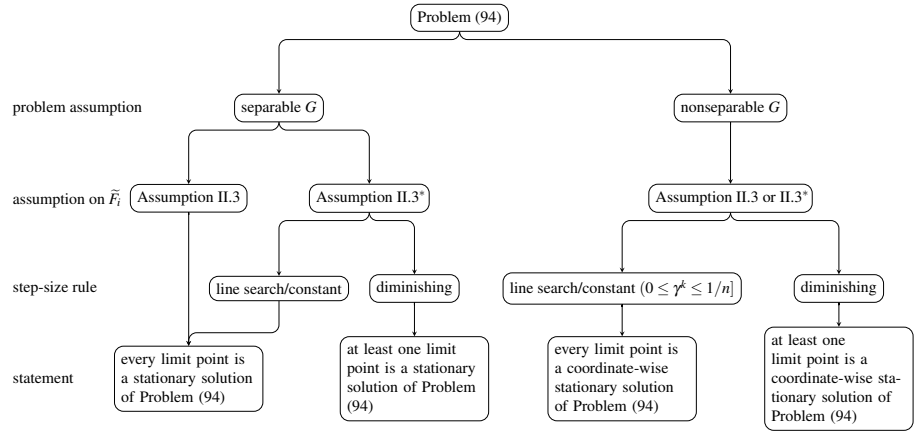

*Fig. II.2: Convergence of FLEXA (Algorithm 6).*

### II.3.1 Discussion on Algorithm 6

The framework described in Algorithm 6 can give rise to very different schemes. We cannot discuss here the entire spectrum of choices; we provide just a few examples of error bound functions $E_i$ and block selection rules $S^k$.

**On the choice of the error bound function** $E_i$. Any function satisfying (110) is a valid candidate for $E_i$. Of course, one can always choose $E_i(\mathbf{x}) = \|\widehat{\mathbf{x}}_i(\mathbf{x}^k) - \mathbf{x}_i^k\|$, corresponding to $\underline{s}_i = \bar{s}_i = 1$ in (110). This is a valuable choice if the computation of $\widehat{\mathbf{x}}_i(\mathbf{x}^k)$ can be easily accomplished. For instance, this is the case in the LASSO problem when the block variables are scalars: $\widehat{x}_i(\mathbf{x}^k)$ can be computed in closed form using the soft-thresholding operator [7]; see Sec. II.5 for details.

In situations where the computation of $\|\widehat{\mathbf{x}}_i(\mathbf{x}^k) - \mathbf{x}_i^k\|$ is not possible or advisable (e.g., when a closed form expression is lacking and the blocks have a large size), one can resort to alternative less expensive metrics satisfying (110). For example, assume momentarily that $G \equiv 0$. Then, it is known [78, Proposition 6.3.1] that, under the stated assumptions, $\|\Pi_{X_i}(\mathbf{x}_i^k - \nabla_{\mathbf{x}_i} F(\mathbf{x}^k)) - \mathbf{x}_i^k\|$ is an error bound for the minimization problem in (96) and therefore it satisfies (110), where $\Pi_{X_i}(\mathbf{y})$ denotes the Euclidean projection of $\mathbf{y}$ onto the closed and convex set $X_i$. In this case, one can choose $E_i(\mathbf{x}^k) = \|\Pi_{X_i}(\mathbf{x}_i^k - \nabla_{\mathbf{x}_i} F(\mathbf{x}^k)) - \mathbf{x}_i^k\|$. If $G(\mathbf{x}) \not\equiv 0$ things become more involved. In several cases of practical interest, adequate error bounds can be derived using [238, Lemma 7].

It is interesting to note that the computation of $E_i$ is only needed if a *partial* update of the (block) variables is performed; otherwise (when $S^k = N$) one can dispense with the computation of $E_i$.

**On the block selection rule $S^k$.** The selection rules satisfying Assumption II.10 are extremely flexible, ranging from deterministic to random-based selection rules. For instance, one can always choose $S^k = N$, resulting in the simultaneous deterministic update of all the (block) variables at each iteration (Algorithm 5). At the other extreme, one can update a single (block) variable per time, thus obtaining a Gauss-Southwell kind of method. Virtually, one can explore all the possibilities "in between", e.g., by choosing properly $S^k$ and leveraging the parameter $\rho$ in (110) to control the desired degree of parallelism or using suitably chosen cyclic-based as well as random-based rules. This flexibility can be coupled with the possibility of computing at each iteration only inexact solutions (Step 3), without affecting the convergence of the resulting scheme (provided that Assumption II.9 is satisfied).

The selection of the most suitable updating rule depends on the specific problem, including the problem scale, computational environment, data acquisition process, as well as the communication among the processors. For instance, versions of Algorithm 6 where all (or most of) the variables are updated at each iteration are particularly amenable to implementation in *distributed* environments (e.g., multi-user communications systems, ad-hoc networks, etc.). In fact, in these settings, not only the calculation of the inexact solutions $\mathbf{z}_i^k$ can be carried out in parallel, but the information that "the $i$-th subproblem" has to exchange with the "other subproblems" in order to compute the next iteration is very limited. A full appreciation of the potentialities of this approach in distributed settings depends however on the specific application under consideration; we discuss some examples in Sec. II.5. The cyclic order has the advantage of being extremely simple to implement. Random selection-based rules are essentially as cheap as cyclic selections while alleviating some of the pitfalls of cyclic updates. They are also relevant in distributed environments wherein data are not available in their entirety, but are acquired either in batches or over a network. In such scenarios, one might be interested in running the optimization at a certain instant even with the limited, randomly available information. A main limitation of random/cyclic selection rules is that they remain disconnected from the status of the optimization process, which instead is exactly the kind of behavior that greedy-based updates try to avoid, in favor of faster convergence, but at the cost of more intensive computation.

We conclude the discussion on the block selection rules providing some specific deterministic and random-based rules that we found effective in our experiments.

• **Deterministic selection**: In addition to the selection rules discussed above, a specific (albeit general) approach is to define first a finite cover $\{S_i\}_{i=1}^{M}$ of $N$ and then update the blocks by selecting the $S_i$'s cyclically. It is also admissible to randomly shuffle the order of the sets $S_i$ before one update cycle.

• **Random-based selection**: The sampling rule $\boldsymbol{S}$ (for notational simplicity the iteration index $k$ will be omitted) is uniquely characterized by the probability mass function
$$\mathbb{P}(S) \triangleq \mathbb{P}(\boldsymbol{S} = S), \quad S \subseteq N,$$
which assign probabilities to the subsets $S$ of $N$. Associated with $\boldsymbol{S}$, define the probabilities $q_j \triangleq \mathbb{P}(|\boldsymbol{S}| = j)$, for $j = 1, \ldots, n$. The following proper sampling rules, proposed in [197] for convex problems with separable $G$, are instances of rules satisfying Assumption II.10.3.

1. *Uniform (U) sampling*: All blocks are selected with the same (non zero) probability:
$$\mathbb{P}(i \in \boldsymbol{S}) = \mathbb{P}(j \in \boldsymbol{S}) = \frac{\mathbb{E}[|\boldsymbol{S}|]}{n}, \quad \forall i \neq j \in N.$$

2. *Doubly Uniform (DU) sampling*: All sets $S$ of equal cardinality are generated with equal probability, i.e., $\mathbb{P}(S) = \mathbb{P}(S')$, for all $S, S' \subseteq N$ such that $|S| = |S'|$. The density function is then
$$\mathbb{P}(S) = \frac{q_{|S|}}{\binom{n}{|S|}}.$$

3. *Nonoverlapping Uniform (NU) sampling*: It is a uniform sampling assigning positive probabilities only to sets forming a partition of $N$. Let $S^1, \ldots, S^p$ be a partition of $N$, with each $|S^i| > 0$, the density function of the NU sampling is:
$$\mathbb{P}(S) = \begin{cases} \dfrac{1}{p}, & \text{if } S \in \{S^1, \ldots, S^p\}; \\ 0 & \text{otherwise;} \end{cases}$$
which corresponds to $\mathbb{P}(i \in \boldsymbol{S}) = 1/p$, for all $i \in N$.

4. *Nice Sampling (NS)*: Given an integer $0 \leq \tau \leq n$, a $\tau$-nice sampling is a DU sampling with $q_\tau = 1$ (i.e., each subset of $\tau$ blocks is chosen with the same probability).

   Using the NS one can control the degree of parallelism of the algorithm by tuning the cardinality $\tau$ of the random sets generated at each iteration, which makes this rule particularly appealing in a multi-core environment. Indeed, one can set $\tau$ equal to the number of available cores/processors, and assign each block coming out from the greedy selection (if implemented) to a dedicated processor/core.

As a final remark, note that the DU/NU rules contain as special cases sequential and fully parallel updates wherein at each iteration a *single* block is updated uniformly at random, or *all* blocks are updated.

5. **Sequential sampling**: It is a DU sampling with $q_1 = 1$, or a NU sampling with $p = n$ and $S^j = j$, for $j = 1, \ldots, p$.
6. **Fully parallel sampling**: It is a DU sampling with $q_n = 1$, or a NU sampling with $p = 1$ and $S^1 = N$.

Other interesting uniform and nonuniform practical rules (still satisfying Assumption II.10) can be found in [196, 197]. Furthermore, see [55, 56] for extensive numerical results comparing the different sampling schemes.

### II.3.2 Convergence analysis of Algorithm 6

In this subsection, we prove convergence of Algorithm 6 (Theorems II.12 and II.13). We consider only deterministic block selection rules (namely Assumption II.10.1 and Assumption II.10.2); the proof under random-based block selection rules follows similar steps and thus is omitted.

### Preliminaries

We first introduce some preliminary technical results that will be used to prove the aforementioned theorems.

**Lemma II.14 (Descent Lemma [14]).** *Let $F : \mathbb{R}^m \to \mathbb{R}$ be continuously differentiable, with L-Lipschitz gradient. Then, there holds:*

$$\left| F(\mathbf{y}) - F(\mathbf{x}) - \nabla F(\mathbf{x})^T (\mathbf{y} - \mathbf{x}) \right| \leq \frac{L}{2} \|\mathbf{y} - \mathbf{x}\|^2, \quad \forall \mathbf{x}, \mathbf{y} \in \mathbb{R}^m. \qquad (112)$$

**Lemma II.15.** *Let $(a_i)_{i=1}^n$ be a n-tuple of nonnegative numbers such that $\sum_{i=1}^n a_i \geq \delta$, with $\delta > 0$. Then, it holds $\sum_{i=1}^n a_i \leq \frac{n}{\delta} \sum_{i=1}^n a_i^2$.*

*Proof.* Define $\mathbf{a} \triangleq [a_1, \ldots, a_n]^T$. The desired result follows readily from $\|\mathbf{a}\|_2^2 \overset{(a)}{\geq} \frac{1}{n} \|\mathbf{a}\|_1^2 \overset{(b)}{\geq} \frac{\delta}{n} \|\mathbf{a}\|_1$, where in (a) we used the Jensen's inequality and $\mathbf{a} \geq 0$ while (b) is due to $\|\mathbf{a}\|_1 \geq \delta$. $\square$

**Lemma II.16** ([17, Lemma 1]). *Let $\{Y^k\}_{k \in \mathbb{N}_+}$, $\{W^k\}_{k \in \mathbb{N}_+}$, and $\{Z^k\}_{k \in \mathbb{N}_+}$ be three sequences such that $W^k$ is nonnegative for all $k$. Assume that*

$$Y^{k+1} \leq Y^k - W^k + Z^k, \quad k = 0, 1, \ldots, \qquad (113)$$

*and that the series $\sum_{k=0}^T Z^k$ converges as $T \to \infty$. Then either $Y^k \to -\infty$, or else $Y^k$ converges to a finite value and $\sum_{k=0}^\infty W^k < \infty$.*

**Lemma II.17.** *Let $\{\mathbf{x}^k\}_{k\in\mathbb{N}_+}$ be the sequence generated by Algorithm 6, with each $\gamma^k \in (0, 1/n)$. For every $k \in \mathbb{N}_+$ and $S^k \subseteq N$, there holds:*

$$G(\mathbf{x}^{k+1}) - G(\mathbf{x}^k) \leq \gamma^k \sum_{i \in S^k} \left( G(\mathbf{z}_i^k, \mathbf{x}_{-i}^k) - G(\mathbf{x}^k) \right), \tag{114}$$

*where $\mathbf{z}_i^k$ is the inexact solution defined in Step 3 of the algorithm.*

*Furthermore, if G is separable, we have: for $\gamma^k \in (0,1)$,*

$$G(\mathbf{x}^{k+1}) - G(\mathbf{x}^k) \leq \gamma^k \sum_{i \in S^k} \left( g_i(\mathbf{z}_i^k) - g_i(\mathbf{x}_i^k) \right). \tag{115}$$

*Proof.* See Appendix–Sec. II.6.2. □

**Lemma II.18.** *Under Assumptions II.1, II.2, and II.9, the inexact solution $\mathbf{z}_i^k$ satisfies*

$$\nabla_{\mathbf{x}_i} F(\mathbf{x}^k)^T (\mathbf{z}_i^k - \mathbf{x}_i^k) + G(\mathbf{z}_i^k, \mathbf{x}_{-i}^k) - G(\mathbf{x}^k) \leq -\frac{\tau_i}{2} \|\mathbf{z}_i^k - \mathbf{x}_i^k\|^2. \tag{116}$$

*Proof.* The proof follows readily from the strong convexity of $\widetilde{F}_i$ and Assumption II.9.2. □

**Lemma II.19.** *Let $S^k$ be selected according to the greedy rule (cf. Assumption II.10.2). Then, there exists a constant $0 < \tilde{c} \leq 1$ such that*

$$\left\| \left( \widehat{\mathbf{x}}(\mathbf{x}^k) - \mathbf{x}^k \right)_{S^k} \right\| \geq \tilde{c} \left\| \widehat{\mathbf{x}}(\mathbf{x}^k) - \mathbf{x}^k \right\|. \tag{117}$$

*Proof.* See Appendix–Sec. II.6.3. □

**Proposition II.20.** *Let $\{\mathbf{x}^k\}_{k\in\mathbb{N}_+}$ be the sequence generated by Algorithm 6, in the setting of Theorem II.12 or Theorem II.13. The following hold [almost surely, if $S^k$ is chosen according to Assumption II.10.3 (random-based rule)]:*

*(a)*
$$\sum_{k=0}^{\infty} \gamma^k \left\| \widehat{\mathbf{z}}^k - \mathbf{x}^k \right\|^2 < +\infty; \tag{118}$$

*(b)*
$$\lim_{k \to \infty} \|\mathbf{x}^{k+1} - \mathbf{x}^k\| = 0. \tag{119}$$

*Proof.* Without loss of generality, we consider next only the case of nonseparable $G$. By the descent lemma (cf. Lemma 114) and Steps 3-4 of the algorithm, we have

$$F(\mathbf{x}^{k+1}) \leq F(\mathbf{x}^k) + \gamma^k \nabla F(\mathbf{x}^k)^T (\widehat{\mathbf{z}}^k - \mathbf{x}^k) + \frac{(\gamma^k)^2 L}{2} \|\widehat{\mathbf{z}}^k - \mathbf{x}^k\|^2.$$

Consider the case of $G$ nonseparable. We have:

$$V(\mathbf{x}^{k+1}) = F(\mathbf{x}^{k+1}) + G(\mathbf{x}^{k+1})$$

$$\stackrel{(114)}{\leq} V(\mathbf{x}^k) + \gamma^k \nabla F(\mathbf{x}^k)^T (\widehat{\mathbf{z}}^k - \mathbf{x}^k) + \frac{(\gamma^k)^2 L}{2} \|\widehat{\mathbf{z}}^k - \mathbf{x}^k\|^2$$

$$+ \gamma^k \sum_{i \in S^k} \left( G(\mathbf{z}_i^k, \mathbf{x}_{-i}^k) - G(\mathbf{x}^k) \right) \tag{120}$$

$$\stackrel{(116)}{\leq} V(\mathbf{x}^k) - \frac{\gamma^k}{2} (c_\tau - \gamma^k L) \|\widehat{\mathbf{z}}^k - \mathbf{x}^k\|^2. \tag{121}$$

If $\gamma^k$ satisfies either Assumption II.11.1 (bounded rule) or Assumption II.11.2 (diminishing rule), statement (a) of the proposition comes readily from Lemma II.16 and Assumption II.1.4.

Consider now the case where $\gamma^k$ is chosen according to Assumption II.11.3 (line search). First of all, we prove that there exists a suitable $\gamma^k \in (0, \gamma^0]$, with $\gamma^0 \in (0, 1/n]$ (if $G$ is separable $\gamma^0 \in (0, 1]$), such that the Armijo-like condition (111) holds. By (120), the line-search condition (111) is satisfied if

$$\gamma^k \left( \nabla F(\mathbf{x}^k)^T (\widehat{\mathbf{z}}^k - \mathbf{x}^k) + \sum_{i \in S^k} \left( G(\mathbf{z}_i^k, \mathbf{x}_{-i}^k) - G(\mathbf{x}^k) \right) \right) + \frac{(\gamma^k)^2 L}{2} \|\widehat{\mathbf{z}}^k - \mathbf{x}^k\|^2$$

$$\leq \alpha \cdot \gamma^k \left( \nabla F(\mathbf{x}^k)^T (\widehat{\mathbf{z}}^k - \mathbf{x}^k) + \sum_{i \in S^k} \left( G(\mathbf{z}_i^k, \mathbf{x}_{-i}^k) - G(\mathbf{x}^k) \right) \right),$$

which, rearranging the terms, yields

$$\frac{\gamma^k \cdot L}{2} \|\widehat{\mathbf{z}}^k - \mathbf{x}^k\|^2$$

$$\leq -(1-\alpha) \left( \nabla F(\mathbf{x}^k)^T (\widehat{\mathbf{z}}^k - \mathbf{x}^k) + \sum_{i \in S^k} \left( G(\mathbf{z}_i^k, \mathbf{x}_{-i}^k) - G(\mathbf{x}^k) \right) \right). \tag{122}$$

Since (cf. Lemma II.18)

$$\|\widehat{\mathbf{z}}^k - \mathbf{x}^k\|^2 \leq -\frac{2}{c_\tau} \cdot \left( \nabla F(\mathbf{x}^k)^T (\widehat{\mathbf{z}}^k - \mathbf{x}^k) + \sum_{i \in S^k} \left( G(\mathbf{z}_i^k, \mathbf{x}_{-i}^k) - G(\mathbf{x}^k) \right) \right),$$

inequality (122) [and thus (111)] is satisfied by any $\gamma^k \leq \min\left\{ \gamma^0, \frac{c_\tau(1-\alpha)}{L} \right\} < +\infty$.

We show next that $\gamma^k$ obtained by (111) is uniformly bounded away from zero. This is equivalent to show $t_k < +\infty$. Without loss of generality we consider $t_k \geq 1$ [otherwise (111) is satisfied by $\gamma^k = \gamma^0$]. Since $t_k$ is the smallest positive integer such that (111) holds with $\gamma^k = \delta^{t_k}$, it must be that the same inequality is not satisfied by $\gamma^k = \delta^{t_k-1}$. Consequently, it must be $\delta^{t_k-1} > \frac{c_\tau(1-\alpha)}{L}$, and thus

$$\gamma^k \geq \min\left\{ \gamma^0, \frac{c_\tau(1-\alpha)}{L} \cdot \delta \right\}. \tag{123}$$

Using (123) in (121), we obtain

$$V(\mathbf{x}^{k+1}) \leq V(\mathbf{x}^k) - \beta_2 \|\widehat{\mathbf{z}}^k - \mathbf{x}^k\|^2, \quad \forall k \in \mathbb{N}_+, \tag{124}$$

where $\beta_2 > 0$ is some finite constant. The rest of the proof follows the same arguments used to prove the statements of the proposition from (121).

We prove now statement (b). By Step 4 of Algorithm 6, it suffices to show that

$$\lim_{k \to \infty} \gamma^k \|\widehat{\mathbf{z}}^k - \mathbf{x}^k\| = 0.$$

Using (118), we have

$$\lim_{k \to \infty} \gamma^k \|\widehat{\mathbf{z}}^k - \mathbf{x}^k\|^2 = 0.$$

Since $\gamma^k \in (0, 1]$, it holds

$$\lim_{k \to \infty} \left(\gamma^k \|\widehat{\mathbf{z}}^k - \mathbf{x}^k\|\right)^2 \leq \lim_{k \to \infty} \gamma^k \|\widehat{\mathbf{z}}^k - \mathbf{x}^k\|^2 = 0. \tag{125}$$

This completes the proof. □

**Proof of Theorem II.12**

We prove (101) for each of the block selection rules in Assumption II.10 separately.

● **Essentially cyclic rule [Assumption II.10.1]:** We start bounding $\|\widehat{\mathbf{x}}(\mathbf{x}^k) - \mathbf{x}^k\|$ as follows:

$$\begin{aligned}
&\|\widehat{\mathbf{x}}(\mathbf{x}^k) - \mathbf{x}^k\| \\
&\leq \sum_{i=1}^n \|\widehat{\mathbf{x}}_i(\mathbf{x}^k) - \mathbf{x}_i^k\| \\
&\leq \sum_{i=1}^n \left( \|\widehat{\mathbf{x}}_i(\mathbf{x}^k) - \widehat{\mathbf{x}}_i(\mathbf{x}^{k+s_i^k})\| + \|\widehat{\mathbf{x}}_i(\mathbf{x}^{k+s_i^k}) - \mathbf{x}_i^{k+s_i^k}\| + \|\mathbf{x}_i^{k+s_i^k} - \mathbf{x}_i^k\| \right)
\end{aligned} \tag{126}$$

where $s_i^k \triangleq \min\{t \in \{1, \ldots, T\} \mid i \in S^{k+t}\}$, so that $k + s_i^k$ is the first time that block $i$ is selected (updated) since iteration $k$. Note that $1 \leq s_i^k \leq T$, for all $i \in N$ and $k \in \mathbb{N}_+$ (due to Assumption II.10.1). We show next that the three terms on the RHS of (126) are asymptotically vanishing, which proves (101).

Since $1 \leq s_i^k \leq T$, we can write

$$\lim_{k \to \infty} \|\mathbf{x}^{k+s_i^k} - \mathbf{x}^k\| \leq \lim_{k \to \infty} \sum_{j=1}^T \|\mathbf{x}^{k+j} - \mathbf{x}^{k+j-1}\| \stackrel{(119)}{=} 0, \tag{127}$$

which, by the continuity of $\widehat{\mathbf{x}}(\bullet)$ (cf. Lemma II.4), leads also to

$$\lim_{k\to\infty} \|\widehat{\mathbf{x}}_i(\mathbf{x}^k) - \widehat{\mathbf{x}}_i(\mathbf{x}^{k+s_i^k})\| = 0. \tag{128}$$

Let $T_i \subseteq \mathbb{N}_+$ be the set of iterations at which block $i$ is updated. It follows from Assumption II.10.1 that $|T_i| = +\infty$, for all $i \in N$. This together with (118) implies

$$\sum_{k \in T_i} \|\mathbf{z}_i^k - \mathbf{x}_i^k\|^2 < +\infty, \quad \forall i \in N,$$

and thus

$$\lim_{T_i \ni k \to \infty} \left\|\mathbf{z}_i^k - \mathbf{x}_i^k\right\| = 0 \quad \Rightarrow \quad \lim_{k\to\infty} \left\|\mathbf{z}_i^{k+s_i^k} - \mathbf{x}_i^{k+s_i^k}\right\| = 0, \quad \forall i \in N.$$

Therefore,

$$\lim_{k\to\infty} \left\|\widehat{\mathbf{x}}_i(\mathbf{x}^{k+s_i^k}) - \mathbf{x}_i^{k+s_i^k}\right\|$$
$$\leq \lim_{k\to\infty} \left\|\widehat{\mathbf{x}}_i(\mathbf{x}^{k+s_i^k}) - \mathbf{z}_i^{k+s_i^k}\right\| + \lim_{k\to\infty} \left\|\mathbf{z}_i^{k+s_i^k} - \mathbf{x}_i^{k+s_i^k}\right\| = 0, \quad \forall i \in N. \tag{129}$$

Combining (126) with (127)-(129) and invoking again the continuity of $\widehat{\mathbf{x}}(\bullet)$, we conclude that $\lim_{k\to\infty} \|\widehat{\mathbf{x}}(\mathbf{x}^k) - \mathbf{x}^k\| = 0$.

• **Greedy rule [Assumption II.10.2]:** We have

$$\tilde{c} \|\widehat{\mathbf{x}}(\mathbf{x}^k) - \mathbf{x}^k\| \overset{(117)}{\leq} \|(\widehat{\mathbf{x}}(\mathbf{x}^k) - \mathbf{x}^k)_{S^k}\|$$
$$\leq \|(\widehat{\mathbf{z}}^k - \mathbf{x}^k)_{S^k}\| + \|(\widehat{\mathbf{x}}(\mathbf{x}^k) - \widehat{\mathbf{z}}^k)_{S^k}\| \tag{130}$$
$$\overset{A.II.9.1}{\leq} \|\widehat{\mathbf{z}}^k - \mathbf{x}^k\| + \sum_{i \in S^k} \varepsilon_i^k \xrightarrow[k\to\infty]{} 0,$$

where the last implication comes from $\lim_{k\to\infty} \varepsilon_i^k = 0$, for all $i \in N$–cf. Assumption II.9.1–and $\lim_{k\to\infty} \|\widehat{\mathbf{z}}^k - \mathbf{x}^k\| = 0$–due to Proposition II.20 and the fact that $\gamma^k$ is bounded away from zero, when the step-size is chosen according to Assumption II.11.1 or Assumption II.11.3 [cf. (123)]. This proves (101).

**Proof of Theorem II.13**

Consider now the diminishing step-size rule [Assumption II.11.2].

**1) Proof of (100):** $\liminf_{k\to\infty} \|\widehat{\mathbf{x}}(\mathbf{x}^k) - \mathbf{x}^k\| = 0$. By Proposition II.20 and the step-size rule, we have $\liminf_{k\to\infty} \|\widehat{\mathbf{z}}^k - \mathbf{x}^k\| = 0$, for all choices of $S^k$. We proceed considering each of the block selection rules in Assumption II.10 separately.

• **Essentially cyclic rule [Assumption II.10.1]:** For notational simplicity, let us assume that $S^k$ is a singleton, that is, $S^k = \{i^k\}$, where $i^k$ denotes the index of the block selected at iteration $k$. The proof can be readily extended to the general case $|S^k| > 1$.

We have

$$\liminf_{k\to\infty} \|\widehat{\mathbf{x}}(\mathbf{x}^k) - \mathbf{x}^k\| \stackrel{(a)}{\leq} \liminf_{r\to\infty} \|\widehat{\mathbf{x}}(\mathbf{x}^{rT}) - \mathbf{x}^{rT}\| \leq \liminf_{r\to\infty} \sum_{i=1}^{n} \|\widehat{\mathbf{x}}_i(\mathbf{x}^{rT}) - \mathbf{x}_i^{rT}\|$$

$$\leq \liminf_{r\to\infty} \sum_{i=1}^{n} \left( \left\|\widehat{\mathbf{x}}_i(\mathbf{x}^{rT+s_i^{rT}}) - \mathbf{x}_i^{rT+s_i^{rT}}\right\| + \left\|\mathbf{x}_i^{rT+s_i^{rT}} - \mathbf{x}_i^{rT}\right\| + \left\|\widehat{\mathbf{x}}_i(\mathbf{x}^{rT}) - \widehat{\mathbf{x}}_i(\mathbf{x}^{rT+s_i^{rT}})\right\| \right)$$

$$= \liminf_{r\to\infty} \sum_{i=1}^{n} \left\|\widehat{\mathbf{x}}_i(\mathbf{x}^{rT+s_i^{rT}}) - \mathbf{x}_i^{rT+s_i^{rT}}\right\|$$

$$+ \underbrace{\lim_{r\to\infty} \sum_{i=1}^{n} \left\|\mathbf{x}_i^{rT+s_i^{rT}} - \mathbf{x}_i^{rT}\right\|}_{\stackrel{(127)}{=}0} + \underbrace{\lim_{r\to\infty} \sum_{i=1}^{n} \left\|\widehat{\mathbf{x}}_i(\mathbf{x}^{rT}) - \widehat{\mathbf{x}}_i(\mathbf{x}^{rT+s_i^{rT}})\right\|}_{\stackrel{(128)}{=}0}$$

$$\leq \liminf_{r\to\infty} \sum_{i=1}^{n} \left\|\mathbf{z}_i^{rT+s_i^{rT}} - \mathbf{x}_i^{rT+s_i^{rT}}\right\| + \underbrace{\lim_{r\to\infty} \sum_{i=1}^{n} \left\|\widehat{\mathbf{x}}_i(\mathbf{x}^{rT+s_i^{rT}}) - \mathbf{z}_i^{rT+s_i^{rT}}\right\|}_{\stackrel{(A.II.9.1)}{\leq} \lim_{r\to\infty} \sum_{i=1}^{n} \varepsilon_i^{rT+s_i^{rT}} \stackrel{(A.II.9.1)}{=} 0}$$

$$\leq \liminf_{r\to\infty} \sum_{k=rT+1}^{(r+1)T} \left\|\mathbf{z}_{i^k}^k - \mathbf{x}_{i^k}^k\right\|, \tag{131}$$

where (a) follows from the fact that the infimum of a subsequence is larger than that of the original sequence.

To complete the proof, we show next that the term on the RHS of (131) is zero. Recalling that if the cyclic block selection rule is implemented, the diminishing step-size $\gamma^k$ is assumed to further satisfy $\eta_1 \leq \gamma^{k+1}/\gamma^k \leq \eta_2$, with $\eta_1 \in (0,1)$ and $\eta_2 \geq 1$ (cf. Assumption II.11.2), we have

$$+\infty \stackrel{(118)}{\geq} \lim_{k\to\infty} \sum_{t=1}^{k} \gamma^t \|\widehat{\mathbf{z}}^t - \mathbf{x}^t\|^2 = \lim_{k\to\infty} \sum_{t=1}^{k} \gamma^t \|\mathbf{z}_{i^t}^t - \mathbf{x}_{i^t}^t\|^2$$

$$= \lim_{k\to\infty} \sum_{r=0}^{k} \sum_{t=rT+1}^{(r+1)T} \gamma^t \|\mathbf{z}_{i^t}^t - \mathbf{x}_{i^t}^t\|^2 \tag{132}$$

$$\geq (\eta_1)^{T-1} \lim_{k\to\infty} \sum_{r=0}^{k} \gamma^{rT+1} \sum_{t=rT+1}^{(r+1)T} \|\mathbf{z}_{i^t}^t - \mathbf{x}_{i^t}^t\|^2,$$

where in the last inequality we used $\gamma^{k+1}/\gamma^k \geq \eta_1$. Since

$$+\infty = \lim_{k\to\infty} \sum_{t=1}^{k} \gamma^t = \lim_{k\to\infty} \sum_{r=0}^{k} \sum_{t=rT+1}^{(r+1)T} \gamma^t \leq T \cdot (\eta_2)^{T-1} \lim_{k\to\infty} \sum_{r=0}^{k} \gamma^{rT+1},$$

it follows from (132) that

$$\liminf_{r\to\infty} \sum_{t=rT+1}^{(r+1)T} \|\mathbf{z}_{i^t}^t - \mathbf{x}_{i^t}^t\|^2 = 0, \tag{133}$$

which, combined with (131), proves the desired result.

• **Greedy rule [Assumption II.10.2]:** Taking the liminf on both sides of (130) leads to the desired result.

**2) Proof of (101):** $\limsup_{k\to\infty} \|\widehat{\mathbf{x}}(\mathbf{x}^k) - \mathbf{x}^k\| = 0$. Recall that in this setting, $\widehat{\mathbf{x}}(\bullet)$ is Lipschitz continuous with constant $\hat{L}$. We prove the result for the essentially cyclic rule and greedy rule separately.

• **Essentially cyclic rule [Assumption II.10.1]:** As in the proof of Theorem II.13, let us assume w.l.o.g. that $S^k = \{i^k\}$, where $i^k$ denotes the index of the block selected at iteration $k$.

By (126)-(129), it is sufficient to prove $\limsup_{k\to\infty} \|\mathbf{z}_i^{k+s_i^k} - \mathbf{x}_i^{k+s_i^k}\| = 0$, for all $i \in N$. Since

$$\limsup_{k\to\infty} \|\mathbf{z}_i^{k+s_i^k} - \mathbf{x}_i^{k+s_i^k}\| \leq \limsup_{k\to\infty} \underbrace{\sum_{t=kT+1}^{(k+1)T} \|\mathbf{z}_{i^t}^t - \mathbf{x}_{i^t}^t\|}_{\triangleq \Delta^k} \tag{134}$$

we prove next $\limsup_{k\to\infty} \Delta^k = 0$.

Assume on the contrary that $\limsup_{k\to\infty} \Delta^k > 0$. Since $\liminf_{k\to\infty} \Delta^k = 0$ [cf. (133)], there exists a $\delta > 0$ such that $\Delta^k < \delta$ for infinitely many $k$ and also $\Delta^k > 2 \cdot T \delta$ for infinitely many $k$. Therefore, there exists a set $K \subseteq \mathbb{N}_+$, with $|K| = \infty$, such that for each $k \in K$, one can find an integer $j_k > k$ such that

$$\Delta^k \geq 2 \cdot \delta \cdot T, \quad \Delta^{j_k} \leq \delta \tag{135}$$
$$\delta < \Delta^\ell < 2 \cdot \delta \cdot T, \quad \text{if } k < \ell < j_k. \tag{136}$$

Define the following quantities: for any $k \in K$, let

$$T_i^k \triangleq \{r \in \{kT+1,\ldots,(k+1)T\} \mid i^r = i\} \quad \text{and} \quad t_i^k \triangleq \min T_i^k. \tag{137}$$

Note that $T_i^k$ (resp. $t_i^k$) is the set of (iteration) indices (resp. the smallest index) within $[kT+1,(k+1)T]$ at which the block index $i$ is selected. Because of Assumption II.10.1, it must be $1 \leq |T_i^k| \leq T$, for all $k$, where $|T_i^k|$ is the number of times block $i$ has been selected in the iteration window $[kT+1,(k+1)T]$. Then we have

$$\delta \cdot T = 2 \cdot \delta \cdot T - \delta \cdot T$$
$$\leq \Delta^k - T \cdot \Delta^{j_k} = \sum_{r=kT+1}^{(k+1)T} \|\mathbf{z}_{i^r}^r - \mathbf{x}_{i^r}^r\| - T \sum_{r=j_k \cdot T+1}^{(j_k+1)T} \|\mathbf{z}_{i^r}^r - \mathbf{x}_{i^r}^r\|$$
$$\leq \sum_{i=1}^n \sum_{r \in T_i^k} \|\mathbf{z}_i^r - \mathbf{x}_i^r\| - \sum_{i=1}^n T \cdot \left\|\mathbf{z}_i^{t_i^{j_k}} - \mathbf{x}_i^{t_i^{j_k}}\right\|$$

$$\leq \sum_{i=1}^{n} \sum_{r \in T_i^k} \|\mathbf{z}_i^r - \mathbf{x}_i^r\| - \sum_{i=1}^{n} |T_i^k| \cdot \left\|\mathbf{z}_i^{t_i^{j_k}} - \mathbf{x}_i^{t_i^{j_k}}\right\| = \sum_{i=1}^{n} \sum_{r \in T_i^k} \left(\|\mathbf{z}_i^r - \mathbf{x}_i^r\| - \left\|\mathbf{z}_i^{t_i^{j_k}} - \mathbf{x}_i^{t_i^{j_k}}\right\|\right)$$

$$\overset{(a)}{\leq} \sum_{i=1}^{n} \sum_{r \in T_i^k} \left(\left\|\widehat{\mathbf{x}}_i(\mathbf{x}^r) - \widehat{\mathbf{x}}_i(\mathbf{x}^{t_i^{j_k}})\right\| + \left\|\mathbf{x}_i^r - \mathbf{x}_i^{t_i^{j_k}}\right\|\right) + \underbrace{\sum_{i=1}^{n} \sum_{r \in T_i^k} (\varepsilon_i^{t_i^{j_k}} + \varepsilon_i^r)}_{\tilde{\varepsilon}_1^k}$$

$$\overset{(b)}{\leq} (1+\hat{L}) \sum_{i=1}^{n} \sum_{r \in T_i^k} \left\|\mathbf{x}^r - \mathbf{x}^{t_i^{j_k}}\right\| + \tilde{\varepsilon}_1^k$$

$$\leq (1+\hat{L}) \sum_{i=1}^{n} \sum_{r \in T_i^k} \sum_{s=r}^{t_i^{j_k}-1} \gamma^s \|\mathbf{z}_{i^s}^s - \mathbf{x}_{i^s}^s\| + \tilde{\varepsilon}_1^k \qquad (138)$$

$$\leq (1+\hat{L}) \sum_{i=1}^{n} |T_i^k| \cdot \sum_{r=kT+1}^{t_i^{j_k}-1} \gamma^r \|\mathbf{z}_{i^r}^r - \mathbf{x}_{i^r}^r\| + \tilde{\varepsilon}_1^k$$

$$\leq (1+\hat{L}) \sum_{i=1}^{n} |T_i^k| \cdot \left(\sum_{r=kT+1}^{j_k \cdot T} \gamma^r \|\mathbf{z}_{i^r}^r - \mathbf{x}_{i^r}^r\| + \sum_{r=j_k \cdot T+1}^{t_i^{j_k}-1} \gamma^r \|\mathbf{z}_{i^r}^r - \mathbf{x}_{i^r}^r\|\right) + \tilde{\varepsilon}_1^k$$

$$\leq (1+\hat{L}) \cdot (n \cdot T) \sum_{r=kT+1}^{j_k \cdot T} \gamma^r \|\mathbf{z}_{i^r}^r - \mathbf{x}_{i^r}^r\| + \underbrace{(1+\hat{L}) \cdot T \sum_{r=j_k \cdot T+1}^{t_i^{j_k}-1} \gamma^r \|\mathbf{z}_{i^r}^r - \mathbf{x}_{i^r}^r\|}_{\tilde{\varepsilon}_2^k} + \tilde{\varepsilon}_1^k$$

$$\leq (1+\hat{L}) \cdot (n \cdot T) \sum_{s=k}^{j_k-1} \sum_{r=sT+1}^{(s+1)T} \gamma^r \|\mathbf{z}_{i^r}^r - \mathbf{x}_{i^r}^r\| + \tilde{\varepsilon}_1^k + \tilde{\varepsilon}_2^k$$

$$\overset{(c)}{\leq} (1+\hat{L}) \cdot (n \cdot T) \sum_{s=k}^{j_k-1} (\eta_2)^{T-1} \gamma^{sT+1} \underbrace{\sum_{r=sT+1}^{(s+1)T} \|\mathbf{z}_{i^r}^r - \mathbf{x}_{i^r}^r\|}_{=\Delta^s} + \tilde{\varepsilon}_1^k + \tilde{\varepsilon}_2^k$$

$$\overset{(135)-(136)}{=} (1+\hat{L}) \cdot (n \cdot T) \cdot (\eta_2)^{T-1} \sum_{s=k}^{j_k-1} \gamma^{sT+1} \sum_{r=sT+1}^{(s+1)T} \frac{\|\mathbf{z}_{i^r}^r - \mathbf{x}_{i^r}^r\|^2}{\|\mathbf{z}_{i^r}^r - \mathbf{x}_{i^r}^r\|} + \tilde{\varepsilon}_1^k + \tilde{\varepsilon}_2^k$$

$$\overset{(d)}{\leq} \frac{T}{\delta} \cdot (1+\hat{L}) \cdot (n \cdot T) \cdot (\eta_2)^{T-1} \underbrace{\sum_{s=k}^{j_k-1} \gamma^{sT+1} \sum_{r=sT+1}^{(s+1)T} \|\mathbf{z}_{i^r}^r - \mathbf{x}_{i^r}^r\|^2}_{\tilde{\varepsilon}_3^k} + \tilde{\varepsilon}_1^k + \tilde{\varepsilon}_2^k,$$

where in (a) we used the reverse triangle inequality and Assumption II.9.1; (b) is due to the Lipschitz continuity of $\widehat{\mathbf{x}}_i$; in (c) we used $\gamma^{k+1}/\gamma^k \leq \eta_2$, with $\eta_2 \geq 1$; and (d) is due to Lemma II.15.

We prove now that $\tilde{\varepsilon}_1^k \downarrow 0$, $\tilde{\varepsilon}_2^k \downarrow 0$, and $\tilde{\varepsilon}_3^k \downarrow 0$. Since $\varepsilon_i^k \downarrow 0$ for all $i \in N$ (cf. Assumption II.9.1), it is not difficult to check that $\tilde{\varepsilon}_1^k \downarrow 0$. The same result for $\tilde{\varepsilon}_2^k$ comes from the following bound:

$$\tilde{\varepsilon}_2^k \leq (1+\hat{L}) \cdot T \sum_{r=j_k \cdot T+1}^{(j_k+1)T} \gamma^r \|\mathbf{z}_{i^r}^r - \mathbf{x}_{i^r}^r\| \leq (1+\hat{L}) \cdot T \cdot (\eta_2)^{T-1} \cdot \gamma^{j_k \cdot T+1} \cdot \Delta^{j_k}$$

$$\stackrel{(135)}{\leq} \left((1+\hat{L}) \cdot T \cdot (\eta_2)^{T-1} \cdot \delta\right) \gamma^{j_k \cdot T+1} \xrightarrow[k \to \infty]{} 0,$$

Finally, it follows from (132) that the series $\sum_{s=0}^\infty \gamma^{sT+1} \sum_{r=sT+1}^{(s+1)T} \|\mathbf{z}_{i^r}^r - \mathbf{x}_{i^r}^r\|^2$ is convergent; the Cauchy convergence criterion implies that

$$\tilde{\varepsilon}_3^k = \sum_{s=k}^{j_k-1} \gamma^{sT+1} \sum_{r=sT+1}^{(s+1)T} \|\mathbf{z}_{i^r}^r - \mathbf{x}_{i^r}^r\|^2 \xrightarrow[k \to \infty]{} 0.$$

By the vanishing properties of $\tilde{\varepsilon}_1^k$, $\tilde{\varepsilon}_2^k$, and $\tilde{\varepsilon}_3^k$, there exists a sufficiently large $k \in K$, say $\bar{k}$, such that

$$\tilde{\varepsilon}_1^k \leq \frac{T\delta}{4}, \quad \tilde{\varepsilon}_2^k \leq \frac{T\delta}{4}, \quad \tilde{\varepsilon}_3^k \leq \frac{\delta^2}{4(1+\hat{L}) \cdot n \cdot T \cdot (\eta_2)^{T-1}}, \quad \forall k \in K, k \geq \bar{k}, \tag{139}$$

which contradicts (138). Therefore, it must be $\limsup_{k \to \infty} \Delta^k = 0$.

• **Greedy rule Assumption II.10.2:** By (130), it is sufficient to prove $\limsup_{k \to +\infty} \widehat{\Delta}^k \triangleq \|\widehat{\mathbf{z}}^k - \mathbf{x}^k\| = 0$. Assume the contrary, that is, $\limsup_{k \to +\infty} \widehat{\Delta}^k > 0$. Since $\liminf_{k \to +\infty} \widehat{\Delta}^k = 0$ (cf. Proposition II.20), there exists a $\delta > 0$ such that $\widehat{\Delta}^k < \delta$ for infinitely many $k$ and also $\widehat{\Delta}^k > 2 \cdot (\delta/\tilde{c})$ for infinitely many $k$, where $0 < \tilde{c} \leq 1$ is the constant defined in Lemma II.19. Therefore, there is a set $\widehat{K} \subseteq \mathbb{N}_+$, with $|\widehat{K}| = \infty$, such that for each $k \in \widehat{K}$, there exists an integer $j_k > k$, such that

$$\widehat{\Delta}^k \geq 2\frac{\delta}{\tilde{c}}, \qquad \widehat{\Delta}^{j_k} \leq \delta, \tag{140}$$

$$\delta < \widehat{\Delta}^t < 2\frac{\delta}{\tilde{c}} \qquad \text{if} \quad k < t < j_k. \tag{141}$$

We have

$$\begin{aligned}
\frac{\delta}{\tilde{c}} &= \frac{2\delta}{\tilde{c}} - \frac{\delta}{\tilde{c}} \\
&\leq \left\|\widehat{\mathbf{z}}^k - \mathbf{x}^k\right\| - \frac{1}{\tilde{c}} \left\|\widehat{\mathbf{z}}^{j_k} - \mathbf{x}^{j_k}\right\| \\
&\stackrel{(a)}{\leq} \|(\widehat{\mathbf{x}}(\mathbf{x}^k) - \mathbf{x}^k)_{S^k}\| - \frac{1}{\tilde{c}} \|(\widehat{\mathbf{x}}(\mathbf{x}^{j_k}) - \mathbf{x}^{j_k})_{S^{j_k}}\| + \underbrace{\sum_{i \in S^k} \varepsilon_i^k + \frac{1}{\tilde{c}} \sum_{i \in S^{j_k}} \varepsilon_i^{j_k}}_{\triangleq \tilde{\varepsilon}^k} \\
&\stackrel{(117)}{\leq} \|\widehat{\mathbf{x}}(\mathbf{x}^k) - \mathbf{x}^k\| - \|\widehat{\mathbf{x}}(\mathbf{x}^{j_k}) - \mathbf{x}^{j_k}\| + \tilde{\varepsilon}^k
\end{aligned}$$

$$
\begin{aligned}
&\leq && \|\widehat{\mathbf{x}}(\mathbf{x}^k) - \widehat{\mathbf{x}}(\mathbf{x}^{j_k})\| + \|\mathbf{x}^k - \mathbf{x}^{j_k}\| + \tilde{\varepsilon}^k \\
&\stackrel{(b)}{\leq} && (1+\hat{L})\|\mathbf{x}^k - \mathbf{x}^{j_k}\| + \tilde{\varepsilon}^k \\
&\leq && (1+\hat{L}) \sum_{t=k}^{j_k-1} \gamma^t \|\widehat{\mathbf{z}}^t - \mathbf{x}^t\| + \tilde{\varepsilon}^k \\
&= && (1+\hat{L}) \sum_{t=k}^{j_k-1} \gamma^t \frac{\|\widehat{\mathbf{z}}^t - \mathbf{x}^t\|^2}{\|\widehat{\mathbf{z}}^t - \mathbf{x}^t\|} + \tilde{\varepsilon}^k \\
&\stackrel{(140)-(141)}{<} && \frac{1+\hat{L}}{\delta} \sum_{t=k}^{j_k-1} \gamma^t \|\widehat{\mathbf{z}}^t - \mathbf{x}^t\|^2 + \tilde{\varepsilon}^k,
\end{aligned} \quad (142)
$$

where in (a) we used the triangle inequality and Assumption II.9.1; and in (b) we used the Lipschitz continuity of $\widehat{\mathbf{x}}$.

Since $\tilde{\varepsilon}^k \downarrow 0$ (due to $\varepsilon_i^k \downarrow 0$, for all $i$; see Assumption II.9.1) and $\sum_{t=k}^{j_k-1} \gamma^t \|\widehat{\mathbf{z}}^t - \mathbf{x}^t\|^2 \xrightarrow[k\to\infty]{} 0$ (due to $\sum_{k=0}^{\infty} \gamma^k \|\widehat{\mathbf{z}}^k - \mathbf{x}^k\|^2 < +\infty$; see Proposition II.20), there exists a sufficiently large $k \in \widehat{K}$, say $\widehat{k}$, such that

$$
\tilde{\varepsilon}^k \leq \frac{\delta}{3 \cdot \tilde{c}} \quad \text{and} \quad \sum_{t=k}^{j_k-1} \gamma^t \|\widehat{\mathbf{z}}^t - \mathbf{x}^t\|^2 \leq \frac{\delta^2}{3 \cdot \tilde{c} \cdot (1+\hat{L})}, \quad \forall k \in \widehat{K}, \quad k \geq \widehat{k}, \quad (143)
$$

which contradicts (142). Therefore it must be $\lim_{k\to\infty} \|\widehat{\mathbf{z}}^k - \mathbf{x}^k\| = 0$.

## II.4. Parallel SCA: Hybrid Schemes

FLEXA (Algorithm 6) and its convergence theory cover fully parallel Jacobi as well as Gauss-Southwell-type methods, and many of their variants. Of course, every block selection rule $S^k$ has advantages and disadvantages, as discussed in Sec. II.3.1. A natural question is whether it is possible to design *hybrid* schemes that inherit the best of the aforementioned plain selection schemes. In the following, we introduce three hybrid strategies that are particularly well suited to parallel optimization on multi-core/processor architectures.

### II.4.1 Random-greedy schemes

In huge-scale optimization problems, the number of variables can be so large that computing the error bound functions $E_i$ for each block $i$ is not computationally affordable. To reduce the computational burden, one can think of a hybrid random/greedy block selection rule that combines random and greedy updates in the following form. First, a random selection is performed—the set $\widehat{S}^k \subseteq N$ is generated. Second, a greedy procedure is invoked to select *within the pool* $\widehat{S}^k$ only the subset of blocks, say $S^k \subseteq \widehat{S}^k$, that are "promising" according to the value of $E_i(\mathbf{x}^k)$. Finally all the blocks in $S^k$ are updated in parallel. This procedure captures both the advantages of random and greedy schemes: the random selection drops off a large proportion of the blocks allowing one to save computation while the greedy selection finds the best candidates that generate the "maximum improvement" on the objective function. The procedure is summarized in Algorithm 7. Convergence of the algorithm was studied in [56] when $\{\gamma^k\}_{k\in\mathbb{N}_+}$ is chosen according to the diminishing step-size rule [Assumption II.11.2]. The proof can be extended to deal with the other step-size rules (Assumption II.11) following similar steps as those in Sec. II.3.2, and thus is omitted.

---

**Algorithm 7: Hybrid Random-Greedy FLEXA**

**Data** : $\mathbf{x}^0 \in X$, $\{\gamma^k \in (0,1]\}_{k\in\mathbb{N}_+}$, $\varepsilon_i^k \geq 0$, for all $i \in N$ and $k \in \mathbb{N}_+$, $\rho \in (0,1]$.
  Set $k = 0$.

(S.1) : If $\mathbf{x}^k$ satisfies a termination criterion: STOP;

(S.2) : Randomly generate a set of blocks $\widehat{S}^k \subseteq N$;

(S.3) : Set $M^k \triangleq \max_{i\in\widehat{S}^k}\{E_i(\mathbf{x}^k)\}$;
  Choose a subset $S^k \subseteq \widehat{S}^k$ that contains at least one index $i$
  for which $E_i(\mathbf{x}^k) \geq \rho M^k$;

(S.4) : For all $i \in S^k$, solve (96) with accuracy $\varepsilon_i^k$:
  find $\mathbf{z}_i^k \in X_i$ s.t. $\|\mathbf{z}_i^k - \widehat{\mathbf{x}}_i(\mathbf{x}^k)\| \leq \varepsilon_i^k$;
  Set $\widehat{\mathbf{z}}_i^k = \mathbf{z}_i^k$ for $i \in S^k$, and $\widehat{\mathbf{z}}_i^k = \mathbf{x}_i^k$ for $i \notin S^k$;

(S.5) : Set $\mathbf{x}^{k+1} \triangleq \mathbf{x}^k + \gamma^k(\widehat{\mathbf{z}}^k - \mathbf{x}^k)$;

(S.6) : $k \leftarrow k+1$, and go to (S.1).

---

### II.4.2 Parallel-cyclic schemes

In a multi-core/processor architecture, another strategy to tackle huge-scale optimization problems is to adopt a hybrid parallel-cyclic strategy whereby the blocks of variables are partitioned among the workers (e.g., cores, processors) and updated *in parallel*, with each worker processing *sequentially* one block at a time.

Specifically, suppose that there are $P$ workers. Let $\{I_p\}_{p=1}^P$ be a partition of $N$ (the indices in $I_p$ follow the natural order). Assign the blocks $\mathbf{x}_i$, with $i \in I_p$, to worker $p$, and write $\mathbf{x}_p \triangleq (\mathbf{x}_{pi})_{i\in I_p}$, where $\mathbf{x}_{pi}$ denotes the $i$-th block assigned to worker $p$; $\mathbf{x}_{-p} \triangleq (\mathbf{x}_i)_{i\notin I_p}$ is the vector of remaining variables, assigned to the other workers. Finally, given $I_p$, we partition $\mathbf{x}_p$ as $\mathbf{x}_p = (\mathbf{x}_{p,i<}, \mathbf{x}_{p,i\geq})$, where $\mathbf{x}_{p,i<} \triangleq (\mathbf{x}_{pj})_{j\in I_p, j<i}$ is the vector containing all the variables in $\mathbf{x}_p$ that appear before $i$ (according to the order in $I_p$) whereas $\mathbf{x}_{p,i\geq} \triangleq (\mathbf{x}_{pj})_{j\in I_p, j\geq i}$ contains the remaining variables.

Once the optimization variables have been assigned to the $P$ workers, one could in principle apply the plain parallel algorithm, described in Algorithm 5, which would lead to the following scheme. All the workers update in parallel their variables $(\mathbf{x}_{pi})_{i\in I_p}$, but each worker can process only one block $\mathbf{x}_{pi}$ at a time. This means that, at iteration $k$, every worker $p$ computes a suitable $\mathbf{z}_{pi}^k$, for each block $i \in I_p$ (one at a time), by keeping all variables but $\mathbf{x}_{pi}$ fixed to $(\mathbf{x}_{pj}^k)_{i\neq j\in I_p}$ and $\mathbf{x}_{-p}^k$. Since we are solving the problems for each group of variables assigned to a worker sequentially but without using the most recent information for the updates, this approach seems a waste of resources. It is much more efficient to use, within each worker, a Gauss-Seidel scheme whereby the *most recent* iterates are used in all subsequent calculations. More specifically, at each iteration $k$, each worker $p$ computes sequentially an inexact version $\mathbf{z}_{pi}^k$ of the solution $\widehat{\mathbf{x}}_{pi}\left(\mathbf{x}_{p,i<}^{k+1}, \mathbf{x}_{p,i\geq}^k, \mathbf{x}_{-p}^k\right)$, with $i \in I_p$, and updates $\mathbf{x}_{pi}$ as $\mathbf{x}_{pi}^{k+1} = \mathbf{x}_{pi}^k + \gamma^k \left(\mathbf{z}_{pi}^k - \mathbf{x}_{pi}^k\right)$. The scheme is summarized in Algorithm 8. Convergence of the algorithm was studied in [79] when $\{\gamma^k\}_{k\in\mathbb{N}_+}$ is chosen according to the diminishing step-size rule [Assumption II.11.2]. The proof can be extended to deal with the other step-size rules (Assumption II.11) following similar steps as those in Sec. II.3.2, and thus is omitted.

---

**Algorithm 8: Hybrid Parallel-Cyclic FLEXA**

---

**Data** : $\mathbf{x}^0 \in X$, $\{\gamma^k \in (0,1]\}_{k\in\mathbb{N}_+}$, $\varepsilon_{pi}^k \geq 0$, for all $i \in I_p$, $p \in \{1,\ldots,P\}$, and $k \in \mathbb{N}_+$, $\rho \in (0,1]$.
    Set $k = 0$.

(S.1) : If $\mathbf{x}^k$ satisfies a termination criterion: STOP;

(S.2) : For all $p \in \{1,\ldots,P\}$ do (in parallel),
    For all $i \in I_p$ do (sequentially)
        a) Find $\mathbf{z}_{pi}^k \in X_i$ s.t.
$$\left\|\mathbf{z}_{pi}^k - \widehat{\mathbf{x}}_{pi}\left(\mathbf{x}_{p,i<}^{k+1}, \mathbf{x}_{p,i\geq}^k, \mathbf{x}_{-p}^k\right)\right\| \leq \varepsilon_{pi}^k;$$
        b) Set $\mathbf{x}_{pi}^{k+1} \triangleq \mathbf{x}_{pi}^k + \gamma^k \left(\mathbf{z}_{pi}^k - \mathbf{x}_{pi}^k\right)$;

(S.3) : $k \leftarrow k+1$, and go to (S.1).

---

### II.4.3 Parallel-Greedy-cyclic schemes

In Algorithm 8, at each iteration $k$, the workers update all their blocks, sequentially over $(\mathbf{x}_{pi})_{i\in I_p}$. However, in some large-scale instances of Problem (94), updating all variables might not always be beneficial or doable. Furthermore, using the latest information in the updates of each worker may require extra calculations (e.g., computing all block gradients) and communication overhead (these aspects are discussed on specific examples in Sec. II.5). It may be of interest then to consider a hybrid parallel-greedy-cyclic scheme, where Algorithm 8 is equipped with a greedy selection rule. More specifically, at each iteration $k$, each worker proceeds as in Algorithm 8 but performing the cyclic sweeping only on a subset $S_p^k$ of its own variables $I_p$, where the subset $S_p^k$ is chosen according to the greedy rule in Assumption II.10. To introduce formally the scheme, we extend the notation used in Algorithm 8 as follows. We reorder the components of $\mathbf{x}_p$ so that the first $|S_p^k|$ variables are those in $S_p^k$ and the remaining variables are those in $I_p \setminus S_p^k$; we write $\mathbf{x}_p = (\mathbf{x}_{S_p^k}, \mathbf{x}_{I_p \setminus S_p^k})$, where $\mathbf{x}_{S_p^k}$ denotes the vector containing the variables indexed in $S_p^k$ while $\mathbf{x}_{I_p \setminus S_p^k}$ contains the remaining variables in $I_p$ but $S_p^k$. Given an index $i \in S_p^k$, we partition $\mathbf{x}_{S_p^k}$ as $\mathbf{x}_{S_p^k} = (\mathbf{x}_{S_p^k i <}, \mathbf{x}_{S_p^k i \geq})$, where $\mathbf{x}_{S_p^k i <}$ is the vector containing all variables in $S_p^k$ that appear before $i$ (in the order assumed in $S_p^k$), while $\mathbf{x}_{S_p^k i \geq}$ are the remaining variables in $S_p^k$. With a slight abuse of notation, we will write $\mathbf{x} = (\mathbf{x}_{S_p^k i <}, \mathbf{x}_{S_p^k i \geq}, \mathbf{x}_{I_p \setminus S_p^k}, \mathbf{x}_{-p})$. The parallel-greedy-cyclic FLEXA is formally described in Algorithm 9.

---

**Algorithm 9: Hybrid Parallel-Greedy-Cyclic FLEXA**

---

**Data** : $\mathbf{x}^0 \in X$, $\{\gamma^k \in (0,1]\}_{k \in \mathbb{N}_+}$, $\varepsilon_{pi}^k \geq 0$, for all $i \in I_p$, $p \in \{1,\ldots,P\}$, and $k \in \mathbb{N}_+$, $\rho \in (0,1]$.

      Set $k = 0$.

(S.1) : If $\mathbf{x}^k$ satisfies a termination criterion: STOP;

(S.2) : Set $M^k \triangleq \max_i \{E_i(\mathbf{x}^k)\}$.

      Choose sets $S_p^k \subseteq I_p$ so that $\cup_{p=1}^P S_p^k$ contains at least one index $i$

      for which $E_i(\mathbf{x}^k) \geq \rho M^k$.

(S.3) : For all $\ell \in P$ do (in parallel),

      For all $i \in S_\ell^k$ do (sequentially)

        a) Find $\mathbf{z}_{\ell i}^k \in X_i$ s.t.

$$\left\| \mathbf{z}_{\ell i}^k - \widehat{\mathbf{x}}_{\ell i}\left( \mathbf{x}_{S_\ell^k i <}^{k+1}, \mathbf{x}_{S_\ell^k i \geq}^k, \mathbf{x}_{I_\ell \setminus S_\ell^k}, \mathbf{x}_{-\ell}^k \right) \right\| \leq \varepsilon_{\ell i}^k;$$

        b) Set $\mathbf{x}_{\ell i}^{k+1} \triangleq \mathbf{x}_{\ell i}^k + \gamma^k \left( \mathbf{z}_{\ell i}^k - \mathbf{x}_{\ell i}^k \right)$

(S.4) : Set $\mathbf{x}_{\ell i}^{k+1} = \mathbf{x}_{\ell i}^k$ for all $i \notin S^k$,

      $k \leftarrow k+1$, and go to (S.1).

---

## II.5. Applications

In this section we apply (Hybrid) FLEXA to some representative convex and nonconvex problems, arising from applications in communications and signal processing/machine learning. More specifically, we consider the following problems: i) the transceiver design in Single-Input-Single-Output (SISO) and Multiple-Input-Multiple-Output (MIMO) multiuser interference systems (cf. Sec. II.5.1); ii) the LASSO problem (cf. Sec. II.5.2); and iii) the logistic regression problem (cf. Sec. II.5.3). For each of the problems above, we show how to choose ad-hoc surrogate functions and provide extensive numerical results.

### II.5.1 Resource allocation in SISO/MIMO multiuser systems

Consider a multiuser system, composed of $I$ transmitter-receiver pairs (users). Each transmitter (resp. receiver) is equipped with $n_T$ (resp. $n_R$) transmit (resp. receive) antennas (w.l.o.g. we assumed that all the transceivers have the same number of antennas). Each transmitter is interested in transmitting its own data stream to its own receiver. No multiple access scheme is fixed a-priori (like OFDMA, TDMA, or CDMA); hence, Multi-User Interference (MUI) is experienced at each receiver. Let $\mathbf{x}_i \in \mathbb{C}^{n_T}$ be the signal transmitted by user $i$. Assuming a linear channel model, the received signal from user $i$ reads

$$\mathbf{y}_i = \underbrace{\mathbf{H}_{ii}\mathbf{x}_i}_{\text{desired signal}} + \underbrace{\sum_{j \neq i} \mathbf{H}_{ij}\mathbf{x}_j}_{\text{multiuser interference}} + \underbrace{\mathbf{n}_i}_{\text{noise}}, \qquad (144)$$

where $\mathbf{H}_{ij} \in \mathbb{C}^{n_R \times n_T}$ is the channel matrix between the transmitter $j$ and the receiver $i$, and $\mathbf{n}_i$ is the additive Gaussian zero-mean noise at the receiver $i$, with variance $\sigma_i^2 > 0$ (the noise is assumed to be white without loss of generality, otherwise one can always pre-whiten the channel matrices). The first term on the RHS of (144) represents the useful signal for user $i$ while the second one is the MUI due to the other users' concurrent transmissions. Note that the system model in (144) captures a fairly general MIMO setup, describing multiuser transmissions over multiple channels, which may represent frequency channels (as in multicarrier systems), time slots (as in time-division multiplexing systems), or spatial channels (as in transmit/receive beamforming systems); each of the aforementioned cases corresponds to a specific structure of the channel matrices $\mathbf{H}_{ij}$.

Denoting by $\mathbf{Q}_i \triangleq \mathbb{E}(\|\mathbf{x}_i\|^2) \succeq \mathbf{0}$ the covariance matrix of the symbols transmitted by agent $i$, each transmitter $i$ is subject to the following general power constraints

$$\mathcal{Q}_i \triangleq \left\{ \mathbf{Q}_i \in \mathbb{C}^{n_T \times n_T} : \mathbf{Q}_i \succeq \mathbf{0}, \quad \text{tr}(\mathbf{Q}_i) \leq P_i^{\text{ave}}, \quad \mathbf{Q}_i \in \mathcal{Z}_i \right\}, \qquad (145)$$

where $\text{tr}(\mathbf{Q}_i) \leq P_i^{\text{ave}}$ is a constraint on the maximum average transmit power, with $P_i^{\text{ave}}$ being the transmit power in unit of energy per transmission; and $\mathcal{Z}_i \subseteq \mathbb{C}^{n_T \times n_T}$ is an arbitrary closed and convex set, which can capture additional power/interference constraints (if any), such as: i) null constraints $\mathbf{U}_i^H \mathbf{Q}_i = \mathbf{0}$, where $\mathbf{U}_i \in \mathbb{C}^{n_T \times r_i}$ is a full rank matrix with $r_i < n_T$, whose columns represent the spatial and/or "frequency" directions along with user $i$ is not allowed to transmit; ii) soft-shaping constraints $\text{tr}\left(\mathbf{G}_i^H \mathbf{Q}_i \mathbf{G}_i\right) \leq I_i^{\text{ave}}$, which permit to control the power radiated (and thus the interference generated) onto the range space of $\mathbf{G}_i \in \mathbb{C}^{n_T \times n_T}$; iii) peak-power constraints $\lambda_{\max}\left(\mathbf{T}_i^H \mathbf{Q}_i \mathbf{T}_i\right) \leq I_i^{\text{peak}}$, which limit the average peak power of transmitter $i$ along the direction spanned by the range space of $\mathbf{T}_i \in \mathbb{C}^{n_T \times n_T}$, with $\lambda_{\max}$ denoting the maximum eigenvalue of the argument matrix; and iv) per-antenna constraints $[\mathbf{Q}_i]_{nn} \leq \alpha_{in}$, which control the maximum average power radiated by each antenna.

Under standard information theoretical assumptions, the maximum achievable rate on each link $i$ can be written as follows: given $\mathbf{Q} \triangleq (\mathbf{Q}_i)_{i=1}^I$,

$$R_i(\mathbf{Q}_i, \mathbf{Q}_{-i}) \triangleq \log\det\left(\mathbf{I} + \mathbf{H}_{ii}^H \mathbf{R}_i(\mathbf{Q}_{-i})^{-1} \mathbf{H}_{ii} \mathbf{Q}_i\right), \qquad (146)$$

where $\det(\bullet)$ is the determinant of the argument matrix; $\mathbf{Q}_{-i} \triangleq (\mathbf{Q}_j)_{j \neq i}$ denotes the tuple of the (complex-valued) covariance matrices of all the transmitters except the $i$-th one; and $\mathbf{R}_i(\mathbf{Q}_{-i}) \triangleq \mathbf{R}_{n_i} + \sum_{j \neq i} \mathbf{H}_{ij} \mathbf{Q}_j \mathbf{H}_{ij}^H$ is the covariance matrix of the multiuser interference plus the thermal noise $\mathbf{R}_{n_i}$ (assumed to be full-rank).

As system design, we consider the maximization of the users' (weighted) sum rate, subject to the power constraint (145), which reads

$$\begin{array}{ll} \underset{\mathbf{Q}_1,\ldots,\mathbf{Q}_I}{\text{maximize}} & \sum_{i=1}^{I} \alpha_i R_i(\mathbf{Q}_i, \mathbf{Q}_{-i}) \\ \text{subject to} & \mathbf{Q}_i \in Q_i, \quad \forall i = 1,\ldots,I, \end{array} \quad (147)$$

where $(\alpha_i)_{i=1}^{I}$ are given positive weights, which one can use to prioritize some user with respect to another. We remark that the proposed algorithmic framework can be applied also to other objective functions involving the rate functions, see [77, 208].

Clearly (147) is an instance of (94) (with $G = 0$ and involving complex variables) and thus we can apply the algorithmic framework described in this lecture. We begin considering the sum-rate maximization problem (147) over SISO frequency selective channels; we then extend the analysis to the more general MIMO case.

**Sum-rate maximization over SISO interference channels**

Given the system model (144), consider SISO frequency selective channels: the channel matrices $\mathbf{H}_{ij}$ are $m \times m$ Toeplitz circulant matrices and $\mathbf{R}_{n_i}$ are $m \times m$ diagonal matrices, with diagonal entries $\sigma_{i1}^2, \ldots, \sigma_{im}^2$ ($\sigma_{i\ell}^2$ is the variance of the noise on channel $\ell$); and $m$ is the length of the transmitted block [note that in (144) it becomes $n_T = n_R = m$]; see, e.g., [246]. The eigendecomposition of each $\mathbf{H}_{ij}$ reads: $\mathbf{H}_{ij} = \mathbf{F} \mathbf{D}_{ij} \mathbf{F}^H$, where $\mathbf{F}$ is the IFFT matrix, i.e., $[\mathbf{F}]_{\ell'\ell} = \exp(j 2\pi (\ell' - 1)(\ell - 1)/N)/\sqrt{N}$, for $\ell', \ell = 1, \ldots N$; and $\mathbf{D}_{ij}$ is the diagonal matrix whose diagonal entries $H_{ij}(1), \ldots, H_{ij}(N)$ are the coefficients of the frequency-response of the channel between the transmitter $j$ and the receiver $i$.

Orthogonal Frequency Division Multiplexing (OFDM) transmissions correspond to the following structure for the covariance matrices: $\mathbf{Q}_i = \mathbf{F} \operatorname{diag}(\mathbf{p}_i) \mathbf{F}^H$, where $\mathbf{p}_i \triangleq (p_{i\ell})_{\ell=1}^{m}$ is the transmit power profile of user $i$ over the $m$ frequency channels. The power constraints read: given $\mathbf{I}_i^{\max} \in \mathbb{R}_+^{q_i}$ and $\mathbf{W}_i \in \mathbb{R}_+^{q_i \times m}$,

$$P_i \triangleq \left\{ \mathbf{p}_i \in \mathbb{R}_+^N : \mathbf{W}_i \mathbf{p}_i \leq \mathbf{I}_i^{\max} \right\}, \quad (148)$$

where the inequality has to be intended component-wise. To avoid redundant constraints, we assume w.l.o.g. that all the columns of $\mathbf{W}_i$ are linearly independent.

The maximum achievable rate on each link $i$ becomes [cf. (146)]

$$r_i(\mathbf{p}_i, \mathbf{p}_{-i}) \triangleq \sum_{\ell=1}^{m} \log \left( 1 + \frac{|H_{ii}(\ell)|^2 p_{i\ell}}{\sigma_{i\ell}^2 + \sum_{j \neq i} |H_{ij}(\ell)|^2 p_{j\ell}} \right), \quad (149)$$

where $\mathbf{p}_{-i} \triangleq (\mathbf{p}_j)_{j \neq i}$ is the power profile of all the users $j \neq i$.

The system design (147) reduces to the following nonconvex optimization problem

$$\underset{\mathbf{p}_1,\ldots,\mathbf{p}_I}{\text{maximize}} \quad \sum_{i=1}^{I} \alpha_i\, r_i(\mathbf{p}_i, \mathbf{p}_{-i}) \tag{150}$$
$$\text{subject to} \quad \mathbf{p}_i \in P_i, \quad \forall i = 1, \ldots, I.$$

We apply next FLEXA (Algorithm 5) to (150); we describe two alternative SCA-decompositions, corresponding to two different choices of the surrogate functions.

**Decomposition #1—Pricing Algorithms:** Since the sum-rate maximization problem (150) is an instance of the problem considered in Example 3 in Sec. II.2.1, a first approach is to use the surrogate (106). Since the rate $r_i(\mathbf{p}_i, \mathbf{p}_{-i})$ is concave in $\mathbf{p}_i$, for any given $\mathbf{p}_{-i} \geq \mathbf{0}$, we have $\tilde{C}_i = \{i\}$ [cf. (105)] and thus $C_i \equiv \tilde{C}_i$, which leads to the following surrogate function: given $\mathbf{p}^k \geq \mathbf{0}$ at iteration $k$,

$$\widetilde{F}(\mathbf{p}\,|\,\mathbf{p}^k) = \sum_{i=1}^{I} \widetilde{F}_i(\mathbf{p}_i\,|\,\mathbf{p}^k),$$

where

$$\widetilde{F}_i(\mathbf{p}_i\,|\,\mathbf{p}^k) \triangleq \alpha_i \cdot r_i(\mathbf{p}_i, \mathbf{p}_{-i}^k) - \boldsymbol{\pi}_i(\mathbf{p}^k)^T (\mathbf{p}_i - \mathbf{p}_i^k) - \frac{\tau_i}{2} \left\| \mathbf{p}_i - \mathbf{p}_i^k \right\|^2,$$

$\tau_i$ is an arbitrary nonnegative constant, and $\boldsymbol{\pi}_i(\mathbf{p}^k) \triangleq (\pi_{i\ell}(\mathbf{p}^k))_{\ell=1}^{m}$ is defined as

$$\pi_{i\ell}(\mathbf{p}^k) \triangleq -\sum_{j \in N_i} \alpha_j |H_{ji}(\ell)|^2 \frac{\mathtt{snr}_{j\ell}^k}{(1 + \mathtt{snr}_{j\ell}^k) \cdot \mathtt{mui}_{j\ell}^k};$$

$N_i$ denotes the set of (out) neighbors of user $i$, i.e., the set of users $j$'s which user $i$ interferers with; and $\mathtt{snr}_{j\ell}^k$ and $\mathtt{mui}_{j\ell}^k$ are the Signal-to-Interference-plus-Noise (SINR) and the multiuser interference-plus-noise power ratios experienced by user $j$ on the frequency $\ell$, generated by the power profile $\mathbf{p}^k$:

$$\mathtt{snr}_{j\ell}^k \triangleq \frac{|H_{jj}(\ell)|^2 p_{j\ell}^k}{\mathtt{mui}_{j\ell}^k}, \quad \text{and} \quad \mathtt{mui}_{j\ell}^k \triangleq \sigma_{j\ell}^2 + \sum_{i \neq j} |H_{ji}(\ell)|^2 p_{i\ell}^k.$$

All the users in parallel will then solve the following strongly concave subproblems: given $\mathbf{p}^k = (\mathbf{p}_i^k)_{i=1}^{I}$,

$$\hat{\mathbf{p}}_i(\mathbf{p}^k) \triangleq \underset{\mathbf{p}_i \in P_i}{\operatorname{argmax}} \left\{ \alpha_i \cdot r_i(\mathbf{p}_i, \mathbf{p}_{-i}^k) - \boldsymbol{\pi}_i(\mathbf{p}^k)^T (\mathbf{p}_i - \mathbf{p}_i^k) - \frac{\tau_i}{2} \left\| \mathbf{p}_i - \mathbf{p}_i^k \right\|^2 \right\}.$$

Note that the best-response $\hat{\mathbf{p}}_i(\mathbf{p}^k)$ can be computed in closed form (up to the multiplies associated with the inequality constraints in $P_i$) according to the following multi-level waterfilling-like expression [209]: setting each $\tau_i > 0$,

$$\hat{\mathbf{p}}_i(\mathbf{p}^k) \triangleq \left[ \frac{1}{2} \mathbf{p}_i^k \circ \left(\mathbf{1} - (\mathbf{snr}_i^k)^{-1}\right) + \right. \\ \left. -\frac{1}{2\tau_i} \left( \tilde{\boldsymbol{\mu}}_i - \sqrt{[\tilde{\boldsymbol{\mu}}_i - \tau_i \mathbf{p}_i^k \circ (\mathbf{1} + (\mathbf{snr}_i^k)^{-1})]^2 + 4\tau_i w_i \mathbf{1}} \right) \right]_+ \quad (151)$$

where $\circ$ denotes the Hadamard product and $[\bullet]_+$ denotes the projection onto the nonnegative orthant $\mathbb{R}_+^m$; $(\mathbf{snr}_i^k)^{-1} \triangleq (1/\mathrm{snr}_{i\ell}^k)_{\ell=1}^m$ and $\tilde{\boldsymbol{\mu}}_i \triangleq \boldsymbol{\pi}_i(\mathbf{p}^k) + \mathbf{W}_i^T \boldsymbol{\mu}_i$, with the multiplier vector $\boldsymbol{\mu}_i$ chosen to satisfy the nonlinear complementarity condition (CC)

$$\mathbf{0} \leq \boldsymbol{\mu}_i \perp \mathbf{I}_i^{\max} - \mathbf{W}_i \hat{\mathbf{p}}_i(\mathbf{p}^k) \geq \mathbf{0}.$$

The optimal $\boldsymbol{\mu}_i$ satisfying the CC can be efficiently computed (in a finite number of steps) using the nested bisection method described in [209, Algorithm 6]; we omit further details here. Note that, in the presence of power budget constraints only, $\boldsymbol{\mu}_i$ reduces to a scalar quantity $\mu_i$ such that $0 \leq \mu_i \perp p_i - \mathbf{1}^T \hat{\mathbf{p}}_i(\mathbf{p}^k) \geq 0$, whose solution can be obtained using the classical bisection algorithm (or the methods in [180]).

Given $\hat{\mathbf{p}}_i(\mathbf{p}^k)$, one can now use, e.g., Algorithm 5, with any of the valid choices for the step-size $\{\gamma^k\}$ [cf. Assumption II.6]. Since there is no coordination among the users as well as no centralized control in network, one is interested in designing distributed algorithms. This naturally suggests the use of a diminishing step-size rule in Algorithm 5. For instance, good candidates are the rules in (108) or (109). Note that the resulting algorithm is fairly distributed. Indeed, given the interference generated by the other users [and thus the MUI coefficients $\mathtt{mui}_{jn}^k$] and the current interference price $\boldsymbol{\pi}_i(\mathbf{p}^k)$, each user can efficiently and locally compute the optimal power allocation $\hat{\mathbf{p}}_i(\mathbf{p}^k)$ via the waterfilling-like expression (151). The estimation of the prices $\pi_{i\ell}(\mathbf{p}^k)$ requires however some signaling among nearby users.

**Decomposition #2−DC Algorithms:** An alternative class of algorithms for the sum-rate maximization problem (150) can be obtained exploring the DC structure of the rate functions (149). By doing so, the sum-rate can be decomposed as the sum of a concave and convex function, namely $U(\mathbf{p}) = f_1(\mathbf{p}) + f_2(\mathbf{p})$, with

$$f_1(\mathbf{p}) \triangleq \sum_{i=1}^I \alpha_i \sum_{\ell=1}^m \log\left( \sigma_{i\ell}^2 + \sum_{j=1}^I |H_{ij}(\ell)|^2 p_{j\ell} \right),$$

$$f_2(\mathbf{p}) \triangleq -\sum_{i=1}^I \alpha_i \sum_{\ell=1}^m \log\left( \sigma_{i\ell}^2 + \sum_{j=1, j\neq i}^I |H_{ij}(\ell)|^2 p_{j\ell} \right).$$

A concave surrogate can be readily obtained from $U(\mathbf{p})$ by linearizing $f_2(\mathbf{p})$ and keeping $f_1(\mathbf{p})$ unaltered. This leads to the following strongly concave subproblem for each agent $i$: given $\mathbf{p}^k \geq \mathbf{0}$,

$$\widetilde{\mathbf{p}}_i(\mathbf{p}^k) \triangleq \underset{\mathbf{p}_i \in P_i}{\mathrm{argmax}} \left\{ f_1(\mathbf{p}_i, \mathbf{p}_{-i}^k) - \boldsymbol{\pi}_i(\mathbf{p}^k)^T (\mathbf{p}_i - \mathbf{p}_i^k) - \frac{\tau_i}{2} \left\| \mathbf{p}_i - \mathbf{p}_i^k \right\|^2 \right\}$$

where $\boldsymbol{\pi}_i(\mathbf{p}^k) \triangleq (\pi_{i\ell}(\mathbf{p}^k))_{\ell=1}^m$, with

$$\pi_{i\ell}(\mathbf{p}^k) \triangleq - \sum_{j \in N_i} \alpha_j |H_{ji}(\ell)|^2 \frac{1}{\mathtt{mui}_{j\ell}^k}. \tag{152}$$

The best-response $\widetilde{\mathbf{p}}_i(\mathbf{p}^k)$ can be efficiently computed using a fixed-point-based procedure, in the same spirit of [181]; we omit further details. Note that the communication overhead to compute the prices (151) and (152) is the same, but the computation of $\widetilde{\mathbf{p}}_i(\mathbf{p}^k)$ requires more (channel state) information exchange than that of $\hat{\mathbf{p}}_i(\mathbf{p}^k)$, since each user $i$ also needs to estimate the cross-channels $\{|H_{ji}(\ell)|^2\}_{j \in N_i}$.

**Numerical example**. We compare now Algorithm 5 based on the best-response $\hat{\mathbf{p}}_i(\mathbf{p}^k)$ in (151) (termed `SR-FLEXA`, SR stands for Sum-Rate), with those proposed in [181] [termed `SCALE` and `SCALE one-step`, the latter being a simplified version of `SCALE` where instead of solving the fixed-point equation (16) in [181], only one iteration of (16) is performed], [206] (which is an instance of the block MM algorithm described in Algorithm 2, and is termed `Block-MM`), and [215] (termed `WMMSE`). Since the algorithms in [181, 206, 215] can only deal with power budget constraints, to run the comparison, we simplified the sum-rate maximization problem (150) considering only power budget constraints (and all $\alpha_i = 1$). We assume the same power budget $P_i^{\text{ave}} = p$, noise variances $\sigma_{i\ell}^2 = \sigma^2$, and $\mathtt{snr} = p/\sigma^2 = 3\text{dB}$ for all the users. We simulated SISO frequency-selective channels with $m = 64$ subcarriers; the channels are generated as FIR filters of order $L = 10$, whose taps are i.i.d. Gaussian random variables with zero mean and variance $1/(d_{ij}^3(L+1)^2)$, where $d_{ij}$ is the distance between the transmitter $j$ and the receiver $i$. All the algorithms are initialized by choosing the uniform power allocation, and are terminated when (the absolute value of) the sum-utility error in two consecutive rounds becomes smaller than $1e$-3. The accuracy in the bisection loops (required by all methods) is set to $1e$-6. In `SR-FLEXA`, we used the rule (108) with $\varepsilon = 1e$-2. In Fig. II.3, we plot the average number of iterations required by the aforementioned algorithms to converge (under the same termination criterion) versus the number of users; the average is taken over 100 independent channel realizations; in Fig. II.3a we set $d_{ij}/d_{ii} = 3$ whereas in Fig. II.3b we have $d_{ij}/d_{ii} = 1$ while in both figures $d_{ij} = d_{ji}$ and $d_{ii} = d_{jj}$, for all $i$ and $j \neq i$; the setting in Fig. II.3a emulates a "low" MUI environment whereas the one in Fig. II.3b a "high" MUI scenario. All the algorithms reach the same average sum-rate. The figures clearly show that the proposed `SR-FLEXA` outperforms all the others (note that `SCALE` and `WMMSE` are also simultaneous-based schemes). For instance, in Fig. II.3a, the gap with `WMMSE` (in terms of number of iterations needed to reach convergence) is about one order of magnitude, for all the network sizes considered in the experiment, which reduces to two times in the "high" interference scenario considered in Fig. II.3b. Such a behavior (requiring less iterations than other methods, with gaps ranging from few times to one order of magnitude) has been observed also for other choices of $d_{ij}/d_{ii}$, termination tolerances, and step-size rules; more experiments can be found in [210, 222]. Note that `SR-FLEXA`, `SCALE one-step`, `WMMSE`, and `Block-MM` have similar per-user computational complexity, whereas `SCALE` is much more demanding and is not appealing for a real-time implementation. Therefore, Fig. II.3 provides also a roughly indication of the per-user cpu time of `SR-FLEXA`, `SCALE one-step`, and `WMMSE`.

**Sum-rate maximization over MIMO interference channels**

Let us focus now on the general MIMO formulation (147).

Similarly to the SISO case, we can invoke the surrogate (106) with $C_i = \{i\}$, corresponding to keeping $R_i$ in (146) unaltered and linearizing the rest of the sum, that is, $\sum_{j \neq i} R_j$. Invoking the Wirtinger calculus (see, e.g., [103, 124, 209]), the sub-problem solved by each agent $i$ at iteration $k$ reads: given $\mathbf{Q}^k = (\mathbf{Q}_i^k)_{i=1}^I$, with each $\mathbf{Q}^k \succeq \mathbf{0}$,

$$\hat{\mathbf{Q}}_i(\mathbf{Q}^k) \triangleq \underset{\mathbf{Q}_i \in Q_i}{\mathrm{argmax}} \left\{ \alpha_i r_i(\mathbf{Q}_i, \mathbf{Q}_{-i}^k) - \left\langle \boldsymbol{\Pi}_i(\mathbf{X}^k), \mathbf{Q}_i \right\rangle - \tau_i \left\| \mathbf{Q}_i - \mathbf{Q}_i^k \right\|_F^2 \right\} \quad (153)$$

where $\langle \mathbf{A}, \mathbf{B} \rangle \triangleq \mathrm{Re}\{\mathrm{tr}(\mathbf{A}^H \mathbf{B})\}$; $\tau_i > 0$,

$$\boldsymbol{\Pi}_i(\mathbf{Q}^k) \triangleq \sum_{j \in N_i} \alpha_j \mathbf{H}_{ji}^H \widetilde{\mathbf{R}}_j(\mathbf{Q}_{-j}^k) \mathbf{H}_{ji},$$

with $N_i$ defined as in the SISO case; and

$$\widetilde{\mathbf{R}}_j(\mathbf{Q}_{-j}^k) \triangleq \mathbf{R}_j(\mathbf{Q}_{-j}^k)^{-1} - (\mathbf{R}_j(\mathbf{Q}_{-j}^k) + \mathbf{H}_{jj} \mathbf{Q}_j^k \mathbf{H}_{jj}^H)^{-1}.$$

Note that, once the price matrix $\boldsymbol{\Pi}_i(\mathbf{Q}^k)$ is given, the best-response $\hat{\mathbf{Q}}_i(\mathbf{Q}^k)$ can be computed locally by each user solving a convex optimization problem. Moreover, for some specific structures of the feasible sets $Q_i$, the case of full-column rank channel matrices $\mathbf{H}_i$, and $\tau_i = 0$, a solution in closed form (up to the multipliers associated with the power budget constraints) is also available [123]; see also [260] for other examples. Given $\hat{\mathbf{Q}}_i(\mathbf{Q}^k)$, one can now use Algorithm 5 (adapted to the complex case), with any of the valid choices for the step-size $\{\gamma^k\}$.

**Complexity analysis and message exchange**. We compare here the computational complexity and signaling (i.e., message exchange) of Algorithm 5 based on the best-response $\hat{\mathbf{Q}}_i(\mathbf{Q}^k)$ (termed `MIMO-SR-FLEXA`) with those of the schemes proposed in the literature for a similar problem, namely the `MIMO-Block-MM` [123, 206], and the `MIMO-WMMSE` [215]. For the purpose of complexity analysis, since all algorithms include a similar bisection step which generally takes few iterations, we will ignore this step in the computation of the complexity. Also, `MIMO-WMMSE` and `MIMO-SR-FLEXA` are simultaneous schemes, while `MIMO-Block-MM` is sequential; we then compare the algorithms by given the *per-round complexity*, where one round means one update from all the users. Recalling that $n_T$ (resp. $n_R$) denotes the number of antennas at each transmitter (resp. receiver), the computational complexity of the algorithms is [210]:

- `MIMO-Block-MM`: $O\left(I^2(n_T n_R^2 + n_T^2 n_R + n_R^3) + I n_T^3\right)$;
- `MIMO-WMMSE`: $O\left(I^2(n_T n_R^2 + n_T^2 n_R + n_T^3) + I n_R^3\right)$ [215];

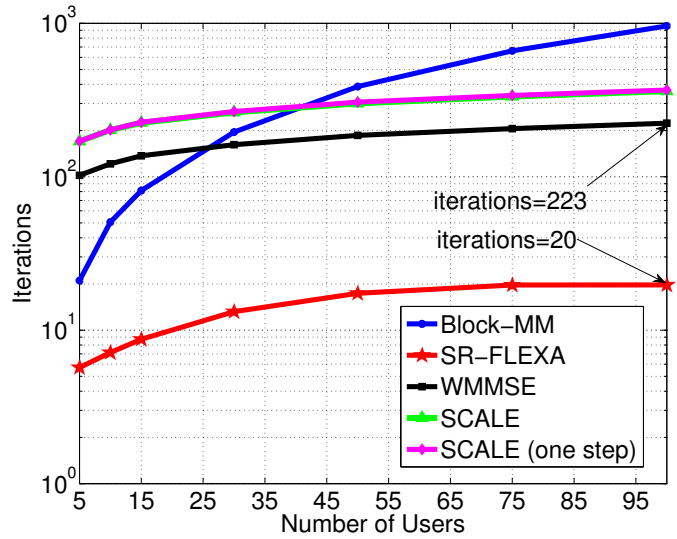

*(a) Low MUI: The proposed method,* `SR-FLEXA`*, is one order of magnitude faster than the WMMSE algorithm.*

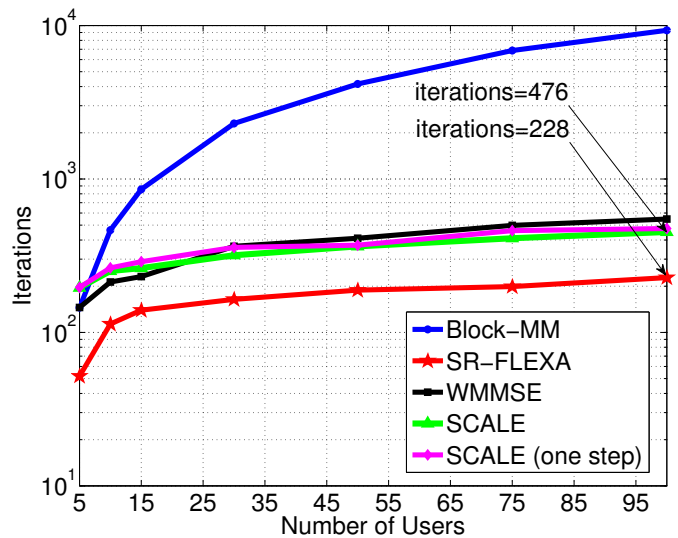

*(b) High MUI: The proposed method,* `SR-FLEXA`*, is two times faster than the WMMSE algorithm.*

*Fig. II.3: Sum-rate maximization problem* (150) *(SISO frequency-selective channels): Average number of iterations versus number of users. Note that all algorithms are simultaneous except* `Block-MM`*, which is sequential. Also, all the algorithms are observed to converge to the same stationary solution of Problem* (150)*. The figures are taken from [210].*

|  | # of users = 10 | | | # of users = 50 | | | # of users = 100 | | |
|---|---|---|---|---|---|---|---|---|---|
|  | d=1 | d=2 | d=3 | d=1 | d=2 | d=3 | d=1 | d=2 | d=3 |
| `MIMO-Block-MM` | 1370.5 | 187 | 54.4 | 4148.5 | 1148 | 348 | 8818 | 1904 | 704 |
| `MIMO-WMMSE` | 169.2 | 68.8 | 53.3 | 138.5 | 115.2 | 76.7 | 154.3 | 126.9 | 103.2 |
| `MIMO-SR-FLEXA` | 169.2 | 24.3 | 6.9 | 115.2 | 34.3 | 9.3 | 114.3 | 28.4 | 9.7 |

*Table II.1: Sum-rate maximization problem* (147) *(MIMO frequency-selective channels): Average number of iterations (termination accuracy=1e-6).*

|  | # of users = 10 | | | # of users = 50 | | | # of users = 100 | | |
|---|---|---|---|---|---|---|---|---|---|
|  | d=1 | d=2 | d=3 | d=1 | d=2 | d=3 | d=1 | d=2 | d=3 |
| `MIMO-Block-MM` | 429.4 | 74.3 | 32.8 | 1739.5 | 465.5 | 202 | 3733 | 882 | 442.6 |
| `MIMO-WMMSE` | 51.6 | 19.2 | 14.7 | 59.6 | 24.9 | 16.3 | 69.8 | 26.0 | 19.2 |
| `MIMO-SR-FLEXA` | 48.6 | 9.4 | 4.0 | 46.9 | 12.6 | 5.1 | 49.7 | 12 | 5.5 |

*Table II.2: Sum-rate maximization problem* (147) *(MIMO frequency-selective channels): Average number of iterations (termination accuracy=1e-3).*

- `MIMO-SR-FLEXA`: $O\left(I^2(n_T n_R^2 + n_T^2 n_R) + I(n_T^3 + n_R^3)\right)$.

The complexity of the three algorithms is very similar, and equivalent in the case in which $n_T = n_R (\triangleq m)$, given by $O(I^2 m^3)$.

In a real system, the MUI covariance matrices $\mathbf{R}_i(\mathbf{Q}_{-i})$ come from an estimation process. It is thus interesting to understand how the complexity changes when the computation of $\mathbf{R}_i(\mathbf{Q}_{-i})$ from $\mathbf{R}_{n_i} + \sum_{j \neq i} \mathbf{H}_{ij} \mathbf{Q}_j \mathbf{H}_{ij}^H$ is not included in the analysis. We obtain the following [210]:

- `MIMO-Block-MM`: $O\left(I^2(n_T n_R^2 + n_T^2 n_R + n_R^3) + I n_T^3\right)$;
- `MIMO-WMMSE`: $O\left(I^2(n_T^2 n_R + n_T^3) + I(n_R^3 + n_T n_R^2)\right)$;
- `MIMO-SR-FLEXA`: $O\left(I^2(n_T n_R^2 + n_T^2 n_R) + I(n_T^3 + n_R^3)\right)$.

Finally, if one is interested in the time necessary to complete one iteration, it can be shown that it is proportional to the above complexity divided by $I$.

As far as the communication overhead is concerned, the same remarks we made about the schemes described in the SISO setting, apply also here for the MIMO case. The only difference is that now the users need to exchange a (pricing) matrix rather than a vector, resulting in $O(I^2 n_R^2)$ amount of message exchange per-iteration for all the algorithms.

**Numerical example #1.** In Tables II.1 and II.2 we compare the `MIMO-SR-FLEXA`, the `MIMO-Block-MM` [123, 206], and the `MIMO-WMMSE` [215], in terms of average number of iterations required to reach convergence, for different number of users, normalized distances $d \triangleq d_{ij}/d_{ii}$ (with $d_{ij} = d_{ji}$ and $d_{ii} = d_{jj}$ for all $i$ and $j \neq i$), and termination accuracy (namely: 1e-3 and 1e-6). All the transmitters/receivers are equipped with 4 antenna; we simulated uncorrelated fading channels, whose coefficients are Gaussian distributed with zero mean and variance $1/d_{ij}^3$ (all the channel matrices are full-column rank); and we set $\mathbf{R}_{n_i} = \sigma^2 \mathbf{I}$ for all $i$, and $\text{snr} \triangleq p/\sigma^2 = 3\text{dB}$. In `MIMO-SR-FLEXA`, we used the step-size rule (108), with $\varepsilon = 1e$-5; in (153) we set $\tau_i = 0$ and computed $\hat{\mathbf{Q}}_i(\mathbf{Q}^k)$ using the closed form solution in [123]. All the algorithms reach the same average sum-rate.

Given the results in Tables II.1 and II.2, the following comments are in order. `MIMO-SR-FLEXA` outperforms the other schemes in terms of iterations, while having similar (or even better) computational complexity. Interestingly, the iteration gap with the other schemes reduces with the distance and the termination accuracy. More specifically: `MIMO-SR-FLEXA` i) seems to be much faster than the other schemes (about one order of magnitude) when $d_{ij}/d_{ii} = 3$ [say low interference scenarios], and just a bit faster (or comparable to `MIMO-WMMSE`) when $d_{ij}/d_{ii} = 1$ [say high interference scenarios]; and ii) it is much faster than the others, if an high termination accuracy is set (see Table II.1). Also, the convergence speed of `MIMO-SR-FLEXA` is not affected too much by the number of users. Finally, in our experiments, we also observed that the performance of `MIMO-SR-FLEXA` is not affected too much by the choice of the parameter $\varepsilon$ in the (108): a change of $\varepsilon$ of many orders of magnitude leads to a difference in the average number of iterations which is within 5%; we refer the reader to [222] for details, where one can also find a comparison of several other step-size rules. We must stress however that `MIMO-Block-MM` and `MIMO-WMMSE` do not need any tuning, which is an advantage with respect to `MIMO-SR-FLEXA`.

**Numerical example #2**: We compare now the `MIMO-WMMSE` [215] and the `MIMO-SR-FLEXA` in a MIMO broadcast cellular system composed of multiple cells, with one Base Station (BS) and multiple randomly generated Mobile Terminals (MTs) in each cell. Each MT experiences both intra-cell and inter-cell interference. We refer to [215] for a detailed description of the system model, the explicit expressions of the BS-MT downlink rates, and the corresponding sum-rate maximization problem.

The setup of our experiments is the following [210]. We simulated seven cells with multiple randomly generated MTs; each BS and MT is equipped with four transmit and receive antennas. Channels are Rayleigh fading, whose path-loss are generated using the 3 GPP(TR 36.814) methodology [1]. We assume white zero-mean Gaussian noise at each mobile receiver, with variance $\sigma^2$, and same power budget $p$ for all the BSs; the SNR is set to $\text{snr} \triangleq p/\sigma^2 = 3$dB. Both algorithms `MIMO-WMMSE` and `MIMO-SR-FLEXA` are initialized by choosing the same feasible randomly generated point, and are terminated when (the absolute value of) the sum-rate error in two consecutive rounds becomes smaller than $1e$-2. In `MIMO-SR-FLEXA`, the step-size rule (108) is used, with $\varepsilon = 1e$-3 and $\gamma^0 = 1$; the unique solution $\hat{\mathbf{Q}}_i(\mathbf{Q}^k)$ of users' subproblems is computed in closed form adapting the procedure in [123]. The experiments were run using Matlab R2012a on a $12 \times 2.40$ GHz Intel Xeon E5645 Processor Cores machine, equipped with 48 GB of memory and 24576 Kbytes of data cache; the operation system is Linux (Red-Hat Enterprise Linux 6.1 2.6.32 Kernel). In Fig. II.4a we plot the average cpu time versus the total number of MTs for the two algorithms under the same termination criterion, whereas in Fig. II.4b we reported the final achieved average sum-rate. The curves are averaged over 1500 channel/topology realizations. It can be observed that `MIMO-SR-FLEXA` significantly outperforms `MIMO-WMMSE` in terms of cpu time when the number of active users is large; moreover `MIMO-SR-FLEXA` also yields better sum-rates. We observed similar results also under different settings (e.g., SNR, number of cells/BSs, etc.); see [222] for more details.

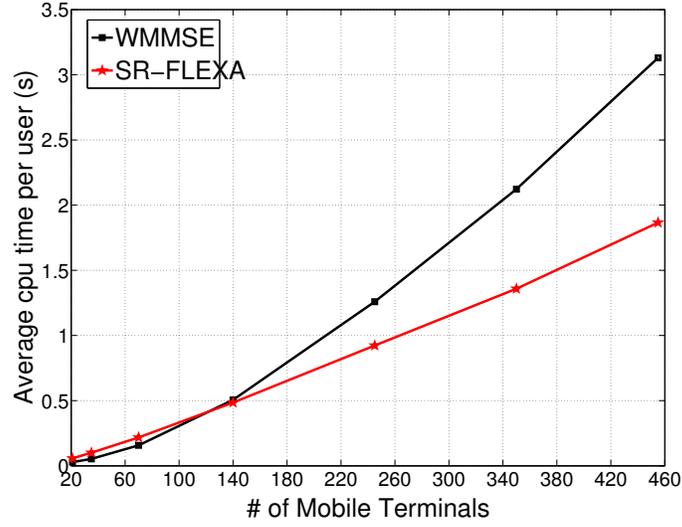

*(a) Average cpu time versus the number of mobile terminals.*

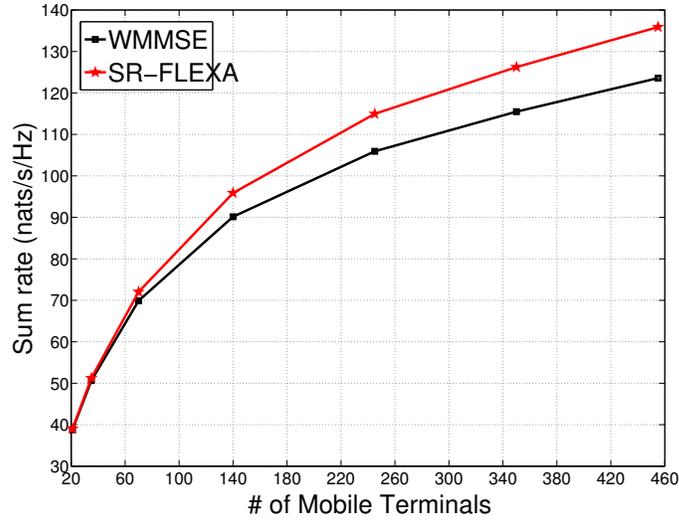

*(b) Average sum-rate versus the number of mobile terminals.*

Fig. II.4: Sum-rate maximization problem over Interference Broadcast Channels: `MIMO-SR-FLEXA` versus `MIMO-WMMSE`. The figures are taken from [210].

### II.5.2 LASSO problem

Consider the LASSO problem in the following form [235] (cf. Sec. II.1.1):

$$\underset{\mathbf{x}}{\text{minimize}}\ V(\mathbf{x}) \triangleq \frac{1}{2}\|\mathbf{z} - \mathbf{A}\mathbf{x}\|^2 + \lambda\|\mathbf{x}\|_1, \qquad (154)$$

where $\mathbf{A} \in \mathbb{R}^{q \times m}$ is the matrix whose columns $\mathbf{a}_i$ are the prediction or feature vectors; $z_i$ is the response variable associated to $\mathbf{a}_i$; and $\lambda > 0$ is the regularization weight.

**FLEXA for LASSO.** Observing that the univariate instance of (154) has a closed form solution, it is convenient to decompose $\mathbf{x}$ in *scalar* components ($m_i = 1$, for all $i \in N$) and update them in parallel. In order to exploit the quadratic structure of $V$ in (154) a natural choice for the surrogate function is (102). Therefore, the subproblem associated with the scalar $x_i$ reads: given $\mathbf{x}^k$,

$$\widehat{x}_i(\mathbf{x}^k) \triangleq \underset{x_i \in \mathbb{R}}{\text{argmin}} \left\{ \frac{1}{2} \left\|\mathbf{r}_i^k - \mathbf{a}_i x_i\right\|^2 + \frac{\tau_i}{2} \cdot (x_i - x_i^k)^2 + \lambda \cdot |x_i| \right\},$$

where the residual $\mathbf{r}_i^k$ is defined as

$$\mathbf{r}_i^k \triangleq \mathbf{z} - \sum_{j \neq i} \mathbf{a}_j x_j^k.$$

Invoking the first order optimality conditions (c.f. Def. I.9) [we write $\widehat{x}_i$ for $\widehat{x}_i(\mathbf{x}^k)$]:

$$-\left(\mathbf{a}_i^T \mathbf{r}_i^k + \tau_i x_i^k\right) + \left(\tau_i + \|\mathbf{a}_i\|^2\right) \widehat{x}_i + \lambda\, \partial|\widehat{x}_i| \ni 0,$$

and the expression of $\partial|x|$ [cf. (9)], one can readily obtain the closed form expression of $\widehat{x}_i(\mathbf{x}^k)$, that is,

$$\widehat{x}_i(\mathbf{x}^k) = \frac{1}{\tau_i + \|\mathbf{a}_i\|^2} \cdot S_\lambda\left(\mathbf{a}_i^T \mathbf{r}_i^k + \tau_i x_i^k\right), \qquad (155)$$

where $S_\lambda(\bullet)$ is the soft-thresholding operator, defined in (54).

We consider the instance of Algorithm 6, with the following choice of the free parameters:

- *Exact solution* $\widehat{x}_i(\mathbf{x}^k)$: In Step 3 we use the best-response $\widehat{x}_i(\mathbf{x}^k)$ as defined in (155), that is, $z_i^k = \widehat{x}_i(\mathbf{x}^k)$ (exact solution).
- *Proximal weights* $\tau_i$: While in the proposed algorithmic framework we considered fixed values of $\tau_i$, varying $\tau_i$ a finite number of times does not affect the theoretical convergence properties of the algorithms. We found that the following choices work well in practice: (i) $\tau_i$ are initially all set to $\tau_i = \text{tr}(\mathbf{A}^T\mathbf{A})/(2m)$, i.e., to half of the mean of the eigenvalues of $\nabla^2 F$; (ii) all $\tau_i$ are doubled if at a certain iteration the objective function does not decrease; and (iii) they are all halved if the objective function decreases for ten consecutive iterations or the relative error on the objective function $\texttt{re}(\mathbf{x})$ is sufficiently small, specifically if

$$\mathrm{re}(\mathbf{x}) \triangleq \frac{V(\mathbf{x}) - V^*}{V^*} \leq 10^{-2}, \tag{156}$$

where $V^*$ is the optimal value of the objective function $V$ (in our experiments on LASSO, $V^*$ is known). In order to avoid increments in the objective function, whenever all $\tau_i$ are doubled, the associated iteration is discarded, and in Step 4 of Algorithm 6 it is set $\mathbf{x}^{k+1} = \mathbf{x}^k$. In any case we limited the number of possible updates of the values of $\tau_i$ to 100.

- *Step-size $\gamma^k$*: The step-size $\gamma^k$ is updated according to the following rule:

$$\gamma^k = \gamma^{k-1}\left(1 - \min\left\{1, \frac{10^{-4}}{\mathrm{re}(\mathbf{x}^k)}\right\} \theta \gamma^{k-1}\right), \quad k = 1, \ldots, \tag{157}$$

with $\gamma^0 = 0.9$ and $\theta = 1e - 7$. The above diminishing rule is based on (108) while guaranteeing that $\gamma^k$ does not become too close to zero before the relative error is sufficiently small.

- *Greedy selection rule $S^k$*: In Step 2, we use the following greedy selection rule (satisfying Assumption II.10.2):

$$S^k = \{i : E_i(\mathbf{x}^k) \geq \sigma \cdot M^k\}, \quad \text{with} \quad E_i(\mathbf{x}^k) = |\widehat{x}_i(\mathbf{x}^k) - x_i^k|.$$

In our tests we consider two options for $\sigma$, namely: i) $\sigma = 0$, which leads to a *fully parallel* scheme wherein at each iteration *all* variables are updated; and ii) $\sigma = 0.5$, which corresponds to updating only a subset of all the variables at each iteration. Note that for both choices of $\sigma$, the resulting set $S^k$ satisfies the requirement in Step 2 of Algorithm 6; indeed, $S^k$ always contains the index $i$ corresponding to the largest $E_i(\mathbf{x}^k)$. We will refer to these two instances of the algorithm as FLEXA $\sigma = 0$ and FLEXA $\sigma = 0.5$.

**Algorithms in the literature**: We compared the above versions of FLEXA with the most competitive parallel and sequential (Block MM) algorithms proposed in the literature to solve the LASSO problem. More specifically, we consider the following schemes.

- FISTA: The Fast Iterative Shrinkage-Thresholding Algorithm (FISTA) proposed in [7] is a first order method and can be regarded as the benchmark algorithm for LASSO problems. Building on the separability of the terms in the objective function $V$, this method can be easily parallelized and thus take advantage of a parallel architecture. We implemented the parallel version that use a backtracking procedure to estimate the Lipschitz constant of $\nabla F$ [7].
- SpaRSA: This is the first order method proposed in [251]; it is a popular spectral projected gradient method that uses a spectral step length together with a nonmonotone line search to enhance convergence. Also this method can be easily parallelized, which is the version implemented in our tests. In all the experiments, we set the parameters of SpaRSA as in [251]: $M = 5$, $\sigma = 0.01$, $\alpha_{\max} = 1e30$, and $\alpha_{\min} = 1e - 30$.
- GRock & Greedy-1BCD: GRock is a parallel algorithm proposed in [187] that performs well on sparse LASSO problems. We tested the instance of GRock where the number of variables simultaneously updated is equal to the number of

Title Suppressed Due to Excessive Length        83the parallel processors. It is important to remark that the theoretical convergence properties of `GRock` are in jeopardy as the number of variables updated in parallel increases; roughly speaking, `GRock` is guaranteed to converge if the columns of the data matrix **A** in the LASSO problem are "almost" orthogonal, a feature that in general is not satisfied by real data. A special instance with convergence guaranteed is the one where only one block per time (chosen in a greedy fashion) is updated; we refer to this special case as `greedy-1BCD`.

• `Parallel ADMM`: This is a classical Alternating Method of Multipliers (ADMM). We implemented the parallel version proposed in [67].

In the implementation of the parallel algorithms, the data matrix **A** of the LASSO problem is generated as follows. Each processor generates a slice of the matrix itself such that $\mathbf{A} = [\mathbf{A}_1 \mathbf{A}_2 \cdots \mathbf{A}_P]$, where $P$ is the number of parallel processors, and each $\mathbf{A}_i$ has $m/P$ columns. Thus the computation of each product $\mathbf{Ax}$ (which is required to evaluate $\nabla F$) and the norm $\|\mathbf{x}\|_1$ (that is $G$) is divided into the parallel jobs of computing $\mathbf{A}_i \mathbf{x}_i$ and $\|\mathbf{x}_i\|_1$, followed by a reducing operation.

**Numerical examples.** We generated six groups of LASSO problems using the random generator proposed by Nesterov [172], which permits to control the sparsity of the solution. For the first five groups, we considered problems with 10,000 variables and matrices **A** with 9,000 rows. The five groups differ in the degree of sparsity of the solution, namely: the percentage of non zeros in the solution is 1%, 10%, 20%, 30%, and 40%, respectively. The last group is formed by instances with 100,000 variables and 5000 rows for **A**, and solutions having 1% of non zero variables. In all experiments and for all the algorithms, the initial point was set to the zero vector.

Results of the experiments for the 10,000 variables groups are reported in Fig. II.5, where we plot the relative error as defined in (156) versus the CPU time; all the curves are obtained using (up to) 40 cores, and averaged over ten independent random realizations. Note that the CPU time includes communication times (for parallel algorithms) and the initial time needed by the methods to perform all pre-iteration computations (this explains why the curves of ADMM start after the others; in fact ADMM requires some nontrivial initializations). For one instance, the one corresponding to 1% of the sparsity of the solution, we plot also the relative error versus iterations [Fig. II.5(a2)]; similar behaviors of the algorithms have been observed also for the other instances, and thus are not reported. Results for the LASSO instance with 100,000 variables are plotted in Fig. II.6. The curves are averaged over five random realizations.

The following comments are in order. On all the tested problems, `FLEXA` $\sigma = 0.5$ outperforms in a consistent manner all the other implemented algorithms. In particular, as the sparsity of the solution decreases, the problems become harder and the selective update operated by `FLEXA` $\sigma = 0.5$ improves over `FLEXA` $\sigma = 0$, where instead all variables are updated at each iteration. `FISTA` is capable to approach relatively fast low accuracy when the solution is not too sparse, but has difficulties in reaching high accuracy. `SpaRSA` seems to be very insensitive to the degree of sparsity of the solution; it behaves well on 10,000 variables problems and not too sparse solutions, but is much less effective on very large-scale problems. The version of `GRock` with $P = 40$ is the closest match to `FLEXA`, but only when the problems are very sparse (but it is not supported by a convergence theory on our test problems). This is consistent with the fact that its convergence properties are at stake when the problems are quite dense. Furthermore, if the problem is very large, updating only 40 variables at each iteration, as `GRock` does, could slow down the convergence, especially when the optimal solution is not very sparse. From this point of view, `FLEXA` $\sigma = 0.5$ seems to strike a good balance between not updating variables that are probably zero at the optimum and nevertheless update a sizeable amount of variables when needed in order to enhance convergence.

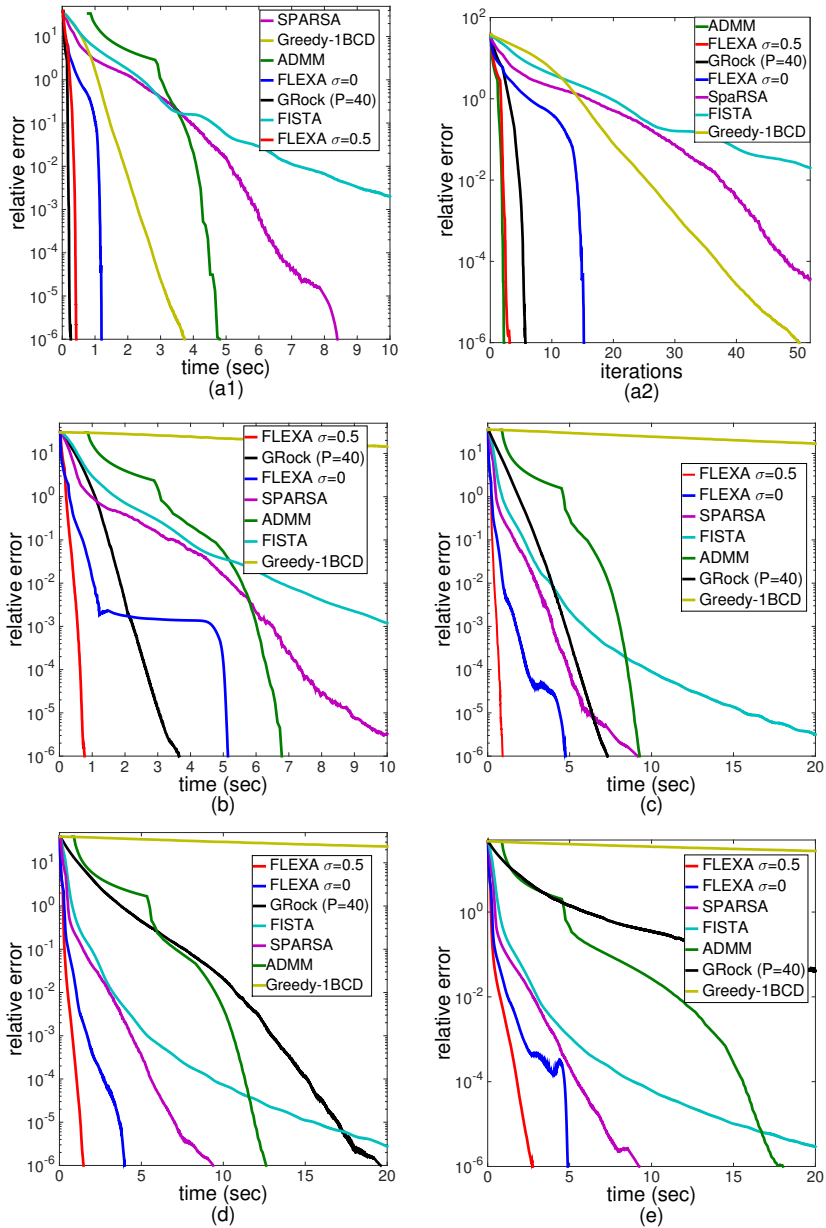

*Fig. II.5: LASSO problem* (154) *with 10,000 variables; relative error vs. time (in seconds) for: (a1) 1% non zeros - (b) 10% non zeros - (c) 20% non zeros - (d) 30% non zeros - (e) 40% non zeros; (a2) relative error vs. iterations for 1% non zeros. The figures are taken from [79].*

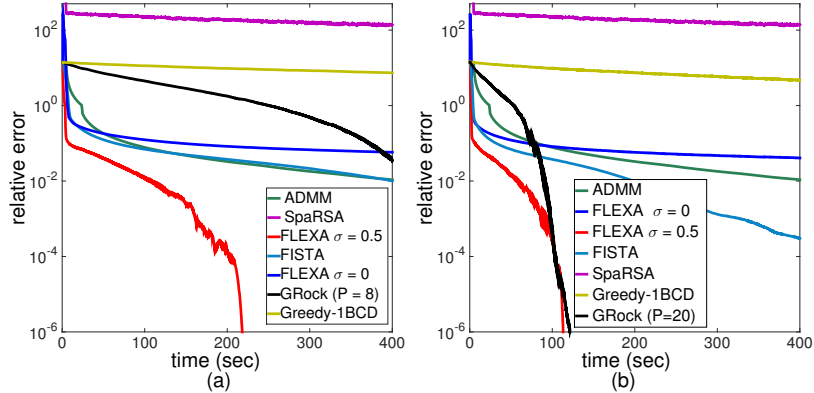

*Fig. II.6: LASSO problem* (154) *with* $10^5$ *variables; Relative error vs. time for: (a) 8 cores - (b) 20 cores. The figures are taken from [79].*

**Remark II.21 (On the parallelism).** *Fig. II.6 shows that* `FLEXA` *seems to exploit well parallelism on LASSO problems. Indeed, when passing from 8 to 20 cores, the running time approximately halves. This kind of behavior has been consistently observed also for smaller problems and different number of cores (not reported here). Note that practical speed-up due to the use of a parallel architecture is given by several factor that are not easily predictable, including communication times among the cores, the data format, etc.. Here we do not pursue a theoretical study of the speed-up but refer to [37] for some initial study. We finally observe that* `GRock` *appears to improve greatly with the number of cores. This is due to the fact that in* `GRock` *the maximum number of variables that is updated in parallel is exactly equal to the number of cores (i.e., the degree of parallelism), and this might become a serious drawback on very large problems (on top of the fact that convergence is in jeopardy). On the contrary, the theory presented in this chapter permits the parallel update of* any *number of variables while guaranteeing convergence.*

**Remark II.22 (On selective updates).** *It is interesting to comment why* `FLEXA` $\sigma = 0.5$ *behaves better than* `FLEXA` $\sigma = 0$. *To understand the reason behind this phenomenon, we first note that Algorithm 6 has the remarkable capability to* iden-tify those variables that will be zero at a solution; *we do not provide here the proof of this statement but only an informal argument. Roughly speaking, it can be shown* *that, for k large enough, those variables that are zero in $\widehat{\mathbf{x}}(\mathbf{x}^k) = (\widehat{x}_i(\mathbf{x}^k))_{i=1}^m$ will be zero also in a limiting solution $\bar{\mathbf{x}}$. Therefore, suppose that k is large enough so that this identification property already takes place (we will say that "we are in the identification phase") and consider an index i such that $\bar{x}_i = 0$. Then, if $x_i^k$ is zero, it is clear, by Steps 3 and 4, that $x_i^{k'}$ will be zero for all indices $k' > k$, independently of whether i belongs to $S^k$ or not. In other words, if a variable that is zero at the solution is already zero when the algorithm enters the identification phase,* that variable will be zero in all subsequent iterations; *this fact, intuitively, should enhance the convergence speed of the algorithm. Conversely, if when we enter the identification phase $x_i^k$ is not zero, the algorithm will have to bring it back to zero iteratively. This explains why updating only variables that we have "strong" reason to believe will be non zero at a solution is a better strategy than updating them all. Of course, there may be a problem dependence and the best value of $\sigma$ can vary from problem to problem. But the above explanation supports the idea that it might be wise to "waste" some calculations and perform only a partial ad-hoc update of the variables.*

### II.5.3 The logistic regression problem

Consider the logistic regression problem in the following form [235]:

$$\underset{\mathbf{x}}{\text{minimize}} \quad V(\mathbf{x}) = \sum_{i=1}^{q} \log(1 + e^{-w_i \cdot \mathbf{z}_i^T \mathbf{x}}) + \lambda \|\mathbf{x}\|_1, \tag{158}$$

where $\mathbf{z}_i \in \mathbb{R}^m$ is the feature vector of sample $i$, with the associated label $w_i \in \{-1, 1\}$.

**FLEXA for logistic regression.** Problem (158) is a highly nonlinear problem involving many exponentials that, notoriously, gives rise to numerical difficulties. Because of these high nonlinearities, a Gauss-Seidel (sequential) approach is expected to be more effective than a pure Jacobi (parallel) method, a fact that was confirmed by the experiments in [79]. For this reason, for the logistic regression problem we tested both `FLEXA` and the hybrid scheme in Algorithm 9, which will term `GJ-FLEXA`. The setting of the free parameters in GJ-FLEXA is essentially the same as the one described for LASSO (cf. Sec. II.5.2), but with the following differences:

- *Exact solution $\widehat{x}_i(\mathbf{x}^k)$:* The surrogate function $\widetilde{F}_i$ is chosen as the second order approximation of the original function $F$: given the current iterate $\mathbf{x}^k$,

$$\widetilde{F}_i(x_i \mid \mathbf{x}^k) = F(\mathbf{x}^k) + \nabla_{x_i} F(\mathbf{x}^k) \cdot (x_i - x_i^k) + \frac{1}{2}(\nabla_{x_i x_i}^2 F(\mathbf{x}^k) + \tau_i) \cdot (x_i - x_i^k)^2 + \lambda \cdot |x_i|,$$

which leads to the following closed form solution for $\widehat{x}_i(\mathbf{x}^k)$:

$$\widehat{x}_i(\mathbf{x}^k) = S_{\lambda \cdot t_i^k}\left(x_i^k - t_i^k \cdot \nabla_{x_i} F(\mathbf{x}^k)\right), \quad \text{with} \quad t_i^k \triangleq \left(\tau_i + \nabla_{x_i x_i}^2 F(\mathbf{x}^k)\right)^{-1},$$

where $S_\lambda(\bullet)$ is the soft-thresholding operator, defined in (54).

- *Proximal weights* $\tau_i$: The initial $\tau_i$ are set to $\mathrm{tr}(\mathbf{Z}^T\mathbf{Z})/(2m)$, for all $i$, where $m$ is the total number of variables and $\mathbf{Z} = [\mathbf{z}_1\,\mathbf{z}_2\,\cdots\,\mathbf{z}_q]^T$.
- *Step-size* $\gamma^k$: We use the step-size rule (157). However, since the optimal value $V^*$ is not known for the logistic regression problem, $\mathrm{re}(\mathbf{x})$ can no longer be computed. We replace $\mathrm{re}(\mathbf{x})$ with the merit function $\|\mathbf{M}(\mathbf{x})\|_\infty$, with

$$\mathbf{M}(\mathbf{x}) \triangleq \nabla F(\mathbf{x}) - \Pi_{[-\lambda,\lambda]^m}\left(\nabla F(\mathbf{x}) - \mathbf{x}\right).$$

Here the projection $\Pi_{[-\lambda,\lambda]^m}(\mathbf{y})$ can be efficiently computed; it acts componentwise on $\mathbf{y}$, since $[-\lambda,\lambda]^m = [-\lambda,\lambda]\times\cdots\times[-\lambda,\lambda]$. Note that $\mathbf{M}(\mathbf{x})$ is a valid optimality measure function; indeed, it is continuous and $\mathbf{M}(\mathbf{x}) = \mathbf{0}$ is equivalent to the standard necessary optimality condition for Problem (94), see [31].

**Algorithms in the literature**: We compared FLEXA ($\sigma = 0.5$) (cf. Sec. II.5.2) and GJ-FLEXA with the other parallel algorithms introduced in Sec. II.5.2 for the LASSO problem (whose tuning of the free parameters is the same as in Fig. II.5 and Fig. II.6), namely: FISTA, SpaRSA, and GRock. For the logistic regression problem, we also tested one more algorithm, that we call CDM. This Coordinate Descent Method is an extremely efficient Gauss-Seidel-type method (customized for logistic regression), and is part of the LIBLINEAR package available at http://www.csie.ntu.edu.tw/~cjlin/.

We tested the aforementioned algorithms on three instances of the logistic regression problem that are widely used in the literature, and whose essential data features are given in Table II.3; we downloaded the data from the LIBSVM repository http://www.csie.ntu.edu.tw/~cjlin/libsvm/, which we refer to for a detailed description of the test problems. In our implementation, the matrix $\mathbf{Z}$ is column-wise partitioned according to $\mathbf{Z} = \begin{bmatrix}\tilde{\mathbf{Z}}_1\,\tilde{\mathbf{Z}}_2\,\cdots\,\tilde{\mathbf{Z}}_P\end{bmatrix}$ and distributively stored across $P$ processors, where $\tilde{\mathbf{Z}}_i$ is the set of columns of $\mathbf{Z}$ owned by processor $i$.

In Fig. II.7, we plotted the relative error vs. the CPU time (the latter defined as in Fig. II.5 and Fig. II.6) achieved by the aforementioned algorithms for the three datasets, and using a different number of cores, namely: 8, 16, 20, 40; for each algorithm but GJ-FLEXA we report only the best performance over the aforementioned numbers of cores. Note that in order to plot the relative error, we had to preliminary estimate $V^*$ (which is not known for logistic regression problems). To do so, we ran GJ-FLEXA until the merit function value $\|\mathbf{M}(\mathbf{x}^k)\|_\infty$ went below $10^{-7}$, and used the corresponding value of the objective function as estimate of $V^*$. We remark that we used this value only to plot the curves. Next to each plot, we also reported the overall FLOPS counted up till reaching the relative errors as indicated in the table. Note that the FLOPS of GRock on real-sim and rcv1 are those counted in 24 hours simulation time; when terminated, the algorithm achieved a relative error that was still very far from the reference values set in our experiment. Specifically, GRock reached 1.16 (instead of $1e-4$) on real-sim and 0.58 (instead of $1e-3$) on rcv1; the counted FLOPS up till those error values are still reported in the tables.

The analysis of the figures shows that, due to the high nonlinearities of the objective function, Gauss-Seidel-type methods outperform the other schemes. In spite of

88	Gesualdo Scutari and Ying Sun

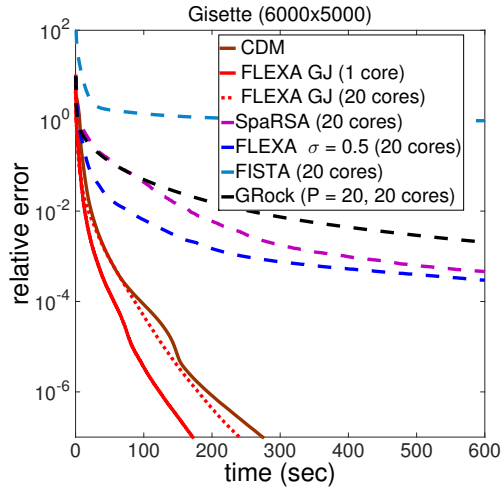

| Algorithms | FLOPS (1e-2/1e-6) |
|---|---|
| GJ-FLEXA (1C) | 1.30e+10/1.23e+11 |
| GJ-FLEXA (20C) | 5.18e+11/5.07e+12 |
| FLEXA $\sigma=0.5$ (20C) | 1.88e+12/4.06e+13 |
| CDM | 2.15e+10/1.68e+11 |
| SpaRSA (20C) | 2.20e+12/5.37e+13 |
| FISTA (20C) | 3.99e+12/5.66e+13 |
| GRock (20C) | 7.18e+12/1.81e+14 |

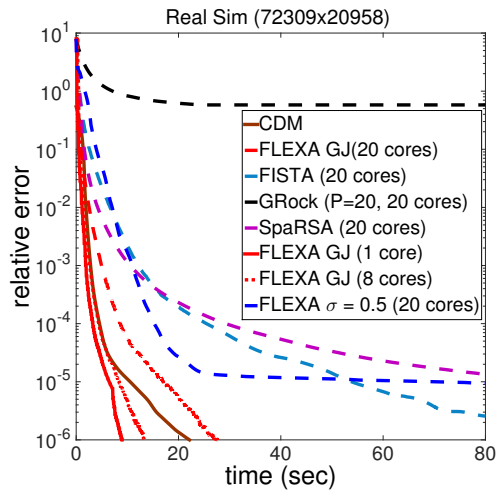

| Algorithms | FLOPS (1e-4/1e-6) |
|---|---|
| GJ-FLEXA (1C) | 2.76e+9/6.60e+9 |
| GJ-FLEXA (20C) | 9.83e+10/2.85e+11 |
| FLEXA $\sigma=0.5$ (20C) | 3.54e+10/4.69e+11 |
| CDM | 4.43e+9/2.18e+10 |
| SpaRSA (20C) | 7.18e+9/1.94e+11 |
| FISTA (20C) | 3.91e+10/1.56e+11 |
| GRock (20C) | 8.30e+14 (after 24h) |

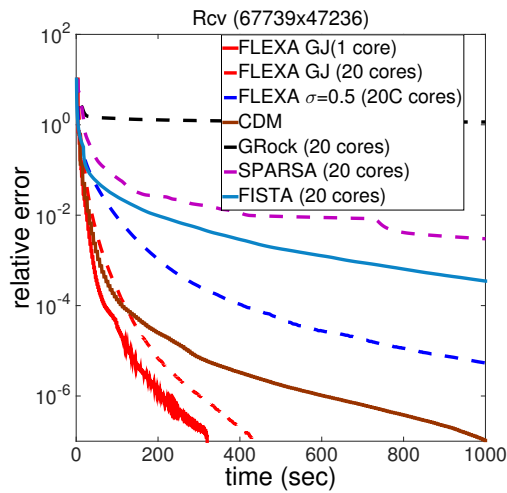

| Algorithms | FLOPS (1e-3/1e-6) |
|---|---|
| GJ-FLEXA (1C) | 3.61e+10/2.43e+11 |
| GJ-FLEXA (20C) | 1.35e+12/6.22e+12 |
| FLEXA $\sigma=0.5$ (20C) | 8.53e+11/7.19e+12 |
| CDM | 5.60e+10/6.00e+11 |
| SpaRSA (20C) | 9.38e+12/7.20e+13 |
| FISTA (20C) | 2.58e+12/2.76e+13 |
| GRock (20C) | 1.72e+15 (after 24h) |

Fig. II.7: Logistic Regression problem (158): Relative error vs. time (in seconds) and FLOPS for i) gisette, ii) real-sim, and iii) rcv. The figures are taken from [79].

| Data set | $q$ | $m$ | $\lambda$ |
|---|---|---|---|
| gisette (scaled) | 6000 | 5000 | 0.25 |
| real-sim | 72309 | 20958 | 4 |
| rcv1 | 677399 | 47236 | 4 |

*Table II.3: Data sets for the logistic regression tests [Problem (158)]*

this, FLEXA still behaves quite well. But GJ-FLEXA with one core, thus a non parallel method, clearly outperforms all other algorithms. The explanation can be the following. GJ-FLEXA with one core is essentially a Gauss-Seidel-type method but with two key differences: the use of a stepsize and more importantly a (greedy) selection rule by which only some variables are updated at each round. As the number of cores increases, the algorithm gets "closer and closer" to a Jacobi-type method, and because of the high nonlinearities, moving along a "Jacobi direction" does not bring improvements. In conclusion, for logistic regression problems, our experiments suggests that while the (opportunistic) selection of variables to update seems useful and brings to improvements even in comparison to the extremely efficient, dedicated CDM algorithm/software, parallelism (at least, in the form embedded in our scheme), does not appear to be beneficial as instead observed for LASSO problems.

## II.6. Appendix

### II.6.1 Proof of Lemma II.4

The continuity of $\widehat{\mathbf{x}}(\bullet)$ follows readily from [199]; see also [104].

We prove next the Lipschitz continuity of $\widehat{\mathbf{x}}(\bullet)$, under the additional assumption that $G$ is separable. Let $\mathbf{x}_i, \mathbf{z}_i \in X_i$. Invoking the optimality conditions of $\widehat{\mathbf{x}}_i(\mathbf{x})$ and $\widehat{\mathbf{x}}_i(\mathbf{z})$, we have

$$(\mathbf{y}_1 - \widehat{\mathbf{x}}_i(\mathbf{x}))^T (\nabla \widetilde{F}_i(\widehat{\mathbf{x}}_i(\mathbf{x}) \,|\, \mathbf{x})) + g_i(\mathbf{y}_1) - g_i(\widehat{\mathbf{x}}_i(\mathbf{x})) \geq 0, \quad \forall \mathbf{y}_1 \in X_i,$$
$$(\mathbf{y}_2 - \widehat{\mathbf{x}}_i(\mathbf{z}))^T (\nabla \widetilde{F}_i(\widehat{\mathbf{x}}_i(\mathbf{z}) \,|\, \mathbf{z})) + g_i(\mathbf{y}_2) - g_i(\widehat{\mathbf{x}}_i(\mathbf{z})) \geq 0, \quad \forall \mathbf{y}_2 \in X_i.$$

Letting $\mathbf{y}_1 = \widehat{\mathbf{x}}_i(\mathbf{z})$ and $\mathbf{y}_2 = \widehat{\mathbf{x}}_i(\mathbf{x})$ and summing the two inequalities above, we obtain

$$(\widehat{\mathbf{x}}_i(\mathbf{z}) - \widehat{\mathbf{x}}_i(\mathbf{x}))^T \left( \nabla \widetilde{F}_i(\widehat{\mathbf{x}}_i(\mathbf{x}) \,|\, \mathbf{x}) - \nabla \widetilde{F}_i(\widehat{\mathbf{x}}_i(\mathbf{z}) \,|\, \mathbf{z}) \right) \geq 0.$$

Adding and subtracting $\nabla \widetilde{F}_i(\widehat{\mathbf{x}}_i(\mathbf{z}) \,|\, \mathbf{x})$ and using the uniform strongly convexity of $\widetilde{F}_i$ with respect to its first argument (cf. Assumption II.2.1) and the Lipschitz continuity of $\nabla \widetilde{F}_i$ with respect to its second argument (cf. Assumption II.3) yield

$$\tau_i \|\widehat{\mathbf{x}}_i(\mathbf{z}) - \widehat{\mathbf{x}}_i(\mathbf{x})\|^2 \leq (\widehat{\mathbf{x}}_i(\mathbf{z}) - \widehat{\mathbf{x}}_i(\mathbf{x}))^T (\nabla \widetilde{F}_i(\widehat{\mathbf{x}}_i(\mathbf{z}) \,|\, \mathbf{x}) - \nabla \widetilde{F}_i(\widehat{\mathbf{x}}_i(\mathbf{z}) \,|\, \mathbf{z}))$$
$$\leq \widetilde{L}_i \|\widehat{\mathbf{x}}_i(\mathbf{z}) - \widehat{\mathbf{x}}_i(\mathbf{x})\| \cdot \|\mathbf{x} - \mathbf{z}\|.$$

Therefore, $\widehat{\mathbf{x}}_i(\bullet)$ is Lipschitz continuous on $X$ with constant $\hat{L}_i \triangleq \widetilde{L}_i/\tau_i$. □

### II.6.2 Proof of Lemma II.17

The proof is adapted by [56, Lemma 10] and reported here for completeness.

With a slight abuse of notation, we will use $(\mathbf{x}_i, \mathbf{x}_j, \mathbf{y}_{-(i,j)})$, with $i < j$, to denote the ordered tuple $(\mathbf{y}_1, \ldots, \mathbf{y}_{i-1}, \mathbf{x}_i, \mathbf{y}_{i+1}, \ldots, \mathbf{y}_{j-1}, \mathbf{x}_j, \mathbf{y}_{j+1}, \ldots, \mathbf{y}_n)$.

Given $k \geq 0$, $S^k \subseteq N$, and $\gamma^k \leq 1/n$, let $\bar{\gamma}^k = \gamma^k n \leq 1$. Define $\check{\mathbf{x}}^k \triangleq (\check{\mathbf{x}}_i^k)_{i \in N}$, with $\check{\mathbf{x}}_i^k = \mathbf{x}_i^k$ if $i \notin S^k$, and

$$\check{\mathbf{x}}_i^k \triangleq \bar{\gamma}^k \widehat{\mathbf{z}}_i^k + (1 - \bar{\gamma}^k) \mathbf{x}_i^k, \tag{159}$$

otherwise. Then $\mathbf{x}^{k+1}$ in Step 4 of the algorithm can be written as

$$\mathbf{x}^{k+1} = \frac{n-1}{n} \mathbf{x}^k + \frac{1}{n} \check{\mathbf{x}}^k. \tag{160}$$

Using (160) and invoking the convexity of $G$, the following recursion holds for all $k$:

$$\begin{aligned}
G(\mathbf{x}^{k+1}) &= G\left(\frac{1}{n}(\check{\mathbf{x}}_1^k, \mathbf{x}_{-1}^k) + \frac{1}{n}(\mathbf{x}_1^k, \check{\mathbf{x}}_{-1}^k) + \frac{n-2}{n}\mathbf{x}^k\right) \\
&= G\left(\frac{1}{n}(\check{\mathbf{x}}_1^k, \mathbf{x}_{-1}^k) + \frac{n-1}{n}\left(\mathbf{x}_1^k, \frac{1}{n-1}\check{\mathbf{x}}_{-1}^k + \frac{n-2}{n-1}\mathbf{x}_{-1}^k\right)\right) \\
&\leq \frac{1}{n} G\left(\check{\mathbf{x}}_1^k, \mathbf{x}_{-1}^k\right) + \frac{n-1}{n} G\left(\mathbf{x}_1^k, \frac{1}{n-1}\check{\mathbf{x}}_{-1}^k + \frac{n-2}{n-1}\mathbf{x}_{-1}^k\right) \\
&= \frac{1}{n} G\left(\check{\mathbf{x}}_1^k, \mathbf{x}_{-1}^k\right) + \frac{n-1}{n} G\left(\frac{1}{n-1}\left(\mathbf{x}_1^k, \check{\mathbf{x}}_{-1}^k\right) + \frac{n-2}{n-1}\mathbf{x}^k\right) \\
&= \frac{1}{n} G\left(\check{\mathbf{x}}_1^k, \mathbf{x}_{-1}^k\right) + \frac{n-1}{n} G\left(\frac{1}{n-1}\left(\check{\mathbf{x}}_2^k, \mathbf{x}_{-2}^k\right)\right. \\
&\qquad\qquad\qquad\qquad\qquad \left.+ \frac{1}{n-1}\left(\mathbf{x}_1^k, \mathbf{x}_2^k, \check{\mathbf{x}}_{-(1,2)}^k\right) + \frac{n-3}{n-1}\mathbf{x}^k\right) \quad (161) \\
&= \frac{1}{n} G\left(\check{\mathbf{x}}_1^k, \mathbf{x}_{-1}^k\right) + \frac{n-1}{n} G\left(\frac{1}{n-1}\left(\check{\mathbf{x}}_2^k, \mathbf{x}_{-2}^k\right)\right. \\
&\qquad\qquad + \left.\frac{n-2}{n-1}\left(\mathbf{x}_1^k, \mathbf{x}_2^k, \frac{1}{n-2}\check{\mathbf{x}}_{-(1,2)}^k + \frac{n-3}{n-2}\mathbf{x}_{-(1,2)}^k\right)\right) \\
&\leq \frac{1}{n} G\left(\check{\mathbf{x}}_1^k, \mathbf{x}_{-1}^k\right) + \frac{1}{n} G\left(\check{\mathbf{x}}_2^k, \mathbf{x}_{-2}^k\right) \\
&\quad + \frac{n-2}{n} G\left(\mathbf{x}_1^k, \mathbf{x}_2^k, \frac{1}{n-2}\check{\mathbf{x}}_{-(1,2)}^k + \frac{n-3}{n-2}\mathbf{x}_{-(1,2)}^k\right) \\
&\leq \quad \cdots \quad \leq \frac{1}{n}\sum_{i \in N} G(\check{\mathbf{x}}_i^k, \mathbf{x}_{-i}^k).
\end{aligned}$$

Using (161), the difference of $G(\mathbf{x}^{k+1})$ and $G(\mathbf{x}^k)$ can be bounded as

$$\begin{aligned}
G(\mathbf{x}^{k+1}) - G(\mathbf{x}^k) &\leq \frac{1}{n} \sum_{i \in N} \left( G(\check{\mathbf{x}}_i^k, \mathbf{x}_{-i}^k) - G(\mathbf{x}^k) \right) \\
&= \frac{1}{n} \sum_{i \in S^k} \left( G(\check{\mathbf{x}}_i^k, \mathbf{x}_{-i}^k) - G(\mathbf{x}^k) \right) \\
&\leq \frac{1}{n} \sum_{i \in S^k} \left( \bar{\gamma}^k G(\widehat{\mathbf{z}}_i^k, \mathbf{x}_{-i}^k) + (1-\bar{\gamma}^k) G(\mathbf{x}^k) - G(\mathbf{x}^k) \right) \\
&= \gamma^k \sum_{i \in S^k} \left( G(\widehat{\mathbf{z}}_i^k, \mathbf{x}_{-i}^k) - G(\mathbf{x}^k) \right).
\end{aligned} \qquad (162)$$

□

### II.6.3 Proof of Lemma II.19

The proof can be found in [79, Lemma 10], and reported here for completeness.

For notational simplicity, we will write $\mathbf{x}_{S^k}$ for $(\mathbf{x})_{S^k}$ [recall that $(\mathbf{x})_{S^k}$ denotes the vector whose block component $i$ is equal to $\mathbf{x}_i$ if $i \in S^k$, and zero otherwise]. Let $j_k$ be an index in $S^k$ such that $E_{j_k}(\mathbf{x}^k) \geq \rho \max_i E_i(\mathbf{x}^k)$ (cf. Assumption II.10.2). Then, by the error bound condition (110) it is easy to check that the following chain of inequalities holds:

$$\begin{aligned}
\bar{s}_{j_k} \|\widehat{\mathbf{x}}_{S^k}(\mathbf{x}^k) - \mathbf{x}_{S^k}^k\| &\geq \bar{s}_{j_k} \|\widehat{\mathbf{x}}_{j_k}(\mathbf{x}^k) - \mathbf{x}_{j_k}^k\| \\
&\geq E_{j_k}(\mathbf{x}^k) \\
&\geq \rho \max_i E_i(\mathbf{x}^k) \\
&\geq \left( \rho \min_i \underline{s}_i \right) \left( \max_i \{\|\widehat{\mathbf{x}}_i(\mathbf{x}^k) - \mathbf{x}_i^k\|\} \right) \\
&\geq \left( \frac{\rho \min_i \underline{s}_i}{n} \right) \|\widehat{\mathbf{x}}(\mathbf{x}^k) - \mathbf{x}^k\|.
\end{aligned}$$

Hence we have for any $k$,

$$\|\widehat{\mathbf{x}}_{S^k}(\mathbf{x}^k) - \mathbf{x}_{S^k}^k\| \geq \left( \frac{\rho \min_i \underline{s}_i}{n \bar{s}_{j_k}} \right) \|\widehat{\mathbf{x}}(\mathbf{x}^k) - \mathbf{x}^k\| \geq \left( \frac{\rho \min_i \underline{s}_i}{n \max_j \bar{s}_j} \right) \|\widehat{\mathbf{x}}(\mathbf{x}^k) - \mathbf{x}^k\|.$$

□

## II.7. Sources and Notes

Although parallel (deterministic and stochastic) block-methods have a long history (mainly for convex problems), recent years have witnessed a revival of such methods and their (probabilistic) analysis; this is mainly due to the current trend towards huge scale optimization and the availability of ever more complex computational architectures that call for efficient, fast, and resilient algorithms. The literature is vast and a comprehensive overview of current methods goes beyond the scope of this commentary. Here we only focus on SCA-related methods and refer to [25, 214, 250] (and references therein) as entry point to other numerical optimization algorithms.

**Parallel SCA-related Methods**: The roots of parallel deterministic SCA schemes (wherein *all* the variables are updated simultaneously) can be traced back at least to the work of Cohen on the so-called auxiliary principle [51, 52] and its related developments, see e.g. [7, 28, 89, 162, 165, 172, 183, 187, 197, 210, 238, 251]. Roughly speaking, these works can be divided in two groups, namely: parallel solution methods for *convex* objective functions [7, 28, 51, 52, 165, 187, 197] and *nonconvex* ones [89, 162, 172, 183, 210, 238, 251]. All methods in the former group (and [89, 162, 172, 238, 251]) are (proximal) gradient schemes; they thus share the classical drawbacks of gradient-like schemes; moreover, by replacing the convex function $F$ with its first order approximation, they do not take any advantage of any structure of $F$ beyond mere differentiability. Exploiting some available structural properties of $F$, instead, has been shown to enhance (practical) convergence speed, see e.g. [210]. Comparing with the second group of works [89, 162, 172, 183, 210, 238, 251], the parallel SCA algorithmic framework introduced in this lecture improves on their convergence properties while adding great flexibility in the selection of the variables to update at each iteration. For instance, with the exception of [68, 143, 187, 204, 238], all the aforementioned works do not allow parallel updates of a *subset* of all variables, a feature that instead, fully explored as we do, can dramatically improve the convergence speed of the algorithm, as shown in Sec. II.5. Moreover, with the exception of [183], they all require an Armijo-type line-search, whereas the scheme in [183] is based on diminishing step-size-rules, but its convergence properties are quite weak: not all the limit points of the sequence generated by this scheme are guaranteed to be stationary solutions of (94).

The SCA-based algorithmic framework introduced in this lecture builds on and generalizes the schemes proposed in [56, 79, 210]. More specifically, Algorithm 5 was proposed in [210] for smooth instances of Problem (94) (i.e., $G = 0$); convergence was established when constant (Assumption II.11.1) or diminishing (Assumption II.11.2) step-sizes are employed (special case of Theorem II.8). In [79], this algorithm was extended to deal with nonsmooth *separable* functions $G$ while incorporating inexact updates (Assumption II.9.1) and the greedy selection rule in Assumption II.10.2 (cf. Algorithm 6) as well as hybrid Jacobi/Gauss-Seidel updates (as described in Algorithm 8); convergence was established under Assumption II.11.2 (diminishing step-size) (special case of Theorem II.12). Finally, in [56], an instance of Algorithm 6 was proposed, to deal with *nonseparable* convex func-tions $G$, and using random block selection rules (Assumption II.10.4) or hybrid random-greedy selection rules (Algorithm 7); convergence was established when a diminishing step-size is employed (special case of Theorem II.13).

While [56, 79, 210] studied some instances of parallel SCA methods in isolation (and only for some block selection rules and step-size rules) the contribution of this lecture is to provide a broader and unified view and analysis of such methods.

**SCA Methods for nonconvex constrained optimization:** (Parallel) SCA methods have been recently extended to deal with nonconvex constraints; state-of-the-art developments can be found in [76, 77, 207] along with their applications to some representative problems in Signal Processing, Communications, and Machine Learning [208]. More specifically, consider the following generalization of Problem (94):

$$\begin{aligned} \underset{\mathbf{x}}{\text{minimize}} \quad & F(\mathbf{x}) + H(\mathbf{x}) \\ \text{s.t.} \quad & \mathbf{x} \in X, \\ & g_j(\mathbf{x}) \leq 0, \quad j = 1, \ldots, J, \end{aligned} \qquad (163)$$

where $g_j(\mathbf{x}) \leq 0$, $j = 1, \ldots, J$, represent *nonconvex nonsmooth* constraints; and $H$ is now a *nonsmooth*, possibly *nonconvex* function.

A natural extension of the SCA idea introduced in this lecture to the general class of nonconvex problems (163) is replacing all the nonconvex functions with suitably chosen convex surrogates, and solve instead the following convexified problems: given $\mathbf{x}^k$,

$$\begin{aligned} \widehat{\mathbf{x}}(\mathbf{x}^k) \triangleq \underset{\mathbf{x}}{\text{argmin}} \quad & \widetilde{F}(\mathbf{x} \,|\, \mathbf{x}^k) + \widetilde{H}(\mathbf{x} \,|\, \mathbf{x}^k) \\ \text{s.t.} \quad & \mathbf{x} \in X, \\ & \widetilde{g}_j(\mathbf{x} \,|\, \mathbf{y}) \leq 0, \quad j = 1, \ldots, J; \end{aligned} \qquad (164)$$

where $\widetilde{F}$, $\widetilde{H}$ and $\widetilde{g}_j$ are (strongly) convex surrogates for $F$, $H$, and $g_j$, respectively. The update of $\mathbf{x}^k$ is then given by

$$\mathbf{x}^{k+1} = \mathbf{x}^{k+1} + \gamma^k \left( \widehat{\mathbf{x}}(\mathbf{x}^k) - \mathbf{x}^k \right).$$

Conditions on the surrogates in (164) and convergence of the resulting SCA algorithms can be found in [207] for the case of smooth $H$ and $g_j$ (or nonsmooth DC), and in [77] for the more general setting of nonsmooth functions; *parallel* and *distributed* implementations are also discussed. Here we only mention that the surrogates $\widetilde{g}_j$ must be a global convex upper bound of the associated nonconvex $g_j$ (as for the MM algorithms–see Lecture I). This condition was removed in [76] where a different structure for the subproblems (164) was proposed. The work [76] also provides a complexity analysis of the SCA-based algorithms.

Other SCA-related methods for nonconvex constrained problem are discussed in [76, 77, 207, 208], which we refer to for details.

**Asynchronous SCA methods:** In the era of data deluge, data-intensive applications give rise to extremely large-scale problems, which naturally call for *asynchronous*, parallel solution methods. In fact, well suited to modern computational architec-tures (e.g., shared memory systems, message passing-based systems, cluster computers, cloud federations), asynchronous methods reduce the idle times of workers, mitigate communication and/or memory-access congestion, and make algorithms more fault-tolerant. Although asynchronous block-methods have a long history (see, e.g., [5, 16, 45, 87, 237]), in the past few years, the study of asynchronous parallel optimization methods has witnessed a revival of interest. Indeed, asynchronous parallelism has been applied to many state-of-the-art optimization algorithms (mainly for convex objective functions and constraints), including stochastic gradient methods [109, 137, 144, 153, 167, 184, 195] and ADMM-like schemes [105, 110, 247]. The asynchronous counterpart of BCD methods has been introduced and studied in the seminal work [146], which motivated and oriented much of subsequent research in the field, see e.g. [60, 61, 147, 185, 186].

*Asynchronous parallel SCA* methods were recently proposed and analyzed in [37, 38] for nonconvex problems in the form (94), with $G$ separable and $X_i$ possibly nonconvex [we refer to such a class of problems as Problem (94)′]. In the asynchronous parallel SCA method [37, 38], workers (e.g., cores, cpus, or machines) continuously and without coordination with each other, update a block-variable by solving a strongly convex block-model of Problem (94)′. More specifically, at iteration $k$, a worker updates a block-variable $\mathbf{x}_{i^k}^k$ of $\mathbf{x}^k$ to $\mathbf{x}_{i^k}^{k+1}$, with $i^k$ in the set $N$, thus generating the vector $\mathbf{x}^{k+1}$. When updating block $i^k$, in general, the worker does not have access to the current vector $\mathbf{x}^k$, but it will use instead the local estimate $\mathbf{x}^{k-\mathbf{d}^k} \triangleq (x_1^{k-d_1^k}, x_2^{k-d_2^k}, \ldots, x_n^{k-d_n^k})$, where $\mathbf{d}^k \triangleq (d_1^k, d_2^k, \ldots, d_n^k)$ is the "vector of delays", whose components $d_i^k$ are nonnegative integers. Note that $\mathbf{x}^{k-\mathbf{d}^k}$ is nothing else but a combination of delayed, block-variables. The way each worker forms its own estimate $\mathbf{x}^{k-\mathbf{d}^k}$ depends on the particular architecture under consideration and it is immaterial to the analysis of the algorithm; see [37]. Given $\mathbf{x}^{k-\mathbf{d}^k}$ and $i^k$, block $\mathbf{x}_{i^k}^k$ is updated by solving the following strongly convex block-approximation of Problem (94)′:

$$\hat{\mathbf{x}}_{i^k}(\mathbf{x}^{k-\mathbf{d}^k}) \triangleq \underset{\mathbf{x}_{i^k} \in \tilde{X}_{i^k}(\mathbf{x}^{k-\mathbf{d}^k})}{\operatorname{argmin}} \tilde{F}_{i^k}(\mathbf{x}_{i^k} \,|\, \mathbf{x}^{k-\mathbf{d}^k}) + g_{i^k}(\mathbf{x}_{i^k}), \tag{165}$$

and then setting

$$\mathbf{x}_{i^k}^{k+1} = \mathbf{x}_{i^k}^k + \gamma \left( \hat{\mathbf{x}}_{i^k}(\mathbf{x}^{k-\mathbf{d}^k}) - \mathbf{x}_{i^k}^k \right). \tag{166}$$

In (165), $\tilde{F}_{i^k}(\bullet \,|\, \mathbf{y})$ represents a strongly convex surrogate of $F$, and $\tilde{X}_{i^k}$ is a convex set obtained replacing the nonconvex functions defining $X_{i^k}$ by suitably chosen upper convex approximations, respectively; both $\tilde{F}_{i^k}$ and $\tilde{X}_{i^k}$ are built using the out-of-sync information $\mathbf{x}^{k-\mathbf{d}^k}$. If the set $X_{i^k}$ is convex, then $\tilde{X}_{i^k} = X_{i^k}$. More details on the choices of $\tilde{F}_{i^k}$ and $\tilde{X}_{i^k}$ can be found in [37].

Almost all modern asynchronous algorithms for convex and nonconvex problems are modeled in a probabilistic way. All current probabilistic models for asynchronous BCD methods are based on the (implicit or explicit) assumption that the random variables $i^k$ and $\mathbf{d}^k$ are *independent*; this greatly simplifies the conver-gence analysis. However, in reality there is a strong dependence of the delays $\mathbf{d}^k$ on the updated block $i^k$; see [37] for a detailed discussion on this issue and several practical examples. Another unrealistic assumption often made in the literature [60, 146, 147, 195] is that the block-indices $i^k$ are selected *uniformly* at random. While this assumption simplifies the convergence analysis, it limits the applicability of the model (see, e.g., Examples 4 and 5 in [37]). In a nutshell, this assumption may be satisfied only if all workers have the same computational power and have access to all variables. In [37] a more general, and sophisticated probabilistic model describing the statistics of $(i^k; \mathbf{d}^k)$ was introduced, and convergence of the asynchronous parallel SCA method (165)-(166) established; theoretical complexity results were also provided, showing nearly ideal linear speedup when the number of workers is not too large. The new model in [37] neither postulates the independence between $i^k$ and $\mathbf{d}^k$ nor requires artificial changes in the algorithm to enforce it (like those recently proposed in the probabilistic models [137, 153, 184] used in stochastic gradient methods); it handles instead the potential dependency among variables directly, fixing thus the theoretical issues that mar most of the aforementioned papers. It also lets one analyze for the first time in a sound way several practically used and effective computing settings and new models of asynchrony. For instance, it is widely accepted that in shared-memory systems, the best performance are obtained by first partitioning the variables among cores, and then letting each core update in an asynchronous fashion their own block-variables, according to some randomized cyclic rule; [37] is the first work proving convergence of such practically effective methods in an asynchronous setting.

Another important feature of the asynchronous algorithm (165)-(166) is its SCA nature, that is, the ability to handle nonconvex objective functions and nonconvex constraints by solving, at each iteration, a strongly convex optimization subproblem. Almost all asynchronous methods cited above can handle only convex optimization problems or, in the case of fixed point problems, nonexpansive mappings. The exceptions are [144, 265] and [60, 61] that studied unconstrained and constrained nonconvex optimization problems, respectively. However, [60, 61] proposed algorithms that require, at each iteration, the global solution of nonconvex subproblems. Except for few cases, the subproblems could be hard to solve and potentially as difficult as the original one. On the other hand, the SCA method [37] needs a feasible initial point and the ability to build approximations $\tilde{X}_i$ satisfying some technical conditions, as given in [37, Assumption D]. The two approaches thus complement each other and may cover different applications.

# Lecture III – Distributed Successive Convex Approximation Methods

This lecture complements the first two, extending the SCA algorithmic framework developed therein to distributed (nonconvex, multi-agent) optimization over networks with *arbitrary*, possibly time-varying, topology.

The SCA methods introduced in Lecture II unlock parallel updates from the workers; however, to perform its update, each worker must have the knowledge of some global information on the optimization problem, such as (part of) the objective function $V$, its gradient, and the current value of the optimization variable of the other agents. This clearly limits the applicability of these methods to network architectures wherein such information can be efficiently acquired (e.g., through suitably defined message-passing protocols and node coordination). Examples of such systems include the so-called multi-layer *hierarchical networks* (HNet); see Fig. III.1. A *HNet* consists of distributed nodes (DNs), cluster heads (CHs) and a master node, each having some local information on the optimization problem. Each CH can communicate with a (possibly dynamically formed) cluster of DNs as well as a higher layer CH, through either deterministic or randomly activated links. The HNet arises in many important applications including sensor networks, cloud-based software defined networks, and shared-memory systems. The HNet is also a generalization of the so-called "star network" (a two-layer HNet) that is commonly adopted in several parallel computing environments; see e.g. the Parameter Server [142] or the popular DiSCO [268] algorithm, just to name a few.

On the other hand, there are networks that lack of a hierarchical structure or "special" nodes; an example is the class of general *mesh networks* (MNet), which consists solely of DNs, and each of them is connected with a subset of neighbors, via possibly time-varying and directional communication links; see Fig. III.1. When the directional links are present, the MNet is referred to as a *digraph*. The MNet has been very popular to model applications such as *ad-hoc* (telecommunication) networks and social networks, where there are no obvious central controllers. Performing the SCA methods introduced in the previous lectures on such networks might incur in a computation/communication inefficient implementation.

The objective of this lecture is to devise distributed algorithms based on SCA techniques that are implementable efficiently on such general network architectures. More specifically, we consider a system of $I$ DNs (we will use interchangeably also the words "workers" or "agents") that can communicate through a network, modeled as a directed graph, possibly time-varying; see Fig. III.2. Agents want to cooperatively solve the following *networked* instance of Problem (94):

$$\underset{\mathbf{x} \in X}{\text{minimize}} \quad V(\mathbf{x}) \triangleq \underbrace{\sum_{i=1}^{I} f_i(\mathbf{x})}_{F(\mathbf{x})} + G(\mathbf{x}), \tag{167}$$

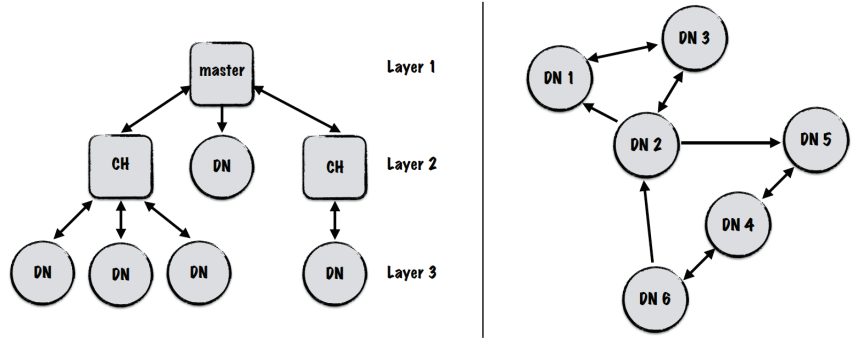

*Fig. III.1: Left: A three-layer **hierarchical network**, with one master node, 2 cluster heads, and 5 distributed nodes. Right: A six-node **mesh network**. The double arrowed (resp. single arrowed) links represent bi-directional (resp. directional) communication links.*

where the objective function $F$ is now the sum of the local cost functions $f_i : O \to \mathbb{R}$ of the agents, assumed to be smooth but possibly nonconvex whereas $G : O \to \mathbb{R}$ is a nonsmooth convex function; $O \supseteq X$ is an open set and $X \subseteq \mathbb{R}^m$ is a convex, closed set. In this networked setting, each agent $i$ knows only its own functions $f_i$ (and $G$ and $X$ as well). The problem and network settings are described in more details in Sec. III.1, along with some motivating applications.

The design of distributed algorithms for Problem (167) faces two challenges, namely: the nonconvexity of $F$ and the lack of full knowledge of $F$ from each agent. To cope with these issues, this lecture builds on the idea of SCA techniques coupled with suitably designed message passing protocols (compatible with the local agent knowledge of the network) aiming at disseminating information among the nodes as well as locally estimating $\nabla F$ from each agent. More specifically, for each agent $i$, a *local copy* $\mathbf{x}_{(i)}$ of the global variable $\mathbf{x}$ is introduced (cf. Fig. III.2). We say that a *consensus* is reached if $\mathbf{x}_{(i)} = \mathbf{x}_{(j)}$, for all $i \neq j$. To solve (167) over a network, *two major steps* are performed iteratively: local computation (to enhance the local solution quality), and local communication (to reach global consensus). In the first step, all the agents in parallel optimize their own variables $\mathbf{x}_{(i)}$ by solving a suitably chosen convex approximation of (167), built using the available local information. In the second step, agents communicate with their neighbors to acquire some new information instrumental to align users' local copies (and thus enforce consensus asymptotically) and update the surrogate function used in their local optimization subproblems. These two steps will be detailed in the rest of the sections of this lecture, as briefly outlined next. Sec. III.2 introduces distributed weighted-averaging algorithms to solve the (unconstrained) consensus problem over both static and time-varying (di-)graphs; a perturbed version of these consensus protocols is also introduced to unlock tracking of time-varying signals over networks. These message-passing protocols constitute the core of the distributed SCA-based algorithms that are discussed in this lecture: they will be used as an underlying mechanism for diffusing the information from one agent to every other agent in the network as well as track locally the gradient of the sum-utility $F$. In Sec. III.3, we build the proposed distributed algorithmic framework combining SCA techniques with the consensus/tracking protocols introduced in Sec. III.2, and study its convergence; a connection with existing (special case) schemes is also discussed. Some numerical results are presented in Sec. III.4. Finally, the main literature on related works is discussed in Sec. III.5 along with some extensions and open problems.

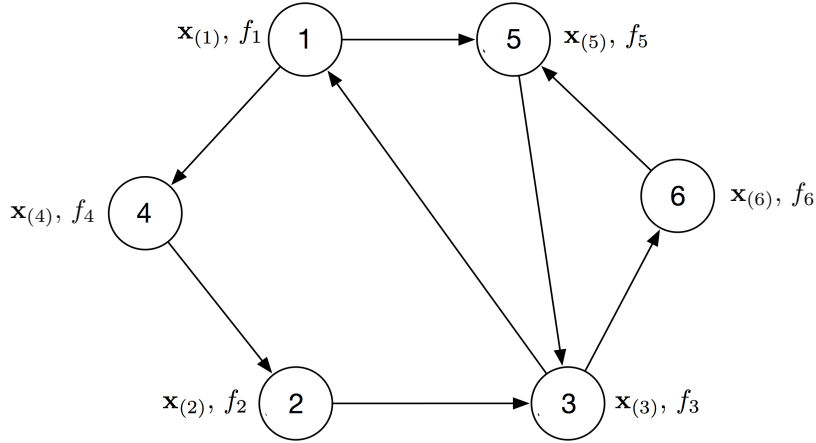

*Fig. III.2: An example of Problem* (167) *over a directed communication network. Each agent i knows only its own function $f_i$. To solve cooperatively* (167), *the agents create a local copy $\mathbf{x}_{(i)}$ of the common set of variables* $\mathbf{x}$. *These local copies are iteratively updated by the owners using only local (neighbor) information, so that asymptotically a consensus among them on a stationary solution of* (167) *is achieved.*

## III.1. Problem Formulation

We study Problem (167) under the following assumptions.

**Assumption III.1.** *Given Problem* (167), *assume that*

1. $\emptyset \neq X \subseteq \mathbb{R}^m$ *is closed and convex;*
2. $f_i : O \to \mathbb{R}$ *is* $C^1$ *on the open set* $O \supseteq X$, *and* $\nabla f_i$ *is* $L_i$-*Lipschitz on X;*
3. $G : O \to \mathbb{R}$ *is convex, possibly nonsmooth;*
4. $V$ *is bounded from below on X.*

Assumption III.1 can be viewed as the distributed counterpart of Assumption II.1 (Lecture II). Furthermore, we make the blanket assumption that each agent $i$ knows only its local function $f_i$, the common regularizer $G$, and the feasible set $X$; therefore, agents must communicate over a network to solve (167). We consider the following network setup.

**Network Model:** Agents are connected through a (communication) network, which is modeled as a graph; the set of agents are the nodes of the graph while the set of edges represents the communication links. We will consider both *static* and *time-varying* graphs, as well as *undirected* and *directed* graphs. We will use the following notation: $G^k = (V, E^k)$ denotes the directed graph that connects the agent at (the discrete) time $k$, where $V \triangleq \{1, \ldots, I\}$ is the set of nodes and $E^k$ is the set of edges (agents' communication links); we use $(i,j) \in E^k$ to indicate that the link is directed from node $i$ to node $j$. The *in*-neighborhood of agent $i$ at time $k$ (including node $i$ itself) is defined as $N_i^{\text{in},k} \triangleq \{j \mid (j,i) \in E^k\} \cup \{i\}$ whereas its *out*-neighborhood (including node $i$ itself) is defined as $N_i^{\text{out},k} \triangleq \{j \mid (i,j) \in E^k\} \cup \{i\}$. Of course, if the graph is undirected, the set of in-neighbors and out-neighbors are identical. These neighbors capture the local view of the network from agent $i$ at time $k$: At the time the communication is performed, agent $i$ can receive information from its current in-neighbors and send information to its current out neighbors. Note that we implicitly assumed that only inter-node (intermittent) communications between single-hop neighbors can be performed. The out-degree of agent $i$ at time $k$ is defined as the cardinality of $N_i^{\text{out},k}$, and is denoted by $d_i^k \triangleq |N_i^{\text{out},k}|$. We will treat static and/or undirected graphs as special cases of the above time-varying directed setting.

An important aspect of graphs is their *connectivity* properties. An *undirected* (static) graph is connected if there is a path connecting every pair of two distinct nodes. A *directed* (static) graph is strongly connected if there is a directed path from any node to any other node in the graph. For time-varying (di-)graphs we will invoke the following "long-term" connectivity property.

**Assumption III.2** (*B-strongly connectivity*). *The digraph sequence $\{G^k\}_{k \in \mathbb{N}_+}$ is B-strongly connected, i.e., there exists an integer $B > 0$ (possibly unknown to the agents) such that the digraph with edge set $\cup_{t=k}^{k+B-1} E^t$ is strongly connected, for all $k \geq 0$.*

Generally speaking, the above assumption permits strong connectivity to occur over a time window of length $B$: the graph obtained by taking the union of any $B$ consecutive graphs is strongly connected. Intuitively, this let information propagate throughout the network. Assumption III.2 is standard and well-accepted in the literature.

### III.1.1 Some motivating applications

Problems in the form (167), under Assumptions III.1 and III.2, have found a wide range of applications in several areas, including network information processing, telecommunications, multi-agent control, and machine learning. In particular, they are a key enabler of many nonconvex in-network "big data" analytic tasks, including nonlinear least squares, dictionary learning, principal/canonical component analysis, low-rank approximation, and matrix completion, just to name a few. Time-varying communications arise, for instance, in mobile wireless networks (e.g., ad-hoc networks) wherein nodes are mobile and/or communicate throughout fading channels. Moreover, since nodes generally transmit at different power and/or communication channels are not symmetric, directed links is a natural assumption. Some illustrative examples are briefly discussed next; see Sec. III.4 for more details and some numerical results.

**Example #1–(Sparse) empirical risk minimization:** In Example #5 in Sec. II.1.1 (Lecture II), we introduced the empirical risk minimization (ERM) problem, which consists in estimating a parameter $\mathbf{x}$ from a given data set $\{D_i\}_{i=1}^{I}$ by minimizing the risk function $F(\mathbf{x}) \triangleq \sum_{i=1}^{I} \ell(h(\mathbf{x}, D_i))$. Consider now the scenario where the data set is not centrally available but split among $I$ agents, connected through a network; agent $i$ only owns the portion $D_i$. All the agents want to collaboratively estimate $\mathbf{x}$, still minimizing $F(\mathbf{x})$. This distributed counterpart of the ERM problem is an instance of (167), with $f_i(\mathbf{x}) \triangleq \ell(h(\mathbf{x}, D_i))$. Many distributed statistical learning problems fall under this umbrella. Examples include: the least squares problem with $f_i(\mathbf{x}) \triangleq \|\mathbf{y}_i - \mathbf{A}_i \mathbf{x}\|^2$, where $D_i \triangleq (\mathbf{y}_i, \mathbf{A}_i)$; the sparse logistic regression with $f_i(\mathbf{x}) \triangleq \sum_{j=1}^{n_i} \log(1 + e^{-w_{ij} \mathbf{y}_{ij}^T \mathbf{x}})$, where $D_i \triangleq \{(w_{ij}, \mathbf{y}_{ij})\}_{j=1}^{n_i}$, and their sparse counterpart with suitable choices of the regularizer $G(\mathbf{x})$ (cf. Table I.1 in Lecture I, Sec. I.4.1). **Example #2–Sparse principal component analysis:** Consider an $m$-dimensional data set $\{\mathbf{d}_i\}_{i=1}^{n}$ with zero mean stored distributively among $I$ agents, each agent $i$ owns $\{\mathbf{d}_j\}_{j \in N_i}$, where $\{N_i\}_{i=1}^{I}$ forms a partition of $\{1, \ldots, n\}$. The problem of sparse principal component analysis is to find a sparse direction $\mathbf{x}$ along which the variance of the data points, measured by $\sum_{i=1}^{n} \|\mathbf{d}_i^T \mathbf{x}\|^2$, is maximized. Construct the matrix $\mathbf{D}_i \in \mathbb{R}^{|N_i| \times m}$ by stacking $\{\mathbf{d}_j\}_{j \in N_i}$ row-wise, the problem can be formulated as an instance of (167) with $f_i(\mathbf{x}) \triangleq -\|\mathbf{D}_i \mathbf{x}\|^2$, $X \triangleq \{\mathbf{x} \mid \|\mathbf{x}\|_2 \leq 1\}$, and $G(\mathbf{x})$ being some sparsity promoting regularizer (cf. Table I.1 in Lecture I, Sec. I.4.1).

**Example #3–Target localization:** Consider the problem of locating $n$ targets using measurements from $I$ sensors, embedded in a network. Each sensor $i$ knows its own position $\mathbf{s}_i$ and $d_{it}$, the latter representing the squared Euclidean distance between the target $t$ and the node. Given the position $\mathbf{x}_t$ of each target $t$, an error measurement of agent $i$ is $e_i(\mathbf{x}_t) \triangleq \sum_{t=1}^{n} p_{it}(d_{it} - \|\mathbf{x}_t - \mathbf{s}_i\|^2)^2$, where $p_{it} \in \{0, 1\}$ is a given binary variable taking value zero if the $i$th agent does not have any measurement related to target $t$. The problem of estimating the locations $\{\mathbf{x}_t\}_{t=1}^{n}$ can be thus formulated as an instance of (167), with $\mathbf{x} \triangleq \{\mathbf{x}_t\}_{t=1}^{n}$, $f_i(\mathbf{x}) \triangleq e_i(\mathbf{x}_t)$, and $X \triangleq \prod_{t=1}^{n} X_t$, where $X_t$ characterizes the region where target $t$ belongs to.

## III.2. Preliminaries: Average Consensus and Tracking

In this section, we introduce some of the building blocks of the distributed algorithmic framework that will be presented in Sec. III.3, namely: i) a consensus algorithm implementable on undirected (Sec. III.2.1) and directed (Sec. III.2.2) time-varying graphs; ii) a dynamic consensus protocol to track the average of time-varying signals over time-varying (directed) graphs (Sec. III.2.3); and iii) a perturbed consensus protocol unifying and generalizing the schemes in i) and ii) (Sec. III.2.4).

### III.2.1 Average consensus over undirected graphs

The consensus problem (also termed agreement problem) is one of the basic problems arising in decentralized coordination and control. Here we are interested in the so-called *average* consensus problem, as introduced next.

Consider a network of $I$ agents, each of which having some initial (vector) variable $\mathbf{u}_i \in \mathbb{R}^m$. The agents are interconnected over a (time-varying) network; the graph modeling the network at time $k$ is denoted by $G^k$ (cf. Sec. III.1). Each agent $i$ controls a local variable $\mathbf{x}_{(i)}$ that is updated at each iteration $k$ using the information of its immediate neighbors $N_i^{\text{in},k}$; we denote by $\mathbf{x}_{(i)}^k$ the value of $\mathbf{x}_{(i)}$ at iteration $k$. The average consensus problem consists in designing a distributed algorithm obeying the communication structure of each graph $G^k$, and enforcing

$$\lim_{k \to \infty} \|\mathbf{x}_{(i)}^k - \bar{\mathbf{u}}\| = 0, \quad \forall i = 1, \ldots, I, \quad \text{with} \quad \bar{\mathbf{u}} \triangleq \frac{1}{I} \sum_{i=1}^I \mathbf{u}_i.$$

One can construct a weighted-averaging protocol that solves the consensus problem as follows. Let each $\mathbf{x}_{(i)}^0 = \mathbf{u}_{(i)}^0$; given the iterate $\mathbf{x}_{(i)}^k$, at time $k+1$, each agent receives values $\mathbf{x}_{(j)}^k$ from its current (in-)neighbors, and updates its variable by setting

$$\mathbf{x}_{(i)}^{k+1} = \sum_{j \in N_i^{\text{in},k}} w_{ij}^k \mathbf{x}_{(j)}^k, \quad (168)$$

where $w_{ij}^k$ are some positive weights, to be properly chosen. For a more compact representation, we define the nonnegative weight-matrix[1] $\mathbf{W}^k \triangleq (w_{ij}^k)_{i,j=1}^I$, whose nonzero pattern is compliant with the topology of the graph $G^k$ (in the sense of Assumption III.3 below). Recall that the set of neighbors $N_i^{\text{in},k}$ contains also agent $i$ (cf. Sec. III.1).

**Assumption III.3.** *Given the graph sequence $\{G^k\}_{k \in \mathbb{N}_+}$, each matrix $\mathbf{W}^k \triangleq (w_{ij})_{i,j=1}^I$ satisfies:*

1. $w_{ij}^k = 0$, *if* $(j,i) \notin E^k$; *and* $w_{ij}^k \geq \kappa$, *if* $(j,i) \in E^k$;
2. $w_{ii}^k \geq \kappa$, *for all* $i = 1, \ldots, I$;

*for some given $\kappa > 0$.*

Using Assumption III.3.1, we can write the consensus protocol (168) in the following equivalent form

$$\mathbf{x}_{(i)}^{k+1} = \sum_{j=1}^I w_{ij}^k \mathbf{x}_{(j)}^k. \quad (169)$$

Note that (169) is compliant with the graph topology: the agents can only exchange information (according to the direction of the edge) if there exists a communication link between them. Also, Assumption III.3.2 states that each agent should include in the update (169) its own current information.

Convergence of $\{\mathbf{x}_{(i)}^k\}_{k \in \mathbb{N}_+}$ in (169) to the average $\bar{\mathbf{u}}$ calls for the following extra assumption.

**Assumption III.4.** *Each $\mathbf{W}^k$ is doubly-stochastic, i.e., $\mathbf{1}^T \mathbf{W}^k = \mathbf{1}^T$ and $\mathbf{W}^k \mathbf{1} = \mathbf{1}$.*

---

[1] Note that, for notational simplicity, here we use reverse links for the weight assignment, that is, each weight $w_{ij}$ is assigned to the directed edge $(j,i) \in E^k$.

*Table III.1: Examples of rules for doubly-stochastic weight matrices compliant to an undirected graph $G = (V,E)$ (or a digraph admitting a double-stochastic matrix). In the Laplacian weight rule, $\lambda$ is a positive constant and $\mathbf{L}$ is the Laplacian of the graph.*

| Rule Name | Weight Expression |
|---|---|
| Metropolis-Hastings | $w_{ij} = \begin{cases} \dfrac{1}{\max\{d_i, d_j\}}, & \text{if } (i,j) \in E, \\ 1 - \sum_{j \neq i} w_{ij}^k, & \text{if } i = j, \\ 0, & \text{if } (i,j) \notin E \text{ and } i \neq j; \end{cases}$ |
| Laplacian | $\mathbf{W} = \mathbf{I} - \lambda \mathbf{L}, \ \lambda > 0$ |
| Maximum-degree | $w_{ij} = \begin{cases} 1/I, & \text{if } (i,j) \in E, \\ 1 - (d_i - 1)/I, & \text{if } i = j, \\ 0, & \text{if } (i,j) \notin E \text{ and } i \neq j. \end{cases}$ |

Assumption III.4 requires **1** being both the left and right eigenvector of $\mathbf{W}^k$ associated to the eigenvalue 1; intuitively, the column stochasticity plays the role of preserving the total sum of the $\mathbf{x}_{(i)}$'s (and thus the average) in the network while the row stochasticity locks consensus.

When the graphs $G^k$ are undirected (or are directed and admits a compliant doubly-stochastic matrix), several rules have been proposed in the literature to build a weight matrix satisfying Assumption III.3 and Assumption III.4. Examples include the Laplacian weight rule [203]; the maximum degree weight, the Metropolis-Hastings, and the least-mean square consensus weight rules [258]. Table III.1 summarizes the aforementioned rules [202], where in the Laplacian weight rule, $\lambda$ is a positive constant and $\mathbf{L}$ is the graph Laplacian, whose the $ij$-th entry $L_{ij}$ is defined as

$$L_{ij} \triangleq \begin{cases} d_i - 1, & \text{if } i = j; \\ -1, & \text{if } (i,j) \in E \text{ and } i \neq j; \\ 0, & \text{otherwise;} \end{cases} \quad (170)$$

where $d_i$ is the degree of node $i$.

Convergence of the average-consensus protocol (169) is stated in the next theorem, whose proof is omitted because it is a special case of that of Theorem III.11 (cf. Sec. III.2.4).

**Theorem III.5.** *Let $\{G^k\}_{k \in \mathbb{N}_+}$ be a sequence of graphs satisfying Assumption III.2. Consider the average-consensus protocol (169), where each $\{\mathbf{W}^k\}_{k \in \mathbb{N}_+}$ is chosen according to Assumptions III.3 and III.4. Then, the sequence $\{\mathbf{x}^k \triangleq (\mathbf{x}_{(i)}^k)_{i=1}^I\}_{k \in \mathbb{N}_+}$ generated by (169) satisfies: for all $k \in \mathbb{N}_+$,*

*(a) Invariance of the average:*

$$\sum_{i=1}^I \mathbf{x}_{(i)}^{k+1} = \sum_{i=1}^I \mathbf{x}_{(i)}^k = \cdots = \sum_{i=1}^I \mathbf{x}_{(i)}^0; \quad (171)$$

*(b) Geometric decay of the consensus error:*

$$\left\| \mathbf{x}_{(i)}^k - \frac{1}{I} \sum_{j=1}^{I} \mathbf{x}_{(j)}^k \right\| \leq c_u \cdot (\rho_u)^k \cdot \|\mathbf{x}^0\|, \quad \forall i = 1, \ldots, I, \tag{172}$$

*where $\mathbf{x}^0 \triangleq (\mathbf{x}_{(i)}^0)_{i=1}^{I}$, and $c_u > 0$ and $\rho_u \in (0,1)$ are constants defined as*

$$c_u \triangleq \frac{2I}{\rho_u} \cdot \frac{2(1 + \kappa^{-(I-1)B})}{1 - \kappa^{(I-1)B}} \quad \text{and} \quad \rho_u \triangleq \left(1 - \kappa^{(I-1)B}\right)^{\frac{1}{(I-1)B}}, \tag{173}$$

*with $B$ and $\kappa$ defined in Assumption III.2 and Assumption III.3, respectively.*

In words, Theorem III.5 states that each $\mathbf{x}_{(i)}^k$ converges to the initial average $(1/I) \cdot \sum_{i=1}^{I} \mathbf{x}_{(i)}^0$ at a geometric rate. Since $\mathbf{x}_{(i)}^0$ is initialized as $\mathbf{x}_{(i)}^0 \triangleq \mathbf{u}_i$, each $\{\mathbf{x}_{(i)}^k\}_{k \in \mathbb{N}_+}$ converges to $\bar{\mathbf{u}}$ geometrically.

**Remark III.6.** *While Theorem III.5 has been stated under Assumption III.4 (because we are mainly interested in the convergence to the average $\bar{\mathbf{u}}$), it is important to remark that the row-stochasticity of each $\mathbf{W}^k$ (rather than doubly-stochasticity) is enough for the sequence $\{\mathbf{x}^k\}_{k \in \mathbb{N}_+}$ generated by the protocol (169) to geometrically reach a consensus, that is, $\lim_{k \to \infty} \|\mathbf{x}_{(i)}^k - \mathbf{x}_{(j)}^k\| = 0$, for all $i, j = 1, \ldots I$ and $i \neq j$. However, the limit point of $\{\mathbf{x}^k\}_{k \in \mathbb{N}_+}$ is no longer the average of the $\mathbf{x}_{(i)}^0$'s.*

### III.2.2 Average consensus over directed graphs

A key assumption for the distributed protocol (169) to reach the *average* consensus is that each matrix $\mathbf{W}^k$, compliant with the graph $G^k$, is doubly-stochastic.

While such constructions exist for networks with bi-directional (possibly time-varying) communication links, they become computationally prohibitive or infeasible for networks with directed links, for several reasons. First of all, not all digraphs admit a compliant (in the sense of Assumption III.3) doubly-stochastic weight matrix; some form of balancedness in the graph is needed [92], which limits the class of networks over which the consensus protocol (169) can be applied. Furthermore, conditions for a digraph to admit such a doubly-stochastic matrix are not easy to be checked in practice; and, even when possible, constructing a doubly-stochastic weight matrix compliant to the digraph calls for computationally intense, generally centralized, algorithms.

To solve the average consensus problem over digraphs that do not admit a doubly-stochastic matrix, a further assumption is needed [101] along with a modification of the protocol (169). Specifically, a standard assumption in the literature is that every agent *i* knows its out-degree $d_i^k$ at each time *k*. This means that, while broadcasting its own message, every agent knows how many nodes will receive it. The problem of computing the out-degree using only local information has been considered in a number of works (see, e.g., [119, 243] and the references therein). Various algorithms have been proposed, mainly based on flooding, which, however, requires significant communication overhead and storage. A less demanding consensus-based approach can be found in [44].

Under the above assumption, the average consensus can be achieved on digraphs using the so-called push-sum protocol [121]. Each agent $i$ controls two local variables, $\mathbf{z}_{(i)} \in \mathbb{R}^m$ and $\phi_i \in \mathbb{R}$, which are updated at each iteration $k$ still using only the information of its immediate neighbors. The push-sum protocols reads: for all $i = 1, \ldots, I$,

$$\begin{aligned} \mathbf{z}_{(i)}^{k+1} &= \sum_{j=1}^{I} a_{ij}^k \mathbf{z}_{(j)}^k, \\ \phi_{(i)}^{k+1} &= \sum_{j=1}^{I} a_{ij}^k \phi_{(j)}^k, \end{aligned} \quad (174)$$

where $\mathbf{z}_{(i)}$ and $\phi_{(i)}$ are initialized as $\mathbf{z}_{(i)}^0 = \mathbf{u}_i$ and $\phi_{(i)}^0 = 1$, respectively; and the coefficient $a_{ij}^k$ are defined as

$$a_{ij}^k \triangleq \begin{cases} \dfrac{1}{d_j^k}, & \text{if } j \in N_i^{\text{in},k}, \\ 0, & \text{otherwise}. \end{cases} \quad (175)$$

Note that the scheme (174) is a broadcast (i.e., one-way) communication protocol: each agent $i$ broadcasts ("pushes out") the values $\mathbf{z}_{(j)}^k/d_j^k$ and $\phi_{(j)}^k/d_j^k$, which are received by its out-neighbors; at the receiver side, every node aggregates the received information according to (174) (i.e., summing the pushed values, which explains the name "push-sum").

Introducing the weight matrices $\mathbf{A}^k \triangleq (a_{ij}^k)_{i,j=1}^I$, it is easy to check that, for general digraphs, $\mathbf{A}^k$ may no longer be row-stochastic (i.e., $\mathbf{A}^k \mathbf{1} = \mathbf{1}$). This means that the $z$- and $\phi$-updates in (174) do not reach a consensus. However, because $\mathbf{1}^T \mathbf{A}^k = \mathbf{1}^T$, the sums of the $z$- and $\phi$-variables are preserved: at every iteration $k \in \mathbb{N}_+$,

$$\begin{aligned} \sum_{i=1}^I \mathbf{z}_{(i)}^{k+1} &= \sum_{i=1}^I \mathbf{z}_{(i)}^k = \cdots = \sum_{i=1}^I \mathbf{z}_{(i)}^0 = \sum_{i=1}^I \mathbf{u}_i, \\ \sum_{i=1}^I \phi_{(i)}^{k+1} &= \sum_{i=1}^I \phi_{(i)}^k = \cdots = \sum_{i=1}^I \phi_{(i)}^0 = I. \end{aligned} \quad (176)$$

This implies that, if the iterates $(\mathbf{z}_{(i)}^k/\phi_{(i)}^k)_{i=1}^I$ converge to a consensus (note that each $\mathbf{z}_{(i)}^k/\phi_{(i)}^k$ is well-defined because the weights $\phi_i^k$ are all positive), then the consensus value must be the average $\bar{\mathbf{u}}$, as shown next. Let $\mathbf{c}^\infty$ be the consensus value, that is, $\lim_{k \to \infty} \mathbf{z}_{(i)}^k/\phi_{(i)}^k = \mathbf{c}^\infty$, for all $i = 1, \ldots, I$. Then, it must be

$$\left\| \sum_{i=1}^I \mathbf{u}_i - I \cdot \mathbf{c}^\infty \right\| \stackrel{(176)}{=} \left\| \sum_{i=1}^I \left( \mathbf{z}_{(i)}^k - \phi_{(i)}^k \cdot \mathbf{c}^\infty \right) \right\| \leq I \cdot \sum_{i=1}^I \left\| \mathbf{z}_{(i)}^k/\phi_{(i)}^k - \mathbf{c}^\infty \right\| \xrightarrow[k \to \infty]{} 0, \quad (177)$$

which shows that $\mathbf{c}^\infty = (1/I) \cdot \sum_{i=1}^I \mathbf{u}_i$. Convergence of $(\mathbf{z}_{(i)}^k/\phi_{(i)}^k)_{i=1}^I$ to the consensus is proved in [11, 121], which we refer to the interested reader.

Here, we study instead an equivalent reformulation of the push-sum algorithm, as given in [211, 231], which is more suitable for the integration with optimization (cf. Sec. III.3). Eliminating the $z$-variables in (174), and considering arbitrary *column-stochastic* weight matrices $\mathbf{A}^k \triangleq (a_{ij}^k)_{i,j=1}^I$, compliant to the graph $G^k$ [not necessarily given by (175)], we have [211, 231]:

$$\phi_{(i)}^{k+1} = \sum_{j=1}^{I} a_{ij}^k \phi_{(j)}^k,$$

$$\mathbf{x}_{(i)}^{k+1} = \frac{1}{\phi_{(i)}^{k+1}} \sum_{j=1}^{I} a_{ij}^k \phi_{(j)}^k \mathbf{x}_{(j)}^k, \qquad (178)$$

where $\mathbf{x}_{(i)}^0$ is set to $\mathbf{x}_{(i)}^0 = \mathbf{u}_i / \phi_{(i)}^0$, for all $i = 1, \ldots, I$; and $\phi_{(i)}^0$ are arbitrary positive scalars such that $\sum_{i=1}^I \phi_{(i)}^0 = I$. For simplicity, hereafter, we tacitly set $\phi_{(i)}^0 = 1$ (implying $\mathbf{x}_{(i)}^0 = \mathbf{u}_i$), for all $i = 1, \ldots, I$. Similarly to (174), in the protocol (178), every agent $i$ controls and updates the variables $\mathbf{x}_{(i)}$ and $\phi_{(i)}$, based on the information $\phi_{(j)}^k$ and $\phi_{(j)}^k \mathbf{x}_{(j)}^k$ received from its current in-neighbors. We will refer to (178) as *condensed push-sum* algorithm.

Combining the weights $a_{ij}^k$ and the $\phi_{(i)}^k$ in the update of $\mathbf{x}_{(i)}^k$ in a single coefficient

$$w_{ij}^k \triangleq \frac{a_{ij}^k \phi_{(j)}^k}{\sum_{j=1}^I a_{ij}^k \phi_{(j)}^k}, \qquad (179)$$

it is not difficult to check that the matrices $\mathbf{W}^k \triangleq (w_{ij}^k)_{i,j=1}^I$ are row-stochastic, that is, $\mathbf{W}^k \mathbf{1} = \mathbf{1}$, and compliant to $G^k$. This means that each $\mathbf{x}_{(i)}^{k+1}$ in (178) is updated performing a convex combination of the variables $(\mathbf{x}_{(j)}^k)_{j \in N_i^{\text{in},k}}$. This is a key property that will be leveraged in Sec. III.3 to build a distributed optimization algorithm for constrained optimization problems wherein the feasibility of the iterates is preserved at each iteration. The above equivalent formulation also sheds light on the role of the $\phi$-variables: they rebuild dynamically the missing row stochasticity of the weights $a_{ij}^k$, thus enforcing the consensus on the $x$-variables.

Since the following quantities are invariants of the dynamics (178) (recall that $\mathbf{A}^k$ are column stochastic), that is, for all $k \in \mathbb{N}_+$,

$$\sum_{i=1}^{I} \phi_{(i)}^{k+1} \mathbf{x}_{(i)}^{k+1} = \sum_{i=1}^{I} \phi_{(i)}^k \mathbf{x}_{(i)}^k = \cdots = \sum_{i=1}^{I} \phi_{(i)}^0 \mathbf{x}_{(i)}^0 = \sum_{i=1}^{I} \mathbf{u}_i,$$

$$\sum_{i=1}^{I} \phi_{(i)}^{k+1} = \sum_{i=1}^{I} \phi_{(i)}^k = \cdots = \sum_{i=1}^{I} \phi_{(i)}^0 = I,$$

by a similar argument used in (177), one can show that, if the $\mathbf{x}_{(i)}^k$ are consensual–$\lim_{k \to \infty} \mathbf{x}_{(i)}^k = \mathbf{x}^\infty$, for all $i = 1, \ldots, I$–it must be $\mathbf{x}^\infty = (1/I) \cdot \sum_{i=1}^I \mathbf{u}_i$. Convergence to the consensus at geometric rate is stated in Theorem III.8 below, under the following assumption on the weight matrices $\mathbf{A}^k$ (the proof of the theorem is omitted because it is a special case of that of Theorem III.11, cf. Sec. III.2.4).

**Assumption III.7.** *Each $\mathbf{A}^k$ is compliant with $G^k$ (i.e., it satisfies Assumption III.3) and it is column stochastic, i.e., $\mathbf{1}^T \mathbf{A}^k = \mathbf{1}^T$.*

**Theorem III.8.** *Let $\{G^k\}_{k \in \mathbb{N}_+}$ be a sequence of graphs satisfying Assumption III.2. Consider the condensed push-sum protocol* (178), *where each $\{\mathbf{A}^k\}_{k \in \mathbb{N}}$ is chosen according to Assumption III.7. Then, the sequence $\{\mathbf{x}^k \triangleq (\mathbf{x}_{(i)}^k)_{i=1}^I\}_{k \in \mathbb{N}_+}$ generated by* (178) *satisfies: for all $k \in \mathbb{N}_+$,*

*(a) Invariance of the weighted-sum:*

$$\sum_{i=1}^I \phi_{(i)}^{k+1} \mathbf{x}_{(i)}^{k+1} = \sum_{i=1}^I \phi_{(i)}^k \mathbf{x}_{(i)}^k = \cdots = \sum_{i=1}^I \phi_{(i)}^0 \mathbf{x}_{(i)}^0; \quad (180)$$

*(b) Geometric decay of the consensus error:*

$$\left\| \mathbf{x}_{(i)}^k - \frac{1}{I} \sum_{j=1}^I \phi_{(j)}^k \mathbf{x}_{(j)}^k \right\| \le c_d \cdot (\rho_d)^k \cdot \|\mathbf{x}^0\|, \quad \forall i = 1, \ldots, I, \quad (181)$$

*where $\mathbf{x}^0 \triangleq (\mathbf{x}_{(i)}^0)_{i=1}^I$, and $c_d > 0$ and $\rho_d \in (0,1)$ are defined as*

$$c_d \triangleq \frac{2I}{\rho} \cdot \frac{2(1+\tilde{\kappa}_d^{-(I-1)B})}{1-\kappa_d^{(I-1)B}} \quad \text{and} \quad \rho_d \triangleq \left(1 - \tilde{\kappa}_d^{(I-1)B}\right)^{\frac{1}{(I-1)B}},$$

*and $\tilde{\kappa}_d \triangleq \kappa^{2(I-1)B+1}/I$, with $B$ and $\kappa$ defined in Assumption III.2 and Assumption III.3, respectively.*

Note that to reach consensus, the condensed push-sum protocol requires that the weight matrices $\mathbf{A}^k$ are column stochastic but not doubly-stochastic. Of course, one can always use the classical push-sum weights as in (175), which is an example of rule satisfying Assumption III.7. Moreover, if the graph $G^k$ is undirected or is directed but admits a compliant doubly-stochastic matrix, one can also choose in (178) a doubly-stochastic $\mathbf{A}^k$. In such a case, the $\phi$-update in (178) indicates that $\phi_{(i)}^k = 1$, for all $i = 1, \ldots, I$ and $k \in \mathbb{N}_+$. Therefore, (178) [and also (174)] reduces to the plain consensus scheme (169).

### III.2.3 Distributed tracking of time-varying signals' average

In this section, we extend the condensed push-sum protocol to the case where the signals $\mathbf{u}_i$ are no longer constant but time-varying; the value of the signal at time $k$ owned by agent $i$ is denoted by $\mathbf{u}_i^k$. The goal becomes designing a distributed algorithm obeying the communication structure of each graph $G^k$ that tracks the average of $(\mathbf{u}_i^k)_{i=1}^I$, i.e.,

$$\lim_{k \to \infty} \|\mathbf{x}_{(i)}^k - \bar{\mathbf{u}}^k\| = 0, \quad \forall i = 1, \ldots, I, \quad \text{with} \quad \bar{\mathbf{u}}^k \triangleq \frac{1}{I} \sum_{i=1}^{I} \mathbf{u}_i^k.$$

We first introduce the algorithm for undirected graphs (or directed ones that admit a doubly-stochastic matrix); we then extend the scheme to the more general setting of arbitrary directed graphs.

**Distributed tracking over undirected graphs**

Consider a (possibly) time-varying network, modeled by a sequence of undirected graphs (or directed graphs admitting a doubly-stochastic matrix) $\{G^k\}_{k \in \mathbb{N}_+}$. As for the average consensus problem, we let each agent $i$ maintain and update a variable $\mathbf{x}_{(i)}^k$ that represents a local estimate of $\bar{\mathbf{u}}^k$; we set $\mathbf{x}_{(i)}^0 \triangleq \mathbf{u}_i^0$. Since $\mathbf{u}_i^k$ is time-varying, a direct application of the protocol (169), developed for the plain average consensus problem, cannot work because it would drive all $\mathbf{x}_{(i)}^k$ to converge to $\bar{\mathbf{u}}^0$. Therefore, we need to modify the vanilla scheme (169) to account for the variability of $\mathbf{u}_i^k$'s. We construct next the distributed tracking algorithm inductively, building on (169).

Recall that, in the average consensus scheme (169), a key property to set the consensus value to the average of the initial values $\mathbf{x}_{(i)}^0$ is the invariance of the average $\sum_{i=1}^{I} \mathbf{x}_{(i)}^k$ throughout the dynamics (169):

$$\sum_{i=1}^{I} \mathbf{x}_{(i)}^{k+1} = \sum_{i=1}^{I} \sum_{j=1}^{I} w_{ij}^k \mathbf{x}_{(j)}^k = \sum_{i=1}^{I} \mathbf{x}_{(i)}^k = \cdots = \sum_{i=1}^{I} \mathbf{x}_{(i)}^0, \tag{182}$$

which is met if the weight matrices $\mathbf{W}^k$ are column stochastic. The row-stochasticity of $\mathbf{W}^k$ enforces asymptotically a consensus. When it comes to solve the tracking problem, it seems then natural to require such an invariance of the (time-varying) average, that is, $\sum_{i=1}^{I} \mathbf{x}_{(i)}^k = \sum_{i=1}^{I} \mathbf{u}_i^k$, for all $k \in \mathbb{N}_+$, while enforcing a consensus on $\mathbf{x}_{(i)}^k$'s. Since $\mathbf{x}_{(i)}^0 \triangleq \mathbf{u}_i^0$ we have $\sum_{i=1}^{I} \mathbf{x}_{(i)}^0 \triangleq \sum_{i=1}^{I} \mathbf{u}_i^0$. Suppose now that, at iteration $k$, we have $\sum_{i=1}^{I} \mathbf{x}_{(i)}^k = \sum_{i=1}^{I} \mathbf{u}_i^k$. In order to satisfy $\sum_{i=1}^{I} \mathbf{x}_{(i)}^{k+1} \triangleq \sum_{i=1}^{I} \mathbf{u}_i^{k+1}$, it must be

$$\begin{aligned}
\sum_{i=1}^{I} \mathbf{x}_{(i)}^{k+1} &= \sum_{i=1}^{I} \mathbf{u}_i^{k+1} \\
&= \sum_{i=1}^{I} \mathbf{u}_i^k + \sum_{i=1}^{I} \left( \mathbf{u}_i^{k+1} - \mathbf{u}_i^k \right) \\
&\stackrel{(a)}{=} \sum_{i=1}^{I} \mathbf{x}_{(i)}^k + \sum_{i=1}^{I} \left( \mathbf{u}_i^{k+1} - \mathbf{u}_i^k \right) \\
&\stackrel{(b)}{=} \sum_{i=1}^{I} \left( \sum_{j=1}^{I} w_{ij}^k \mathbf{x}_{(j)}^k + \left( \mathbf{u}_i^{k+1} - \mathbf{u}_i^k \right) \right),
\end{aligned} \tag{183}$$

where in (a) we used $\sum_{i=1}^{I} \mathbf{x}_{(i)}^{k} = \sum_{i=1}^{I} \mathbf{u}_{i}^{k}$; and (b) follows from the column stochasticity of $\mathbf{W}^{k}$ [cf. (182)]. This naturally suggests the following modification of the protocol (169): for all $i = 1, \ldots, I$,

$$\mathbf{x}_{(i)}^{k+1} = \sum_{j=1}^{I} w_{ij}^{k} \mathbf{x}_{(j)}^{k} + \left( \mathbf{u}_{i}^{k+1} - \mathbf{u}_{i}^{k} \right), \tag{184}$$

with $\mathbf{x}_{(i)}^{0} = \mathbf{u}_{i}^{0}$, $i = 1, \ldots, I$. Generally speaking, (184) has the following interpretation: by averaging neighbors' information, the first term $\sum_{j=1}^{I} w_{ij}^{k} \mathbf{x}_{(j)}^{k}$ aims at enforcing a consensus among the $\mathbf{x}_{(i)}^{k}$'s while the second term is a perturbation that bias the current sum $\sum_{i=1}^{I} \mathbf{x}_{(i)}^{k}$ towards $\sum_{i=1}^{I} \mathbf{u}_{i}^{k}$. By (183) and a similar argument as in (177), it is not difficult to check that if the $\mathbf{x}_{(i)}^{k}$ are consensual it must be $\lim_{k \to \infty} \|\mathbf{x}_{(i)}^{k} - \bar{\mathbf{u}}^{k}\| = 0$, for all $i = 1, \ldots, I$. Theorem III.10 proves this result (as a special case of the tracking scheme over arbitrary directed graphs).

**Distributed tracking over arbitrary directed graphs**

Consider now the case of digraphs $\{G^{k}\}_{k \in \mathbb{N}_{+}}$ with arbitrary topology. With the results established in Sec. III.2.2, the tracking mechanism (184) can be naturally generalized to this more general setting by "perturbing" the condensed push-sum scheme (178) as follows: for $i = 1, \ldots, I$,

$$\begin{aligned}
\phi_{(i)}^{k+1} &= \sum_{j=1}^{I} a_{ij}^{k} \phi_{(j)}^{k}, \\
\mathbf{x}_{(i)}^{k+1} &= \frac{1}{\phi_{(i)}^{k+1}} \sum_{j=1}^{I} a_{ij}^{k} \phi_{(j)}^{k} \mathbf{x}_{(j)}^{k} + \frac{1}{\phi_{(i)}^{k+1}} \cdot \left( \mathbf{u}_{i}^{k+1} - \mathbf{u}_{i}^{k} \right),
\end{aligned} \tag{185}$$

where $\mathbf{x}_{(i)}^{0} = \mathbf{u}_{i}^{0}/\phi_{(i)}^{0}$, $i = 1, \ldots, I$; and $(\phi_{(i)}^{0})_{i=1}^{I}$ are arbitrary positive scalars such that $\sum_{i=1}^{I} \phi_{(i)}^{k} = I$. For simplicity, hereafter, we tacitly set $\phi_{(i)}^{0} = 1$, for all $i = 1, \ldots, I$. Note that, differently from (184), in (185), we have scaled the perturbation by $(\phi_{(i)}^{k+1})^{-1}$ so that the weighted average is preserved, that is, $\sum_{i=1}^{I} \phi_{(i)}^{k} \mathbf{x}_{(i)}^{k} = \sum_{i=1}^{I} \mathbf{u}_{i}^{k}$, for all $k \in \mathbb{N}_{+}$.

Convergence of the tracking scheme (185) is stated in Theorem III.10 below, whose proof is omitted because is a special case of the more general result stated in Theorem III.11 (cf. Sec. III.2.4).

**Assumption III.9.** *Let $\{\mathbf{u}_{i}^{k}\}_{k \in \mathbb{N}_{+}}$ be such that*

$$\lim_{k \to \infty} \|\mathbf{u}_{i}^{k+1} - \mathbf{u}_{i}^{k}\| = 0, \quad \forall i = 1, \ldots, I.$$

**Theorem III.10.** *Let $\{G^{k}\}_{k \in \mathbb{N}_{+}}$ be a sequence of graphs satisfying Assumption III.2. Consider the distributed tracking protocol* (185)*, where each $\mathbf{A}^{k}$ is chosen according* to Assumption III.7. Then, the sequence $\{\mathbf{x}^k \triangleq (\mathbf{x}^k_{(i)})^I_{i=1}\}_{k\in\mathbb{N}_+}$ generated by (185) satisfies:

*(a) Invariance of the weighted-sum:*

$$\sum_{i=1}^I \phi^k_{(i)} \mathbf{x}^k_{(i)} = \sum_{i=1}^I \mathbf{u}^k_i, \quad \forall k \in \mathbb{N}_+; \tag{186}$$

*(b) Asymptotic consensus: if, in addition, the sequence of signals* $\{(\mathbf{u}^k_i)^I_{i=1}\}_{k\in\mathbb{N}_+}$ *satisfies Assumption III.9, then*

$$\lim_{k\to\infty} \left\| \mathbf{x}^k_{(i)} - \bar{\mathbf{u}}^k \right\| = 0, \quad \forall i = 1,\ldots,I. \tag{187}$$

### III.2.4 Perturbed condensed push-sum

In this section, we provide a unified proof of Theorems III.5, III.8, and III.10 by interpreting the consensus and tracking schemes introduced in the previous sections [namely: (169), (178) and (185)] as special instances of the following *perturbed condensed push-sum* protocol: for all $i = 1,\ldots,I$, and $k \in \mathbb{N}_+$,

$$\begin{aligned} \phi^{k+1}_{(i)} &= \sum_{j=1}^I a^k_{ij} \phi^k_{(j)}, \\ \mathbf{x}^{k+1}_{(i)} &= \frac{1}{\phi^{k+1}_{(i)}} \sum_{j=1}^I a^k_{ij} \phi^k_{(j)} \mathbf{x}^k_{(j)} + \boldsymbol{\varepsilon}^k_i, \end{aligned} \tag{188}$$

with given $\mathbf{x}^0 \triangleq (\mathbf{x}^0_{(i)})^I_{i=1}$, and $\phi^0_{(i)} = 1$, for all $i = 1,\ldots,I$; where $\boldsymbol{\varepsilon}^k_i$ models a perturbation locally injected into the system by agent $i$ at iteration $k$. Clearly, the schemes (178) and (185) are special case of (188), obtained setting $\boldsymbol{\varepsilon}^k_i \triangleq \mathbf{0}$ and $\boldsymbol{\varepsilon}^k_i \triangleq (\phi^{k+1}_{(i)})^{-1}(\mathbf{u}^{k+1}_i - \mathbf{u}^k_i)$, respectively. Furthermore, if the matrices $\mathbf{A}^k$ are doubly-stochastic, we obtain as special cases the schemes (169) and (184), respectively.

Convergence of the scheme (188) is given in the theorem below.

**Theorem III.11.** *Let* $\{G^k\}_{k\in\mathbb{N}_+}$ *be a sequence of graphs satisfying Assumption III.2. Consider the perturbed condensed push-sum protocol* (188)*, where each* $\mathbf{A}^k$ *is chosen according to Assumption III.7. Then, the sequences* $\{\mathbf{x}^k \triangleq (\mathbf{x}^k_{(i)})^I_{i=1}\}_{k\in\mathbb{N}_+}$ *and* $\{\boldsymbol{\phi}^k \triangleq (\phi^k_{(i)})^I_{i=1}\}_{k\in\mathbb{N}_+}$ *generated by* (188) *satisfy:*

*(a) Bounded* $\{\boldsymbol{\phi}^k\}_{k\in\mathbb{N}_+}$:

$$\begin{aligned} \phi_{lb} &\triangleq \inf_{k\in\mathbb{N}_+} \min_{1\leq i \leq I} \phi^k_{(i)} \geq \kappa^{2(I-1)B}, \\ \phi_{ub} &\triangleq \sup_{k\in\mathbb{N}_+} \max_{1\leq i \leq I} \phi^k_{(i)} \leq I - \kappa^{2(I-1)B}; \end{aligned} \tag{189}$$

*(b) Bounded consensus error:* for all $k \in \mathbb{N}_+$,

$$\left\| \mathbf{x}_{(i)}^k - \frac{1}{I} \sum_{j=1}^{I} \phi_{(j)}^k \mathbf{x}_{(j)}^k \right\| \leq c \cdot \left( (\rho)^k \|\mathbf{x}^0\| + \sum_{t=0}^{k-1} (\rho)^{k-1-t} \|\boldsymbol{\varepsilon}^t\| \right), \quad \forall i = i, \ldots, I, \tag{190}$$

where $\boldsymbol{\varepsilon}^k \triangleq (\boldsymbol{\varepsilon}_i^k)_{i=1}^I$, and $c > 0$ and $\rho \in (0,1)$ are constants defined as

$$c \triangleq \frac{2I}{\rho} \cdot \frac{2(1 + \tilde{\kappa}^{-(I-1)B})}{1 - \tilde{\kappa}^{(I-1)B}} \quad \text{and} \quad \rho \triangleq \left(1 - \tilde{\kappa}^{(I-1)B}\right)^{\frac{1}{(I-1)B}},$$

with $\tilde{\kappa} \triangleq \kappa \cdot (\phi_{lb}/\phi_{ub})$, and $B$ and $\kappa$ defined in Assumption III.2 and Assumption III.3, respectively.

*Proof.* See Sec. III.2.5. □

**Discussion**

Let us apply now Theorem III.11 to study the impact on the consensus value of the following three different perturbation errors:

1) **Error free**: $\boldsymbol{\varepsilon}_i^k = \mathbf{0}$, for all $i = 1, \ldots, I$, and $k \in \mathbb{N}_+$;

2) **Vanishing error**: $\lim_{k \to \infty} \|\boldsymbol{\varepsilon}_i^k\| = 0$, for all $i = 1, \ldots, I$, and $k \in \mathbb{N}_+$;

3) **Bounded error**: There exists a constant $0 \leq M < +\infty$ such that $\|\boldsymbol{\varepsilon}_i^k\| \leq M$, for all $i = 1, \ldots, I$, and $k \in \mathbb{N}_+$.

**Case 1: Error free**. Since $\boldsymbol{\varepsilon}_i^k = \mathbf{0}$, $i = 1, \ldots, I$, the perturbed consensus protocol (188) reduces to the condensed push-sum consensus scheme (178). According to Theorem III.11, each $\mathbf{x}_{(i)}^k$ converges to the weighted average $(1/I) \cdot \sum_{i=1}^I \phi_{(i)}^k \mathbf{x}_{(i)}^k$ geometrically. Since $(1/I) \cdot \sum_{i=1}^I \phi_{(i)}^k \mathbf{x}_{(i)}^k = (1/I) \cdot \sum_{i=1}^I \phi_{(i)}^0 \mathbf{x}_{(i)}^0 = \bar{\mathbf{u}}$, for all $k \in \mathbb{N}_+$, each $\|\mathbf{x}_{(i)}^k - \bar{\mathbf{u}}\|$ vanishes geometrically, which proves Theorem III.8.

If, in addition, the weighted matrices $\mathbf{A}^k$ are also row-stochastic (and thus doubly-stochastic), we have $\phi_{(i)}^k = 1$, for all $k \in \mathbb{N}_+$. Therefore, (188) reduces to the vanilla average consensus protocol (169); and in Theorem III.11 we have $\tilde{\kappa} = \kappa$, which proves Theorem III.5 as a special case.

**Case 2: Vanishing error**. To study the consensus value in this setting, we need the following lemma on the convergence of product of sequences.

**Lemma III.12 ([148, Lemma 7]).** *Let $0 < \lambda < 1$, and let $\{\beta^k\}_{k \in \mathbb{N}_+}$ and $\{v^k\}_{k \in \mathbb{N}_+}$ be two positive scalar sequences. Then, the following hold:*

*(a) If $\lim_{k \to \infty} \beta^k = 0$, then*

$$\lim_{k \to \infty} \sum_{t=1}^{k} (\lambda)^{k-t} \beta^t = 0.$$

(b) If $\sum_{k=1}^{\infty} (\beta^k)^2 < \infty$ and $\sum_{k=1}^{\infty} (v^k)^2 < \infty$, then

1) $\lim_{k\to\infty} \sum_{t=1}^{k} \sum_{l=1}^{t} (\lambda)^{t-l} (\beta^l)^2 < \infty$;

2) $\lim_{k\to\infty} \sum_{t=1}^{k} \sum_{l=1}^{t} (\lambda)^{t-l} \beta^t v^l < \infty$.

$\square$

Consider (190): Since $\rho \in (0,1)$, invoking Lemma III.12(a) and $\lim_{k\to\infty} \|\boldsymbol{\varepsilon}^k\| = 0$, we have

$$\lim_{k\to\infty} \sum_{t=0}^{k-1} (\rho)^{k-1-t} \|\boldsymbol{\varepsilon}^t\| = 0; \qquad (191)$$

hence $\lim_{k\to\infty} \|\mathbf{x}_{(i)}^k - (1/I) \cdot \sum_{j=1}^{I} \phi_{(j)}^k \mathbf{x}_{(j)}^k\| = 0$. Note that, in this case, the rate of convergence to the weighted average may not be geometric.

Consider, as special case, the tracking algorithm (185). i.e., set in (188) $\boldsymbol{\varepsilon}_i^k = (\phi_{(i)}^{k+1})^{-1}(\mathbf{u}_i^{k+1} - \mathbf{u}_i^k)$. We have $(1/I) \cdot \sum_{i=1}^{I} \phi_{(i)}^k \mathbf{x}_{(i)}^k = \bar{\mathbf{u}}^k$, which proves Theorem III.10.

**Case 3: Bounded error**. Since $\|\boldsymbol{\varepsilon}_i^k\| \leq M$, $i = 1, \ldots, I$, the consensus error in (190) can be bounded as

$$\begin{aligned}\left\| \mathbf{x}_{(i)}^k - \frac{1}{I} \sum_{j=1}^{I} \phi_{(j)}^k \mathbf{x}_{(j)}^k \right\| &\leq c \left( (\rho)^k \|\mathbf{x}^0\| + M \cdot \sum_{t=0}^{k-1} (\rho)^{k-1-t} \right) \\ &\leq c \left( (\rho)^k \|\mathbf{x}^0\| + M \cdot \frac{1-(\rho)^{k-1}}{1-\rho} \right).\end{aligned} \qquad (192)$$

### III.2.5 Proof of Theorem III.11

We prove now Theorem III.11, following the analysis in [211]. We first introduce some intermediate results and useful notation.

#### Preliminaries

It is convenient to rewrite the perturbed consensus protocol (188) in a more compact form. To do so, let us introduce the following notation: given the matrix $\mathbf{A}^k$ compliant to the graph $G^k$ (cf. Assumption III.3) and $\mathbf{W}^k$ defined in (179), let

$$\mathbf{x}^k \triangleq [\mathbf{x}_{(1)}^{kT}, \ldots, \mathbf{x}_{(I)}^{kT}]^T, \qquad (193)$$

$$\boldsymbol{\varepsilon}^k \triangleq [\boldsymbol{\varepsilon}_1^{kT}, \ldots, \boldsymbol{\varepsilon}_I^{kT}]^T, \qquad (194)$$

$$\boldsymbol{\phi}^k \triangleq [\phi_{(1)}^k, \ldots, \phi_{(I)}^k]^T, \qquad (195)$$

$$\boldsymbol{\Phi}^k \triangleq \text{Diag}(\boldsymbol{\phi}^k), \qquad (196)$$

$$\widehat{\boldsymbol{\Phi}}^k \triangleq \boldsymbol{\Phi}^k \otimes \mathbf{I}_m, \qquad (197)$$

$$\widehat{\mathbf{A}}^k \triangleq \mathbf{A}^k \otimes \mathbf{I}_m, \tag{198}$$

$$\widehat{\mathbf{W}}^k \triangleq \mathbf{W}^k \otimes \mathbf{I}_m, \tag{199}$$

where $\text{Diag}(\bullet)$ denotes a diagonal matrix whose diagonal entries are the elements of the vector argument. Under the column stochasticity of $\mathbf{A}^k$ (Assumption III.7), it is not difficult to check that the following relationship exists between $\mathbf{W}^k$ and $\mathbf{A}^k$ (and $\widehat{\mathbf{W}}^k$ and $\widehat{\mathbf{A}}^k$):

$$\mathbf{W}^k = \left(\boldsymbol{\Phi}^{k+1}\right)^{-1} \mathbf{A}^k \boldsymbol{\Phi}^k \quad \text{and} \quad \widehat{\mathbf{W}}^k = \left(\widehat{\boldsymbol{\Phi}}^{k+1}\right)^{-1} \widehat{\mathbf{A}}^k \widehat{\boldsymbol{\Phi}}^k. \tag{200}$$

Using the above notation, the perturbed consensus protocol (188) can be rewritten in matrix form as

$$\boldsymbol{\phi}^{k+1} = \mathbf{A}^k \boldsymbol{\phi}^k \quad \text{and} \quad \mathbf{x}^{k+1} = \widehat{\mathbf{W}}^k \mathbf{x}^k + \boldsymbol{\varepsilon}^k. \tag{201}$$

To study the dynamics of the consensus error in (201), let us introduce the matrix products: given $k,t \in \mathbb{N}_+$, with $k \geq t$,

$$\mathbf{A}^{k:t} \triangleq \begin{cases} \mathbf{A}^k \mathbf{A}^{k-1} \cdots \mathbf{A}^t, & \text{if } k > t \geq 0, \\ \mathbf{A}^k, & \text{if } k = t \geq 0, \end{cases} \tag{202}$$

$$\mathbf{W}^{k:t} \triangleq \begin{cases} \mathbf{W}^k \mathbf{W}^{k-1} \cdots \mathbf{W}^t, & \text{if } k > t \geq 0, \\ \mathbf{W}^k, & \text{if } k = t \geq 0, \end{cases} \tag{203}$$

and

$$\widehat{\mathbf{A}}^{k:t} \triangleq \mathbf{A}^{k:t} \otimes \mathbf{I}_m, \tag{204}$$

$$\widehat{\mathbf{W}}^{k:t} \triangleq \mathbf{W}^{k:t} \otimes \mathbf{I}_m. \tag{205}$$

Define the weight-averaging matrix as

$$\mathbf{J}_{\boldsymbol{\phi}^k} \triangleq \frac{1}{I} \left(\mathbf{1} (\boldsymbol{\phi}^k)^T\right) \otimes \mathbf{I}_m, \tag{206}$$

so that

$$\mathbf{J}_{\boldsymbol{\phi}^k} \mathbf{x}^k = \mathbf{1} \otimes \frac{1}{I} \sum_{i=1}^{I} \phi_{(i)}^k \mathbf{x}_{(i)}^k,$$

i.e., $\mathbf{J}_{\boldsymbol{\phi}^k} \mathbf{x}^k$ computes the average of $\mathbf{x}_{(i)}^k$ weighted by $(\phi_{(i)}^k)_{i=1}^I$ and stacks it $I$ times in a column vector. Under the column stochasticity of $\mathbf{A}^k$ (Assumption III.7), it is not difficult to check that the following property holds between $\mathbf{J}_{\boldsymbol{\phi}^k}$ and $\widehat{\mathbf{W}}^{k:t}$: for all $k,t \in \mathbb{N}_+$, with $k \geq t$,

$$\mathbf{J}_{\boldsymbol{\phi}^{k+1}} \widehat{\mathbf{W}}^{k:t} = \mathbf{J}_{\boldsymbol{\phi}^t} = \widehat{\mathbf{W}}^{k:t} \mathbf{J}_{\boldsymbol{\phi}^t}. \tag{207}$$

Convergence of the perturbed consensus protocol boils down to studying the dynamics of $\|\widehat{\mathbf{W}}^{k:t} - \mathbf{J}_{\boldsymbol{\phi}^t}\|$ (this will be more clear in the next subsection). The follow-

ing lemma shows that, in the setting of Theorem III.11, $\widehat{\mathbf{W}}^{k:t}$ converges geometrically to $\mathbf{J}_{\boldsymbol{\phi}^t}$; the proof of the lemma is omitted and can be found in [211, Lemma 2].

**Lemma III.13.** *Let $\{G^k\}_{k\in\mathbb{N}_+}$ be a sequence of graphs satisfying Assumption III.2; let $\{\mathbf{A}^k\}_{k\in\mathbb{N}_+}$ be a sequence of weight matrices satisfying Assumption III.7; and let $\{\mathbf{W}^k\}_{k\in\mathbb{N}_+}$ be the sequence of (row-stochastic) matrices with $\mathbf{W}^k$ related to $\mathbf{A}^k$ by (200). Then, the following holds: for all $k, t \in \mathbb{N}_+$, with $k \geq t$,*

$$\left\|\widehat{\mathbf{W}}^{k:t} - \mathbf{J}_{\boldsymbol{\phi}^t}\right\|_2 \leq c \cdot (\rho)^{k-t+1}, \tag{208}$$

*where $c > 0$ and $\rho \in (0,1)$ are defined in Theorem III.11.*

**Proof of Theorem III.11**

We are now ready to prove the theorem. We start rewriting the dynamics of the consensus error $\mathbf{x}^k - \mathbf{J}_{\boldsymbol{\phi}^k}\mathbf{x}^k$ in a form that permits the application of Lemma III.13. Applying the $x$-update in (201) recursively, we have

$$\begin{aligned}
\mathbf{x}^k &= \widehat{\mathbf{W}}^{k-1}\mathbf{x}^{k-1} + \boldsymbol{\varepsilon}^{k-1} \\
&= \widehat{\mathbf{W}}^{k-1}\left(\widehat{\mathbf{W}}^{k-2}\mathbf{x}^{k-2} + \boldsymbol{\varepsilon}^{k-2}\right) + \boldsymbol{\varepsilon}^{k-1} \\
&= \cdots = \widehat{\mathbf{W}}^{k-1:0}\mathbf{x}^0 + \sum_{t=1}^{k-1}\widehat{\mathbf{W}}^{k-1:t}\boldsymbol{\varepsilon}^{t-1} + \boldsymbol{\varepsilon}^{k-1}.
\end{aligned} \tag{209}$$

Therefore, the weighted average $\mathbf{J}_{\boldsymbol{\phi}^k}\mathbf{x}^k$ can be written as

$$\begin{aligned}
\mathbf{J}_{\boldsymbol{\phi}^k}\mathbf{x}^k &\stackrel{(209)}{=} \mathbf{J}_{\boldsymbol{\phi}^k}\widehat{\mathbf{W}}^{k-1:0}\mathbf{x}^0 + \mathbf{J}_{\boldsymbol{\phi}^k}\sum_{t=1}^{k-1}\widehat{\mathbf{W}}^{k-1:t}\boldsymbol{\varepsilon}^{t-1} + \mathbf{J}_{\boldsymbol{\phi}^k}\boldsymbol{\varepsilon}^{k-1} \\
&\stackrel{(207)}{=} \mathbf{J}_{\boldsymbol{\phi}^0}\mathbf{x}^0 + \sum_{t=1}^{k-1}\mathbf{J}_{\boldsymbol{\phi}^t}\boldsymbol{\varepsilon}^{t-1} + \mathbf{J}_{\boldsymbol{\phi}^k}\boldsymbol{\varepsilon}^{k-1}.
\end{aligned} \tag{210}$$

Using (209) and (210) and Lemma III.13, the consensus error can be bounded as

$$\begin{aligned}
&\left\|\mathbf{x}^k - \mathbf{J}_{\boldsymbol{\phi}^k}\mathbf{x}^k\right\| \\
&= \left\|\left(\widehat{\mathbf{W}}^{k-1:0} - \mathbf{J}_{\boldsymbol{\phi}^0}\right)\mathbf{x}^0 + \sum_{t=1}^{k-1}\left(\widehat{\mathbf{W}}^{k-1:t} - \mathbf{J}_{\boldsymbol{\phi}^t}\right)\boldsymbol{\varepsilon}^{t-1} + \left(\mathbf{I} - \mathbf{J}_{\boldsymbol{\phi}^k}\right)\boldsymbol{\varepsilon}^{k-1}\right\| \\
&\stackrel{(208)}{\leq} c \cdot (\rho)^k \|\mathbf{x}^0\| + c \cdot \sum_{t=1}^{k-1}(\rho)^{k-t}\|\boldsymbol{\varepsilon}^{t-1}\| + \underbrace{\left\|\mathbf{I} - \mathbf{J}_{\boldsymbol{\phi}^k}\right\|_2}_{\leq \sqrt{2I}}\|\boldsymbol{\varepsilon}^{k-1}\| \\
&\leq c \cdot \left((\rho)^k\|\mathbf{x}^0\| + \sum_{t=0}^{k-1}(\rho)^{k-t-1}\|\boldsymbol{\varepsilon}^t\|\right),
\end{aligned} \tag{211}$$

where in the last inequality we used $c > \sqrt{2I}$. The inequality $\|\mathbf{I} - \mathbf{J}_{\boldsymbol{\phi}^k}\|_2 \leq \sqrt{2I}$ can be proved as follows. Let $\mathbf{z} \in \mathbb{R}^{I \cdot m}$ be an arbitrary vector; partition $\mathbf{z}$ as $\mathbf{z} = (\mathbf{z}_i)_{i=1}^I$, with each $\mathbf{z}_i \in \mathbb{R}^m$. Then,

$$\left\| \mathbf{z} - \mathbf{J}_{\boldsymbol{\phi}^k}\mathbf{z} \right\| \leq \left\| \mathbf{z} - \mathbf{J_1}\mathbf{z} \right\| + \left\| \mathbf{J_1}\mathbf{z} - \mathbf{J}_{\boldsymbol{\phi}^k}\mathbf{z} \right\| \stackrel{(a)}{\leq} \|\mathbf{z}\| + \frac{\sqrt{I}}{I} \left\| \sum_{i=1}^I \mathbf{z}_i - \sum_{i=1}^I \phi_i^k \mathbf{z}_i \right\|$$

$$\leq \|\mathbf{z}\| + \frac{\sqrt{I}}{I} \sqrt{I^2 - I} \, \|\mathbf{z}\| \leq \sqrt{2I} \, \|\mathbf{z}\|,$$

where in (a) we used $\|\mathbf{I} - \mathbf{J_1}\| = 1$ (note that $\mathbf{I} - \mathbf{J_1}$ is a Toeplitz matrix, with diagonal entries equal to $1 - 1/I$ and off-diagonal entries all equal to $-1/I$; therefore, its eigenspectrum is given by $\{0, 1, \ldots, 1\}$).

The inequality (211) proves the theorem. □

### III.3. Distributed SCA over Time-Varying Digraphs

We are now ready to introduce the proposed distributed algorithmic framework to solve Problem (167), which combines SCA techniques (introduced in the previous lectures) with the consensus/tracking protocols described in Sec. III.2. We consider the optimization over time-varying (B-strongly connected) digraphs (cf. Assumption III.2); distributed algorithms for undirected or time invariant networks can be obtained as special cases.

As already anticipated, each agent $i$ maintains and updates iteratively a *local copy* $\mathbf{x}_{(i)}$ of the global variable $\mathbf{x}$, along with an auxiliary variable $\mathbf{y}_{(i)} \in \mathbb{R}^m$, whose goal is to track locally the average of the gradients $(1/I) \cdot \sum_{i=1}^I \nabla f_i$ (the importance of this extra variable will be clarified shortly), an information that is not available at the agent's side; let $\mathbf{x}_{(i)}^k$ and $\mathbf{y}_{(i)}^k$ denote the values of $\mathbf{x}_{(i)}$ and $\mathbf{y}_{(i)}$ at iteration $k$, respectively. To update these variables, *two major steps* are performed iteratively, namely:

**Step 1–Local SCA (optimization):** Given $\mathbf{x}_{(i)}^k$ and $\mathbf{y}_{(i)}^k$, each agent $i$ solves a convexification of Problem (167), wherein $V$ is replaced by a suitably chosen strongly convex surrogate, which is built using only the available local information $\mathbf{x}_{(i)}^k$ and $\mathbf{y}_{(i)}^k$;

**Step 2–Communication:** All the agents broadcast the solutions computed in Step 1, and update their own variables $\mathbf{x}_{(i)}^k \to \mathbf{x}_{(i)}^{k+1}$ and $\mathbf{y}_{(i)}^k \to \mathbf{y}_{(i)}^{k+1}$, based on the information received from their neighbors.

The two steps above need to be designed so that: i) all the $\mathbf{x}_{(i)}^k$ will be asymptotically consensual, that is, $\lim_{k \to \infty} \|\mathbf{x}_{(i)}^k - (1/I) \cdot \sum_{j=1}^I \mathbf{x}_{(j)}^k\| = 0$, for all $i$; and ii) every limit point of $(1/I) \cdot \sum_{j=1}^I \mathbf{x}_{(j)}^k$ is a stationary solution of Problem (167). We describe next the above two steps in detail.

**Step 1: Local SCA (optimization).** Each agent $i$ faces two issues to solve Problem (167), namely: $f_i$ is not convex and $\sum_{j \neq i} f_j$ is not known. To cope with the first issue, we leverage the SCA techniques introduced in the previous lectures. More specifically, at iteration $k$, agent $i$ solves a convexification of $V$ in (167) having the following form

$$\widehat{\mathbf{x}}_i(\mathbf{x}_{(i)}^k) \triangleq \underset{\mathbf{x}_{(i)} \in X}{\operatorname{argmin}} \; \widetilde{F}_i(\mathbf{x}_{(i)} \mid \mathbf{x}_{(i)}^k) + G(\mathbf{x}_{(i)}), \tag{212}$$

where $\widetilde{F}_i : O \times O \to \mathbb{R}$ is a suitably chosen surrogate of $F$. To guarantee that a fixed point of $\widehat{\mathbf{x}}_i(\bullet)$ is a stationary solution of (167), a naive application of the SCA theory developed in Lecture II, would call $\widetilde{F}_i$ to satisfy the following gradient consistency condition (cf. Assumption II.2, Sec. II.2): $\widetilde{F}_i$ is $C^1$ on $O$ and

$$\nabla \widetilde{F}_i(\mathbf{x}_{(i)}^k \mid \mathbf{x}_{(i)}^k) = \nabla F(\mathbf{x}_{(i)}^k) = \nabla f_i(\mathbf{x}_{(i)}^k) + \sum_{j \neq i} \nabla f_j(\mathbf{x}_{(i)}^k), \quad k \in \mathbb{N}_+. \tag{213}$$

For example, a surrogate function satisfying the above condition would be:

$$\widetilde{F}_i(\mathbf{x}_{(i)} \mid \mathbf{x}_{(i)}^k) = \widetilde{f}_i(\mathbf{x}_{(i)} \mid \mathbf{x}_{(i)}^k) + \left( \sum_{j \neq i} \nabla f_j(\mathbf{x}_{(i)}^k) \right)^T (\mathbf{x}_{(i)} - \mathbf{x}_{(i)}^k), \tag{214}$$

where $\widetilde{f}_i(\bullet \mid \mathbf{x}_{(i)}^k) : O \to \mathbb{R}$ is a strongly convex surrogate of $f_i$ on $X$, consistent to $f_i$, in the following sense (cf. Assumption II.2).

**Assumption III.14.** *Each function $\widetilde{f}_i : O \times O \to \mathbb{R}$ satisfies the following conditions:*

1. $\widetilde{f}_i(\bullet \mid \mathbf{x})$ *is $\tau_i$-strongly convex on $X$, for all $\mathbf{x} \in X$;*
2. $\widetilde{f}_i(\bullet \mid \mathbf{x})$ *is $C^1$ on $O$ and $\nabla \widetilde{f}_i(\mathbf{x} \mid \mathbf{x}) = \nabla f_i(\mathbf{x})$, for all $\mathbf{x} \in X$;*
3. $\nabla \widetilde{f}_i(\mathbf{x} \mid \bullet)$ *is $\widetilde{L}_i$-Lipschitz on $X$, for all $\mathbf{x} \in X$.*

Unfortunately, the surrogate function in (214) cannot be used by agent $i$, because of the lack of knowledge of $\sum_{j \neq i} \nabla f_j(\mathbf{x}_{(i)}^k)$; hence a gradient consistency condition in the form (213) cannot be enforced in such a distributed setting.

To cope with this issue, the idea, first proposed in [148] and further developed in [211], is to replace $\sum_{j \neq i} \nabla f_j(\mathbf{x}_{(i)}^k)$ in (214) [and in (213)] with a *local*, *asymptotically consistent*, approximation, so that condition (213) will be satisfied in the limit, as $k \to \infty$. This can be accomplished, e.g., using the following surrogate function:

$$\widetilde{F}_i\left(\mathbf{x}_{(i)} \mid \mathbf{x}_{(i)}^k, \mathbf{y}_{(i)}^k\right) = \widetilde{f}_i\left(\mathbf{x}_{(i)} \mid \mathbf{x}_{(i)}^k\right) + \left(I \cdot \mathbf{y}_{(i)}^k - \nabla f_i(\mathbf{x}_{(i)}^k)\right)^T \left(\mathbf{x}_{(i)} - \mathbf{x}_{(i)}^k\right), \tag{215}$$

where $\widetilde{f}_i$ is defined as in (214); and $\mathbf{y}_{(i)}^k$ in (215) is an auxiliary variable controlled by agent $i$, aiming at tracking locally the average of the gradients $(1/I) \cdot \sum_{j=1}^I \nabla f_j(\mathbf{x}_{(i)}^k)$, that is, $\lim_{k \to \infty} \|\mathbf{y}_{(i)}^k - (1/I) \cdot \sum_{j=1}^I \nabla f_j(\mathbf{x}_{(i)}^k)\| = 0$. This explains the role of the linear term in (215): under the claimed tracking property of $\mathbf{y}_{(i)}^k$, we have

$$\lim_{k \to \infty} \left\| \left( I \cdot \mathbf{y}_{(i)}^k - \nabla f_i(\mathbf{x}_{(i)}^k) \right) - \sum_{j \neq i} \nabla f_j(\mathbf{x}_{(i)}^k) \right\| = 0, \tag{216}$$

which would guarantee that the gradient consistency condition (213), with now $\widetilde{F}_i(\bullet \mid \mathbf{x}_{(i)}^k)$ replaced by $\widetilde{F}_i(\bullet \mid \mathbf{x}_{(i)}^k, \mathbf{y}_{(i)}^k)$ in (215), will be *asymptotically* satisfied, that is,

$$\lim_{k \to \infty} \left\| \nabla \widetilde{F}_i(\mathbf{x}_{(i)}^k \mid \mathbf{x}_{(i)}^k, \mathbf{y}_{(i)}^k) - \nabla F(\mathbf{x}_{(i)}^k) \right\| = 0.$$

As it will be shown later, this relaxed condition is in fact enough to prove that, if convergence and consensus are asymptotically achieved, the limit point of all the local variables $\mathbf{x}_{(i)}^k$ is a stationary solution of Problem (167). Leveraging the distributed tracking protocol (185) (cf. Sec. III.2.3), in Step 2 below, we show how to devise an update for the $\mathbf{y}_{(i)}^k$ variables that uses only local information and such that (216) asymptotically holds.

Using $\widetilde{F}_i(\bullet \mid \mathbf{x}_{(i)}^k, \mathbf{y}_{(i)}^k)$ defined in (215), the local optimization step performing by each agent $i$ consists then in solving the following strongly convex problem:

$$\widetilde{\mathbf{x}}_i^k \triangleq \underset{\mathbf{x}_{(i)} \in X}{\operatorname{argmin}} \ \widetilde{F}_i \left( \mathbf{x}_{(i)} \mid \mathbf{x}_{(i)}^k, \mathbf{y}_{(i)}^k \right) + G \left( \mathbf{x}_{(i)} \right), \tag{217}$$

followed by the step-size update

$$\mathbf{x}_{(i)}^{k+1/2} = \mathbf{x}_{(i)}^k + \gamma^k \left( \widetilde{\mathbf{x}}_i^k - \mathbf{x}_{(i)}^k \right), \tag{218}$$

where $\gamma^k \in (0, 1]$ is the step-size, to be properly chosen.

**Step 2–Communication.** Given $\mathbf{x}_{(i)}^{k+1/2}$ and $\mathbf{y}_{(i)}^k$, each agent $i$ communicates with its current neighbors in order to achieve asymptotic consensus on $\mathbf{x}_{(i)}$'s as well as track $(1/I) \cdot \sum_{j=1}^I \nabla f_j(\mathbf{x}_{(i)}^k)$ by $\mathbf{y}_{(i)}^k$. Both goals can be accomplished using (two instances of) the condensed perturbed push-sum protocol (188), introduced in Sec. III.2.4.

More specifically, after obtaining $\mathbf{x}_{(j)}^{k+1/2}$ from its neighbors, each agent $i$ updates its own local estimate $\mathbf{x}_{(i)}$ employing:

$$\phi_{(i)}^{k+1} = \sum_{j=1}^I a_{ij}^k \phi_{(j)}^k; \tag{219}$$

$$\mathbf{x}_{(i)}^{k+1} = \frac{1}{\phi_{(i)}^{k+1}} \sum_{j=1}^I a_{ij}^k \phi_{(j)}^k \mathbf{x}_{(j)}^{k+1/2}, \tag{220}$$

where the weights $\mathbf{A}^k \triangleq (a_{ij}^k)_{i=1}^I$ are chosen according to Assumption III.7; and the $\phi_{(i)}^k$ variables are initialized to $\phi_{(i)}^0 = 1$. This update can be clearly implemented locally: All agents i) send their local variables $\phi_{(j)}^k$ and $\phi_{(j)}^k \mathbf{x}_{(j)}^{k+1/2}$ to their out-neighbors; and ii) linearly combine with coefficients $a_{ij}^k$ the information coming from their in-neighbors.

A local update for the $\mathbf{y}_{(i)}^k$ variable aiming at tracking $(1/I) \cdot \sum_{j=1}^I \nabla f_j(\mathbf{x}_{(i)}^k)$ can be readily obtained invoking the distributed tracking protocol (185), and setting therein

$\mathbf{u}_i^k \triangleq \nabla f_i(\mathbf{x}_{(i)}^k)$. This leads to

$$\mathbf{y}_{(i)}^{k+1} = \frac{1}{\phi_{(i)}^{k+1}} \sum_{j=1}^{I} a_{ij}^k \phi_{(j)}^k \mathbf{y}_{(j)}^k + \frac{1}{\phi_{(i)}^{k+1}} \left( \nabla f_i(\mathbf{x}_{(i)}^{k+1}) - \nabla f_i(\mathbf{x}_{(i)}^k) \right), \qquad (221)$$

where $\phi_{(i)}^{k+1}$ is defined in (219), and $\mathbf{y}_{(i)}^0 = \nabla f_i(\mathbf{x}_{(i)}^0)$. Note that, as for $\mathbf{x}_{(i)}^{k+1}$, the update of $\mathbf{y}_{(i)}^k$ is performed by agent $i$ using only the information coming from its neighbors, with the same signaling as for (219)-(220).

The described distributed SCA method (Step 1 and Step 2) is summarized in Algorithm 10, and termed *distributed Successive cONvex Approximation algorithm over Time-varying digrAphs* (SONATA).

---

**Algorithm 10: SCA over Time-varying Digraphs (SONATA)**

---

**Data** : $\mathbf{x}_{(i)}^0 \in X$, $\phi_{(i)}^0 = 1$, $\mathbf{y}_{(i)}^0 = \nabla f_i(\mathbf{x}_{(i)}^0)$, for all $i = 1, \ldots, I$; $\{\gamma^k \in (0,1]\}_{k \in \mathbb{N}_+}$.
Set $k = 0$.
(S.1) : If $\mathbf{x}^k$ satisfies a termination criterion: STOP;
(S.2) : **Local SCA.** Each agent $i$ computes

$$\widetilde{\mathbf{x}}_i^k \triangleq \underset{\mathbf{x}_{(i)} \in X}{\operatorname{argmin}} \; \widetilde{f}_i \left( \mathbf{x}_{(i)} \,|\, \mathbf{x}_{(i)}^k \right) + \left( I \cdot \mathbf{y}_{(i)}^k - \nabla f_i(\mathbf{x}_{(i)}^k) \right)^T \left( \mathbf{x}_{(i)} - \mathbf{x}_{(i)}^k \right) + G \left( \mathbf{x}_{(i)} \right),$$

$$\mathbf{x}_{(i)}^{k+1/2} = \mathbf{x}_{(i)}^k + \gamma^k \left( \widetilde{\mathbf{x}}_i^k - \mathbf{x}_{(i)}^k \right);$$

(S.3) : **Averaging and gradient tracking.**
Each agent $i$ sends out its local variables $\phi_{(i)}^k$, $\mathbf{x}_{(i)}^{k+1/2}$ and $\mathbf{y}_{(i)}^k$, and receives $\phi_{(j)}^k$, $\mathbf{x}_{(j)}^{k+1/2}$ and $\mathbf{y}_{(j)}^k$, with $j \in N_i^{\text{in},k} \setminus \{i\}$. Then, it updates:

$$\phi_{(i)}^{k+1} = \sum_{j=1}^{I} a_{ij}^k \phi_{(j)}^k,$$

$$\mathbf{x}_{(i)}^{k+1} = \frac{1}{\phi_{(i)}^{k+1}} \sum_{j=1}^{I} a_{ij}^k \phi_{(j)}^k \mathbf{x}_{(j)}^{k+1/2},$$

$$\mathbf{y}_{(i)}^{k+1} = \frac{1}{\phi_{(i)}^{k+1}} \sum_{j=1}^{I} a_{ij}^k \phi_{(j)}^k \mathbf{y}_{(j)}^k + \frac{1}{\phi_{(i)}^{k+1}} \left( \nabla f_i(\mathbf{x}_{(i)}^{k+1}) - \nabla f_i(\mathbf{x}_{(i)}^k) \right);$$

(S.4) : $k \leftarrow k+1$, and go to (S.1).

---

Convergence of Algorithm 10 is stated in Theorem III.16 below. We first introduce a standard condition on the steps-size $\gamma^k$ and a proper merit function assessing the convergence of the algorithm.

**Assumption III.15.** *The step-size $\gamma^k \in (0,1]$ satisfies the standard diminishing rule: $\lim_{k \to \infty} \gamma^k = 0$ and $\sum_{k=0}^{\infty} \gamma^k = +\infty$.*

Given $\{\mathbf{x}^k \triangleq (\mathbf{x}_{(i)}^k)_{i=1}^I\}_{k \in \mathbb{N}_+}$ generated by Algorithm 10, convergence of the algorithm is stated measuring the distance of the average sequence $\bar{\mathbf{x}}^k \triangleq (1/I) \cdot \sum_{i=1}^I \mathbf{x}_{(i)}^k$ from optimality and well as the consensus disagreement among the local variables $\mathbf{x}_{(i)}^k$'s. More specifically, let us introduce the following function as a measure of optimality:

$$J(\bar{\mathbf{x}}^k) \triangleq \left\| \bar{\mathbf{x}}^k - \operatorname*{argmin}_{\mathbf{y} \in X} \left\{ \nabla F(\bar{\mathbf{x}}^k)^T (\mathbf{y} - \bar{\mathbf{x}}^k) + \frac{1}{2} \|\mathbf{y} - \bar{\mathbf{x}}^k\|_2^2 + G(\mathbf{y}) \right\} \right\|. \quad (222)$$

Note that $J$ is a valid measure of stationarity because it is continuous and $J(\bar{\mathbf{x}}^\infty) = 0$ if and only if $\bar{\mathbf{x}}^\infty$ is a d-stationary solution of Problem (167). The consensus disagreement at iteration $k$ is defined as

$$D(\mathbf{x}^k) \triangleq \|\mathbf{x}^k - \mathbf{1}_I \otimes \bar{\mathbf{x}}^k\|_2.$$

Note that $D$ is equal to 0 if and only if all the $\mathbf{x}_{(i)}^k$'s are consensual. We combine the metrics $J$ and $D$ in a single merit function, defined as

$$M(\mathbf{x}^k) \triangleq \max \left\{ J(\bar{\mathbf{x}}^k)^2, D(\mathbf{x}^k)^2 \right\},$$

which captures the progresses of the algorithm towards optimality and consensus.

We are now ready to state the main convergence results for Algorithm 10.

**Theorem III.16.** *Consider Problem* (167) *under Assumption III.1; and let $\{G^k\}_{k \in \mathbb{N}_+}$ be a sequence of graphs satisfying Assumption III.2. Let $\{\mathbf{x}^k \triangleq (\mathbf{x}_{(i)}^k)_{i=1}^I\}_{k \in \mathbb{N}_+}$ be the sequence generated by Algorithm 10 under Assumptions III.7 and III.14; and let $\bar{\mathbf{x}}^k = (1/I) \cdot \sum_{i=1}^I \mathbf{x}_{(i)}^k$ be the average sequence. Furthermore, suppose that either one of the following is satisfied.*

*(a) (`diminishing step-size`): The step-size $\gamma^k$ satisfies Assumption III.15;*

*(b) (`constant step-size`): The step-size $\gamma^k$ is fixed–$\gamma^k = \gamma$, for all $k \in \mathbb{N}_+$– and it is sufficiently small (see [211, Theorem 5] for the specific expression of the upper bound on $\gamma$).*

*Then, there holds*

$$\lim_{k \to \infty} M(\mathbf{x}^k) = 0. \quad (223)$$

*Proof.* See Sec. III.3.2. □

Theorem III.17 below provides an upper bound on the number of iterations needed to decrease $M_U(\mathbf{x}^k)$ below a given accuracy $\varepsilon > 0$; we omit the proof and refer to [211, Theorem 6] for more details.

**Theorem III.17.** *Consider Problem* (167) *and Algorithm 10 in the setting of Theorem III.16. Given $\varepsilon > 0$, let $T_\varepsilon$ be the first iteration $k$ such that $M(\mathbf{x}^k) \leq \varepsilon$.*

(a) (`diminishing step-size`): *Suppose that the step-size $\gamma^k$ satisfies Assumption III.15. Then,*

$$T_\varepsilon \leq \inf\left\{k \in \mathbb{N}_+ : \sum_{t=0}^{k} \gamma^t \geq \frac{B_0}{\varepsilon}\right\},$$

*where $B_0 > 0$ is a constant independent on $\varepsilon$ [211];*

(b) (`constant step-size`): *Suppose that the step-size $\gamma^k$ is fixed, $\gamma^k = \gamma$, for all $k \in \mathbb{N}_+$. Then, there exists a sufficiently small $\bar{\gamma} \in (0,1]$–independent of $\varepsilon$ (see [211, Theorem 6] for the specific expression of $\bar{\gamma}$)–such that, if $\gamma \in (0, \bar{\gamma}]$, then it holds*

$$T_\varepsilon = O\left(\frac{1}{\varepsilon}\right).$$

**Discussion and generalizations**

**On the convergence:** Stating convergence for constrained (nonsmooth) optimization problems in the form (167), Theorem III.16 and Theorems III.17 (proved in our work [211]) significantly enlarge the class of convex and nonconvex problems which distributed algorithms can be applied to with convergence guarantees. We remark that convergence is established without requiring that the (sub)gradients of $F$ or $G$ is bounded; this is a major improvement with respect to current distributed methods for nonconvex problems [20, 148, 232, 245] and nonsmooth convex ones [168].

We remark that convergence (as stated in the above theorems) can also be established weakening the assumption on the strongly convexity of the surrogates $\widetilde{f}_i$ (Assumption III.14) to just convexity, as long as the feasible set $X$ is compact. Also, with mild additional assumptions on $G$–see Lecture II–convergence can be also proved in the case wherein agents solve their subproblems (217) inexactly.

**ATC- versus CAA updates.** As a final remark, we note that variants of SONATA wherein the order of the consensus, tracking, and local updates are differently combined, are still convergent, in the sense of Theorems III.16 and III.17. We briefly elaborate on this issue next.

Using a jargon well established in the literature [202], the update of the $x$-variables in Step 3 of Algorithm 10 is in the form of so-called *Adapt-Then-Combine (ATC)* strategy: eliminating the intermediate variable $\mathbf{x}_{(i)}^{k+1/2}$, each $\mathbf{x}_i^k$ follows the dynamic

$$\mathbf{x}_{(i)}^{k+1} = \frac{1}{\phi_{(i)}^{k+1}} \sum_{j=1}^{I} a_{ij}^k \phi_{(j)}^k \left(\mathbf{x}_{(j)}^k + \gamma^k \left(\widetilde{\mathbf{x}}_j^k - \mathbf{x}_{(j)}^k\right)\right). \tag{224}$$

The name ATC comes from the form of (224): each agent $i$ first "adapts" its local copy $\mathbf{x}_{(i)}^k$ moving along the direction $\widetilde{\mathbf{x}}_i^k - \mathbf{x}_{(i)}^k$, that is, $\mathbf{x}_{(i)}^k \to \mathbf{x}_{(i)}^k + \gamma^k (\widetilde{\mathbf{x}}_i^k - \mathbf{x}_{(i)}^k)$; and then it "combines" its new update with that of its in-neighbors.

As an alternative to (224), one can employ the so-called *Combine-And-Adapt (CAA)* update (also termed "consensus strategy" in [202]), which reads

$$\mathbf{x}_{(i)}^{k+1} = \frac{1}{\phi_{(i)}^{k+1}} \sum_{j=1}^{I} a_{ij}^k \phi_{(j)}^k \mathbf{x}_{(j)}^k + \frac{\phi_{(i)}^k}{\phi_{(i)}^{k+1}} \cdot \gamma^k \left( \widetilde{\mathbf{x}}_i^k - \mathbf{x}_{(i)}^k \right). \tag{225}$$

According to this protocol, each agent $i$ first "combines" (by weighted-averaging) its current $\mathbf{x}_{(i)}^k$ with those of its neighbors, and then "adapt" the resulting update moving along the direction $\widetilde{\mathbf{x}}_i^k - \mathbf{x}_{(i)}^k$. Note that, when dealing with constraint optimization, the CAA update in general does not preserve the feasibility of the iterates while the ATC protocol deos.

The ATC and CAA protocols can be interchangeably used also in the update of the tracking variables $\mathbf{y}_{(i)}^k$ in Step 3 of SONATA. While the $y$-update as stated in the algorithm [cf. (221)] is in the CAA form, one can also use the ATC-based update, which reads

$$\mathbf{y}_{(i)}^{k+1} = \frac{1}{\phi_{(i)}^{k+1}} \sum_{j=1}^{I} a_{ij}^k \left( \phi_{(j)}^k \mathbf{y}_{(j)}^k + \nabla f_j(\mathbf{x}_{(j)}^{k+1}) - \nabla f_j(\mathbf{x}_{(j)}^k) \right).$$

One can show that the above versions of SONATA are all convergent.

### III.3.1 SONATA (Algorithm 10) and special cases

The SONATA framework represents a gamut of algorithms, each of them corresponding to a specific choice of the surrogate functions, step-size, and weight matrices. In this section, we focus on recent proposals in the literature that built on the idea of distributed gradient tracking [69, 70, 148, 170, 190, 257, 259], and we show that all of them are in fact special cases of SONATA. A more detailed analysis of the state of the art can be found in Sec. III.5.

The idea of tracking the gradient averages through the use of consensus coupled with distributed optimization was independently introduced in the NEXT & SONATA framework [69, 70, 148] & [211, 231] for the general class of (convex) constrained nonsmooth nonconvex problems (167) and in [259] for the special case of strongly convex unconstrained smooth optimization. The algorithmic framework in [69, 70, 148] is applicable to optimization problems over time-varying graphs, but requires the use of doubly stochastic matrices. This assumption was removed in SONATA [211, 231] by using column-stochastic matrices, which are more general than the degree-based column-stochastic matrices of distributed push-sum-based methods. The scheme in [259] is implementable only over undirected fixed graphs. A convergence rate analysis of the scheme in [259] in the case of strongly convex smooth unconstrained optimization problems was later developed in [170, 190] for undirected graphs and in [170] for time-varying directed graphs. Complexity results for NEXT and SONATA for (strongly) convex and nonconvex constrained optimization problems over (time-varying) digraphs can be found in [211, 229]; differently from [170, 190], the analysis in [211, 229] applies to general surrogate functions (satisfying Assumption III.14).

We establish next a formal connection between SONATA and all these schemes.

**Preliminaries**

Since all the aforementioned works but [69, 70, 148] are applicable only to the special instance of Problem (167) wherein $X = \mathbb{R}^m$ (unconstrained), $G = 0$ (only smooth objectives), and $F$ is strongly convex, throughout this section, for a fair comparison, we only consider this setting. We begin customizing Algorithm 10 to this special instance of (167) as follows. Choose each surrogate function $\widetilde{f}_i$ in (215) as

$$\widetilde{f}_i(\mathbf{x}_{(i)} \,|\, \mathbf{x}_{(i)}^k) = f_i(\mathbf{x}_{(i)}^k) + \nabla f_i(\mathbf{x}_{(i)}^k)^T (\mathbf{x}_{(i)} - \mathbf{x}_{(i)}^k) + \frac{I}{2} \left\| \mathbf{x}_{(i)} - \mathbf{x}_{(i)}^k \right\|^2.$$

This leads to the following closed form expression for $\widetilde{\mathbf{x}}_i^k$ in (217) (recall that $X = \mathbb{R}^m$ and $G = 0$):

$$\begin{aligned}
\widetilde{\mathbf{x}}_i^k &= \underset{\mathbf{x}_i}{\mathrm{argmin}} \, \left( I \cdot \mathbf{y}_{(i)}^k \right)^\top \left( \mathbf{x}_{(i)} - \mathbf{x}_{(i)}^k \right) + \frac{I}{2} \left\| \mathbf{x}_{(i)} - \mathbf{x}_{(i)}^k \right\|^2 \\
&= \mathbf{x}_{(i)}^k - \mathbf{y}_{(i)}^k.
\end{aligned} \quad (226)$$

Define now $\mathbf{g}_i^k \triangleq \nabla f_i(\mathbf{x}_{(i)}^k)$, $\mathbf{g}^k = [\mathbf{g}_1^{kT}, \ldots, \mathbf{y}_I^{kT}]^T$, and $\mathbf{y}^k = [\mathbf{y}_{(1)}^{kT}, \ldots, \mathbf{g}_{(I)}^{kT}]^T$; and recall the definitions of $\boldsymbol{\phi}^k$, $\boldsymbol{\Phi}^k$, $\mathbf{A}^k$, $\mathbf{W}^k$ and $\widehat{\mathbf{W}}^k$ introduced in Sec. III.2.5. Using (226), Algorithm 10 under the ATC or CAA updates can be rewritten in vector/matrix form as: for all $k \in \mathbb{N}_+$,

$$\begin{aligned}
\boldsymbol{\phi}^{k+1} &= \mathbf{A}^k \boldsymbol{\phi}^k, \\
\mathbf{W}^k &= (\boldsymbol{\Phi}^{k+1})^{-1} \mathbf{A}^k \boldsymbol{\Phi}^k, \\
\mathbf{x}^{k+1} &= \begin{cases} \widehat{\mathbf{W}}^k \left( \mathbf{x}^k - \gamma^k \mathbf{y}^k \right), & \text{if ATC is employed;} \\ \widehat{\mathbf{W}}^k \mathbf{x}^k - \gamma^k (\widehat{\boldsymbol{\Phi}}^{k+1})^{-1} \widehat{\boldsymbol{\Phi}}^k \mathbf{y}^k, & \text{if CAA is employed;} \end{cases} \\
\mathbf{y}^{k+1} &= \begin{cases} \widehat{\mathbf{W}}^k \left( \mathbf{y}^k + (\widehat{\boldsymbol{\Phi}}^k)^{-1} (\mathbf{g}^{k+1} - \mathbf{g}^k) \right), & \text{if ATC is employed;} \\ \widehat{\mathbf{W}}^k \mathbf{y}^k + (\widehat{\boldsymbol{\Phi}}^{k+1})^{-1} \left( \mathbf{g}^{k+1} - \mathbf{g}^k \right), & \text{if CAA is employed;} \end{cases}
\end{aligned} \quad (227)$$

which we will refer to as *ATC/CAA-SONATA-L* (L stands for "linearized").

In the special case where all $\mathbf{A}^k$ are doubly-stochastic matrices, we have $\mathbf{W}^k = \mathbf{A}^k$ and $\widehat{\mathbf{W}}^k \triangleq \mathbf{A}^k \otimes \mathbf{I}_m$; ATC/CAA-SONATA-L reduces to

$$\begin{aligned}
\mathbf{x}^{k+1} &= \begin{cases} \widehat{\mathbf{W}}^k \left( \mathbf{x}^k - \gamma^k \mathbf{y}^k \right), & \text{for the ATC update;} \\ \widehat{\mathbf{W}}^k \mathbf{x}^k - \gamma^k \mathbf{y}^k, & \text{for the CAA update;} \end{cases} \\
\mathbf{y}^{k+1} &= \begin{cases} \widehat{\mathbf{W}}^k \left( \mathbf{y}^k + \mathbf{g}^{k+1} - \mathbf{g}^k \right), & \text{for the ATC update;} \\ \widehat{\mathbf{W}}^k \mathbf{y}^k + \left( \mathbf{g}^{k+1} - \mathbf{g}^k \right), & \text{for the CAA update;} \end{cases}
\end{aligned} \quad (228)$$

which is referred to as *ATC/CAA-SONATA-NEXT-L* (because the algorithm becomes an instance of NEXT [69, 70, 148]).

**Connection with current algorithms**

We are now in the position to show that the algorithms in [170, 190, 257, 259] are all special cases of SONATA and the earlier proposal NEXT [69, 70, 148].

**Aug-DGM [259] and Algorithm in [190].** Introduced in [259] for *undirected, time-invariant* graphs, the Aug-DGM algorithm reads

$$\begin{aligned}\mathbf{x}^{k+1} &= \widehat{\mathbf{W}}\left(\mathbf{x}^k - \text{Diag}\left(\boldsymbol{\gamma}\otimes\mathbf{1}_m\right)\mathbf{y}^k\right),\\ \mathbf{y}^{k+1} &= \widehat{\mathbf{W}}\left(\mathbf{y}^k + \mathbf{g}^{k+1} - \mathbf{g}^k\right),\end{aligned} \quad (229)$$

where $\widehat{\mathbf{W}} \triangleq \mathbf{W}\otimes\mathbf{I}_m$, $\mathbf{W}$ is a doubly-stochastic matrix compliant with the graph $G$ (cf. Assumption III.3), and $\boldsymbol{\gamma} \triangleq (\gamma_i)_{i=1}^I$ is the vector of agents' step-sizes.

A similar algorithm was proposed independently in [190] (in the same network setting of [259]), which reads

$$\begin{aligned}\mathbf{x}^{k+1} &= \widehat{\mathbf{W}}\left(\mathbf{x}^k - \gamma\mathbf{y}^k\right),\\ \mathbf{y}^{k+1} &= \widehat{\mathbf{W}}\mathbf{y}^k + \mathbf{g}^{k+1} - \mathbf{g}^k.\end{aligned} \quad (230)$$

Clearly Aug-DGM [259] in (229) and Algorithm [190] in (230) are both special cases of ATC-SONATA-NEXT-L [cf. (228)].

**(Push-)DIGing [170].** Appeared in [170] and applicable to *time-varying undirected graphs*, the DIGing Algorithm reads

$$\begin{aligned}\mathbf{x}^{k+1} &= \widehat{\mathbf{W}}^k\mathbf{x}^k - \gamma\mathbf{y}^k,\\ \mathbf{y}^{k+1} &= \widehat{\mathbf{W}}^k\mathbf{y}^k + \mathbf{g}^{k+1} - \mathbf{g}^k,\end{aligned} \quad (231)$$

where $\widehat{\mathbf{W}}^k \triangleq \mathbf{W}^k \otimes \mathbf{I}_m$ and $\mathbf{W}^k$ is a doubly-stochastic matrix compliant with the graph $G^k$. Clearly, DIGing coincides with CAA-SONATA-NEXT-L [cf. (228)], earlier proposed in [69, 70, 148].

The push-DIGing algorithm [170] extends DIGing to time-varying *digraphs*, and it is an instance of ATC-SONATA-L [cf. (227)], with $a_{ij}^k = 1/d_j^k$, $i,j = 1,\ldots I$.

**ADD-OPT [257].** Finally, we mention the ADD-OPT algorithm, proposed in [257] for *static digraphs*, which takes the following form:

$$\begin{aligned}
\mathbf{z}^{k+1} &= \widehat{\mathbf{A}}\mathbf{z}^k - \gamma\widetilde{\mathbf{y}}^k, \\
\boldsymbol{\phi}^{k+1} &= \mathbf{A}\boldsymbol{\phi}^k, \\
\mathbf{x}^{k+1} &= (\widehat{\boldsymbol{\Phi}}^{k+1})^{-1}\mathbf{z}^{k+1}, \\
\widetilde{\mathbf{y}}^{k+1} &= \widehat{\mathbf{A}}\widetilde{\mathbf{y}}^k + \mathbf{g}^{k+1} - \mathbf{g}^k.
\end{aligned} \quad (232)$$

Introducing the transformation $\mathbf{y}^k = (\widehat{\boldsymbol{\Phi}}^k)^{-1}\widetilde{\mathbf{y}}^k$, it is not difficult to check that (232) can be rewritten as

$$\begin{aligned}
\boldsymbol{\phi}^{k+1} &= \mathbf{A}\boldsymbol{\phi}^k, \\
\mathbf{W}^k &= (\boldsymbol{\Phi}^{k+1})^{-1}\mathbf{A}\boldsymbol{\Phi}^k, \\
\mathbf{x}^{k+1} &= \widehat{\mathbf{W}}^k\mathbf{x}^k - \gamma(\widehat{\boldsymbol{\Phi}}^{k+1})^{-1}\widehat{\boldsymbol{\Phi}}^k\mathbf{y}^k, \\
\mathbf{y}^{k+1} &= \widehat{\mathbf{W}}^k\mathbf{y}^k + (\widehat{\boldsymbol{\Phi}}^{k+1})^{-1}\left(\mathbf{g}^{k+1} - \mathbf{g}^k\right),
\end{aligned} \quad (233)$$

where $\widehat{\mathbf{W}}^k \triangleq \mathbf{W}^k \otimes \mathbf{I}_m$. Comparing (227) with (233), one can readily see that ADD-OPT is an instance of CAA-SONATA-L.

We summarize the connections between the different versions of SONATA(-NEXT) and its special cases in Table III.2.

*Table III.2: Connection of SONATA [211, 231] with current algorithms employing gradient tracking.*

| Algorithms | Special cases of | Instance of Problem (167) | Graph topology |
|---|---|---|---|
| NEXT [69, 148] | SONATA | $F$ nonconvex<br>$G \neq 0$<br>$X \subseteq \mathbb{R}^m$ | time-varying<br>doubly-stochasticable digraph |
| Aug-DGM [190, 259] | ATC-SONATA-NEXT-L ($\gamma = \gamma\mathbf{1}_I$) (228) | $F$ convex<br>$G = 0$<br>$X = \mathbb{R}^m$ | static undirected graph |
| DIGing [170] | CAA-SONATA-NEXT-L (228) | $F$ convex<br>$G = 0$<br>$X = \mathbb{R}^m$ | time-varying<br>doubly-stochasticable digraph |
| push-DIGing [170] | ATC-SONATA-L (227) | $F$ convex<br>$G = 0$<br>$X = \mathbb{R}^m$ | time-varying<br>digraph |
| ADD-OPT [257] | CAA-SONATA-L (227) | $F$ convex<br>$G = 0$<br>$X = \mathbb{R}^m$ | static digraph |

### III.3.2 Proof of Theorem III.16

The proof of Theorem III.16 is quite involved and can be found in [211]. Here we provide a simplified version, under the extra Assumption III.18 on Problem (167) (stated below) and the use of a square-summable (and thus diminishing) step-size in the algorithm.

**Assumption III.18.** *Given Problem* (167), *in addition to Assumption III.1, suppose that*

1. *The gradient of $F$ is bounded on $X$, i.e., there exists a constant $0 < L_F < +\infty$ such that $\|\nabla F(\mathbf{x})\| \leq L_F$, $\forall \mathbf{x} \in X$;*
2. *The subgradient of $G$ is bounded on $X$, i.e., there exists a constant $0 < L_G < +\infty$ such that $\|\partial G(\mathbf{x})\| \leq L_G$, $\forall \mathbf{x} \in X$.*

**Assumption III.19.** *The step size $\gamma^k \in (0, 1]$ satisfies the diminishing rule: $\sum_{k=0}^{\infty} \gamma^k = +\infty$ and $\sum_{k=0}^{\infty} (\gamma^k)^2 < +\infty$.*

Next, we prove separately

$$\lim_{k \to \infty} D(\mathbf{x}^k) = 0, \tag{234}$$

and

$$\lim_{k \to \infty} J(\mathbf{x}^k) = 0, \tag{235}$$

which imply (223).

### Technical preliminaries and sketch of the proof

We introduce here some preliminary definitions and results along with a sketch of the proof of the theorem.

**Weighted averages $\bar{\mathbf{x}}_\phi^k$ and $\bar{\mathbf{y}}_\phi^k$:** Define the weighted averages for the local copies $\mathbf{x}_{(i)}$ and the tracking variables $\mathbf{y}_{(i)}$:

$$\bar{\mathbf{x}}_\phi^k \triangleq \frac{1}{I} \sum_{i=1}^{I} \phi_{(i)}^k \mathbf{x}_{(i)}^k \quad \text{and} \quad \bar{\mathbf{y}}_\phi^k \triangleq \frac{1}{I} \sum_{i=1}^{I} \phi_{(i)}^k \mathbf{y}_{(i)}^k. \tag{236}$$

Using (218), (220) and (221), the dynamics of $\{\bar{\mathbf{x}}_\phi^k\}_{k \in \mathbb{N}_+}$ and $\{\bar{\mathbf{y}}_\phi^k\}_{k \in \mathbb{N}_+}$ generated by Algorithm 10 read: for all $k \in \mathbb{N}_+$,

$$\bar{\mathbf{x}}_\phi^{k+1} = \bar{\mathbf{x}}_\phi^k + \frac{\gamma^k}{I} \sum_{i=1}^{I} \phi_{(i)}^k \left( \widetilde{\mathbf{x}}_i^k - \bar{\mathbf{x}}_\phi^k \right), \tag{237}$$

and

$$\bar{\mathbf{y}}_\phi^{k+1} = \bar{\mathbf{y}}_\phi^k + \frac{1}{I} \sum_{i=1}^{I} \left( \mathbf{u}_i^{k+1} - \mathbf{u}_i^k \right), \quad \text{with} \quad \mathbf{u}_i^k \triangleq \nabla f_i(\mathbf{x}_{(i)}^k), \tag{238}$$

respectively. Let $\mathbf{u}^k \triangleq (\mathbf{u}_i^k)_{i=1}^{I}$. Note that, since each $\mathbf{y}_{(i)}^0 = \mathbf{u}_i^0 = \nabla f_i(\mathbf{x}_{(i)}^0)$ and $\phi_i^0 = 1$, we have [cf. Theorem III.10(a)]: for all $k \in \mathbb{N}_+$,

$$\bar{\mathbf{y}}_\phi^k = \frac{1}{I} \sum_{i=1}^{I} \mathbf{u}_i^k. \qquad (239)$$

The average quantities $\bar{\mathbf{x}}_\phi^k$ and $\bar{\mathbf{y}}_\phi^k$ will play a key role in proving asymptotic consensus and tracking. In fact, by (239), tracking is asymptotically achieved if $\lim_{k\to\infty} \|\mathbf{y}_{(i)}^k - \bar{\mathbf{y}}_\phi^k\| = 0$, for all $i = 1, \ldots, I$; and by

$$\left\| \mathbf{x}_{(i)}^k - \bar{\mathbf{x}}^k \right\| \leq \left\| \mathbf{x}_{(i)}^k - \bar{\mathbf{x}}_\phi^k \right\| + \left\| \frac{1}{I} \sum_{j=1}^{I} \left( \mathbf{x}_{(j)}^k - \bar{\mathbf{x}}_\phi^k \right) \right\| \leq B_1 \sum_{j=1}^{I} \left\| \mathbf{x}_{(j)}^k - \bar{\mathbf{x}}_\phi^k \right\|, \quad i = 1, \ldots I, \qquad (240)$$

with $B_1$ being a finite positive constant, it follows that consensus is asymptotically reached if $\lim_{k\to\infty} \|\mathbf{x}_{(i)}^k - \bar{\mathbf{x}}_\phi^k\| = 0$, for all $i = 1, \ldots, I$. These facts will will be proved in Step 1 of the proof, implying (234).

**Properties of the best-response $\widetilde{\mathbf{x}}_i^k$ and associated quantities:** We study here the connection between agents' best-responses $\widetilde{\mathbf{x}}_i^k$, defined in (217), and the "ideal" best-response $\widehat{\mathbf{x}}_i(\mathbf{x}_{(i)}^k)$ defined in (212) (not computable locally by the agents), with $\widetilde{F}_i$ given by (214), which we rewrite here for convenience: given $\mathbf{z} \in X$,

$$\widehat{\mathbf{x}}_i(\mathbf{z}) \triangleq \operatorname*{argmin}_{\mathbf{x}_{(i)} \in X} \widetilde{f}_i(\mathbf{x}_{(i)} \mid \mathbf{z}) + \left( \nabla F(\mathbf{z}) - \nabla f_i(\mathbf{z}) \right)^T (\mathbf{x}_{(i)} - \mathbf{z}) + G(\mathbf{x}_{(i)}). \qquad (241)$$

As observed in Sec. III.3, $\widetilde{\mathbf{x}}_i^k$ can be interpreted as a locally computable proxy of $\widehat{\mathbf{x}}_i(\mathbf{x}_{(i)}^k)$. We establish next the following connection among $\widetilde{\mathbf{x}}_i^k$, $\widehat{\mathbf{x}}_i(\bullet)$, and the stationary solutions of Problem (167): i) every fixed point of $\widehat{\mathbf{x}}_i(\bullet)$ is a stationary solution of Problem (167) (cf. Lemma III.20); and ii) the distance between these two mappings, $\|\widetilde{\mathbf{x}}_i^k - \widehat{\mathbf{x}}_i(\mathbf{x}_{(i)}^k)\|$, asymptotically vanishes, if consensus and tracking are achieved (cf. Lemma III.21). This also establishes the desired link between the limit points of $\mathbf{x}_i^k$ and the fixed points of $\widehat{\mathbf{x}}_i(\bullet)$ [and thus the stationary solutions of (167)].

**Lemma III.20.** *In the setting of Theorem III.16, the best-response map $X \ni \mathbf{z} \mapsto \widehat{\mathbf{x}}_i(\mathbf{z})$, defined in (241), with $i = 1, \ldots I$, enjoys the following properties:*

*(a) $\widehat{\mathbf{x}}_i(\bullet)$ is $\hat{L}_i$-Lipschitz continuous on $X$;*

*(b) The set of fixed points of $\widehat{\mathbf{x}}_i(\bullet)$ coincides with the set of stationary solutions of Problem (167). Therefore, $\widehat{\mathbf{x}}_i(\bullet)$ has a fixed point.*

*Proof.* The proof follows the same steps of that of Lemma II.4 and Lemma II.5 (Lecture II), and thus is omitted. □

**Lemma III.21.** *Let $\{\mathbf{x}^k = (\mathbf{x}_{(i)}^k)_{i=1}^{I}\}_{k \in \mathbb{N}_+}$ and $\{\mathbf{y}^k = (\mathbf{y}_{(i)}^k)_{i=1}^{I}\}_{k \in \mathbb{N}_+}$, be the sequence generated by Algorithm 10. Given $\widehat{\mathbf{x}}_i(\bullet)$ and $\widetilde{\mathbf{x}}_i^k$ in the setting of Theorem III.16, with $i = 1, \ldots, I$, the following holds: for every $k \in \mathbb{N}_+$ and $i = 1, \ldots, I$,*

$$\left\| \widetilde{\mathbf{x}}_i^k - \widehat{\mathbf{x}}_i\big(\mathbf{x}_{(i)}^k\big) \right\| \leq B_2 \left( \left\| \mathbf{y}_{(i)}^k - \bar{\mathbf{y}}_\phi^k \right\| + \sum_{j=1}^{I} \left\| \mathbf{x}_{(j)}^k - \bar{\mathbf{x}}_\phi^k \right\| \right), \tag{242}$$

where $B_2$ is some positive, finite, constant.

*Proof.* For notational simplicity, let us define $\widehat{\mathbf{x}}_i^k \triangleq \widehat{\mathbf{x}}_i\big(\mathbf{x}_{(i)}^k\big)$. We will also use the following shorthand: $\pm a \triangleq +a - a$, with $a \in \mathbb{R}$. Invoking the first order optimality condition for $\widetilde{\mathbf{x}}_i^k$ and $\widehat{\mathbf{x}}_i^k$, we have

$$\left( \widehat{\mathbf{x}}_i^k - \widetilde{\mathbf{x}}_i^k \right)^T \left( \nabla \widetilde{f}_i(\widetilde{\mathbf{x}}_i^k \mid \mathbf{x}_{(i)}^k) + I \cdot \mathbf{y}_{(i)}^k - \nabla f_i(\mathbf{x}_{(i)}^k) \right) + G(\widehat{\mathbf{x}}_i^k) - G(\widetilde{\mathbf{x}}_i^k) \geq 0$$

and

$$\left( \widetilde{\mathbf{x}}_i^k - \widehat{\mathbf{x}}_i^k \right)^T \left( \nabla \widetilde{f}_i(\widehat{\mathbf{x}}_i^k \mid \mathbf{x}_{(i)}^k) + \nabla F\big(\mathbf{x}_{(i)}^k\big) - \nabla f_i\big(\mathbf{x}_{(i)}^k\big) \right) + G(\widetilde{\mathbf{x}}_i^k) - G(\widehat{\mathbf{x}}_i) \geq 0,$$

respectively. Summing the two inequalities and using the strongly convexity of $\widetilde{f}_i(\bullet \mid \mathbf{x}_{(i)}^k)$ [cf. Assumption III.14.1], leads to the desired result

$$\tau_i \left\| \widetilde{\mathbf{x}}_i^k - \widehat{\mathbf{x}}_i^k \right\| \leq \left\| \nabla F\big(\mathbf{x}_{(i)}^k\big) - I \cdot \mathbf{y}_{(i)}^k \pm I \cdot \bar{\mathbf{y}}_\phi^k \right\|$$

$$\stackrel{(239), A.III.1.2}{\leq} I \cdot \left\| \mathbf{y}_{(i)}^k - \bar{\mathbf{y}}_\phi^k \right\| + \sum_{j=1}^{I} L_j \left\| \mathbf{x}_{(i)}^k - \mathbf{x}_{(j)}^k \right\|.$$

$\square$

**Structure of the proof.** The rest of the proof is organized in the following steps:

- **Step 1:** We first study the dynamics of the consensus and tracking errors, proving, among other results, that asymptotic consensus and tracking are achieved, that is, $\lim_{k\to\infty} \|\mathbf{x}_{(i)}^k - \bar{\mathbf{x}}_\phi^k\| = 0$ and $\lim_{k\to\infty} \|\mathbf{y}_{(i)}^k - \bar{\mathbf{y}}_\phi^k\| = 0$, for all $i = 1, \ldots, I$. By (240), this also proves (234);
- **Step 2:** We proceed studying the descent properties of $\{V(\bar{\mathbf{x}}_\phi^k)\}_{k \in \mathbb{N}_+}$, from which we will infer $\lim_{k\to\infty} \|\widetilde{\mathbf{x}}_i^k - \bar{\mathbf{x}}_\phi^k\| = 0$, for all $i = 1, \ldots, I$;
- **Step 3:** Finally, using the above results, we prove (235).

**Step 1: Asymptotic consensus and tracking**

We begin observing that the dynamics $\{\mathbf{x}^k \triangleq (\mathbf{x}_{(i)}^k)_{i=1}^I\}_{k \in \mathbb{N}_+}$ [cf. (220)] and $\{\mathbf{y}^k \triangleq (\mathbf{y}_{(i)}^k)_{i=1}^I\}_{k \in \mathbb{N}_+}$ [cf. (221)] generated by Algorithm 10 are instances of the perturbed condensed push-sum protocol (188), with errors

$$\boldsymbol{\varepsilon}_i^k = \gamma^k (\phi_{(i)}^{k+1})^{-1} \sum_{j=1}^{I} a_{ij}^k \phi_{(j)}^k \left( \widetilde{\mathbf{x}}_j^k - \mathbf{x}_{(j)}^k \right) \tag{243}$$

and
$$\boldsymbol{\varepsilon}_i^k = (\phi_{(i)}^{k+1})^{-1}\left(\mathbf{u}_i^{k+1} - \mathbf{u}_i^k\right), \tag{244}$$

respectively. We can then leverage the convergence results introduced in Sec. III.2.4 to prove the desired asymptotic consensus and tracking. To do so, we first show that some related quantities are bounded.

**Lemma III.22.** *Let* $\{\mathbf{x}^k = (\mathbf{x}_{(i)}^k)_{i=1}^I\}_{k \in \mathbb{N}_+}$, $\{\mathbf{y}^k = (\mathbf{y}_{(i)}^k)_{i=1}^I\}_{k \in \mathbb{N}_+}$, *and* $\{\boldsymbol{\phi}^k = (\phi_{(i)}^k)_{i=1}^I\}_{k \in \mathbb{N}_+}$ *be the sequence generated by Algorithm 10, in the setting of Theorem III.16 and under the extra Assumptions III.18-III.19. Then, the following hold: for all* $i = 1, \ldots, I$, *(a)* $\{\boldsymbol{\phi}^k\}_{k \in \mathbb{N}_+}$ *is uniformly bounded:*

$$\phi_{lb} \cdot \mathbf{1} \le \boldsymbol{\phi}^k \le \phi_{ub} \cdot \mathbf{1}, \quad \forall k \in \mathbb{N}_+, \tag{245}$$

*where* $\phi_{lb}$ *and* $\phi_{ub}$ *are finite, positive constants, defined in (189);*
*(b)*
$$\sup_{k \in \mathbb{N}_+} \left\| \mathbf{y}_{(i)}^k - \bar{\mathbf{y}}_\phi^k \right\| < \infty; \tag{246}$$

*(c)*
$$\sup_{k \in \mathbb{N}_+} \left\| \mathbf{x}_{(i)}^k - \widetilde{\mathbf{x}}_i^k \right\| < \infty. \tag{247}$$

*Proof.* Statement (a) is a consequence of Theorem III.11(a). Statement (b) follows readily from (192), observing that the errors $\boldsymbol{\varepsilon}_i^k = (\phi_{(i)}^{k+1})^{-1}(\mathbf{u}_i^{k+1} - \mathbf{u}_i^k)$ are all uniformly bounded, due to Assumption III.18.1 and (245).

We prove now statement (c). Since $\widetilde{\mathbf{x}}_{(i)}^k$ is the unique optimal solution of Problem (217), invoking the first order optimality condition, we have

$$\left(\mathbf{x}_{(i)}^k - \widetilde{\mathbf{x}}_i^k\right)^T \left(\nabla \widetilde{f}_i(\widetilde{\mathbf{x}}_i^k \mid \mathbf{x}_{(i)}^k) + I \cdot \mathbf{y}_{(i)}^k - \nabla f_i(\mathbf{x}_{(i)}^k) + \boldsymbol{\xi}^k\right) \ge 0,$$

where $\boldsymbol{\xi}^k \in \partial G(\widetilde{\mathbf{x}}_i^k)$, and $\|\boldsymbol{\xi}^k\| \le L_G$ (cf. Assumption III.18.2). Since $\widetilde{f}_i(\bullet \mid \mathbf{x}_i^k)$ is $\tau_i$-strongly convex (cf. Assumption III.14), we have

$$\begin{aligned}
\left\|\mathbf{x}_{(i)}^k - \widetilde{\mathbf{x}}_i^k\right\| &\le \frac{I}{\tau_i} \left\|\mathbf{y}_{(i)}^k\right\| + \frac{L_G}{\tau_i} \\
&\le \frac{I}{\tau_i} \left\|\mathbf{y}_{(i)}^k - \bar{\mathbf{y}}_\phi^k\right\| + \left\|\bar{\mathbf{y}}_\phi^k\right\| + \frac{L_G}{\tau_i} \\
&\stackrel{(239),(246)}{\le} B_3 + \frac{1}{I}\sum_{i=1}^I \left\|\mathbf{u}_i^k\right\| \le B_4,
\end{aligned} \tag{248}$$

for all $k \in \mathbb{N}_+$, where the last inequality follows from Assumption III.18.1, and $B_3$ and $B_4$ are some positive, finite, constants. □

We can now study the dynamics of the consensus error.

**Proposition III.23.** *Let $\{\mathbf{x}^k = (\mathbf{x}_{(i)}^k)_{i=1}^I\}_{k\in\mathbb{N}_+}$ be the sequence generated by Algorithm 10, in the setting of Theorem III.16, and under the extra Assumptions III.18-III.19; and let $\{\bar{\mathbf{x}}_\phi^k\}_{k\in\mathbb{N}_+}$, with $\bar{\mathbf{x}}_\phi^k$ defined in (237). Then, all the $\mathbf{x}_{(i)}^k$ are asymptotically consensual, that is,*

$$\lim_{k\to\infty} \left\|\mathbf{x}_{(i)}^k - \bar{\mathbf{x}}_\phi^k\right\| = 0, \quad \forall i = 1,\ldots,I. \tag{249}$$

*Furthermore, there hold: for all $i = 1,\ldots,I$,*

$$\sum_{k=0}^\infty \gamma^k \left\|\mathbf{x}_{(i)}^k - \bar{\mathbf{x}}_\phi^k\right\| < \infty, \tag{250}$$

$$\sum_{k=0}^\infty \left\|\mathbf{x}_{(i)}^k - \bar{\mathbf{x}}_\phi^k\right\|^2 < \infty. \tag{251}$$

*Proof.* We use again the connection between (220) and the perturbed condensed push-sum protocol (188), with error $\boldsymbol{\varepsilon}_i^k$ given by (243). Invoking Theorem III.11 and Lemma III.12(a) [cf. (191)], to prove (249), it is then sufficient to show that all the errors $\boldsymbol{\varepsilon}_i^k$ are asymptotically vanishing. This follows readily from the following facts: $\gamma^k \downarrow 0$ (Assumption III.19); $0 < \phi_{lb} \leq \phi_{(i)}^k \leq \phi_{ub} < \infty$, for all $k \in \mathbb{N}_+$ and $i = 1,\ldots,I$ [cf. Theorem III.11(a)]; and $\sup_{k\in\mathbb{N}_+} \|\mathbf{x}_{(i)}^k - \widetilde{\mathbf{x}}_i^k\| < \infty$ [Lemma III.22(c)].

We prove now (250). Invoking Theorem III.11(b), we can write

$$\lim_{k\to\infty}\sum_{t=1}^k \gamma^t \left\|\mathbf{x}_{(i)}^t - \bar{\mathbf{x}}_\phi^t\right\| \stackrel{(190)}{\leq} c \cdot \lim_{k\to\infty}\sum_{t=1}^k \gamma^t \left((\rho)^t \|\mathbf{x}^0\| + \sum_{s=0}^{t-1}(\rho)^{t-1-s}\|\boldsymbol{\varepsilon}^s\|\right)$$

$$\stackrel{(a)}{\leq} B_5 \cdot \lim_{k\to\infty}\sum_{t=1}^k (\rho)^t + B_6 \cdot \underbrace{\lim_{k\to\infty}\sum_{t=1}^k\sum_{s=0}^{t-1}(\rho)^{t-1-s}\gamma^s\gamma^t}_{<\infty \quad [\text{Lemma } III.12(b2)]} < \infty,$$

where $B_5$ and $B_6$ are some positive, finite, constants, and in (a) we used $\|\boldsymbol{\varepsilon}^t\| \leq B_7\gamma^t$, for some positive, finite $B_7$ [which follows from the same arguments used to prove (249)].

To prove (251), we can use similar steps and write:

$$\lim_{k\to\infty}\sum_{t=1}^k \left\|\mathbf{x}_{(i)}^t - \bar{\mathbf{x}}_\phi^t\right\|^2 \leq \lim_{k\to\infty}\sum_{t=1}^k \left(B_5(\rho)^t + B_6 \cdot \sum_{s=0}^{t-1}(\rho)^{t-1-s}\gamma^s\right)^2$$

$$\leq 2\cdot B_5^2 \lim_{k\to\infty}\sum_{t=0}^k (\rho)^{2\cdot t} + 2B_6^2 \cdot \lim_{k\to\infty}\sum_{t=0}^k\sum_{s=0}^{t-1}\sum_{s'=0}^{t-1}(\rho)^{t-1-s}(\rho)^{t-1-s'}\gamma^s\gamma^{s'} \stackrel{(a)}{<} \infty,$$

where in (a) we used

$$\lim_{k\to\infty}\sum_{t=0}^k\sum_{s=0}^{t-1}\sum_{s'=0}^{t-1}(\rho)^{t-1-s}(\rho)^{t-1-s'}\gamma^s\gamma^{s'}$$

$$\leq \lim_{k\to\infty} \sum_{t=0}^{k} \sum_{s=0}^{t-1} \sum_{s'=0}^{t-1} (\rho)^{t-1-s} (\rho)^{t-1-s'} (\gamma^s)^2$$

$$= \lim_{k\to\infty} \sum_{t=0}^{k} \underbrace{\left(\sum_{s'=0}^{t-1}(\rho)^{t-1-s'}\right)}_{\leq \frac{1}{1-\rho}} \sum_{s=0}^{t-1}(\rho)^{t-1-s}(\gamma^s)^2 < \infty,$$

where the last implication is a consequence of Lemma III.12(b1). □

We prove next a similar result for the tracking error.

**Proposition III.24.** *Let $\{\mathbf{y}^k = (\mathbf{y}_{(i)}^k)_{i=1}^I\}_{k\in\mathbb{N}_+}$ be the sequence generated by Algorithm 10, in the setting of Theorem III.16, and under the extra Assumptions III.18-III.19; and let $\{\bar{\mathbf{y}}_\phi^k\}_{k\in\mathbb{N}_+}$, with $\bar{\mathbf{y}}_\phi^k$ defined in (238). Then, all the $\mathbf{y}_{(i)}^k$ asymptotically track the gradient average $(1/I)\cdot\sum_{i=1}^I \nabla f_i(\mathbf{x}_{(i)}^k)$, that is,*

$$\lim_{k\to\infty}\left\|\mathbf{y}_{(i)}^k - \bar{\mathbf{y}}_\phi^k\right\| = 0, \quad \forall i = 1,\ldots,I. \tag{252}$$

*Furthermore, there holds*

$$\sum_{k=0}^{\infty} \gamma^k \left\|\mathbf{y}_{(i)}^k - \bar{\mathbf{y}}_\phi^k\right\| < \infty, \quad \forall i = 1,\ldots,I. \tag{253}$$

*Proof.* Invoking Theorem III.10 and using (239), to prove (252), it is sufficient to show that Assumption III.9 holds. By Assumption III.1.2, this reduces to prove $\lim_{k\to\infty}\|\mathbf{x}_{(i)}^{k+1} - \mathbf{x}_{(i)}^k\| = 0$, for all $i=1,\ldots,I$. We write:

$$\lim_{k\to\infty}\left\|\mathbf{x}_{(i)}^{k+1} - \mathbf{x}_{(i)}^k\right\| \leq \underbrace{\lim_{k\to\infty}\left\|\mathbf{x}_{(i)}^{k+1} - \bar{\mathbf{x}}_\phi^{k+1}\right\|}_{\stackrel{(249)}{=}0} + \underbrace{\lim_{k\to\infty}\left\|\mathbf{x}_{(i)}^k - \bar{\mathbf{x}}_\phi^k\right\|}_{\stackrel{(249)}{=}0}$$
$$+ \underbrace{\lim_{k\to\infty}\left\|\frac{\gamma^k}{I}\sum_{i=1}^I \phi_{(i)}^k\left(\widetilde{\mathbf{x}}_i^k - \bar{\mathbf{x}}_\phi^k\right)\right\|}_{\stackrel{(247)}{\leq}\gamma^k\cdot B_8} = 0, \tag{254}$$

where $B_8$ is some positive, finite, constant.

We prove next (253). Invoking Theorem III.11(b), with $\boldsymbol{\varepsilon}_i^t$ defined in (244), we have

$$\lim_{k\to\infty}\sum_{t=1}^k \gamma^t \left\|\mathbf{y}_{(i)}^t - \bar{\mathbf{y}}_\phi^t\right\|$$
$$\stackrel{(190)}{\leq} c\cdot\lim_{k\to\infty}\sum_{t=1}^k \gamma^t \left((\rho)^t\|\mathbf{u}^0\| + \sum_{s=0}^{t-1}(\rho)^{t-1-s}\|\boldsymbol{\varepsilon}^s\|\right)$$

$$\stackrel{(a)}{\leq} B_9 \cdot \lim_{k\to\infty} \sum_{t=1}^{k} (\rho)^t + B_{10} \cdot \lim_{k\to\infty} \sum_{i=1}^{I} \sum_{t=1}^{k} \sum_{s=0}^{t-1} (\rho)^{t-1-s} \gamma^t \left\| \mathbf{x}_{(i)}^{s+1} - \mathbf{x}_{(i)}^{s} \right\|$$

$$\stackrel{(b)}{\leq} B_9 \cdot \lim_{k\to\infty} \sum_{t=1}^{k} (\rho)^t + B_{10} \cdot B_8 \cdot I \cdot \underbrace{\lim_{k\to\infty} \sum_{t=1}^{k} \sum_{s=0}^{t-1} (\rho)^{t-1-s} \gamma^t \gamma^s}_{<\infty \quad [\text{Lemma } III.12(b2)]}$$

$$+ B_{10} \cdot \underbrace{\lim_{k\to\infty} \sum_{i=1}^{I} \sum_{t=1}^{k} \sum_{s=0}^{t-1} (\rho)^{t-1-s} \gamma^t \left( \left\| \mathbf{x}_{(i)}^{s+1} - \bar{\mathbf{x}}_{\phi}^{s+1} \right\| + \left\| \mathbf{x}_{(i)}^{s} - \bar{\mathbf{x}}_{\phi}^{s} \right\| \right)}_{<\infty \quad [(251) \& \text{ Lemma } III.12(b2)]} < \infty,$$

where $B_9$ and $B_{10}$ are some positive, finite, constants; in (a) we used the definition of $\boldsymbol{\varepsilon}_i^t$ [cf. (244)] along with Assumption III.1.2; and in (b) we invoked the same upper bound used in (254).

$\square$

**Step 2: Descent of** $\{V(\bar{\mathbf{x}}_\phi^k)\}_{k\in\mathbb{N}_+}$

We study now the descent of $V$ along $\{\bar{\mathbf{x}}_\phi^k\}_{k\in\mathbb{N}_+}$. For notational simplicity, define

$$\Delta \mathbf{x}_{(i)}^k \triangleq \mathbf{x}_{(i)}^k - \bar{\mathbf{x}}_\phi^k \quad \text{and} \quad \Delta \widetilde{\mathbf{x}}_{(i)}^k \triangleq \widetilde{\mathbf{x}}_i^k - \bar{\mathbf{x}}_\phi^k, \quad i=1,\ldots,I, \tag{255}$$

where we recall that $\widetilde{\mathbf{x}}_i^k$ is defined in (217). Invoking the first order optimality condition of $\widetilde{\mathbf{x}}_i^k$, it is not difficult to check that the following property holds for $\Delta \widetilde{\mathbf{x}}_{(i)}^k$:

$$\left( \Delta \widetilde{\mathbf{x}}_{(i)}^k \right)^T \left( \nabla \widetilde{f}_i(\bar{\mathbf{x}}_\phi^k \mid \mathbf{x}_{(i)}^k) - \nabla f_i(\mathbf{x}_{(i)}^k) + I \cdot \mathbf{y}_{(i)}^k \right) + G(\widetilde{\mathbf{x}}_i^k) - G(\bar{\mathbf{x}}_\phi^k) \leq -\tau_i \left\| \Delta \widetilde{\mathbf{x}}_{(i)}^k \right\|^2. \tag{256}$$

We begin studying the dynamics of $F$ along $\{\bar{\mathbf{x}}_\phi^k\}_{k\in\mathbb{N}_+}$. Using (237) and invoking the descent lemma [cf. Lemma II.14], with $L \triangleq \sum_{i=1}^{I} L_i$, we have:

$$F(\bar{\mathbf{x}}_\phi^{k+1})$$
$$\leq F(\bar{\mathbf{x}}_\phi^k) + \frac{\gamma^k}{I} \sum_{i=1}^{I} \phi_{(i)}^k \nabla F(\bar{\mathbf{x}}_\phi^k)^T \Delta \widetilde{\mathbf{x}}_{(i)}^k + \frac{L}{2} \frac{(\gamma^k)^2}{I} \sum_{i=1}^{I} (\phi_{(i)}^k)^2 \left\| \Delta \widetilde{\mathbf{x}}_{(i)}^k \right\|^2$$
$$\stackrel{(256)}{\leq} F(\bar{\mathbf{x}}_\phi^k) - \frac{\gamma^k}{I} \sum_{i=1}^{I} \phi_{(i)}^k \left( \tau_i \left\| \Delta \widetilde{\mathbf{x}}_{(i)}^k \right\|^2 + G(\widetilde{\mathbf{x}}_{(i)}^k) - G(\bar{\mathbf{x}}_\phi^k) \right)$$
$$+ \frac{\gamma^k}{I} \sum_{i=1}^{I} \phi_{(i)}^k \left( \nabla f_i(\mathbf{x}_{(i)}^k) - \nabla \widetilde{f}_i(\bar{\mathbf{x}}_\phi^k \mid \mathbf{x}_{(i)}^k) \pm \nabla \widetilde{f}_i(\bar{\mathbf{x}}_\phi^k \mid \bar{\mathbf{x}}_\phi^k) \right)^T \Delta \widetilde{\mathbf{x}}_{(i)}^k$$
$$+ \frac{\gamma^k}{I} \sum_{i=1}^{I} \phi_{(i)}^k \left( \nabla F(\bar{\mathbf{x}}_\phi^k) - I \cdot \mathbf{y}_{(i)}^k \pm (I \cdot \bar{\mathbf{y}}_\phi^k) \right)^T \Delta \widetilde{\mathbf{x}}_{(i)}^k + \frac{L}{2} \frac{(\gamma^k)^2}{I} \sum_{i=1}^{I} (\phi_{(i)}^k)^2 \left\| \Delta \widetilde{\mathbf{x}}_{(i)}^k \right\|^2$$

Title Suppressed Due to Excessive Length 131$$\stackrel{(a)}{\leq} F(\bar{\mathbf{x}}_\phi^k) - \frac{\gamma^k}{I}\sum_{i=1}^{I} \phi_{(i)}^k \left( \tau_i \left\| \Delta \widetilde{\mathbf{x}}_{(i)}^k \right\|^2 + G(\widetilde{\mathbf{x}}_{(i)}^k) - G(\bar{\mathbf{x}}_\phi^k) \right)$$

$$+ \frac{\gamma^k}{I}\sum_{i=1}^{I} \phi_{(i)}^k (L_i + \tilde{L}_i) \left\| \Delta \mathbf{x}_{(i)}^k \right\| \left\| \Delta \widetilde{\mathbf{x}}_{(i)}^k \right\|$$

$$+ \frac{\gamma^k}{I}\sum_{i=1}^{I} \phi_{(i)}^k \left( \sum_{j=1}^{I} L_j \left\| \Delta \mathbf{x}_{(j)}^k \right\| \right) \left\| \Delta \widetilde{\mathbf{x}}_{(i)}^k \right\| + \gamma^k \sum_{i=1}^{I} \phi_{(i)}^k \left\| \mathbf{y}_{(i)}^k - \bar{\mathbf{y}}_\phi^k \right\| \left\| \Delta \widetilde{\mathbf{x}}_{(i)}^k \right\|$$

$$+ \frac{L}{2}\frac{(\gamma^k)^2}{I}\sum_{i=1}^{I} (\phi_{(i)}^k)^2 \left\| \Delta \widetilde{\mathbf{x}}_{(i)}^k \right\|^2, \tag{257}$$

where in (a) we used (239), Assumption III.14.2, and Assumption III.1.2.

Invoking the convexity of $G$ and (237), we can write

$$G(\bar{\mathbf{x}}_\phi^{k+1}) \leq G(\bar{\mathbf{x}}_\phi^k) - \frac{\gamma^k}{I}\sum_{i=1}^{I} \phi_{(i)}^k \left( G(\bar{\mathbf{x}}_\phi^k) - G(\widetilde{\mathbf{x}}_i^k) \right). \tag{258}$$

Substituting (258) in (257) and using (245) [cf. Lemma III.22], $\tau_{\min} \triangleq \min_i \tau_i > 0$, and (249), we can write: for sufficiently large $k \in \mathbb{N}_+$,

$$V(\bar{\mathbf{x}}_\phi^{k+1}) \leq V(\bar{\mathbf{x}}_\phi^k) - \gamma^k \left( \frac{\tau_{\min} \phi_{lb}}{I} - \frac{L\phi_{ub}^2 \gamma^k}{2I} \right) \sum_{i=1}^{I} \left\| \Delta \widetilde{\mathbf{x}}_{(i)}^k \right\|^2$$

$$+ B_{11} \cdot \gamma^k \sum_{i=1}^{I} \left( \left\| \Delta \mathbf{x}_{(i)}^k \right\| + \sum_{j=1}^{I} \left\| \Delta \mathbf{x}_{(j)}^k \right\| + \left\| \mathbf{y}_{(i)}^k - \bar{\mathbf{y}}_\phi^k \right\| \right) \left\| \Delta \widetilde{\mathbf{x}}_{(i)}^k \right\|$$

$$\leq V(\bar{\mathbf{x}}_\phi^k) - B_{12}\, \gamma^k \sum_{i=1}^{I} \left\| \Delta \widetilde{\mathbf{x}}_{(i)}^k \right\|^2 + \underbrace{B_{13} \cdot \gamma^k \sum_{i=1}^{I} \left( \left\| \mathbf{y}_{(i)}^k - \bar{\mathbf{y}}_\phi^k \right\| + \left\| \Delta \mathbf{x}_{(i)}^k \right\| \right)}_{\triangleq T^k} \tag{259}$$

where $B_{11}, B_{12}, B_{13}$ are some finite, positive, constants, and in the last inequality we used the fact that $\left\| \Delta \widetilde{\mathbf{x}}_{(i)}^k \right\| \leq B_{14}$, with $i = 1, \ldots, I$, for sufficiently large $k \in \mathbb{N}_+$ and some finite, positive, constant $B_{14}$, which is a direct consequence of (247) (cf. Lemma III.22) and (249) (cf. Proposition III.23).

Note that, $\sum_{k=0}^{\infty} T^k < \infty$, due to Proposition III.23 [cf. (250)] and Proposition III.24 [cf (253)]. We can then apply Lemma II.16 to (259), with the following identifications:

$$Y^k = V(\bar{\mathbf{x}}_\phi^k), \quad Z^k = T^k, \quad \text{and} \quad W^k = B_{12}\, \gamma^k \sum_{i=1}^{I} \left\| \Delta \widetilde{\mathbf{x}}_{(i)}^k \right\|^2,$$

which, using Assumption II.1.4, yields i) $\lim_{k\to\infty} V(\bar{\mathbf{x}}_\phi^k) = V^\infty$, with $V^\infty$ finite; and ii)

$$\sum_{k=0}^{\infty} \gamma^k \sum_{i=1}^{I} \left\| \Delta \widetilde{\mathbf{x}}_{(i)}^k \right\|^2 < \infty. \tag{260}$$

Since $\sum_{k=0}^{\infty} \gamma^k = \infty$, (260) implies $\liminf_{k \to \infty} \sum_{i=1}^{I} \left\| \Delta \widetilde{\mathbf{x}}_{(i)}^k \right\| = 0$.

We prove next that $\limsup_{k \to \infty} \sum_{i=1}^{I} \left\| \Delta \widetilde{\mathbf{x}}_{(i)}^k \right\| = 0$, for all $i = 1, \ldots, I$. Assume on the contrary that $\limsup_{k \to \infty} \sum_{i=1}^{I} \left\| \Delta \widetilde{\mathbf{x}}_{(i)}^k \right\| > 0$. Since $\liminf_{k \to \infty} \sum_{i=1}^{I} \left\| \Delta \widetilde{\mathbf{x}}_{(i)}^k \right\| = 0$, there exists an infinite set $K \subseteq \mathbb{N}_+$, such that for all $k \in K$, one can find an integer $t_k > k$ such that

$$\sum_{i=1}^{I} \| \Delta \widetilde{\mathbf{x}}_{(i)}^k \| < \delta, \qquad \sum_{i=1}^{I} \| \Delta \widetilde{\mathbf{x}}_{(i)}^{t_k} \| > 2\delta, \tag{261}$$

$$\delta \le \sum_{i=1}^{I} \| \Delta \widetilde{\mathbf{x}}_{(i)}^j \| \le 2\delta, \qquad \text{if} \quad k < j < t_k. \tag{262}$$

Therefore, for all $k \in K$,

$$\delta < \sum_{i=1}^{I} \left\| \Delta \widetilde{\mathbf{x}}_{(i)}^{t_k} \right\| - \sum_{i=1}^{I} \left\| \Delta \widetilde{\mathbf{x}}_{(i)}^k \right\| \le \sum_{i=1}^{I} \left\| \Delta \widetilde{\mathbf{x}}_{(i)}^{t_k} - \Delta \widetilde{\mathbf{x}}_{(i)}^k \pm \widehat{\mathbf{x}}_i(\bar{\mathbf{x}}_\phi^k) \pm \widehat{\mathbf{x}}_i(\bar{\mathbf{x}}_\phi^{t_k}) \right\|$$

$$\overset{(a)}{\le} \sum_{i=1}^{I} (1+\hat{L}_i) \left\| \bar{\mathbf{x}}_\phi^{t_k} - \bar{\mathbf{x}}_\phi^k \right\| + \underbrace{\sum_{i=1}^{I} \left\| \widetilde{\mathbf{x}}_{(i)}^{t_k} - \widehat{\mathbf{x}}_i(\bar{\mathbf{x}}_\phi^{t_k}) \right\| + \sum_{i=1}^{I} \left\| \widetilde{\mathbf{x}}_{(i)}^k - \widehat{\mathbf{x}}_i(\bar{\mathbf{x}}_\phi^k) \right\|}_{\triangleq e_1^k}$$

$$\overset{(245)}{\le} \sum_{i=1}^{I} (1+\hat{L}_i) \frac{\phi_{ub}}{I} \sum_{j=k}^{t_k-1} \gamma^j \sum_{\ell=1}^{I} \left\| \Delta \widetilde{\mathbf{x}}_{(\ell)}^j \right\| + e_1^k$$

$$\le \sum_{i=1}^{I} (1+\hat{L}_i) \frac{\phi_{ub}}{I} \sum_{j=k+1}^{t_k-1} \gamma^\tau \sum_{\ell=1}^{I} \left\| \Delta \widetilde{\mathbf{x}}_{(\ell)}^j \right\| + e_1^k + \underbrace{\sum_{i=1}^{I} (1+\hat{L}_i) \frac{\phi_{ub}}{I} \gamma^k \sum_{\ell=1}^{I} \left\| \Delta \widetilde{\mathbf{x}}_{(\ell)}^k \right\|}_{\triangleq e_2^k}$$

$$\overset{(b)}{\le} \sum_{i=1}^{I} (1+\hat{L}_i) \frac{\phi_{ub}}{\delta} \sum_{j=k+1}^{t_k-1} \gamma^j \sum_{\ell=1}^{I} \left\| \Delta \widetilde{\mathbf{x}}_{(\ell)}^j \right\|^2 + e_1^k + e_2^k,$$
$$\tag{263}$$

where in (a) we used Lemma III.20(a); and in (b) we used (262) and Lemma II.15.

Note that

(i) $\lim_{k \to \infty} e_1^k = 0$, due to Lemma III.21 [cf. (242)], Proposition III.23 [cf. (249)], and Proposition III.24 [cf. (252)];

(ii) $\lim_{k \to \infty} e_2^k = 0$, due to (260), which implies $\lim_{k \to \infty} \gamma^k \| \Delta \widetilde{\mathbf{x}}_{(i)}^k \|^2 = 0$, and thus $\lim_{k \to \infty} \gamma^k \| \Delta \widetilde{\mathbf{x}}_{(i)}^k \| = 0$ (recall that $\gamma^k \in (0,1]$); and

(iii)
$$\lim_{k \to \infty} \sum_{j=k+1}^{t_k-1} \gamma^j \sum_{i=1}^{I} \left\| \Delta \widetilde{\mathbf{x}}_{(i)}^j \right\|^2 = 0,$$

due to (260).

This however contradicts (263). Therefore, it must be

$$\lim_{k\to\infty} \left\| \Delta \widetilde{\mathbf{x}}_{(i)}^k \right\| = 0, \qquad i = 1, \ldots, I. \tag{264}$$

**Step 3:** $\lim_{k\to\infty} J(\bar{\mathbf{x}}^k) = 0$

Recalling the definition of $J(\bar{\mathbf{x}}^k) = \|\bar{\mathbf{x}}^k - \bar{\mathbf{x}}(\bar{\mathbf{x}}^k)\|$ [cf.(222)], where for notational simplicity we introduced

$$\bar{\mathbf{x}}(\bar{\mathbf{x}}^k) \triangleq \operatorname*{argmin}_{\mathbf{y} \in X} \left\{ \nabla F(\bar{\mathbf{x}}^k)^T (\mathbf{y} - \bar{\mathbf{x}}^k) + \frac{1}{2}\|\mathbf{y} - \bar{\mathbf{x}}^k\|_2^2 + G(\mathbf{y}) \right\}, \tag{265}$$

we start bounding $J$ as follows

$$J(\bar{\mathbf{x}}^k) = \left\| \bar{\mathbf{x}}^k - \bar{\mathbf{x}}(\bar{\mathbf{x}}^k) \pm \widehat{\mathbf{x}}_i(\bar{\mathbf{x}}^k) \right\| \leq \underbrace{\left\| \widehat{\mathbf{x}}_i(\bar{\mathbf{x}}^k) - \bar{\mathbf{x}}^k \right\|}_{\text{term I}} + \underbrace{\left\| \bar{\mathbf{x}}(\bar{\mathbf{x}}^k) - \widehat{\mathbf{x}}_i(\bar{\mathbf{x}}^k) \right\|}_{\text{term II}}, \tag{266}$$

where $\widehat{\mathbf{x}}_i(\bullet)$ is defined in (241). To prove $\lim_{k\to\infty} J(\bar{\mathbf{x}}^k) = 0$, it is then sufficient to show that both Term I and Term II in (266) are asymptotically vanishing. We study the two terms separately.

• `Term I`: We prove $\lim_{k\to\infty} \|\widehat{\mathbf{x}}_i(\bar{\mathbf{x}}^k) - \bar{\mathbf{x}}^k\| = 0$.

We begin bounding $\|\widehat{\mathbf{x}}_i(\bar{\mathbf{x}}^k) - \bar{\mathbf{x}}^k\|$ as

$$\begin{aligned}\left\|\widehat{\mathbf{x}}_i(\bar{\mathbf{x}}^k) - \bar{\mathbf{x}}^k\right\| &= \left\| \widehat{\mathbf{x}}_i(\bar{\mathbf{x}}^k) \pm \widehat{\mathbf{x}}_i(\bar{\mathbf{x}}_\phi^k) \pm \bar{\mathbf{x}}_\phi^k - \bar{\mathbf{x}}^k \right\| \\ &\leq \left\| \widehat{\mathbf{x}}_i(\bar{\mathbf{x}}_\phi^k) - \bar{\mathbf{x}}_\phi^k \right\| + (1+\hat{L}_i)\left\| \bar{\mathbf{x}}_\phi^k - \bar{\mathbf{x}}^k \right\|. \end{aligned} \tag{267}$$

By (240) and Proposition III.23 [cf. (249)] it follows that

$$\lim_{k\to\infty} \left\| \bar{\mathbf{x}}_\phi^k - \bar{\mathbf{x}}^k \right\| = 0. \tag{268}$$

Therefore, to prove $\lim_{k\to\infty} \|\widehat{\mathbf{x}}_i(\bar{\mathbf{x}}^k) - \bar{\mathbf{x}}^k\| = 0$, it is sufficient to show that $\|\widehat{\mathbf{x}}_i(\bar{\mathbf{x}}_\phi^k) - \bar{\mathbf{x}}_\phi^k\|$ is asymptotically vanishing, as proved next.

We bound $\|\widehat{\mathbf{x}}_i(\bar{\mathbf{x}}_\phi^k) - \bar{\mathbf{x}}_\phi^k\|$ as:

$$\begin{aligned}\left\| \widehat{\mathbf{x}}_i(\bar{\mathbf{x}}_\phi^k) - \bar{\mathbf{x}}_\phi^k \pm \widetilde{\mathbf{x}}_{(i)}^k \pm \widehat{\mathbf{x}}_i(\mathbf{x}_{(i)}^k) \right\| &\leq \left\| \Delta\widetilde{\mathbf{x}}_{(i)}^k \right\| + \left\| \widehat{\mathbf{x}}_i(\mathbf{x}_{(i)}^k) - \widetilde{\mathbf{x}}_{(i)}^k \right\| + \left\| \widehat{\mathbf{x}}_i(\mathbf{x}_{(i)}^k) - \widehat{\mathbf{x}}_i(\bar{\mathbf{x}}_\phi^k) \right\| \\ &\leq \left\| \Delta\widetilde{\mathbf{x}}_{(i)}^k \right\| + B_{15} \sum_{j=1}^I \left\| \mathbf{x}_{(j)}^k - \bar{\mathbf{x}}_\phi^k \right\| + B_{16} \left\| \mathbf{y}_{(i)}^k - \bar{\mathbf{y}}_\phi^k \right\|, \end{aligned} \tag{269}$$

where $B_{15}$ and $B_{16}$ are some positive, finite, constants; and in the last inequality we used Lemma III.20 and Lemma III.21 [cf. (242)]. Invoking

i) $\lim_{k\to\infty} \left\| \Delta \widetilde{\mathbf{x}}_{(i)}^k \right\| = 0$ [cf. (264)];
ii) $\lim_{k\to\infty} \left\| \mathbf{x}_{(i)}^k - \bar{\mathbf{x}}_\phi^k \right\| = 0$, for all $i = 1, \ldots, I$, [cf. (249)]; and
iii) $\lim_{k\to\infty} \left\| \mathbf{y}_{(i)}^k - \bar{\mathbf{y}}_\phi^k \right\| = 0$, for all $i = 1, \ldots, I$, [cf. (252)],

the desired result, $\lim_{k\to\infty} \left\| \widehat{\mathbf{x}}_i(\bar{\mathbf{x}}_\phi^k) - \bar{\mathbf{x}}_\phi^k \right\| = 0$, follows readily. This together with (268) and (267) prove

$$\lim_{k\to\infty} \left\| \widehat{\mathbf{x}}_i(\bar{\mathbf{x}}^k) - \bar{\mathbf{x}}^k \right\| = 0. \qquad (270)$$

- `Term II`: We prove $\lim_{k\to\infty} \|\bar{\mathbf{x}}(\bar{\mathbf{x}}^k) - \widehat{\mathbf{x}}_i(\bar{\mathbf{x}}^k)\| = 0$.

We begin deriving a proper upper bound of $\|\bar{\mathbf{x}}(\bar{\mathbf{x}}^k) - \widehat{\mathbf{x}}_i(\bar{\mathbf{x}}^k)\|$. By the first order optimality condition of $\bar{\mathbf{x}}(\bar{\mathbf{x}}^k)$ and $\widehat{\mathbf{x}}_i(\bar{\mathbf{x}}^k)$ we have

$$\left(\widehat{\mathbf{x}}_i(\bar{\mathbf{x}}^k) - \bar{\mathbf{x}}(\bar{\mathbf{x}}^k)\right)^T \left(\nabla F(\bar{\mathbf{x}}^k) + \bar{\mathbf{x}}(\bar{\mathbf{x}}^k) - \bar{\mathbf{x}}^k\right) + G\left(\widehat{\mathbf{x}}_i(\bar{\mathbf{x}}^k)\right) - G\left(\bar{\mathbf{x}}(\bar{\mathbf{x}}^k)\right) \geq 0,$$

$$\left(\bar{\mathbf{x}}(\bar{\mathbf{x}}^k) - \widehat{\mathbf{x}}_i(\bar{\mathbf{x}}^k)\right)^T \left(\nabla \widetilde{f}_i(\widehat{\mathbf{x}}_i(\bar{\mathbf{x}}^k) \,|\, \bar{\mathbf{x}}^k) + \nabla F(\bar{\mathbf{x}}^k) - \nabla f_i(\bar{\mathbf{x}}^k)\right)$$
$$+ G\left(\bar{\mathbf{x}}(\bar{\mathbf{x}}^k)\right) - G\left(\widehat{\mathbf{x}}_i(\bar{\mathbf{x}}^k)\right) \geq 0,$$

which yields

$$\begin{aligned}
\left\| \bar{\mathbf{x}}(\bar{\mathbf{x}}^k) - \widehat{\mathbf{x}}_i(\bar{\mathbf{x}}^k) \right\| &\leq \left\| \nabla \widetilde{f}_i(\widehat{\mathbf{x}}_i(\bar{\mathbf{x}}^k) \,|\, \bar{\mathbf{x}}^k) - \nabla f_i(\bar{\mathbf{x}}^k) - \widehat{\mathbf{x}}_i(\bar{\mathbf{x}}^k) + \bar{\mathbf{x}}^k \right\| \\
&\leq \left\| \nabla \widetilde{f}_i(\widehat{\mathbf{x}}_i(\bar{\mathbf{x}}^k) \,|\, \bar{\mathbf{x}}^k) \pm \nabla f_i(\widehat{\mathbf{x}}_i(\bar{\mathbf{x}}^k)) - \nabla f_i(\bar{\mathbf{x}}^k) \right\| + \left\| \widehat{\mathbf{x}}_i(\bar{\mathbf{x}}^k) - \bar{\mathbf{x}}^k \right\| \\
&\leq B_{17} \left\| \widehat{\mathbf{x}}_i(\bar{\mathbf{x}}^k) - \bar{\mathbf{x}}^k \right\|,
\end{aligned} \qquad (271)$$

where $B_{17}$ is some positive, finite, constant; and in the last inequality we used the Lipschitz continuity of $\nabla f_i(\bullet)$ (cf. Assumption III.1) and $\nabla \widetilde{f}_i(\mathbf{x} \,|\, \bullet)$ (cf. Assumption III.14).

Using (270) and (271), we finally have

$$\lim_{k\to\infty} \left\| \bar{\mathbf{x}}(\bar{\mathbf{x}}^k) - \widehat{\mathbf{x}}_i(\bar{\mathbf{x}}^k) \right\| = 0. \qquad (272)$$

Therefore, we conclude

$$\lim_{k\to\infty} J(\bar{\mathbf{x}}^k) \overset{(266)}{\leq} \underbrace{\lim_{k\to\infty} \left\| \widehat{\mathbf{x}}_i(\bar{\mathbf{x}}^k) - \bar{\mathbf{x}}^k \right\|}_{\overset{(270)}{=}0} + \underbrace{\lim_{k\to\infty} \left\| \bar{\mathbf{x}}(\bar{\mathbf{x}}^k) - \widehat{\mathbf{x}}_i(\bar{\mathbf{x}}^k) \right\|}_{\overset{(272)}{=}0} = 0. \qquad (273)$$

This completes the proof of Step 3, and also the proof of Theorem III.16. □

## III.4. Applications

In this section, we test the performance of SONATA (Algorithm 10) on both convex and nonconvex problems. More specifically, as convex instance of Problem (167), we consider a distributed linear regression problem (cf. Sec. III.4.1) whereas as nonconvex instances of (167), we study a target localization problem (cf. Sec. III.4.2). Other applications and numerical results can be found in [211, 230, 231].

### III.4.1 Distributed robust regression

The distributed robust linear regression problem is an instance of the empirical risk minimization considered in Example 1 in Sec. III.1.1: Agents in the network want to cooperatively estimate a common parameter $\mathbf{x}$ of a linear model from a set of distributed measures corrupted by noise and outliers; let $D_i \triangleq \{d_{i1}, \ldots, d_{in_i}\}$ be the set of $n_i$ measurements taken by agent $i$. To be robust to the heavy-tailed errors or outliers in the response, a widely explored approach is to use the Huber's criterion as loss function, which leads to the following formulation

$$\underset{\mathbf{x}}{\text{minimize}} \quad F(\mathbf{x}) \triangleq \sum_{i=1}^{I} \sum_{j=1}^{n_i} H\left(\mathbf{b}_{ij}^T \mathbf{x} - d_{ij}\right), \tag{274}$$

where $\mathbf{b}_{ij} \in \mathbb{R}^m$ is the vector of features (or predictors) associated with the response $d_{ij}$, owned by agent $i$, with $j = 1, \ldots, n_i$ and $i = 1, \ldots, I$; and $H : \mathbb{R} \to \mathbb{R}$ is the Huber loss function, defined as:

$$H(r) = \begin{cases} r^2, & \text{if } |r| \leq \alpha, \\ \alpha \cdot (2|r| - \alpha), & \text{otherwise;} \end{cases}$$

for some given $\alpha > 0$. This function is quadratic for small values of the residual $r$ (like the least square loss) but grows linearly (like the absolute distance loss) for large values of $r$. The cut-off parameter $\alpha$ describes where the transition from quadratic to linear takes place. Note that $H$ is convex (but not strongly convex) and differentiable, with derivative:

$$H'(r) = \begin{cases} -2\alpha, & \text{if } r < -\alpha, \\ 2r, & \text{if } r \in [-\alpha, \alpha], \\ 2\alpha, & \text{if } r > \alpha. \end{cases} \tag{275}$$

Introducing

$$f_i(\mathbf{x}) \triangleq \sum_{j=1}^{n_i} H\left(\mathbf{b}_{ij}^T \mathbf{x} - d_{ij}\right),$$

(274) is clearly an instance of Problem (167), with $F = \sum_{i=1}^{I} f_i$ and $G = 0$. It is not difficult to check that Assumption III.1 is satisfied.

We apply SONATA to Problem (274) considering two alternative choices of the surrogate functions $\widetilde{f}_i$. The first choice is the linear approximation of $f_i$ (plus a proximal regularization): given the local copy $\mathbf{x}_{(i)}^k$,

$$\widetilde{f}_i\left(\mathbf{x}_{(i)} \mid \mathbf{x}_{(i)}^k\right) = f_i\left(\mathbf{x}_{(i)}^k\right) + \nabla f_i\left(\mathbf{x}_{(i)}^k\right)^T \left(\mathbf{x}_{(i)} - \mathbf{x}_{(i)}^k\right) + \frac{\tau_i}{2} \left\|\mathbf{x}_{(i)} - \mathbf{x}_{(i)}^k\right\|^2; \quad (276)$$

with

$$\nabla f_i\left(\mathbf{x}_{(i)}^k\right) = \sum_{j=1}^{n_i} \mathbf{b}_{ij} \cdot H'\left(\mathbf{b}_{ij}^T \mathbf{x}_{(i)}^k - d_{ij}\right),$$

where $H'(\bullet)$ is defined in (275). An alternative choice for $\widetilde{f}_i$ is a quadratic approximation of $f_i$ at $\mathbf{x}_{(i)}^k$:

$$\widetilde{f}_i\left(\mathbf{x}_{(i)} \mid \mathbf{x}_{(i)}^k\right) = \sum_{j=1}^{n_i} \widetilde{H}_{ij}\left(\mathbf{x}_{(i)} \mid \mathbf{x}_{(i)}^k\right) + \frac{\tau_i}{2} \|\mathbf{x}_{(i)} - \mathbf{x}_{(i)}^k\|^2, \quad (277)$$

where $\widetilde{H}_{ij}$ is given by

$$\widetilde{H}_{ij}\left(\mathbf{x}_{(i)} \mid \mathbf{x}_{(i)}^k\right) = \begin{cases} \dfrac{\alpha}{r_{ij}^k} \cdot \left(\mathbf{b}_{ij}^T \mathbf{x}_{(i)} - d_{ij}\right)^2, & \text{if } \left|r_{ij}^k\right| \geq \alpha, \\ \left(\mathbf{b}_{ij}^T \mathbf{x}_{(i)} - d_{ij}\right)^2, & \text{if } \left|r_{ij}^k\right| < \alpha; \end{cases}$$

with $r_{ij}^k \triangleq \mathbf{b}_{ij}^T \mathbf{x}_{(i)}^k - d_{ij}$.

Note that the subproblems in Step 2 of Algorithm 10, with $\widetilde{f}_i$ given by (276) or (277), have a closed form solution. In particular, $\widetilde{\mathbf{x}}_{(i)}^k$ associated with the surrogate in (277) is given by

$$\widetilde{\mathbf{x}}_i^k = \left(\tau_i \mathbf{I} + 2\mathbf{B}_i^T \mathbf{D}_i^k \mathbf{B}_i\right)^{-1} \left(\tau_i \mathbf{x}^k - (I \cdot \mathbf{y}_{(i)}^k - \nabla f_i(\mathbf{x}_{(i)}^k)) + 2\mathbf{B}_i^T \mathbf{D}_i \mathbf{d}_i\right),$$

where $\mathbf{B}_i$ is the matrix whose $j$-th row is $\mathbf{b}_{ij}^T$; $\mathbf{d}_i$ is the vector whose $j$-th component is $d_{ij}$; and $\mathbf{D}_i^k$ is a diagonal matrix whose $j$-th diagonal entry is $\min\{1, \alpha/r_{ij}^k\}$.

**Numerical example.** We simulate a network of $I = 30$ agents, modeled as a directed time-varying graph. At each iteration, the graph is generated according to the following procedure. All the agents are connected through a time-varying ring; at each iteration the order of the agents in the ring is randomly permuted. In addition to the ring topology, each agent has one out-neighbor, selected at each iteration uniformly at random. In (274), the parameter $\mathbf{x}$ to estimate has dimension 200, with *i.i.d.* uniformly distributed entries in $[-1,1]$; $n_i = 20$ (number of measures per agent); and the elements of the vectors $\mathbf{b}_{ij}$ are *i.i.d.* Gaussian distributed and then normalized so that $\|\mathbf{b}_{ij}\| = 1$; the noise affecting the measurements is generated according to a Gaussian distribution, with standard deviation $\sigma = 0.1$, and each agent has one mea-surement corrupted by one outlier following a Gaussian distribution, with standard deviation $5 \cdot \sigma$; and finally, the cut-off parameter $\alpha$ is set to be $\alpha = 3 \cdot \sigma$.

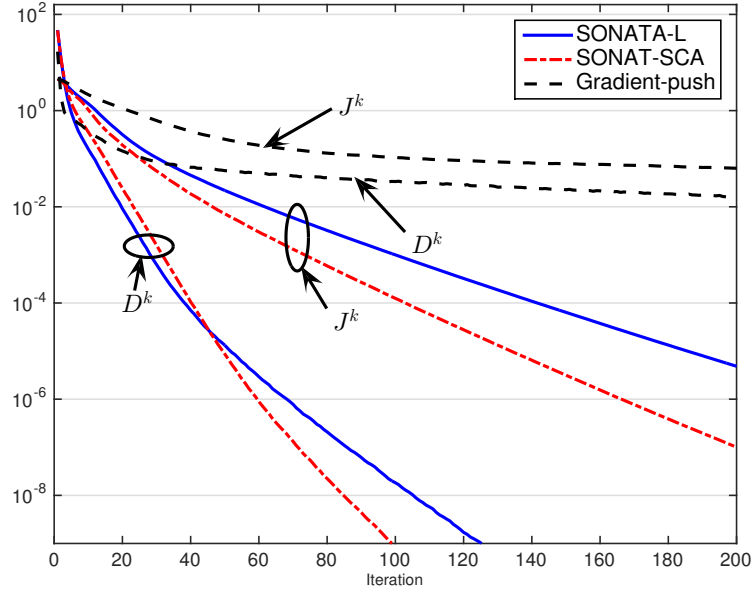

*Fig. III.3: Robust linear regression* (274)*: Distance from optimality $J^k$ and consensus disagreement $D^k$ versus the number of message exchanges. Both instances of SONATA significantly outperform the (sub)gradient-push [168].*

The free parameters in SONATA (Algorithm 10) are tuned as follows. The proximal parameters $\tau_i$ in (276) and (277) are set to $\tau_i = 2$ and $\tau_i = 1.5$, respectively, for all $i = 1, \ldots, I$. In both instances of the algorithm, the step-size $\gamma^k$ is chosen according to the rule (108), with parameters $\gamma^0 = 0.1$ and $\mu = 0.01$. In the consensus and tracking updates (Step 3), the matrices $\mathbf{A}^k$ are chosen according to the push-sum rule: $a_{ij}^k = 1/d_j^k$. We will refer to the two instances of SONATA as *SONATA-L* ("L" stands for linearized) for the one associated with the surrogate in (276), and *SONATA-SCA* for the one associated with the surrogate (277). We compare SONATA with the subgradient-push algorithm proposed in [168]. The step-size in the algorithm is chosen according to (108), with parameters $\gamma^0 = 0.5$, $\mu = 0.01$ (this choice resulted in the best practical performance among all the tested tunings).

To compare the two algorithms, as merit functions, we used $J^k \triangleq \|\nabla F(\bar{\mathbf{x}}^k)\|_\infty$ and $D^k \triangleq (1/I) \cdot \sum_{i=1}^{I} \|\mathbf{x}_{(i)}^k - \bar{\mathbf{x}}^k\|^2$; the former measures the distance of the average of the iterates $\bar{\mathbf{x}}^k$ from stationarity whereas the latter is a measure of the consensus disagreement. In Fig. III.3 we plot $J^k$ and $D^k$ for all the algorithms versus the number of messages exchanged at each iteration $k$. For the (sub)gradient-push algorithm, this number coincides with the iteration index $k$ whereas for the two instances of SONATA it is $2 \cdot k$ (recall that SONATA employs two communications at each iteration). The curves are averaged over 100 Monte-Carlo simulations: in each trial, the parameter **x** does not change while the noise and graph connectivity are randomly generated (as described above). The analysis of the figure clearly shows that both instances of SONATA lock consensus and converge much faster than the (sub)gradient-push. Furthermore, SONATA-SCA outperforms SONATA-L, since it better exploits the convexity of the loss function.

### III.4.2 Target localization

Consider the target localization problem described in Example 3 of Sec. III.1.1:

$$\begin{aligned} \underset{\mathbf{x} \triangleq (\mathbf{x}_t)_{t=1}^n}{\text{minimize}} \quad & F(\mathbf{x}) \triangleq \sum_{i=1}^{I} \sum_{t=1}^{n} p_{it} \left( d_{it} - \|\mathbf{x}_t - \mathbf{s}_i\|^2 \right)^2 \\ \text{subject to} \quad & \mathbf{x}_{(i)} \in X \subseteq \mathbb{R}^m, \quad \forall i = 1, \ldots, I, \end{aligned} \quad (278)$$

where we recall that $n$ is the number of targets; $\mathbf{s}_i$ is the vector of the coordinates of sensor $i$'s location; $d_{it}$ is the squared Euclidean distance between the position of sensor $i$ and that of target $t$; $\mathbf{x}_t$ is an estimate (to be found) of the location of target $t$; $X$ is a closed convex set; and $p_{it} \in \{0, 1\}$ is a scalar taking value zero if the $i$th agent does not have any measurement related to target $t$. Clearly (278) is an instance of Problem (167), with $F = \sum_{i=1}^{I} f_i$, $G = 0$, and

$$f_i(\mathbf{x}) = \sum_{t=1}^{n} p_{it} \left( d_{it} - \|\mathbf{x}_t - \mathbf{s}_i\|^2 \right)^2. \quad (279)$$

We consider next two instances of SONATA, associated with two alternative surrogate functions $\widetilde{f}_i$. A first choice is using the first order approximation of $f_i$ in (279), that is,

$$\widetilde{f}_i \left( \mathbf{x}_{(i)} \mid \mathbf{x}_{(i)}^k \right) = f_i \left( \mathbf{x}_{(i)}^k \right) + \nabla f_i \left( \mathbf{x}_{(i)}^k \right)^T \left( \mathbf{x}_{(i)} - \mathbf{x}_{(i)}^k \right) + \frac{\tau_i}{2} \left\| \mathbf{x}_{(i)} - \mathbf{x}_{(i)}^k \right\|^2, \quad (280)$$

where $\mathbf{x}_{(i)}^k \triangleq (\mathbf{x}_{(i),t}^k)_{t=1}^n$, $\nabla f_i(\mathbf{x}_{(i)}^k) = [\nabla_1 f_i(\mathbf{x}_{(i),1}^k)^T, \ldots, \nabla_n f_i(\mathbf{x}_{(i),n}^k)^T]^T$, and

$$\nabla_t f_i(\mathbf{x}_{(i),t}^k) = -4 \cdot \left( d_{it} - \|\mathbf{x}_{(i),t}^k - \mathbf{s}_i\|^2 \right) \cdot \left( \mathbf{x}_{(i),t}^k - \mathbf{s}_i \right).$$

A second choice for $\widetilde{f}_i$ is motivated by the observation that (279) is a fourth-order polynomial in each $\mathbf{x}_t$; one may want to preserve the "partial" convexity in $f_i$ by keeping the first and second order (convex) terms in each summand of (279) unaltered and linearizing the higher order terms. This leads to the following:

$$\widetilde{f}_i \left( \mathbf{x}_{(i)} \mid \mathbf{x}_{(i)}^k \right) = \sum_{t=1}^{n} p_{it} \left( \widetilde{f}_{it} \left( \mathbf{x}_{(i),t} \mid \mathbf{x}_{(i)}^k \right) + \frac{\tau_i}{2} \left\| \mathbf{x}_{(i),t} - \mathbf{x}_{(i),t}^k \right\|^2 \right), \quad (281)$$

where $\mathbf{x}_{(i)} \triangleq (\mathbf{x}_{(i),t})_{t=1}^n$,

$$\widetilde{f}_{it}\left(\mathbf{x}_{(i),t} \mid \mathbf{x}^k\right) = \mathbf{x}_{(i),t}^T \mathbf{S}_i \mathbf{x}_{(i),t} - \mathbf{b}_{it}^{kT}\left(\mathbf{x}_{(i),t} - \mathbf{x}_{(i),t}^k\right),$$

with $\mathbf{S}_i = 4 \cdot \mathbf{s}_i \mathbf{s}_i^T + 2 \cdot \|\mathbf{s}_i\|^2 \mathbf{I}$, and

$$\mathbf{b}_{it}^k = 4\|\mathbf{s}_i\|^2 \mathbf{s}_i - 4\left(\|\mathbf{x}_{(i),t}^k\|^2 - d_{it}\right)\left(\mathbf{x}_{(i),t}^k - \mathbf{s}_i\right) + 8\left(\mathbf{s}_i^T \mathbf{x}_{(i),t}^k\right)\mathbf{x}_{(i),t}^k.$$

**Numerical example.** We simulate a network of $I = 30$ sensors, modeled as a directed time-varying graph; and we deploy $n = 5$ targets. The graph is generated as described in the previous example (cf. Sec. III.4.1). The locations of the sensors and targets are generated uniformly at random in $[0,1]^2$. The parameters $p_{it}$ take values $\{0,1\}$, with equal probability. The distances $d_{it}$ are corrupted by *i.i.d.* zero-mean Gaussian noise, with standard deviation set to be the minimum pairwise distance between sensors and targets. We let the constraint set $X = \mathbb{R}^2$.

The free parameters in SONATA (Algorithm 10) are tuned as follows. The proximal parameters $\tau_i$ in (280) and (281) are set to $\tau_i = 7$ and $\tau_i = 5$, respectively, for all $i = 1, \ldots, I$. In both instances of the algorithm, the step-size $\gamma^k$ is chosen according to the rule (108), with parameters $\gamma^0 = 0.1$ and $\mu = 10^{-4}$. In the consensus and tracking updates (Step 3), the matrices $\mathbf{A}^k$ are chosen according to the push-sum rule: $a_{ij}^k = 1/d_j^k$. We will refer to the two instances of SONATA as *SONATA-L* ("L" stands for linearized) for the one associated with the surrogate in (280), and *SONATA-SCA* for the one associated with the surrogate (281). We compare SONATA with the distributed gradient [232] for unconstrained optimization (adapted here to account for the constraints; note that for such a version there is no formal proof of convergence in the literature). The step-size in the algorithm [232] is chosen according to (108), with parameters $\gamma^0 = 0.5$, $\mu = 0.01$ (this choice resulted in the best practical performance among all the tested tunings).

A comparison of the algorithms is given in Fig. III.4, where we plotted the merit functions $J^k \triangleq \|\nabla F(\bar{\mathbf{x}}^k)\|_\infty$ and $D^k \triangleq (1/I) \cdot \sum_{i=1}^{I} \|\mathbf{x}_{(i)}^k - \bar{\mathbf{x}}^k\|^2$ versus the number of messages exchanged at each iteration $k$ (cf. Sec. III.4.1). The curves are averaged over 100 Monte-Carlo simulations; in each trial, the position of the sensors and targets were kept fixed, while the noise and graph connectivity were randomly generated every trial. Fig. III.4 shows that both instances of SONATA converge much faster than the benchmark distributed gradient algorithm. This can be explained by the fact that SONATA better exploits the problem's structure and employs the tracking of the average of agents' gradients .

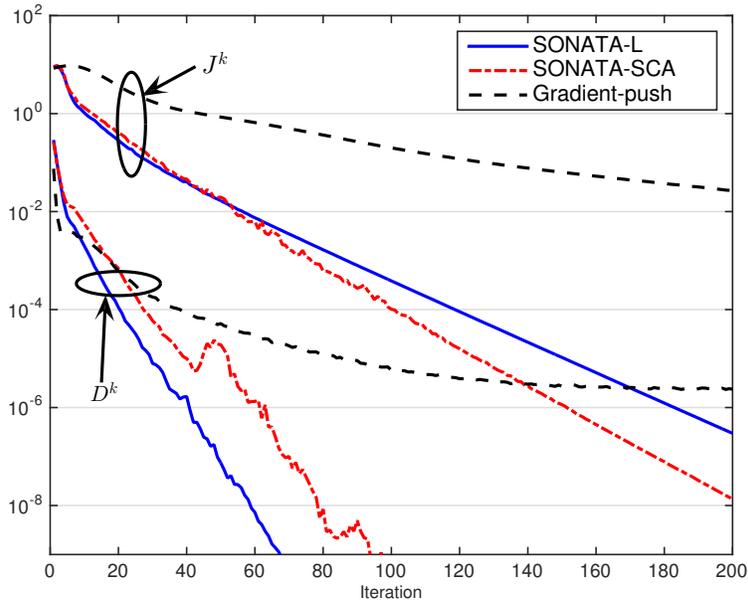

*Fig. III.4: Target localization problem* (278)*: Distance from optimality $J^k$ and consensus disagreement $D^k$ versus the number of message exchanges. Both instances of SONATA significantly outperform the distributed gradient algorithm [232].*

## III.5. Sources and Notes

Distributed optimization has a rich literature, dating back to the 1980s [240]. Since then, due to the emergence of large-scale networks, advances in computing and storage, data acquisition, and communications, the development of distributed algorithms for solving optimization problems over networks has received significant attention. Besides the distributed approaches we are going to summarize below, many efforts have been made to solve optimization problems over networks with *specific topology*, including multilayer HNets (e.g., master-slave architectures and shared-memory systems)–see Fig. III.1. However, all these schemes require some degree of centralization to be implemented, due to the use of master nodes. Since the focus of this lecture is on distributed algorithms implementable on general network architectures, in the following we omit to discuss the literature on algorithms for HNets. Although the lecture considered mainly nonconvex optimization, we begin overviewing the much abundant literature of distributed algorithms for convex problems. In fact, some of the methods and tools introduced in this lecture for nonconvex problems can be applied and are useful also for convex instances.

**Distributed convex optimization**

The literature on distributed solution methods for convex problems is vast and a comprehensive overview of current methods goes beyond the scope of this commentary. Here we briefly focus on the methods somehow related to the proposed SCA approach and distributed gradient tracking. Generally speaking, we roughly divide the literature in two groups, namely: primal and dual methods.

**Primal methods**. While substantially different, primal methods can be generically abstracted as a combination of a local (sub)gradient step and a subsequent consensus-like update (or multiple consensus updates); examples include [114, 168, 169, 216, 217]. Algorithms for adaptation and learning tasks based on in-network diffusion techniques were proposed in [39, 46, 201]. Schemes in [39, 46, 169, 201, 217] are applicable only to *undirected* graphs; [39, 169, 217] require the consensus matrices to be *doubly-stochastic* whereas [46, 201] use *row-stochastic* matrices but are applicable only to strongly convex agents' cost functions having a *common* minimizer. When the graph is *directed*, doubly-stochastic weight matrices compliant to the graph may not exist; and, when they exist, their construction is computationally prohibitive for networks with arbitrary topology [92]. Recently, some new algorithms have been developed that do not require doubly-stochastic weights. Examples, include [168], where the authors combined the sub-gradient algorithm [169] with the push-sum consensus protocol [121]; and [256, 257]. However, the schemes in [169, 256, 257] considered only *unconstrained* optimization problems; in fact, up until the very recent works [211, 230, 231], it was not clear how to leverage push-sum protocols to deal with constraints (e.g., to preserve the feasibility of the iterates).

In all the aforementioned algorithms but [211, 231], agents perform their local optimization using only the gradient of their *own* function, while neglecting the rest of the sum-cost function $F$ [cf. (167)]. This makes the algorithms simple but usu-ally slow: Even if the objective functions are differentiable and strongly convex, these methods still need to use a diminishing step-size to converge to a consensual solution, and thus converge at a sub-linear rate. On the other hand, with a constant (sufficiently small) step-size, these distributed methods can be faster, but they only converge to a neighborhood of the solution set. This phenomenon creates an exactness-speed dilemma. In addition, convergence of these algorithms was proved only under the assumption that the objective function has *bounded (sub)gradient*, which limits the range of applicability of these schemes.

A different class of distributed (primal) approaches that bypasses the above issues is the one based on the idea of *gradient tracking*: each agent updates its own local variables along a surrogate direction that tracks the gradient average $\sum_i \nabla f_i$. This idea was proposed independently in [69, 70, 148] for the general class of constrained nonsmooth nonconvex problems (167), and in [259] for the special case of strongly convex, unconstrained, smooth, optimization. The works [211, 231] extended the algorithms to (possibly) time-varying digraphs (still in the nonconvex setting of [69, 70, 148]). A convergence rate analysis of the scheme [259] was later developed in [170, 190], with [170] considering time-varying (directed) graphs. The benefit of these schemes is threefold: i) a constant step-size can be employed, instead of a diminishing one; ii) the assumption that the gradient is bounded is no longer needed; and iii) the convergence speed of the algorithms significantly improves– linear convergence was in fact proved in [170, 190] for unconstrained, smooth, strongly convex $\sum_i f_i$, and more recently in [229] for the more general class of constrained, nonsmooth, problems (167) (with $F$ strongly convex).

Finally, it is worth mentioning that a difference algorithmic structure was proposed in [216] to cancel the steady state error in the aforementioned plain decentralized gradient descent methods and converge to a consensual solution using a fixed (sufficiently small) step-size. Convergence at a sub-linear rate was proved when the objective function is smooth and convex (and the problem unconstrained), and a Q-linear rate when the objective function is (smooth) strongly convex (still unconstrained optimization).

**Dual-based methods**: This class of algorithms is based on a different approach: slack variables are first introduced to decouple the sum-utility function while forcing consistency among these local copies by adding consensus equality constraints (matching the graph topology). Lagrangian dual variables are then introduced to deal with such coupling constraints. The resulting algorithms build on primal-dual updates, aiming at converging to a saddle point of the (augmented) Lagrangian function. Examples of such algorithms include ADMM-like methods [42, 113, 217, 247] as well as inexact primal-dual instances wherein a first order [43] or second-order [163, 164] approximation of the Lagrangian function is minimized at each iteration. All these algorithms can handle only *static and undirected* graphs. Dual methods cannot be trivially extended to time-varying graphs or digraphs, as it is not clear how to enforce consensus via equality constraints over time-varying or directed networks. Furthermore, all the above schemes but [42, 217] require the objective function to be *smooth* and the optimization problem to be *unconstrained*.

|  |  | Proj-DGM [20] | NEXT [148] | Push-sum DGM [232] | Prox-PDA [106] | DeFW [245] | SONATA [211] (this chapter) |
|---|---|---|---|---|---|---|---|
| **nonsmooth $G^+$** | |  | ✓ |  |  |  | ✓ |
| **constraints** | | ✓ | ✓ |  |  | $X$ compact | ✓ |
| **unbounded gradient** | |  |  |  | ✓ |  | ✓ |
| **network topology** | time-varying |  | ✓ | ✓ |  |  | ✓ |
|  | digraph |  | restricted | ✓ |  |  | ✓ |
| **step-size:** | constant |  |  |  | ✓ |  | ✓ |
|  | diminishing | ✓ | ✓ | ✓ | ✓ | ✓ | ✓ |
| **complexity** | |  |  | ✓ | ✓ | ✓ | ✓ |

*Table III.3: Overview of current distributed methods for nonconvex optimization over networks.*

### III.5.1 Distributed nonconvex optimization

For the general nonconvex problem (167), results are scarse. We remark that the parallel algorithm described in Lecture II as well as those discussed in Sec. II.7 are not applicable to the distributed setting (167): either they require the graph to be complete (i.e., at each iteration, agents must be able to exchange information with *all* the others) or they assume full knowledge of $F$ from the agents. There are only a few distributed algorithms designed for (167), including primal methods [20, 148, 232, 245] and dual methods [106, 275]. The key features of these algorithms are summarized in Table III.3 and briefly discussed next.

**Primal methods**: The distributed stochastic projection scheme [20] implements random gossip between agents and uses diminishing step-size, and can handle *smooth* objective functions over *undirected static* graphs. In [232], the authors showed that the distributed push-sum gradient algorithm with diminishing step-size, earlier proposed for convex objectives in [168], converges also when applied to nonconvex *smooth unconstrained* problems. The first provably convergent distributed scheme for (167), with $G \neq 0$ and (convex) constraints $X$, over time-varying digraphs is NEXT [69, 70, 148]. A special instance of NEXT, applied to smooth $V$ over *undirected static* graphs, was later proved in [245] to have sub-linear convergence rate. However, the consensus protocol employed in [148, 245] uses *doubly-stochastic* weight matrices, and thus is limited to special digraphs (or undirected graphs). Moreover, all the algorithms discussed above require that the (sub)gradient of $V$ is *bounded* on $X$ (or $\mathbb{R}^m$).

This issue was resolved in SONATA (Algorithm 10, Sec. III.3), proposed in [211, 231], which combined a judiciously designed perturbed push-sum consensus scheme with SCA techniques. Based on a new line of analysis whereby optimization and consensus dynamics are *jointly* analyzed by introducing a novel Lyapunov function, SONATA was shown to converge sub-linearly to a stationary solution of (167), without requiring the (sub)gradient of $V$ being bounded. SONATA is also the first provably convergent algorithm that i) can deal with constrained, nonsmooth, (convex and nonconvex) problems in the form (167), and time-varying digraphs; and ii) employs either a constant or diminishing step-size.

**Dual-based methods**: As their convex counterparts, dual-based schemes [106, 275] cannot deal with time-varying or directed graphs. The algorithm proposed in [275] is an early attempt of distributed method for nonconvex problems. However, it calls for the solution of nonconvex subproblems and converges to (stationary) solutions of an auxiliary problem, which are not necessarily stationary for the original problem. The prox-PDA algorithm in [106], proposed for linearly constrained problems, can be applied to Problem (167), but it is limited to the special classes where the functions $f_i$ are *smooth* and $X = \mathbb{R}^m$. The dual-based algorithms discussed above but [106] assume that the (sub)gradient of $V$ is *bounded*.

**Practical implementation issues and open questions.** The majority of the methods discussed above assume an "ideal" distributed computing environment wherein all the communication links operate reliably and noise-free, and the agents communicate and update their own variables in a synchronous fashion. There are however many practical issues that arise for these models, including difficulties associated with maintaining time-synchronization in face of computational and communication delays, reaching deadlock situations requiring re-initializations of the update process due to node/link failures.

Asynchronous methods based on randomization techniques were proposed in [18–20, 22, 110, 128, 140, 166, 174, 224, 247, 272]. These algorithms assume that agents are activated (wakeup) randomly (with no coordination), run some local computation, and then wake up (some or all) their neighbors to pass them their most recent update. Both random (e.g., based on gossip consensus schemes) [18–20, 110, 140, 166, 224, 247, 272] and deterministic (e.g., based on the push-sum protocol) activations [22, 128, 174] of the neighbors were proposed. While asynchronous in the activation of the agents performing the update, these schemes still impose a coordination among the agents, cannot tolerate (arbitrary) delays in the communications (packet losses are considered in [22, 128, 272]) or optimization (out-of-sync information), and require the agents to stay all the time in listening mode. A first attempt to deal with this issue was proposed in [255]: a dual-based scheme is employed wherein agents are randomly activated, and update their local copies using possibly outdated information. The probabilistic model for asynchrony requires that the random variables triggering the activation of the agents and the delay vector used by the agent to performs its update are *independent*. While this greatly simplifies the convergence analysis, it makes the model not realistic–see Sec. II.7 in Lec. II for a discussion on this issue. Furthermore, [255] is applicable only to convex instances of (167), and over fixed undirected graphs. At the time of this writing, designing distributed fully asynchronous algorithms for the nonconvex problem (167)–in the sense that i) agents update their own variables with no coordination; and ii) can handle arbitrary (deterministic or random) delay profiles–remains an open problem. A first attempt towards this direction can be found in [234].

Despite of asynchrony, another question of practical interest is how to design and implement distributed algorithms under low communication requirements, and what trade-offs are involved in such designs/implementations. For instance, in some sensor networks, broadcast communications can be too expensive in terms of consumed power. Also, communications are subject to rate constraints (systems are bandwidth limited). Furthermore, in the emerging distributed big-data applications, optimization problems are huge-scale: computing at each iteration the gradient with respect to all the optimization variables, communicating at each iteration, and sending to the neighbors the entire set of variables would incur in an unaffordable computational cost and communication overhead or is just infeasible. Some initial investigations along this lines were presented in the literature, following alternative, scattered, paths. For instance, in [175, 176], a distributed algorithm for convex and nonconvex instances of (167) was proposed, aiming at reducing both communication and computation costs by optimizing and transmitting at each iteration only a subset of the entire optimization variables. While asymptotic convergence of the method was proved, the tradeoff between communication and computation has not been studied yet. Some insight on the impact of multiple communications steps (without optimizing) or multiple optimization steps (without communicating) can be found in [13, 239]. Finally, some works proposed distributed methods for convex optimization problems implementing quantized communications. Of course, much more research needs to be done on these important topics.

Title Suppressed Due to Excessive Length 151